\documentclass[12pt]{article}
\usepackage[authoryear,round]{natbib}
\usepackage{graphicx}
\usepackage{latexsym}
\usepackage[left=1.5in,top=1.5in,right=1.5in,bottom=1.5in]{geometry}\usepackage{amsmath}
\usepackage{booktabs}
\usepackage{color}
\usepackage{tikz}
\usetikzlibrary{shapes}
\usetikzlibrary{arrows}
\definecolor{darkblue}{rgb}{0,0.4,0.9}
\definecolor{gray10}{rgb}{0.1,0.1,0.1}
\definecolor{gray20}{rgb}{0.2,0.2,0.2}
\definecolor{gray30}{rgb}{0.3,0.3,0.3}
\definecolor{gray40}{rgb}{0.4,0.4,0.4}
\definecolor{gray60}{rgb}{0.6,0.6,0.6}
\definecolor{gray80}{rgb}{0.8,0.8,0.8}
\definecolor{gray90}{rgb}{0.9,0.9,.9}
\definecolor{gray95}{rgb}{0.95,0.95,.95}
\definecolor{gray96}{rgb}{0.96,0.96,.96}
\definecolor{lgreen} {RGB}{180,210,100}
\definecolor{dblue}  {RGB}{20,66,129}
\definecolor{ddblue} {RGB}{11,36,69}
\definecolor{lred}   {RGB}{220,0,0}
\definecolor{nred}   {RGB}{224,0,0}
\definecolor{norange}{RGB}{230,120,20}
\definecolor{nyellow}{RGB}{255,221,0}
\definecolor{ngreen} {RGB}{98,158,31}
\definecolor{dgreen} {RGB}{78,138,21}
\definecolor{nblue}  {RGB}{28,130,185}
\definecolor{jblue}  {RGB}{20,50,100}
\definecolor{nnyellow}{RGB}{235,200,0}
\definecolor{purple}{RGB}{150, 0, 120}
\definecolor{sgGreen} {RGB}{20, 180, 50}
\definecolor{revised}{rgb}{0,0,0.9}
\newtheorem{definition}{Definition}

\newtheorem{theorem}{Theorem}
\newtheorem{corollary}{Corollary}
\newtheorem{lemma}{Lemma}

\newtheorem{assumption}{Assumption}
\newcommand{\nl}{\newline}
\newcommand{\pl}{\parallel}
\newcommand{\openr}{\hbox{${\rm I\kern-.2em R}$}}
\newcommand{\openn}{\hbox{${\rm I\kern-.2em N}$}}

\usepackage{authblk}

\title{Higher Order Spline Highly Adaptive Lasso Estimators of Functional Parameters: Pointwise Asymptotic Normality and Uniform Convergence Rates}

\date{\today}
\author{Mark van der Laan}
\affil{University of California, Berkeley}
 
\begin{document}
\maketitle
\begin{abstract} 

We consider estimation of a functional of the data distribution based on observing a sample of i.i.d. observations. We assume the target function can be defined as the minimizer of the expectation of a loss function over a class of $d$-variate real valued cadlag functions that  have finite sectional variation norm. Examples of loss functions are the squared error loss for regression, or the log-likelihood loss for the density. 

For all $k=0,1,\ldots$, we define a  $k$-th order smoothness class of functions as $d$-variate functions on the unit cube for which each of a sequentially defined $k$-th  order Radon-Nikodym derivative w.r.t. Lebesgue measure  is cadlag  and of bounded variation. For a target function in this $k$-th order smoothness class we provide a representation of the target function as an infinite linear combination of tensor products of $\leq k$-th order spline basis functions indexed by a knot-point, where the lower (than $k$) order spline basis functions are used to represent the function at the $0$-edges.  The $L_1$-norm of the coefficients in this representation represents the sum of the variation norms across all the $k$-th order derivatives, which is called the $k$-th order sectional variation norm of the target function. This generalizes the previous results for the zero order spline representation  of cadlag functions with bounded sectional variation norm to higher order smoothness classes. We use this $k$-th order spline representation of a function to define the $k$-th order spline sieve minimum loss estimator (MLE), Highly Adaptive Lasso (HAL) MLE, and Relax HAL-MLE.    For first and higher order smoothness classes, in this article we analyze these three classes of estimators and establish pointwise asymptotic normality and uniform convergence at dimension free rate $n^{-k^*/(2k*+1)}$ up till a power of $\log n$ depending on the dimension, where $k^*=k+1$, assuming appropriate undersmoothing is used in selecting the $L_1$-norm. We also establish asymptotic linearity of plug-in estimators of pathwise differentiable features of the target function. 
\end{abstract}

\section{Introduction}
We consider estimation of a functional of the data distribution based on observing a sample of i.i.d. Euclidean valued observations. We assume the functional can be defined as the minimizer over cadlag functions \citep{Neuhaus71} with finite sectional variation norm of the expectation of a loss function, such as the squared error loss for regression, or the log-likelihood loss for the data density, over a nonparametric function class. 

Statistical inference for multivariate real valued functions in nonparametric or, more generally, infinite dimensional models is a highly challenging area of research. These functions are generally non-pathwise differentiable functionals of the data density so that there is no efficiency theory guiding the construction of asymptotically linear estimators
\citep{Bickeletal97,vanderLaan&Robins03,vanderLaan&Rubin06,vanderLaan&Rose11,vanderLaan&Rose18}.
 The typical kernel or local averaging type estimators often allow for a formal asymptotic normality result but suffer enormously from the curse of dimensionality. This is reflected by their rate of convergence strongly depending on the dimension of the function, so that a large degree of smoothness is required to still achieve reasonable rates of convergence. Many typical machine learning algorithms do not allow for statistical inference due to their ad hoc nature resulting in these estimators not solving score equations, the basis for establishing asymptotic normality. The sieve maximum likelihood estimators are more suitable for being  asymptotically  normally distributed by solving score equations 
 but the choice of sieve heavily affects the rate of convergence, and controlling bias by selection of tuning parameters are generally a big challenge, explaining the  lack of formal asymptotic normality theorems for these type of estimators. Some general asymptotic normality theory for smooth (pathwise differentiable) features of plug-in sieve MLE has been developed in the 
literature \citep{Shen97,Shen07,Qiu&Luedtke&Carone21,Newey94,vanderLaan&Benkeser&Cai19}, but the current article is focussing on the more challenging asymptotic normality of the function estimator itself. 
 
 In \citep{vanderLaan15,Benkeser&vanderLaan16,vanderLaan17} we proposed an MLE over the class of $d$-variate real valued cadlag functions with a bound on the sectional variation norm \citep{Gill&vanderLaan&Wellner95}.
It was shown that this MLE corresponds with minimizing the empirical risk over linear combinations of tensor products of zero-order spline basis functions under an $L_1$-norm constraint, where the $L_1$-norm actually equals the sectional variation norm of the function.  Therefore it was named the highly adaptive lasso MLE (HAL-MLE). We also proposed an explicit non-informative sufficient initial set of knot-point of size $\sim n (2^d-1)$ on which one can run the HAL-MLE, where this knot-point set is  immediately implied by the observed data points. 
It was shown that this HAL-MLE (which has  at most $n-1$ non-zero  knot-point specific coefficients) converges in loss-based dissimilarity at a rate $n^{-2/3}(\log n)^{d}$, which generally implies an $L^2$-rate of convergence of $n^{-1/3}(\log n)^{d/2}$ to the  true target function \citep{Bibaut&vanderLaan19,Fangetal19,Gao13}. This rate improved on the initial rate established in \citep{vanderLaan17}. Thus this rate of convergence is only weakly related to the dimension of the target function through a power of the $\log n$-factor. 
In addition \citep{vanderLaan&Benkeser&Cai19} showed that an undersmoothed HAL-MLE estimates smooth features of the target function efficiently (i.e., asymptotically normal limit distribution with minimal asymptotic variance). Moreover, this asymptotic efficiency applies  uniformly over a large class of smooth features, such as the class of target features whose efficient influence curve are contained in set of  cadlag functions with a universal bound on the sectional variation norm. 
Moreover our work on higher order TMLE \citep{vanderLaan&Wang21} shows that the plug-in undersmoothed HAL-MLE for a pathwise differentiable target feature of true function  can be represented as  a higher order TMLE itself that approximately solves higher order efficient influence curve equations that  reduce the general second order remainder of a first order TMLE  to a higher order difference. 

Given these powerful results on HAL-MLE w.r.t. $L^2$-rate of convergence and global estimation of pathwise differentiable target features of the target function, one wonders if the HAL-MLE or an higher order spline generalization, when evaluated at a point, is asymptotically normally distributed as well, and, if so,  at what rate?
 In addition, one wonders what is the supremum norm rate of convergence for the HAL-MLE or its higher order analogue? These are the questions addressed in this article.

In this article we will define a general $k$-th order spline representation of a $k$-th order smooth function whose $k$-th order derivatives are cadlag and have bounded variation norm, thereby generalizing the representation theorem for cadlag functions ($k=0$). This then naturally defines the generalization of the zero-order spline HAL-MLE to a $k$-th order spline HAL-MLE, where the $L_1$-norm now corresponds with the sum of variation norms of $k$-th order derivatives.

The lasso defines a data adaptive working model (i.e., linear combination of $\leq k$-th order spline basis functions with non-zero coefficients), and the HAL-MLE acts as a regularized MLE for that data adaptive working model. To avoid these two complexities initially, we also  define an analogue  $k$-th order spline sieve MLE indexed by an a-priori selected family of growing sets of $\leq k$-th order spline basis functions, and use a regular MLE instead of $L_1$-regularized MLE, while selecting the specific set of basis functions with cross-validation (or undersmoother). In addition, we consider the relax HAL-MLE that first runs the HAL-MLE and then refits the working model defined by the non-zero coefficients with regular MLE. Our understanding of this $k$-th order spline sieve MLE will also  guide our theoretical results for the $k$-th order spline HAL-MLE and the corresponding relax HAL-MLE.
This article is concerned with analyzing these three classes of $k$-th order spline sieve and HAL-MLEs.

Beyond the independent interest in the $k$-th order spline sieve  MLE, a maximal version of our working models for this sieve MLE can be inputted into the $k$-th order spline Highly Adaptive Lasso as initial working model, chosen so that the (non-penalized) MLE using this working model is already converging at the optimal rate, while selecting the $L_1$-norm with cross-validation. In this manner, the resulting $k$-th order spline HAL-MLE becomes more scalable,  and we still have the benefit of adaptively selected knot-points that take into account the outcome data.

 \subsection{Organization of article}
 
 In the introductory Section \ref{section2} we first formulate the estimation problem and thereby define the target function $Q_0=Q(P_0)$ we have to  learn from the i.i.d. observations from $P_0$, and review the zero-order spline representation of a cadlag function of finite sectional variation norm, an additive model across all subsets of $\{1,\ldots,d\}$, and the corresponding zero-order spline HAL-MLE. 
 In Appendix \ref{AppendixA} we present some key results from the current literature on approximating $d$-variate cumulative distribution functions with linear combinations of zero-order splines that proves that there is a set of ${\cal R}(d,J)$ of $J$ knot-points so that the resulting  $L^2(\mu)$-approximation error of a $d$-variate cumulative distribution function 
 is $O(r(d,J))$, where $r(d,J)\sim O(1/J)$ up till a power of $\log J$-factor depending on the dimension \citep{Fangetal19}. Specifically, throughout this article ${\cal R}(d,J)$ satisfies this zero order spline $L^2$-approximation.
 The set ${\cal R}(m,J)$ also defines the analogue set ${\cal R}(s,J)$ for modeling an $m$-dimensional function of $(x(j): j\in s)$ for any $s$ of size $m=\mid s\mid$. 
 We use these sets $({\cal R}(s,J):J)$ of $J$ knot-points across different subsets $s$ of the $d$-components of size $m$  to define the total set  ${\cal R}^0(d,{\bf J})=\cup_{s\subset\{1,\ldots,d\}}\{(s,u): u\in {\cal R}(s,{\bf J}(s)\}$ of zero order spline basis functions indexed by subset $s$ and knot-point $u$, with ${\cal R}(s,J)$ the set ${\cal R}(\mid s\mid,J)$ but applied to modeling a function with coordinates $(x(j): j\in s)$. Subsets ${\cal R}^{0,*}(d,{\bf J})\subset{\cal R}^0(d,{\bf J})$ define a class of additive submodels of interest. This set ${\cal R}^0(d,{\bf J})$  implies a  corresponding working model $D^{(0)}({\cal R}^0(d,{\bf J}))=\{\sum_{(s,u)}\beta(s,u)\phi_u^0: (s,u)\in {\cal R}^0(d,{\bf J})\}$,  where ${\bf J}=({\bf J}(s): s\subset\{1,\ldots,d\})$ is a vector of sizes, one for each subset $s$ of $d$ components representing the number of zero order spline to model this $s$-specific function in the zero order spline representation of $Q_0$. 
 
 In Section \ref{section4} this same family $({\cal R}(s,J):J,s)$ of knot-point sets in $J$ and subset $s$ will  define in the similar manner our total set ${\cal R}^k(d,{\bf J})=\cup_{\bar{s}(k+1)}\{(\bar{s}(k+1),u): u\in {\cal R}(s_{k+1},{\bf J}(\bar{s}(k+1))\}$  of $\leq k$-th order spline basis functions for the $k$-th order spline sieve MLE, $k=1,\ldots$. In this case, the basis functions are indexed by a nested sequence of subsets $\bar{s}(k+1)=(s_1,\ldots,s_{k+1})$, $s_{k+1}\subset s_k\ldots\subset s_1$,  and a knot-point $u$ in a corresponding index set ${\cal R}(s_{k+1},{\bf J}(\bar{s}(k+1))$ of size ${\bf J}(\bar{s}(k+1)$. 
 Given this family of zero order splines  ${\cal R}^{0}(d,{\bf J})$ and corresponding zero-order spline working models $D^{(0)}({\cal R}^0(d,{\bf J}))$, indexed by ${\bf J}$, we define the corresponding sieve zero-order spline-MLE, the HAL and relax HAL-MLE, its score equations,  and demonstrate it with the key examples. It will be clear that these examples have trivial analogues for the $k$-th order spline sieves. We also lay out our general (outline of)  proof for establishing that these estimators centered at the data dependent working-model best approximation of true target function are asymptotically normal at rate $(n/d_n)^{-1/2}$ with $d_n$ the number of non-zero coefficients in the fit. It follows that the rate of convergence w.r.t. target function is completely driven by the uniform approximation error of the sieve, by making sure the bias shrinks to zero faster than $(n/d_n)^{-1/2}$. At this point the reader already understands the asymptotic normality  and should be motivated for seeing the $k$-th order spline representations and their uniform approximation error as a function of $d_n$ and smoothness of  the true target function. 

 In Section \ref{section3} we provide a general $k$-th order spline representation of a $k$-th order differentiable (w.r.t. Lebesgue) multivariate real valued function, $k=1,\ldots,$ generalizing the zero-order spline representation. This defines a smoothness class $D^{(k)}([0,1]^d)$ and its representation as a large linear model in $\leq k$-th order splines. 
 The total set of basis functions for modeling $D^{(k)}([0,1]^d)$ is denoted with ${\cal R}^k(d)$. We also consider various submodels $D^{(k)}({\cal R}^{k,*}(d))\subset D^{(k)}({\cal R}^k(d))\equiv D^{(k)}([0,1]^d)$ of interest of this general $k$-th order spline model indexed by infinite collections ${\cal R}^{k,*}(d)\subset {\cal R}^k(d)$ of $\leq k$-th order splines, demonstrating its enormous flexibility in modeling the true target function. As with the zero order spline representation, it is clear one can think of the model  for a  target function $Q_0\in D^{(k)}([0,1]^d)$ as an additive model $\sum_j Q_j$ and submodels are obtained by restricting the set of functions in the additive model, and possibly by restricting the set of splines (i.e., knot-points) modeling a particular function $Q_j$ in this additive model. The index set for $j$ depends on and grows with $k$ since it is represented by the set of vectors $\bar{s}(k+1)=(s_1,\ldots,s_{k+1})$ of nested subsets $s_1\supset s_2\ldots\supset s_{k+1}$ of $\{1,\ldots,d\}$. So the zero order spline model is an additive model over functions indexed by $s_1\subset\{1,\ldots,d\}$, while the first order spline model is a larger additive model over functions indexed by nested pairs of subsets $(s_1,s_2)$ with $s_1\subset s_2$, and so on.  A large variety of additive submodels are obtained by restricting these choices $\bar{s}(k+1)$. 
 
In Section \ref{section4} we define the finite knot-point set ${\cal R}^k(d,{\bf J})=\cup_{\bar{s}(k+1)}\{ (\bar{s}(k+1),u): u\in {\cal R}(s_{k+1},{\bf J}(\bar{s}(k+1))\} \subset {\cal R}^k(d)$ of size $d^k(J)=O(J)$ and corresponding finite dimensional linear $k$-th order spline working model
 $D^{(k)}({\cal R}^k(d,{\bf J}))\subset D^{(k)}([0,1]^d)$, in terms of $({\cal  R}(s,J):s,J)$ satisfying  the $O(r(\mid s\mid,J))$-$L^2$-approximation property. $J$ represents the average number ${\bf J}(\bar{s}(k+1))$  of basis functions for modeling each $\bar{s}(k+1)$-specific function in the additive model.
 
Subsequently, in Section \ref{section5} we express the uniform approximation error  of this ${\cal R}^k(d,{\bf J})$-specific $k$-th  order spline working model $D^{(k)}({\cal R}^k(d,{\bf J}))$, w.r.t.  a $k$-th order smooth target function in $D^{(k)}([0,1]^d)$, in terms of  the approximation error of a linear combination of ${\bf J}(\bar{s}(k+1))$ $k$-th order splines for approximating a   $k$-th order primitive function in coordinates $(x(j):j\in s_{k+1})$. This result teaches us that it remains to understand how well a $J$-dimensional $k$-th order spline model with knot-point set ${\cal R}(d,J)$ uniformly approximates a $d$-variate $k$-th order primitive function. In Appendix  \ref{AppendixI} we establish that the latter uniform approximation error is $O(r(d,J)^{k+1})$, where as mentioned  $r(d,J)\approx J^{-1}$ is  known rate for approximating a cumulative distribution function with a linear combination of $J$ zero-order splines w.r.t. $L^2$-norm. This result proves that there exists an element in our ${\cal R}^k(d,{\bf J})$-specific working model $D^{(k)}({\cal R}^k(d,{\bf J}))$ of linear combinations $Q_{{\bf J},\beta}=\sum_{j\in {\cal R}^k(d,{\bf J})}\beta(j)\phi_j$ of  $\leq k$-th order splines  $\phi_j$ that provides a supremum norm approximation of the $k$-th order smooth true function $Q_0$ with approximation error $O(r(d,J)^{k+1})$, where the size of ${\cal R}(d,{\bf J})$ is proportional to $J$. In Section \ref{section6}  we prove that this result also implies that the minimizer  $Q_{0,{\bf J}}=\arg\min_{Q\in D^{(k)}({\cal R}^k(d,{\bf J}))} P_0L(Q)$ of the true risk over this  ${\cal R}^k(d,{\bf J})$-specific working model $D^{(k)}({\cal R}^k(d,{\bf J}))$ approximates $Q_0$ at the same rate $O(r(d,J)^{k+1})$ w.r.t. supremum norm. This finalizes the bias analysis for our working models $D^{(k)}({\cal R}^d(d,{\bf J}))$ w.r.t $Q_0$ in the sense that we now understand $\pl Q_{0,{\bf J}}-Q_0\pl_{\infty}=O(r(d,J)^{k+1})$, assuming ${\bf J}(\bar{s}(k+1))$ is proportional to $J$ across all $\bar{s}(k+1)$ in our assumed index set. 


In Section \ref{section7} we define the $k$-th order spline-MLE $Q_{n,{\bf J}}=\arg\min_{Q\in D^{(k)}({\cal R}^k(d,{\bf J}))}P_n L(Q)$ over $D^{(k)}({\cal R}^k(d,{\bf J}))$ indexed by knot-point set ${\cal R}^k(d,{\bf J})$ and $k$-th order spline HAL-MLE $Q_{n,C}=\arg\min_{Q_{{\bf J},\beta}\in D^{(k)}({\cal R}^k(d,{\bf J}_{max})),\pl \beta\pl_1\leq C}P_nL(Q)$ indexed by an $L_1$-norm $C$, and an initial  knot-point set ${\cal R}^k(d,{\bf J}_{max,n})$. 
We also define the $k$-th order spline relax HAL-MLE that refits the model selected by HAL-MLE with an unpenalized MLE. Each of these estimators, when selecting the tuning parameters with a data adaptive selector, end up with a set of $\leq k$-th order spline basis functions with non-zero coefficients corresponding with an index set ${\cal R}_n$, such as ${\cal R}^k(d,{\bf J}_n)$ for the sieve MLE. Moreover, in this section \ref{section7}  we establish a rate of convergence in loss-based dissimilarity for these estimators as estimators of $Q_0$, which shows convergence in $L^2$-norm at the remarkably fast rate $n^{-k^*/(2k^*+1)}$ up till $\log n$-factors, where $k^*=k+1$.

In Sections \ref{section8} and \ref{section10} we establish pointwise asymptotic normality of these three estimators of $Q_0$, at this same rate for logistic and least squares linear regression and general loss functions, respectively. Each of these MLEs solve the score equations of the linear working model $D^{(k)}({\cal R}_n)$ or $D^{(k)}_{C_n}({\cal R}_n)$ (with $L_1$-norm constraint $C_n$) implied by the set ${\cal R}_n$ of basis functions with non-zero coefficients, and thus also the linear span of these score equations. Our general proof outlined in Section \ref{section2}  establishes asymptotic normality of an MLE for a fixed user supplied working model $D^{(k)}({\cal R}_{0,n})$ as an estimator of its loss-based projection $Q_{0,{\cal R}_{0,n}}=\arg\min_{Q\in D^{(k)}({\cal R}_{0,n})}P_0L(Q)$, where (say) ${\cal R}_{0,n}={\cal R}^k(d,{\bf J}_{0,n})$ for fixed sequence ${\bf J}_{0,n}$.  This immediately yields the asymptotic normality for  the sieve MLE $Q_{n,{\bf J}_{0,n}}$ for a fixed non-informative selector ${\bf J}_{0,n}$, and combined with our uniform approximation result, it shows that the optimal rate w.r.t $Q_0$ is achieved by setting ${\bf J}_{0,n}=C n^{1/(2k^*+1)}$ (for any $C$) up till $\log n$-factors across its components. and thereby proves  asymptotic normality of $Q_{n,{\bf J}_{0,n}}$ w.r.t. target function $Q_0$ at rate $n^{-k^*/(2k^*+1)}$ up till $\log n$-factors. However, even though this is of theoretical interest, and such a choice ${\bf J}_{0,n}$ can play the role of ${\bf J}_{max,n}$ in definition of HAL-MLEs, to make this sieve MLE practical the tuning parameters ${\bf J}$ need to be chosen data adaptively, such as with the cross-validation selector. 

In order to deal with the challenge of an informatively selected working model $D^{(k)}({\cal R}_n)$ (especially for the HAL-MLEs), we define ${\cal R}_{0,n}$ as the algorithm ${\cal R}(P_n)$ but applied to the empirical measure $P_n^{\#}$ of an  independent sample from $P_0$. An option is to shrink ${\cal R}_{0,n}$ relative to ${\cal R}(P_n^{\#})$ but our results appear to show that this is generally not needed. 
We then assume that the algorithm ${\cal R}(P_n)$ is chosen so that  these MLEs $Q_n$  also solve the  score equations  of the corresponding fixed working model $D^{(k)}({\cal R}_{0,n})$, up till small enough error, where the projection $Q_{0,{\cal R}_{0,n}}$ onto the working model satisfies $\pl Q_{0,{\cal R}_{0,n}}-Q_0\pl_{\infty}=O(r(d,J_{0,n})^{k+1})$. 
As a consequence of the estimators acting as an MLE for a fixed working model $D^{(k)}({\cal R}_{0,n})$, our asymptotic normality proof applies. 

In Sections \ref{section8} 
 we consider the weighted squared error loss function for mean outcome regression for a continuous outcome and the log-likelihood for the conditional distribution of a binary outcome using the logit-link. 
In this section we  prove  pointwise convergence in distribution of the standardized (w.r.t. its projection $Q_{0,{\cal R}_{0,n}}$) $k$-th-order spline sieve and HAL-MLE at rate $O((d^k(J_{0,n}))/n)^{-1/2})$  to a normal mean zero limit distribution.  Combined with our bias result $Q_{0,{\cal R}_{0,n}}-Q_0=O(r(d,J_{0,n})^{k+1})$, this then also provides us asymptotic normality of the standardized (w.r.t. $Q_0(x)$)  sieve and HAL-MLE when undersmoothing   relative to the cross-validation selector,  and corresponding pointwise confidence intervals for $Q_0(x)$, at the remarkable rate $n^{-k^*/(2k^*+1)}$ up till a power of $\log n$-factor with $k^*=k+1$.  In addition, we also establish a simultaneous   confidence interval for $Q_0$.  

In Section \ref{section10}  we generalize these results for the regression case to  general loss functions $L(Q)$ with $Q_0=\arg\min_fP_0L(Q)$.

In Sections \ref{section9} and \ref{section11} we utilize our  empirical mean approximation of $(Q_n-Q_{0,{\cal R}_n})(x)$, generalized to actual set ${\cal R}_n$ giving nicer second order remainder,  with known influence curve to 
establish asymptotic linearity of the plug-in estimator $\Phi(Q_n)$ of a pathwise differentiable parameter $\Phi(Q_0)$. We use this empirical mean expansion, combined with the delta-method, to analyze $\Phi(Q_n)-\Phi(Q_{0,{\cal R}_n})$, while showing that the bias $\Phi(Q_{0,{\cal R}_n})-\Phi(Q_0)$ is second order. Specifically, our results show that if ${\cal R}_n$ is adaptive and thereby converges to a set ${\cal R}_0$ (depending on true $Q_0$)  rich enough so that $Q_0\in D^{(k)}({\cal R}_0)$, then the resulting estimator will still be asymptoticallly linear but with a super-efficient influence curve.  By undersmoothing one can arrange that the estimators $\Phi(Q_n)$ are asymptotically efficient and regular as well, but even without this our asymptotic normality provides asymptotically valid confidence intervals.  

 We conclude with a discussion in Section \ref{section12}. 
 In Section \ref{section13} before the start of the Appendix we provide a notation index. 
Our Appendix contains technical proofs of our results, and various extra results of interest.
Overall main ideas and proofs are presented in the main part of paper so that it reads as a self-contained journey towards the final results.

 \section{The zero order spline sieve and HAL-MLEs; key examples; natural generalization to $k$-th order sieve; and general asymptotic normality proof}\label{section2}
 We observe $n$ i.i.d. copies $O_1,\ldots,O_n$ with common probability distribution $P_0$.
 We are interested in estimating a functional target parameter $Q:{\cal M}\rightarrow D^{(0)}_{C}([0,1]^d)$ with parameter space being contained in the set $D^{(0)}_C([0,1]^d)$ of $d$-variate real valued cadlag functions with a universal bound $C$ on its sectional variation norm. Let $D^{(0)}([0,1]^d)$ be the space of cadlag functions only assuming that the sectional variation norm is bounded.  
 We let $L:D^{(0)}([0,1]^d)\rightarrow L^2(P_0)$ be a loss function so that $Q(P)=\arg\min_{Q\in D^{(0)}([0,1]^d)} P L(Q)$. 
 
  From the literature on HAL \citep{vanderLaan17} we have that any function $Q\in D^{(0)}([0,1]^d)$ can be represented as $Q(x)=\int_{[0,x]} dQ(u)=\int_{u\in [0,1]^d}\phi_u(x)dQ(u)$ with $\phi_u(x)\equiv I(u\leq x)=\prod_{j=1}^d I(u_j\leq x_j)$  being a tensor product of zero-order splines. Here $[0,1]^d=\cup_{s\subset\{1,\ldots,d\}}E_s$ is a union of disjoint lower-dimensional $0$-edges of $[0,1]^d$ defined by  $E_s=\{(x(s),0(-s)):x\in [0,1]^d\}$, beyond the full  inner set $(0,1]^d$ corresponding with $s=\{1,\ldots,d\}$. 
 As a convention we include the empty subset $s$ in which case $E_s=\{0\}$ is the singleton $0$. In addition, in this integral representation we have that $dQ(u)=\sum_s I(u\in E_s) dQ(0(-s),du(s))$ is a sum measure with disjoint complementary support: it represents the measure  $dQ(u)$   for $u$ in the inner cube $E_{\{1,\ldots,d\}}=(0,1]^d$, while for the real subsets $s\subset\{1,\ldots,d\}$ it represents the measure $Q(0(-s),du(s))$ implied by its section $Q_s: x(s)\rightarrow Q(x(s),0(-s))$ on the lower dimensional zero-edge $E_s$, and, finally, for the empty subset $s=\emptyset$, $dQ(\{0\})=Q(0)$ is a pointmass on the singleton $\{0\}$.   That is, the integral over the zero-edges $E_s$ is  defined by setting the coordinates outside the subset $s\subset\{1,\ldots,d\}$ equal to zero and defining the integral w.r.t  $dQ_s(u(s))=Q(du(s),0(-s))$. Therefore, our short-hand notation truly represents the following representation: \[
 \int_{[0,x]}dQ(u)\equiv Q(0)+\sum_{s\subset\{1,\ldots,d\}}\int_{(0(s),x(s)]}Q(du(s),0(-s)).\] The  sectional variation norm of $Q$ is defined by the sum over $s$ of the variation norms of the measures generated by the sections $Q_s$: \[
 \pl Q\pl_v^*\equiv \int_{[0,1]^d} \mid dQ(u)\mid =
 \mid Q(0)\mid+\sum_{s\subset\{1,\ldots,d\}}\int_{(0(s),1(s)]}\mid Q(du(s),0(-s))\mid.
 \]
 Let ${\cal S}^0(d)$ be the set of all subsets $s\subset\{1,\ldots,d\}$, including the empty subset. 
 For each subset $s$, let ${\cal R}^0(s)\equiv (0(s),1(s)]$. 
 Let ${\cal R}^0(d)=\{(s,u): s\in {\cal S}^0(d), u\in {\cal R}^0(s)=(0(s),1(s)]\}$, where for the empty subset $s$, only $u=0$ is a knot-point. This represents the set of all basis functions needed to span $D^{(0)}([0,1]^d)$, so that we can represent $D^{(0)}([0,1]^d)$ as
 \begin{equation}\label{Rd0}
 D^{(0)}({\cal R}^0(d))\equiv \left\{\sum_{s\in {\cal S}^0(d)}\int_{u\in {\cal R}^0(s)}\phi_u^0(\cdot)dQ(u): Q\in D^{(0)}([0,1]^d\right\}.
 \end{equation}
 This representation now suggests a general class of submodels of $D^{(0)}([0,1]^d)$ defined by 
  a subset ${\cal R}^{0,*}(d)=\{(s,u): s\in {\cal S}^{0,*}(d)\subset {\cal S}^{0,*}(d),u\in {\cal R}^{0,*}(s)\subset(0(s),1(s)]\}$ of ${\cal R}^0(d)$. Then, we can define the submodel
 \[
 D^{(0)}({\cal R}^{0,*}(d))\equiv \left\{\sum_{s\in {\cal S}^{0,*}(d)}\int_{u\in {\cal R}^{0,*}(s)}\phi_u^0 dQ(u): Q\in D^{(0)}([0,1]^d)\right\}.\]
 So this is the closure of the linear span of all zero order spline basis functions $\phi_{s,u}^0\equiv I(u\in E_s)\phi_u^0$ indexed by $s\in {\cal S}^{0,*}(d)$ and $u\in {\cal R}^0(s)$. This representation of these subspaces indexed by a joint $(s,u)$ instead of just $u\in [0,1]^d$ is now easily generalized to the higher order spaces $D^{(k)}([0,1]^d)$ as we will see. 
 We assume that the parameter space $Q({\cal M})=\{Q(P):P\in {\cal M}\}$ implied by the statistical model equals 
  \begin{equation}\label{D0}
  D^{(0)}_{C^u}({\cal R}^{0,*}(d))=\left\{\sum_{s\in {\cal S}^{0,*}(d)}\int_{u\in {\cal R}^{0,*}(s)\subset (0(s),1(s)]}\phi_u^0(\cdot) dQ(u):Q\in D([0,1]^d), \pl Q\pl_v^*<C^u\right\}.\end{equation}
Note that this parameter space corresponds with a general additive model $\sum_{s\subset{\cal R}^{0,*}(d)} Q_s$ with each $Q_s$ modelled with the linear span of $\phi_{u,s}$ with $u\in {\cal R}^0(s)$. 
  Our $k$-th order spline representation theorems generalize to $k$-th order smooth functions in $D^{(0)}({\cal R}^{0,*}(d))$, which just induces a restriction on the resulting set of $k$-th order spline basis functions. 
  Our $k$-th order spline representation for a $k$-th order smooth function in $D^{(0)}({\cal R}^{0,*}(d))$ is then defined by restricting $s_1$ in the vector  $\bar{s}(k+1)=(s_1,\ldots,s_{k+1})$ of nested subsets $s_{k+1}\subset s_k\ldots\subset s_1$  to be an element of ${\cal S}^{0,*}(d)$, and, by definition all other subsets in $\bar{s}(k+1)$ are then subsets of these $s_1$ and thereby automatically restricted as well. 
For notational convenience, we will focus our presentation on the nonparametric case  ${\cal R}^0(d)$, but we also showcase the many submodels one can construct by subsetting the set of potential  basis functions, including generating $k$-th order smooth submodels of $D^{(0)}({\cal R}^{0,*}(d))$.

Throughout this article we assume that the loss function is uniformly bounded  and has a quadratic loss-based dissimilarity in the sense that 
\begin{equation}\label{boundloss}
\sup_{Q \in D^{(0)}({\cal R}^0(d))}\pl L(Q)\pl_{\infty}<\infty \mbox{ and }\sup_{Q,Q_0\in D^{(0)}({\cal R}^0(d)) }\frac{P_0(L(Q)-L(Q_0))^2}{d_0(Q,Q_0)}<\infty .
\end{equation}
 This represents a standard assumption for establishing the rate of convergence of the MLE w.r.t. loss-based dissimilarity and establishing the asymptotic equivalence of the cross-validation selector with the oracle selector \citep{vanderLaan&Dudoit03,vanderVaart&Dudoit&vanderLaan06,vanderLaan&Dudoit&vanderVaart06}.

 \subsection{Construction of sequence of  sieve zero-order spline-MLE of increasing complexity}
 In Corollary \ref{defknots}  we define ${\cal R}(m,J)\subset (0,1]^m$ as any (non-informative) knot-point set of size $J$ whose linear combination of corresponding zero order splines  provides an $L^2(\mu)$-approximation of size $O(C(M)r(m,J))$ uniformly in the set ${\cal F}^0_M((0,1]^d)$ of $m$-variate cumulative distribution functions on $[0,1]^m$ bounded from above by $M$, a result derived from the current literature in Appendix \ref{AppendixA}.
 \begin{definition}
Let 
\begin{equation}\label{defrdJMm}
r(d,J,M)\equiv \arg\min_{u_1,\ldots,u_J\in (0,1]^d}\sup_{Q\in {\cal F}^{(0)}_M((0,1]^d)}\inf_{\beta}
\pl \sum_{j=1}^J\beta(j)\phi_{u_j}^0-Q\pl_{L^2(\mu)}.
\end{equation}
In Appendix \ref{AppendixA} we show $r(d,J,M)=O( C(M)r(d,J))$ with $r(d,J)=J^{-1}(\log J)^{2(d-1)}$ and $C(M)=M$. 
\end{definition}
 
 \begin{definition}
 Recall the definition of  $r(d,J,M)$ above and the bound $r(d,J,M)\sim C(M)r(d,J)$ with $r(d,J)=J^{-1}(\log J)^{2(d-1)}$ and $C(M)=M$.
In the remaining of this article ${\cal R}(d,J)$ represents a set of $J$ knot-points that satisfies \begin{equation}\label{Rbm}
\sup_{Q\in {\cal F}^0_M((0,1]^d)}\inf_{\alpha}\pl \sum_{v\in {\cal R}(d,J)}\alpha(v)\phi_v^0-Q\pl_{\mu}=O(C(M) r(d,J)).
\end{equation}
\end{definition}
 For each subset $s\subset\{1,\ldots,d\}$ we can then define a set ${\cal R}^0(s,{\bf J}(s))\subset {\cal R}^0(s)\}$ of ${\bf J}(s)$ knot-points, and define the total index set as 
 \[
 {\cal R}^0(d,{\bf J})\equiv \{(s,u): s\subset{\cal S}^0(s), u\in {\cal R}^0(s,{\bf J}(s))\},\] where ${\bf J}$ denotes the vector of all  sizes ${\bf J}(s)$ across the subsets $s$.  For $s$ the empty set, we define ${\cal R}(s,{\bf J}(s))=\{0\}$ and ${\bf J}(s)=1$.  
Similarly, we can define ${\cal R}^{0,*}(d,{\bf J})\subset {\cal R}^{0,*}(d)$. Generally speaking, we would make sure that, for all $s\in {\cal S}^{0,*}(d)$ allowed by our model $D^{(0)}({\cal R}^{0,*}(d))$, for some $\delta>0$ and $M<\infty$ $\delta J<{\bf J}(s)<M J $. In that way  all sizes behave as $\sim J$ for an integer $J$, so that each $s$-specific section in the additive model is approximated at rate $r(\mid s\mid, J)$. Let $d^0(J)=\sum_{s}{\bf J}(s)$ denote the size of ${\cal R}^0(d,{\bf J})$.
According to our representation (\ref{D0}) of $D^{(0)}({\cal R}^{0}(d))$, the finite dimensional working model $D^{(0)}({\cal R}^{0}(d,{\bf J}))$ provides an $O(r(d,J))$-$L^2$-approximation of $D^{(0)}({\cal R}^0(d))=D^{(0)}([0,1]^d)$. 

 Let $\beta_{n,{\bf J}}^0=\arg\min_{\beta\in \openr^{d^0(J)} } P_n L(\sum_{j\in {\cal R}^0(d,{\bf J})}  \beta(u)\phi_j^0)$ be the MLE for this zero-order spline working model $D^{(0)}({\cal R}^0(d,{\bf J}))=\{\sum_{(s,u)\in {\cal R}^0(d,{\bf J})} \beta(s,u)\phi_{s,u}^0:\beta\in \openr^{d^0(J)} \}\subset D^{(0)}({\cal R}^0(d))$. 
 To emphasize $\beta_{n,{\bf J}}^0$ as an estimator we  would write $\beta_{n,{\bf J}}^0=\hat{\beta}_{\bf J}^0(P_n)$.

 
 
 Given the specification of ${\cal R}^0(d,{\bf J})$, and resulting MLE $Q^0_{{\bf J},n}=\hat{f}^0_{\bf J}(P_n)=\sum_{j\in {\cal R}^0(d,{\bf J})}\hat{\beta}_{\bf J}^0(P_n)(j)\phi_j^0$, we may select ${\bf J}$ with the cross-validation selector ${\bf J}_{n,cv}=\arg\min_{{\bf J}}1/V\sum_{v=1}^V P_{n,v}^1 L(\hat{f}_{{\bf J}}^0(P_{n,v}))$ using $V$-fold cross-validation, where $V$ is the number of sample splits; $P_{n,v}^1$ is the empirical measure of the validation sample and $P_{n,v}$ is the empirical measure of the training sample, across the sample splits $v=1,\ldots,V$.
  The  sieve zero-order spline estimator, using cross-validation to select the knot-point set, is then $Q_n^0=\hat{f}^0(P_n)$ defined by $\hat{f}^0(P_n)=\hat{f}^0_{{\bf J}_{n,cv} }(P_n)$. 
  
  One may  also consider the case that we use a data adaptive selector ${\bf J}_n$ whose components are slightly larger than  ${\bf J}_{n,cv}$. We will refer to this as undersmoothing which would be necessary to guarantee that the plug-in estimator of smooth features of $Q_n$ will be asymptotically efficient \citep{vanderLaan&Benkeser&Cai19}. 

  For notational convenience, let ${\cal R}_n^0$ play the role of ${\cal R}^0(d,{\bf J}_n)$.
Let $D^{(0)}({\cal R}_n^0)=\{\sum_{j\in {\cal R}_n^0}\beta(j)\phi_j^0:\beta\in\openr^{d^0(J_n)}\}$ be the finite dimensional zero-order spline working model identified by ${\cal R}_n^0$. Even though the sieve MLE is defined as above, it is important to note that this sieve MLE is  also an MLE for the data adaptive working model $D^{(0)}({\cal R}_n^0)$, since that determines the set of score equations the sieve MLE  solves. 
 Specifically,  $Q_n^0=\arg\min_{Q\in D^{(0)}({\cal R}_n^0)}P_n L(Q)$ is the MLE over $D^{(0)}({\cal R}_n^0)$,  while $Q_{0,n}^0=\arg\min_{Q\in D^{(0)}({\cal R}_n^0) }P_0 L(Q)$ represents the $P_0$-MLE and thereby best approximation of the true target function $Q_0$ in the working model $D^{(0)}({\cal R}_n^0)$.   Let $\beta_n^0=\beta^0_{J_n}(P_n)$ the corresponding vector of coefficients so that $Q_n^0=\sum_{j\in {\cal R}_n^0}\beta_n^0(j)\phi_j^0$.  

\subsection{Zero-order spline HAL-MLE and relax HAL-MLE}
Let ${\bf J}_{max}$ be an initial specification, chosen rich enough, typically with ${\bf J}_{max}(s)$ is constant across all subsets $s$. 
A knot-point set ${\cal R}^0(d,{\bf J}_{max})$ defines a working model $D^{(0)}({\cal R}^{0}(d,{\bf J}_{max}))=\{\sum_{j\in {\cal R}^0(d,{\bf J}_{max})}\beta(j)\phi_j^0:\beta\}$ and the $L_1$-restricted working model $D^{(0)}_C({\cal R}^0(d,{\bf J}_{max}))=\{\sum_{j\in {\cal R}^0(d,{\bf J}_{max})}\beta(j)\phi_j^0:\beta,\pl\beta\pl_1<C\}$. 
A corresponding zero-order spline HAL-MLE could be defined by defining $Q_{n,{\bf J}_{max},C}^0=\arg\min_{Q\in D^{(0)}_C({\cal R}^{0}(d,{\bf J}_{max}))}P_n L(Q)$, 
and, one could then select $({\bf J}_{max},C)$ with a cross-validation selector or undersmoothed version of that, or, one could keep ${\bf J}_{max}$ fixed. We can also define the relax zero-order spline HAL-MLE $Q_{n,{\bf J}_{max},C,r}^0$ by first computing $Q_{n,{\bf J}_{max},C}^0$, identify the set ${\cal R}_n^0({\bf J}_{max},C)$ of basis functions with non-zero coefficients, and thereby its corresponding working model $D^{(0)}_C({\cal R}_n^0(d,{\bf J}_{max}))$, and then refit this model with  the unpenalized MLE. Again, we can then select $({\bf J}_{max},C)$ with the cross-validation selector w.r.t. the $({\bf J}_{\max},C)$-specific relax lasso $Q_{n,{\bf J},C,r}^0$,  or an undersmoother such as the cross-validation selector for the HAL-MLE.

 
 Similarly, one notes that the HAL-MLE and relax HAL-MLE are $L_1$-regularized MLEs or unpenalized MLE for the working model $D^{(0)}({\cal R}_n^0)$, where  ${\cal R}_n^0$ be the set of knot-points with non-zero coefficients in $Q_n^0=Q^0_{n,{\bf J}_{max,n},C_n}$ and $Q_{n,r}^0=Q^0_{n,{\bf J}_{max,n},C_n,r}$, respectively.

 \subsection{Generalization to $k$-th order spline sieve and HAL-MLEs}
 In the next section we will define a $k$-th order smoothness class $D^{(k)}_M([0,1]
 ^d)\subset D^{(0)}([0,1]^d)$, $k=0,1,\ldots$, consisting of functions for which $k$-th order derivatives are cadlag and have a  combined sectional variation norm bounded by $M$.
 Our representation theorem for functions in this class shows that any submodel of this linear space  is indexed by an index set ${\cal R}^k(d)$ indicating that these are functions in the closure of the linear span of $\leq $-kth order spline basis functions identified by ${\cal R}^k(d)$, so that this model will also be denoted with $D^{(k)}({\cal R}^k(d))$.
 More precisely, these submodels of $D^{(k)}([0,1]^d)$ are additive models $Q=\sum_{\bar{s}(k+1)}Q_{\bar{s}(k+1)}$ over functions $Q_{\bar{s}(k+1)}$ indexed by vectors of nested subsets $\bar{s}(k+1)$, where each of these functions is modelled by its own set of basis functions $\phi_{\bar{s}(k+1),u}$, $u\in {\cal R}^0(s_{k+1})=(0(s_{k+1}),1(s_{k+1})]$. We will provide finite dimensional submodels $D^{(k)}_M({\cal R}^k(d,{\bf J})$ indexed by a vector ${\bf J}$ of sizes, and dimension $d(J)=\sum_{\bar{s}(k+1)}{\bf J}(\bar{s}(k+1))$ of order $J$ (representing average size), which describe linear combinations of $d(J)$ $\leq k$-th order spline basis functions whose vector of coefficients have an $L_1$-norm bounded by $M$. It will be shown that these finite dimensional submodels uniformly  approximate $D^{(k)}({\cal R}^k(d))$ at rate $r(d,J)^{k+1}$, which is $(1/J)^{k+1}$ up till $\log J$-factors. Analogue to above, we can define an MLE over $D^{(k)}({\cal R}(d,{\bf J}))$ for each ${\bf J}$, and define a $L_1$-constrained MLE over $D^{(k)}({\cal R}(d,{\bf J}_{max}))$, and define the corresponding sieve-MLE involving data adaptive selection of ${\bf J}$, and HAL-MLE involving data adaptive selection of $C$ (and possibly ${\bf J}_max$, and corresponding relax HAL-MLE. 
Similarly as for $D^{(0)}([0,1]^d)$,  one notes that the HAL-MLE and relax HAL-MLE are $L_1$-regularized MLEs or unpenalized MLE for the working model $D^{(k)}({\cal R}_n^k)$, where  ${\cal R}_n^k$ indicates the set of non-zero coefficients. 
We could denote the sieve MLE with $Q_{n,{\bf J}_n}^k$, HAL-MLE with  $Q^k_{n,{\bf J}_{max,n},C_n}$ and relax HAL-MLE with $Q^k_{n,{\bf J}_{max,n},C_n,r}$.
Any of these estimators are of form $Q_n^k=\sum_{j\in {\cal R}_n^k}\beta_n^k(j)\phi_j^k$ with $\phi_j^k$ representing tensor products of $\leq k$-th order spline basis functions indexed by $j=(\bar{s}(k+1),u)$. 
In this manner our examples  in the next subsections have immediate analogues to $Q_n^k$ for any $k=0,\ldots$, while below we describe them for the zero order spline estimator $Q_n^0$.  
 
 {\bf Extending $k$-th order spline sieve and HAL-MLEs to handle binary covariates:}
 The zero-order sieve and HAL-MLEs naturally handle binary covariates. A binary covariate $X_j\in \{0,1\}$ only generates zero-order splines $I(X_j\geq 1)$ with a single knot-point. Similarly, the $k$-th order spline sieve and HAL-MLEs are easily extended to handle binary covariates. For that purpose we should assume $Q_0\in D^{(k)}(\{0,1\}^m\times [0,1]^d)$, where $m$ represents the number of binary components, where this $k$-th order smoothness class simply assumes that $Q_0(b,\cdot)\in D^{(k)}([0,1]^d)$ for all $b\in \{0,1\}^m$. Our $k$-th order spline representation of a $Q\in D^{(k)}([0,1]^d)$ implies the $k$-th order spline representation for $Q\in D^{(k)}(\{0,1\}^m\times [0,1]^d)$ by applying it to each $Q(b,\cdot)\in D^{(k)}([0,1]^d)$. From a practical perspective, it just means that we multiply all  the basis functions that model $D^{(k)}([0,1]^d)$ with $I(b\geq u)=\prod_{j=1,u_j=1}I(b(j)\geq 1)$ for $u\in \{0,1\}^m$.

\subsection{Score equations for sieve zero-order spline-MLE, oracle MLE, and HAL analogues}
Since  $Q_n^0$ is an MLE over $D^{(0)}({\cal R}_n^0)$ it solves the following score equations: \[
 0=P_n \frac{d}{d\beta_n^0}L\left(\sum_{j\in {\cal R}_n^0}\beta_n^0(j)\phi_j^0\right)=
 P_n \frac{d}{dQ_n^0}L(Q_n^0)(\phi^0),\]
 where $\phi^0=(\phi_j^0:j\in {\cal R}_n^0)$. Sometimes, we will also denote this vector of basis functions with $\phi_n^0$ to acknowledge its dependence on ${\cal R}_n^0$.
 
The $P_0$-MLE $Q_{0,n}^0=\arg\min_{Q\in D^{(0)}({\cal R}_n^0) }P_0 L(Q)$ solves the same score equations under $P_0$: 
 \[
 0= P_0 \frac{d}{d\beta_{0,n}^0}L\left(\sum_{j\in {\cal R}_n^0}\beta_{0,n}^0(j)\phi_j^0\right)=
 P_0\frac{d}{dQ_{0,n}^0}L(Q_{0,n}^0)(\phi_n^0).\]
  Let $S_{Q}\equiv \frac{d}{dQ}L(Q)(\phi_n^0)$ so that we have $P_n S_{Q_n^0}=P_0 S_{Q_{0,n}^0}=0$. Later we will also use notation $S_{Q}^0(\phi)=\frac{d}{dQ}L(Q)(\phi)$.
  Alternatively, one can represent $S_{Q}^0$ as function in $\beta$ and denote it with $S_{\beta}^0$ in which case we can write
  $P_n S_{\beta_n}=P_0 S_{\beta_{0,n}}=0$. 
  Similarly, the relax HAL-MLE solves these score equations exactly as well. 
  
  {\bf HAL-MLE solves score equations at desired approximation:}
  As shown in Appendix \ref{AppendixJ}, the HAL-MLE $Q_n^0=Q^{0}_{n,{\bf J}_n,C_n}$ solves these  same score equations $P_n S_{Q_n^0}\approx 0$ up till a specified rate that depends in a known manner on  the amount of undersmoothing $C_n$ used. That is, the score equations solved by HAL-MLE along paths $(1+\delta h(j))\beta_n^0(j)$ with $h\perp \mid \beta_n\mid$ that keep the $L_1$-norm constant for small $\delta$ approximate these score equations.   
  As shown in \citep{vanderLaan&Benkeser&Cai19}, it also follows that, even without undersmoothing, $P_n S_{Q_n^0}(\phi_u^0) =o_P(n^{-1/2})$ under a very weak undersmoothing condition,  due to 
  \[
  P_n S_{Q_n^0}(\phi_j^0)=(P_n-P_0)S_{Q_n^0}(\phi_j^0-\tilde{\phi}_j^0)+P_0 (S_{Q_n^0}-S_{Q_{0,n}^0})(\phi_j^0-\tilde{\phi}_j^0),\]
  where $\tilde{\phi}_j^0$ is a linear combination of $\{\phi_j^0:j\in {\cal R}_n^0\}$ satisfying that $P_n S_{Q_n^0}(\tilde{\phi}_j^0)=0$ and chosen to approximate $\phi_j$, where the latter approximation will shrink with $C_n$. This identity uses that $P_0 S_{Q_{0,n}^0}(\phi_j^0)=P_0 S_{Q_{0,n}^0}(\tilde{\phi}_j^0)=0$. This shows that $P_n S_{Q_n^0}(\phi_j^0)$ is a second order remainder in $d_0(Q_n^0,Q_{0,n}^0)^{1/2}$ and an $L^2$-norm $\pl \phi_j^0-\tilde{\phi}_j^0\pl$.  The first term can be bounded by empirical process theory (e.g, $O_P(n^{-2/3}(\log n)^d)$) and the second term can generally be bounded with Cauchy-Schwarz inequality. In essence, one wants that $\phi_j^0-\tilde{\phi}_j^0$ converges to zero at rate faster than $n^{-1/6}$. 
  Therefore, we can conclude that, possibly using minor undersmoothing,  the $L_1$-regularization does not harm solving the desired score equations at the desired precision $o_P((n/d_n)^{-1/2})$ in our asymptotic normality proofs. 
  
 \subsection{An overview of our asymptotic normality proof for sieve and HAL-MLEs.}
 In the second part of this article we analyze $(Q_n^k-Q_{0,n}^k)(x)$ with the modification that $Q_{0,n}^k$ is the oracle MLE over a working model $D^{(k)}({\cal R}_{0,n}^k)$ with ${\cal R}_{0,n}^k$ being the same algorithm that maps $P_n$ into ${\cal R}_n^k$ (in order to obtain a clean CLT not having to deal with the fact that the oracle MLE $Q_{0,n}$ still depends on the data through ${\cal R}_n$), but applied to an independent i.i.d. sample $P_n^{\#}$ from $P_0$. For now, we ignore this modification, in order to communicate the  main idea of the analysis, and how it is  particularly powerful for likelihood behaving loss functions. This analysis of $(Q_n^k-Q_{0,n}^k)(x)$ combined with establishing $\pl Q_{0,n}^k-Q_0\pl_{\infty}=O(r(d,d_n)^{k+1})$,  as established in subsequent sections, this provides an overall analysis of $(Q_n^k-Q_0)(x)$ at a rate defined by trading off
 $(n/d_n)^{-1/2}$ with $r(d,d_n)^{k+1}$, where $d_n$ is the size of ${\cal R}_n^k$, which is proportional to $J_n$ (representing an average size across the $\bar{s}(k+1)$-specific sets of $k$-th order spline basis functions).
 
For simplicity, we assume that $Q_{0,n}^k$ is fixed/non -informative (e.g., consider the MLE over $D^{(k)}({\cal R}(d,{\bf J}_{0,n}))$ for a non-informative selector ${\bf J}_{0,n}$), and that $P_n S_{Q_n^k}(\phi^k_j)=0$, $j\in {\cal R}_n^k$, exactly as for the relax HAL-MLE and sieve MLE.  The proof below for asymptotic normality of $(n/d_n)^{1/2}(Q_n^k-Q_{0,n}^k)(x)$ applies equally to $k=0$.
Firstly, we note that for each $j\in {\cal R}_n^k$
 \[
 P_0\{ S_{Q_n^k}(\phi_j^k)-S_{Q_{0,n}^k}(\phi_j^k)\}=-(P_n-P_0) S_{Q_n^k}(\phi_j^k).\]
 The left-hand side can be approximated with a Tailor expansion so that we obtain 
 \[
 P_0\frac{d}{dQ_{0,n}^k}S_{Q_{0,n}^k}(\phi_j^k)(Q_n^k-Q_{0,n}^k)\approx -(P_n-P_0)S_{Q_n^k}(\phi_j^k).\]
 By empirical process theory one approximates $(P_n-P_0)S_{Q_n^k}(\phi_j^k)$ with $(P_n-P_0)S_{Q_{0,n}^k}(\phi_j^k)$ so that we obtain
 \[
 -P_0\frac{d}{dQ_{0,n}^k}S_{Q_{0,n}^k}(\phi_j^k)(Q_n^k-Q_{0,n}^k)\approx (P_n-P_0)S_{Q_{0,n}^k }(\phi_j^k).\] 
 The left-hand side is a linear operator on $D^{(k)}({\cal R}_n^k)$ applied to $(Q_n^k-Q_{0,n}^k)\in D^{(k)}({\cal R}_n^k)$, so  that by the Riesz-representation theorem it allows an inner product representation for any given inner product on $D^{(k)}({\cal R}_n^k)$. 
 However, for likelihood behaving loss functions we have a particularly nice inner product representation:
 \[
 -P_0 \frac{d}{dQ_0}S_{Q_0}(\phi)(Q_n^k-Q_0)=P_0 S_{Q_0}(\phi)S_{Q_0}(Q_n^k-Q_0).\]
 Therefore, up till a second order remainder we also have
 \begin{equation}\label{likelihoodlossprop}
 -P_0\frac{d}{dQ_{0,n}^k}S_{Q_{0,n}^k}(\phi)(Q_n^k-Q_{0,n}^k)\approx P_0 S_{Q_{0,n}^k}(\phi)S_{Q_{0,n}^k}(Q_n^k-Q_{0,n}^k).\end{equation}
 Thus, we then obtain
 \[
 P_0 S_{Q_{0,n}^k}(\phi_j^k)S_{Q_{0,n}^k}(Q_n^k-Q_{0,n}^k)\approx (P_n-P_0)S_{Q_{0,n}^k}(\phi_j^k).\]
 The left-hand side represents an inner product $\langle \phi_j^k,(Q_n^k-Q_{0,n}^k)\rangle_{{\cal R}_n^k}$  on $D^{(k)}({\cal R}_n^k)$. 
 We can now choose an orthonormal basis $\{\phi_j^*: j\in {\cal R}_n^k\}$ of $\{\phi_j^k:j\in {\cal R}_n^k\}$, w.r.t. to this inner product, where each $\phi_j^*$ is a particular linear combination of $\{\phi_j^k:j\in {\cal R}_n^k\}$.  Due to linearity of above terms in $\phi_j^k$, we then obtain
 \begin{equation}\label{key21}
 \langle \phi_j^*,Q_n^k-Q_{0,n}^k\rangle_{{\cal R}_n^k}\approx (P_n-P_0)S_{Q_{0,n}^k}(\phi_j^*).\end{equation}
 Since $(Q_n^k-Q_{0,n}^k)=\sum_{j\in {\cal R}_n^k}\langle \phi_j^*,Q_n^k-Q_{0,n}^k\rangle_{{\cal R}_n^k} \phi_j^*$, this then shows
 \[
 (Q_n^k-Q_{0,n}^k)(x)\approx (P_n-P_0)D_{Q_{0,n}^k,x},\]
 where \[
 D_{Q_{0,n}^k,x}=\sum_{j\in {\cal R}_n^k}S_{Q_{0,n}^k}(\phi_j^*)\phi_j^*(x).\]
 However, due to the orthonormality of $\phi_j^*$ w.r.t this covariance of score inner product $\langle \cdot,\cdot\rangle_{{\cal R}_n^k}$, we have $P_0 S_{Q_{0,n}^k}(\phi_{j_1}^*)S_{Q_{0,n}}(\phi_{j_2}^*)=0$ if $j_1\not =j_2$ and it equals $1$ if $j_1=j_2$. 
 This thus proves that with $d_n=\mid {\cal R}_n^k\mid$
 \[
 P_0\{D_{Q_{0,n}^k,x}\}^2=\sum_{j\in {\cal R}_n^k} \{\phi_j^*(x)\}^2,\]
 so that we obtain, by the CLT, 
 \[
 (n/d_n)^{1/2}(Q_n^k-Q_{0,n}^k)(x)/\tilde{\sigma}_n(x)\Rightarrow_d N(0,1),\]
 where 
 \[
 \tilde{\sigma}^2_{0,n}(x)\equiv \frac{1}{d_n}\sum_{j\in {\cal R}_n^k}\{\phi_j^*(x)\}^2.\]
 This particularly powerful proof establishing pointwise asymptotic normality of $Q_n^k$ at rate $(n/d_n)^{-1/2}$ was driven by the approximate identity (\ref{likelihoodlossprop}), known to hold for log-likelihood loss, but more generally, we will refer to likelihood behaving loss functions as loss functions for which (\ref{likelihoodlossprop}) holds.  
 In the examples below this property can be explicitly verified. 
 
 Let's now consider a non-likelihood behaving loss function. We can still define as inner product $\langle \phi_j,Q_n^k-Q_{0,n}^k\rangle_{{\cal R}_n^k}\equiv -P_0 \frac{d}{dQ_{0,n}^k}S_{Q_{0,n}^k}(\phi)(Q_n^k-Q_{0,n}^k)$  (so that we obtain an exact Riesz-representation of the right-hand side) and define $\{\phi_j^*:j\in {\cal R}_n^k\}$ as the orthonormal basis for $\{\phi_j^k:j\in {\cal R}_n^k\}$ and thus $D^{(k)}({\cal R}_n^k)$ w.r.t. this inner product.
 We then still obtain (\ref{key21}) and thus the same linear expansion $(Q_n^k-Q_{0,n}^k)(x)\approx (P_n-P_0)D_{Q_{0,n}^k,x}$ with $D_{Q_{0,n},x}=\sum_{j\in {\cal R}_n^k}S_{Q_{0,n}^k}(\phi_j^*)\phi_j^*(x)$. The only difference is that the asymptotic variance is now given by 
\[
\tilde{\sigma}^2_{0,n}=\frac{1}{d_n}\sum_{j_1,j_2}\Sigma_{0,n}(j_1,j_2)\phi_{j_1}^*(x)\phi_{j_2}^*(x),\]
 where $\Sigma_{0,n}(j_1,j_2)=P_0 S_{Q_{0,n}}(\phi_{j_1}^*)S_{Q_{0,n}}(\phi_{j_2}^*)$ is now not an identity or diagonal matrix. Lemma \ref{lemmaboundedvar} shows that this will be bounded by the maximal eigenvalue of $\Sigma_{0,n}$, and we suggest that that will not converge faster to infinity (if at all) than a $\log n$-factor. 
 
  \subsection{Least squares regression for continuous outcome}
 Consider the special  case of least squares regression with  $O=(X,Y)$; $Q_0(X)=E_0(Y\mid X)$; $L(Q)(X,Y)=(Y-Q(X))^2$, $S_{Q}=d/dQ L(Q)(\phi_n^0)$ is given by $S_{Q}(\phi_n^0)=\phi_n^0(X)(Y-Q(X))$. Now $Q_n^0$ is a least squares regression estimator based on the parametric linear regression model $D^{(0)}({\cal R}_n^0)$.
 We also consider the case that $L_{\sigma^2_{0}}(Q)=(Y-Q(X))^2/\sigma^2_{0}(X)$ is  a weighted least squares loss function, where $\sigma^2_{0}(X)=E_0(Y-Q_{0}(X))^2\mid X)$ represents a nuisance parameter that needs to be estimated as well. In that case, $S_{Q}^0=\phi_n^0(X)/\sigma^2_{0}(X)(Y-Q(X))$.
 The advantage of the latter loss function is that it behaves as a log-likelihood loss in the sense that
 \begin{eqnarray*}
 P_0 d/dQ_{0}S_{Q_{0}}(\phi)(Q_n^0-Q_{0})&=&-P_0 \phi/\sigma^2_{0} (Q_n^0-Q_{0})\\
 &=&-P_0 S_{Q_{0}}(\phi)S_{Q_{0}}(Q_n^0-Q_{0}),
\end{eqnarray*}
so that the above general proof of asymptotic normality applies with diagonal covariance matrix $\Sigma_{0,n}$.
 
 \subsection{Logistic regression for binary outcome}
 Let $\mbox{Logit}x=x/(1-x)$ and $\mbox{expit}(x)=1/(1+\exp(-x))$.
 Let $O=(X,Y)$; $Y\in \{0,1\}$; $Q_0(X)=\mbox{Logit}P_0(Y=1\mid X)$; \[
 L(Q)(X,Y)=-\left\{ Y\log \mbox{expit}(X)+(1-Y)\log (1-\mbox{expit} (X))\right\} ;\]
  $S_{Q}(\phi_n^0)=d/dQL(Q)(\phi_n^0)=\phi_n^0(X)(Y-Q(X))$. Now, $Q_n^0=\sum_{j\in {\cal R}_n^0}\beta_n^0(j)\phi_j^0$ is a logistic linear regression estimator based on the parametric logistic linear regression working model $P(Y=1\mid X)\sim \mbox{expit}(\sum_{j\in {\cal R}_n^0}\beta(j)\phi_j^0(X))$.
  Due to $L(Q)$ being a log-likelihood loss we have the identity
  \begin{eqnarray*}
  P_0 \frac{d}{dQ_{0}}S_{Q_{0}}(\phi)(Q_n^0-Q_0)&=&
  -P_0 \phi \frac{d}{dQ_{0}}m_{Q_{0}}(Q_n^0-Q_0)\\
  &=&
  -P_0\phi m_{Q_0}(1-m_{Q_0})(Q_n^0-Q_0)\\
  &=&
  -P_0 S_{Q_0}(\phi)S_{Q_0}(Q_n^0-Q_0).
  \end{eqnarray*}

\subsection{Marginal  univariate density estimation}
 Another  example is that $O$ is a real valued random variable and the target function is defined as  the hazard $\lambda_0=\exp(Q_0)$ using exponential link so that  the data density is parametrized as $p_{Q}(o)=\exp(Q(o))\exp(-\int_0^{o}\exp(Q(x))dx )$. 
Then, the log-likelihood loss is given by $L(Q)(o)=-\log p_Q(o)=-Q(o)+\int_0^{o}\exp(Q(x))  dx$, while $\frac{d}{dQ}L(Q)(h)(o)=
-h(o)+\int_0^{o}\exp(Q(x)) h(x) dx$. Thus the MLE $Q_n^0=\sum_{j\in {\cal R}_n^0}\beta_n^0(j)\phi_j^0$ over $D^{(0)}({\cal R}_n^0)$ solves  the score equations $P_n \frac{d}{dQ_n^0}L(Q_n^0)(\phi_n^0)=0$, where, for each $j\in {\cal R}_n^0$, we have
\[
S_{Q_n^0}(\phi_j^0)(O)=-\phi_j^0(O)+\int_0^{O}\exp(Q_n^0(x)) \phi_j^0(x) dx.\]
Similarly, the oracle MLE $Q_{0,n}^0=\sum_{j\in {\cal R}_n^0}\beta_{0,n}^0(j)\phi_j^0$ solves $0=P_0 S_{Q_{0,n}^0}^0(\phi_j^0)=0$.
Note that $Q_n^0$ is now the partial likelihood estimator for this parametric working model $D^{(0)}({\cal R}_n^0)$ for its linear form $\log \lambda_0$.
Again, we have the identity \[
-P_0 \frac{d}{dQ_0}S_{Q_0}(\phi)(Q_n^0-Q_0)=
P_0 S_{Q_0}(\phi)S_{Q_0}(Q_n^0-Q_0).\]


\subsection{Conditional density estimation}
Consider now the case that $O=(T,X)$ and that our target function is the conditional hazard  $\lambda_0(t\mid X)=\exp(Q_0(t,X))$ 
again parametrized through the exponential link. Then, the conditional density of $T$, given $X$ is modeled as $p_{Q}(t\mid X)=\exp(Q(t,X))\exp(-\int_0^{\cdot}\exp(Q(t,X))dt )$.
The log-likelihood loss of the conditional density is given by $L(Q)(t,x)=-\log p_Q(t\mid x)=-Q(t,x)+\int_0^{t}\exp(Q(s,X))ds$, while the scores are given by \[
\frac{d}{dQ}L(Q)(h)(t,x)=
-h(t,x)+\int_0^{t}\exp(Q(s,X)) h(s,X) ds.\] Thus, the sieve spline-MLE $Q_n^0=\arg\min_{Q\in D^{(0)}({\cal R}_n^0)}P_n L(Q)$ implied by the index sets ${\cal R}_n^0$ has the form $Q_n^0=Q_{\beta_n^0}^0=\sum_{j\in {\cal R}_n^0}\beta_n^0(j)\phi_j^0$ with $\beta_n^0$ solving the score equations $P_n \frac{d}{dQ_n^0}L(Q_n^0)(\phi_j^0)=0$, where, for each $j\in {\cal R}_n^0$, we have
\[
S_{Q_n,j}(T,X)=\frac{d}{dQ_n^0}L(Q_n^0)(\phi_j^0)(T,X)=- \phi_j^0(T,X)+\int_0^{T}\exp(Q_n^0(s,X)) \phi_j^0(s,X) ds.\]
Similarly, $P_0 S_{Q_{0,n}^0}=0$ is solved by the oracle MLE $Q_{0,n}^0=\arg\min_{Q\in D^{(0)}({\cal R}_n^0)}P_0L(Q)$.
One can verify that $-P_0\frac{d}{dQ_0}S_{Q_0}(\phi)(Q_n^0-Q_0)=P_0 S_{Q_0}(\phi)S_{Q_0}(Q_n^0-Q_0)$.

\subsection{Examples covered in this article.}
The uniform approximation error results for our proposed finite dimensional $k$-th order splines working models are not specific to the loss function and therefore apply to all examples. 
The asymptotic pointwise asymptotic normality, and asymptotic linearity for smooth features of $Q_0$,  of the linear least squares regression and logistic linear regression sieve and HAL-ML estimators are presented in Sections \ref{section8} and \ref{section10}, respectively. The asymptotic pointwise asymptotic normality and asymptotic linearity for smooth features of $Q_0$, of the sieve and HAL-MLEs of the univariate density and conditional density  based on log-likelihood loss are covered by our general asymptotic normality results for log-likelihood behaving loss functions in Section \ref{section9} and \ref{section11}, respectively.

\section{General $k$-th order spline representations of $k$-th order differentiable multivariate real valued functions}\label{section3}
In the next definition we define a $k$-th order sequentially defined Lebesgue Radon-Nikodym derivative $Q^{(k)}_{\bar{s}(k+1)}$ for a $k+1$-dimensional vector of nested subsets $\bar{s}(k+1)=(s_1,\ldots,s_{k+1})$, $s_{k+1}\subset s_{k}\ldots\subset s_1\subset \{1,\ldots,d\}$. 
\begin{definition}
For a given $\bar{s}(k+1)$ we define $m(\bar{s}(k+1))$ as the smallest $m
in \{1,\ldots,k\}$ for which $s_{m+1}$ is the empty set. 
\end{definition}
\begin{definition}
For a given subset $s\subset\{1,\ldots,d\}$, we define $\mu_s$ as the Lebesgue measure on $(0(s),1(s)]$, so that $d\mu_s(u(s))=\prod_{j\in s}du(j)$.
\end{definition}

\begin{definition}{\bf (Defining $Q^{(k)}_{\bar{s}(k+1)}$)}
Let $Q:[0,1]^d\rightarrow\openr$ be a given real valued function on the $d$ dimensional unit cube. It is assumed that $Q$ is $k$-times differentiable in the sense that, for all $\bar{s}(k+1)$, $Q_{\bar{s}(k+1)}^{(k)}$ as defined below exists. 
Let $k\in \{0,1,2,\ldots\}$ be given.
Let $\bar{s}(k+1)=(s_1,\ldots,s_{k+1})$ be a vector of nested subsets of $\{1,\ldots,d\}$ so that
$s_{k+1}\subset s_{k}\ldots\subset s_1\subset\{1,\ldots,d\}$. For each $j$, we include the empty subset $s_j$ as one of the possible subsets of $s_{j-1}$. 
Let ${\cal S}^{k+1}(d)$ be the set of such vectors $\bar{s}(k+1)$ of $k+1$ nested subsets of $\{1,\ldots,d\}$. 
For a given $\bar{s}(k+1)\in {\cal S}^{k+1}(d)$ we define $Q^{(k)}_{\bar{s}(k+1)}$  (and $Q^{(j)}_{\bar{s}(j)}$, $Q^{(j)}_{\bar{s}(j+1)}$, $j=1,\ldots,k$)  recursively as follows.
\begin{itemize}
\item ($j=0$) For a given $s_1$ we define $Q_{s_1}^{(0)}(x(s_1))=Q(x(s_1),0(-s_1))$ as the $s_1$-specific section of $Q^{(0)}=Q$, and, for empty subset $s_1$ we have $Q_{s_1}^{(0)}=Q(0)$, and $Q^{(k)}_{\bar{s}(k+1)}=Q(0)$ for the $\bar{s}(k+1)$ with $s_1$ empty subset and $k=0,1,\ldots$.
For non-empty subsets $s_1$,  let $Q^{(1)}_{s_1}=\frac{dQ_{s_1}^{(0)}}{d\mu_{s_1}}$ be the Radon-Nikodym derivative of section $Q_{s_1}^{(0)}$ w.r.t. Lebesgue measure $\mu_{s_1}$ with $d\mu_{s_1}(x(s_1))=\prod_{j\in s_1}dx(j)$.
\item ($j=1$) 
For non-empty $s_1$, for a given $\bar{s}(2)=(s_1,s_2)$, given $Q^{(1)}_{s_1}$ defined above, we define $Q_{s_1,s_2}^{(1)}(x(s_2))=Q^{(1)}_{s_1}(x(s_2),0(s_1/s_2))$ and for empty subset $s_2$, we have $Q_{s_1,s_2}^{(1)}=Q^{(1)}_{s_1}(0(s_1))$, and $Q^{(l)}_{\bar{s}(l+1)}=Q^{(1)}_{s_1}(0(s_1))$ for all $l=2,\ldots,k$ and $s_2=\ldots=s_{l+1}$.  For non-empty $s_1,s_2$, we then define $Q_{s_1,s_2}^{(2)}(x(s_2))=\frac{dQ_{s_1,s_2}^{(1)}}{d\mu_{s_2}}$. 

\item (general $j$ up till $j=k$)
For $j=1,\ldots,k$,  given $(s_1,\ldots,s_j)$ are non-empty (i.e., $s_j$ is non-empty) and given  $Q_{\bar{s}(j)}^{(j)}=Q_{s_1,\ldots,s_{j}}^{(j)}$  as a function of $x(s_j)$, for a non-empty subset $s_{j+1}\subset s_j$,  define $Q_{s_1,\ldots,s_{j},s_{j+1}}^{(j)}(x(s_{j+1}))=Q_{\bar{s}(j)}^{(j)}(x(s_{j+1}),0(s_{j}/s_{j+1}))$. For empty subset $s_{j+1}$ we have $Q_{\bar{s}(j),s_{j+1}}^{(j)}=Q^{(j)}_{\bar{s}(j)}(0(s_j))$, and $Q^{(l)}_{\bar{s}(l+1)}=Q^{(j)}_{\bar{s}(j)}(0(s_j))$ for all  $l=j+1,\ldots,k$ and $s_{j+1}=\ldots=s_{l+1}$.  If $j<k$, for non-empty subset $s_{j+1}$,  let $Q_{\bar{s}(j+1)}^{(j+1)}=\frac{d Q_{\bar{s}_{j+1}}^{(j)}} {d\mu_{s_{j+1} } }$, and iterate to next $j=j+1$.
At $j=k$, given $Q_{\bar{s}(k)}^{(k)}$, this then defines $Q_{\bar{s}(k+1)}^{(k)}$.
\end{itemize}
\end{definition}
\begin{lemma}
We note that for $\bar{s}(k+1)\in {\cal S}^{k+1}(d)$ for which there is a smallest $m=m(\bar{s}(k+1))\in \{1,\ldots,k\}$, for which $s_{m+1}$ is empty set, we have
\[
Q^{(k)}_{\bar{s}(k+1)}=Q^{(m)}_{\bar{s}(m)}(0(s_m))\] 
is a constant function equal to $Q^{(m)}_{\bar{s}(m)}(0(s_m))$.
\end{lemma}

This allows us now to define the $k$-th order smoothness class $D^{(k)}([0,1]^d)$ for $k=1,\ldots$, and  a corresponding $k$-th order sectional variation norm $\pl Q\pl_{v,k}^*$ of a function $Q\in D^{(k)}([0,1]^d)$, which represent a sum of variation norms of all the $k$-th order derivatives $Q^{(k)}_{\bar{s}(k+1)}$, where the variation norm of a constant function is just the absolute value of the constant. 

\begin{definition}{\bf (Defining $k$-th order smoothness class $D^{(k)}([0,1]^d)$ and $k$-th order sectional variation norm)}
Recall  $D^{(0)}([0,1]^d)$ is the set of $d$-variate real valued cadlag functions $Q:[0,1]^d\rightarrow\openr$ and $D^{(0)}_M([0,1]^d)$ is the set of $d$-variate real valued cadlag functions with sectional variation norm $\pl Q^{(0)}\pl_v^*=\sum_{s_1}\pl Q_{s_1}^{(0)}\pl_v\leq M$ bounded by $M$, where, for $s_1$ the empty subset we define $\pl Q_{s_1}^{(0)}\pl_v=\mid Q^{(0)}(0)\mid$.

Let $k\in \{1,\ldots\}$ be an integer. Let $D^{(k)}([0,1]^d)$ be the set of functions $Q:[0,1]^d\rightarrow\openr$ so that for each $\bar{s}(k+1)\in {\cal S}^{k+1}(d)$, the $k$-th order derivative $x(s_k)\rightarrow Q_{\bar{s}(k+1)}^{(k)}(x(s_{k}))$ exists and is an $\mid s_{k}\mid$-dimensional real valued cadlag function with finite  variation norm, $\pl Q_{\bar{s}(k+1)}^{(k)}\pl_v<\infty$, where,  for an $\bar{s}(k+1)$  with $s_{k+1}$ the empty subset we have $\pl Q_{\bar{s}(k+1)}^{(k)}\pl_v=\mid Q^{(m)}_{\bar{s}(m)}(0(s_m))\mid$ with $m=m(\bar{s}(k+1)) $  the smallest integer for which $s_{m+1}$ is the empty set. 
We define the $k$-th order sectional variation norm as 
\begin{eqnarray*}
\pl Q\pl_{v,k}^*&\equiv& \sum_{\bar{s}(k+1)\in {\cal S}^{k+1}(d)}\pl Q_{\bar{s}(k+1)}^{(k)}\pl_v\\
&=& \sum_{\bar{s}(k)\in {\cal S}^{k}(d)}\pl Q_{\bar{s}(k)}^{(k)}\pl_v^* ,
\end{eqnarray*}
 which equals the sum of  the   variation norms of $Q_{\bar{s}(k+1)}^{(k)}$ across all $\bar{s}(k+1)\in {\cal S}^k(d)$, or, equivalently, is the sum of sectional variation norms of $Q_{\bar{s}(k)}^{(k)}$ over all $\bar{s}(k)$, where the variation or sectional variation norm of a constant is the absolute value of the constant.   Note $\pl Q\pl_v^*=\pl Q\pl_{v,k=0}^*$.
We define $D^{(k)}_M([0,1]^d)=\{Q\in D^{(k)}_M([0,1]^d): \pl Q\pl_{v,k}^*\leq M\}$.
\end{definition}

The following theorem proves that if $Q\in D^{(k)}([0,1]^d)$, then $Q$ can be represented as an infinite linear combination of tensor products of $\leq k$-th order spline basis functions with knot-point in interior $[0,1]^d$. Moreover, it shows that the $L_1$-norm of its coefficients in this representation equals $\pl Q\pl_{k,v}^*$. Before we can present the theorem we need to define $k$-th order spline basis functions and $k$-th order primitives. 

Let's first define the $k$-th order spline basis functions and a particular tensor product $\bar{\phi}_{\bar{s}(k+1)}$ of $\leq k$-th order spline basis functions that appears in our representation. 
\begin{definition}{\bf ( Defining $k$-th order splines and $\bar{\phi}_{\bar{s}(k+1)}$)}
Let $d\mu(u)=\prod_{j=1}^d du(j)$ be the Lebesgue measure. 
Recall $\phi_u^0(x)=I(u\leq x)$ is the zero-order spline with knot-point $u\in [0,1]^d$.
For $j=1,\ldots$ and for a knot-point $u\in [0,1]^d$, we define
\[
\phi_{u}^j(x)=\int_{(u,x]}\phi^{j-1}_{u_1}(x) d\mu(u_1) .\]
  We refer to $\phi_u^j$ as the $j$-th order spline basis function with knot-point $u$.

For $\bar{s}(k+1)\in {\cal S}^{k+1}(d)$, we also define the function
 \begin{equation}\label{barphi}
 \bar{\phi}_{\bar{s}(k+1)}(x)\equiv \prod_{j=1,\mid s_j\mid>0}^k \phi^j_0(x(s_j/s_{j+1}))=\prod_{j=1}^{\min(k,m(\bar{s}(k+1))}\phi^j_0(x(s_j/s_{j+1})).
 \end{equation} 
 \end{definition}
 If $s_{k}$ is non-empty, then $\bar{\phi}_{\bar{s}(k+1)}$ is a function of $x$ through $x(s_1/s_{k+1})$. 
If  $m=m(\bar{s}_{k+1})<k$, then
 \[
 \bar{\phi}_{\bar{s}(k+1)}(x)=\bar{\phi}_{\bar{s}(m+1)}(x)=\prod_{j=1}^{m} \phi_0^j(x(s_j/s_{j+1})),\]
 which is now a function of $x(s_1/s_{m+1}))$. 

{\bf Examples of higher order splines:}
For example, if $d=1$, then $\phi_u^0(x)=I(x\geq u)$ while the first order spine equals $
\phi_u^1(x)=\int_{(u,x]}\phi^0_{u_1}(x)du_1=\int_{(u,x]}I(u_1\leq x) du_1=
(x-u)I(x\geq u)$. For $d=1$, the second order spline is given by $
\phi_u^2(x)=\int_{(u,x]}\phi_{u_1}^1(x)du_1=
\int_{(u,x]}I(x\geq u_1)(x-u_1) du_1$, which thus equals $\phi_u^2(x)=1/2(x-u)^2I(x\geq u)$,  and so on. Similarly, for $d=2$, the first order spline equals $\phi_u^1(x)=\int_{(u,x]}I(u_1\leq x) du_1=
(x_1-u_1)(x_2-u_2)I(x_1\geq u_1,x_2\geq u_2)$ and the second order spline equals
$\phi_u^2(x)=1/4(x_1-u_1)^2(x_2-u_2)^2I(x_1\geq u_1,x_2\geq u_2)$.  

So we have the following trivial observation.
\begin{lemma}
For general dimensions $d$, the $k$-th order spline $\phi_u^k(x)=\prod_{j=1}^d\phi_{u_j}^k(x_j)$, where $
\phi_{u_j}^k(x_j)$ is the $k$-th order univariate spline basis function with knot-point $u_j$. That is, the $k$-th order spline basis functions are tensor products of the $k$-th order spline basis functions for dimension $d=1$.
\end{lemma}


We also want to emphasize an equivalent representation of $Q\in D^{(k)}([0,1]^d)$ in terms of $k$-th order primitives, since that most naturally yields our uniform approximation error results. We first define what we mean with a $k$-th order primitive function.

\begin{definition}{\bf (Defining $k$-th order primitive function and corresponding class)}
Let $\mu$ be the Lebesgue measure. Let $k$ be an integer in $\{1,\ldots\}$.
A $d$-variate $k$-th order primitive  function is defined as  function $Q:(0,1]^d\rightarrow\openr$ such that $Q^{(m)}=dQ^{(m-1)}/d\mu$ exists as $d$-variate function $Q^{(m)}:(0,1]^d\rightarrow\openr$, for $m=1,\ldots,k$.  A first order primitive function satisfies $Q(x)=\int_{(0,x]}Q^{(1)}(y)d\mu(y)$;  a second order primitive function satisfies $Q^{(2)}(x)=\int_{(0,x]}\int_{(0,y_1]}Q^{(2)}(y_2)d\mu(y_2)d\mu(y_1)$;  and, in general, a $k$-th order primitive function  can be represented as \[
Q(x)=\int_{(0,x]}\int_{(0,y_1]}\ldots\int_{(0,y_{k-1}]} Q^{(k)}(y_k) \prod_{j=k}^{1} d\mu(y_j).\]  
We also say that $Q$ is the $k$-th order primitive of $Q^{(k)}$ and we use notation $Q=\mu^k(Q^{(k)})$. 
For a general function $Q$, $\mu^1(Q)=\mu(Q)=\int_{(0,x]}fd\mu$ and for general $j=2,\ldots$, we have that $\mu^j(Q)$ is defined by the recursive relation $\mu^j(Q)=\int_{(0,x]}\mu^{j-1}(Q) d\mu$. 
We note that if $Q$ is only a function of  $y$ through $y(s)$, then $\mu(Q)(x(s))=\int_{(0(s),x(s)]}Q(y(s))d\mu(y(s))$.
For given $k\in \{1,2,\ldots\}$, let ${\cal F}^{(k)}((0,1]^d)=\{\mu^k(Q): Q\in {\cal F}^{(0)}((0,1]^d)\}$ be the class of $d$-dimensional $k$-th order primitive functions on $(0,1]^d$ whose $k$-th order density $Q^{(k)}$ has bounded variation norm, $\pl Q^{(k)}\pl_v<\infty$.
We also define ${\cal F}^{(k)}_{M}((0,1]^d)=\{Q\in {\cal F}^{(k)}(0,1]^d: \pl Q^{(k)}\pl_{v}<M\}$. Recall definition ${\cal F}^{(0)}((0,1]^d)$ and ${\cal F}^{(0)}_M((0,1]^d)$ as cadlag functions $Q$ for which $Q(x)=\int_{(0,x]} dQ(u)$ with $\pl Q\pl_v<\infty$ and $\pl Q\pl_v\leq M$, respectively. 
\end{definition}

The following lemma shows that the $k$-th order primitive of the zero-order spline basis function equals the $k$-th order spline basis function. 
\begin{lemma}\label{lemmamuphij}
We have
\[
\mu^k(\phi_u^0)=\phi_u^k.\]
\end{lemma}
{\bf Proof:}
For example, $\mu(\phi_u^0)(x)=\int_{(0,x]}I(u\leq y)d\mu(y)=\int_{[u,x]}d\mu(y)=\phi_u^1(x)$.
Similarly, $\mu(\phi_u^1)(x)=\int_{(0,x]} \phi_u^1(y) d\mu(y)=\int_{[u,x]}\phi_u^1(y)d\mu(y)=\phi_u^2(x)$. This argument iterates to obtain $\mu^k(\phi_u^0)=\phi_u^k$.
$\Box$.

More generally, we have for $\tilde{Q}\in D^{(0)}([0,1]^d)$,  $\mu^k(\tilde{Q})=\int_u\phi_u^k(\cdot)d\tilde{Q}(u)$.
\begin{lemma}\label{lemmahandyz1}
If $\tilde{Q}\in D^{(0)}([0,1]^d)$, then 
$\mu^k(\tilde{Q})
=\int_u \phi_u^k(\cdot)d\tilde{Q}(u)$ is a linear combination of $k$-th order splines with ''$L_1$''-norm the sectional variation norm of $\tilde{Q}$ (which also equals the $k$-th order sectional variation norm of $\mu^k(\tilde{Q})$). 
\end{lemma}
{\bf Proof:}
If $\tilde{Q}\in D^{(0)}([0,1]^d)$, then  $\tilde{Q}(x)=\int_{[0,x]}d\tilde{Q}(u)=\int \phi_u^0(x)d\tilde{Q}(u)$, so that, using that $\mu^k(\tilde{Q})$ is a linear operator in $\tilde{Q}$ and Lemma \ref{lemmamuphij},  we have
\begin{eqnarray*}
\mu^k(\tilde{Q})&=&\mu^k\left(\int \phi_u^0(\cdot)d\tilde{Q}(u)\right)=\int_u \mu^k(\phi_u^0)(\cdot) d\tilde{Q}(u)\\
&=&\int_u \phi_u^k(\cdot)d\tilde{Q}(u),
\end{eqnarray*}
which is thus a linear combination of $k$-th order splines. $\Box$

The latter result allows us to switch between  a linear combination of $k$-th order splines of form $\int_u \phi_u^k(\cdot)d\tilde{Q}(u)$ and $\mu^k(\tilde{Q})$, which will correspond with two equivalent representations for our $k$-th order smoothness class given in next theorem.  
We are now ready to present for a function $Q\in D^{(k)}([0,1]^d)$, its  $k$-th order spline representation and the (by Lemma \ref{lemmahandyz1}) corresponding representation in terms of $k$-th order primitives.

\begin{theorem}\label{theoremkthorderspline}
{\bf (Exact representation of $k$-th order smooth function as linear combination of $\leq k$-th order splines)}\ \newline
For $\bar{s}(k+1)\in {\cal S}^{k+1}(d)$ with $\mid s_{k+1}\mid >0$, we define
 $\tilde{Q}^{(k)}_{\bar{s}(k+1)}(u(s_{k+1}))\equiv  \int_{(0(s_{k+1}),u(s_{k+1})]} Q^{(k)}_{\bar{s}(k+1)}(dy(s_{k+1}))$.
\newline
 {\bf Conventions:} We make the convention that $\prod_{j=1}^m a_j =1$ for $m=0$. If
$s_{k+1}$ is  the empty set, and  for $m=m(\bar{s}(k+1))$, $s_{m+1}$ is the  first empty set in the sequence $\bar{s}(k+1)$, then 
\begin{eqnarray}
\mu^k(\tilde{Q}^{(k)}_{\bar{s}(k+1)})&\equiv& Q^{(m)}_{\bar{s}(m)}(0(s_m)),\label{conventiona}
\end{eqnarray}
and, thus, by Lemma \ref{lemmahandyz1}, we also have \[
\bar{\phi}_{\bar{s}(k+1)}\int \phi^k_{u_k(s_{k+1})}(x(s_{k+1})) Q^{(k)}_{\bar{s}(k+1)}(du_k(s_{k+1}) )\equiv
\bar{\phi}_{\bar{s}(k+1)} Q^{(m)}_{\bar{s}(m)}(0(s_m)).\]
Let $Q\in D^{(k)}([0,1]^d)$. With the above convention (\ref{conventiona}), we have
\begin{eqnarray}
Q(x)&=&\sum_{\bar{s}(k+1)\in {\cal S}^{k+1}(d)}\bar{\phi}_{\bar{s}(k+1)}(x )\int_{(0(s_{k+1}),x(s_{k+1})]}\phi^k_{u_k(s_{k+1})}(x(s_{k+1})) Q^{(k)}_{\bar{s}(k+1)}(du_k(s_{k+1}) )\nonumber \\
&=& \sum_{\bar{s}(k+1)\subset{\cal S}^{k+1}(d)}\bar{\phi}_{\bar{s}(k+1)}(x )\int \phi^k_{u_k(s_{k+1})}(x(s_{k+1})) Q^{(k)}_{\bar{s}(k+1)}(du_k(s_{k+1}) )\nonumber \\
&=&  \sum_{\bar{s}(k+1), s_{k+1}=\emptyset} \bar{\phi}_{\bar{s}(m+1)}(x )Q^{(m)}_{\bar{s}(m)}(0(s_m)) \mbox{ (with $m=m(\bar{s}(k+1))$ in each term)}\nonumber \\
&&+\sum_{\bar{s}(k+1), s_{k+1}\not =\emptyset}  \bar{\phi}_{\bar{s}(k+1)}(x )\int \phi^k_{u_k(s_{k+1})}(x(s_{k+1})) Q^{(k)}_{\bar{s}(k+1)}(du_k(s_{k+1}) ). \label{kthorderrepr}
\end{eqnarray}

We also have the alternative representation in terms of $k$-th order primitives of $\tilde{Q}^{(k)}_{\bar{s}(k+1)}$ across all $\bar{s}(k+1)\subset{\cal S}^{k+1}(d)$: \begin{eqnarray*}
Q(x)&=&
\sum_{\bar{s}(k+1)\subset{\cal S}^{k+1}(d)} \bar{\phi}_{\bar{s}(k+1)}(x) 
\mu^{k}\left ( \tilde{Q}^{(k)}_{\bar{s}(k+1)}\right)\\
 &=&\sum_{\bar{s}(k+1), s_{k+1}=\emptyset} \bar{\phi}_{\bar{s}(m+1)}(x )Q^{(m)}_{\bar{s}(m)}(0(s_m)) \mbox{ (with $m=m(\bar{s}(k+1))$ in each term)}\\
 &&+\sum_{\bar{s}(k+1),  s_{k+1}\not = \emptyset} \bar{\phi}_{\bar{s}(k+1)}(x) 
\mu^{k}\left ( \tilde{Q}^{(k)}_{\bar{s}(k+1)}\right).
\end{eqnarray*}
\end{theorem}
The conventions are due to our definition of $Q^{(k)}_{\bar{s}(k+1)}$ for $\bar{s}(k+1)$ with $s_{k+1}$ being the empty set.
Note that indeed
\[
\begin{array}{l}
\pl Q\pl_{v,k}^*= \sum_{\bar{s}(k+1), s_{k+1}=\emptyset}\mid Q^{(m)}_{\bar{s}(m)}(0(s_m)) \mid 
+
 \sum_{\bar{s}(k+1), s_{k+1}\not =\emptyset}  \int \mid Q^{(k)}_{\bar{s}(k+1)}(du_k(s_{k+1})) \mid 
 \end{array}
 \]
 is precisely the $L_1$-norm of the coefficient ''vector''  in the linear $k$-th order spline representation (\ref{kthorderrepr}) of $Q$.

\subsection{Example}
Consider the case that $d=3$ and $k=1$. 
Then,
\begin{eqnarray*}
Q(x)&=&\sum_{s_1\subset\{1,2,3\}} \phi^1_0(x(s_1)) Q^{(1)}_{s_1}(0(s_1))\\
&&+\sum_{\bar{s}(2), \mid s_2\mid>0}\bar{\phi}_{\bar{s}(2)}(x)\int \phi^1_{u(s_2)}(x(s_2))Q^{(1)}_{\bar{s}(2)}(du(s_2)),
\end{eqnarray*}
where $\bar{\phi}_{\bar{s}(2)}(x)=\phi_0^1(x(s_1/s_2))=\prod_{l\in s_1/s_2}x(l)$.
Note that for $k=1$, the representation is a linear combination of first order splines, just also including knot-points at zero. 
Consider the case that $d=3$ and $k=2$. 
Then,
\begin{eqnarray*}
Q(x)&=&\sum_{s_1\subset\{1,2,3\}} \phi^1_0(x(s_1)) Q^{(1)}_{s_1}(0(s_1))\\
&&+\sum_{\bar{s}_2}\phi^1_0(x(s_1/s_2))\phi^2_0(x(s_2))Q^{(2)}_{\bar{s}(2)}(0(s_2))\\
&&+\sum_{\bar{s}(3), \mid s_3\mid>0}\phi^1_0(x(s_1/s_2))\phi_0^2(x(s_2/s_3))\int \phi^2_{u(s_3)}(x(s_3)) Q^{(2)}_{\bar{s}(3)}(du(s_3))
\end{eqnarray*}
In this case, we see that the representation includes linear combinations of first order splines at the zero edges, beyond the second order splines.

\subsection{General submodels of the $k$-th order spline smoothness class}
Note that by defining ${\cal R}(\bar{s}(k+1))\equiv {\cal R}(s_{k+1})=(0(s_{k+1}),1(s_{k+1})]$ when $s(k+1)$ is non-empty,  we  can represent the target function $Q\in D^{(k)}([0,1]^d)$ as
\begin{eqnarray*}
Q&=&\sum_{\bar{s}(k+1)\in {\cal S}^{k+1}(d),\mid s_{k+1}\mid>0} \int_{u\in {\cal R}({s}(k+1))}\bar{\phi}_{\bar{s}(k+1)}  \phi_u^k Q^{(k)}_{\bar{s}(k+1)}(du)\\
&&+\sum_{\bar{s}(k+1)\in {\cal S}^{k+1}(d),\mid s_{k+1}\mid =0} \bar{\phi}_{\bar{s}(m+1)}Q^{(m)}_{\bar{s}(m)}(0(s_m)),
\end{eqnarray*}
where $m=m(\bar{s}(k+1))$
We could define the total index set as \[
{\cal R}^k(d)=\{(\bar{s}(k+1),u):\bar{s}(k+1),\mid s_{k+1}\mid>0,u\in {\cal R}(\bar{s}(k+1))\}
\cup \{\bar{s}(k+1):\mid s_{k+1}\mid =0\},\]
where $\bar{s}(k+1)\in {\cal S}^{k+1}(d)$.
 By our Theorem \ref{theoremkthorderspline} $D^{(k)}_M([0,1]^d)$ is just a family of infinite linear combination of basis functions $x\rightarrow\phi_{\bar{s}(k+1),u}\equiv \bar{\phi}_{\bar{s}(k+1)}(x)\phi_u^k(x)$ indexed by $(\bar{s}(k+1),u)\in {\cal R}^k(d)$ with $\mid s_{k+1}\mid>0$, and $x\rightarrow \bar{\phi}_{\bar{s}(m+1)}$ for $\bar{s}(k+1)\in {\cal R}^k(d)$ with $\mid s_{k+1}\mid =0$ with $m=m(\bar{s}(k+1))$,  where the  $L_1$-norm represents the $k$-th order sectional variation norm $\pl Q\pl_{k,v}^*\leq M$.
 Theorem \ref{theoremkthorderspline} states that  $D^{(k)}_M({\cal R}^k(d))=D^{(k)}_M([0,1]^d)$.

This $k$-th order spline representation of functions in $D^{(k)}([0,1]^d)$ and $D^{(k)}_M([0,1]^d)$ also suggest a general class of submodels of this smoothness class $D^{(k)}([0,1]^d)$ obtained by replacing the unrestricted sets ${\cal S}^{k+1}(d)$ for $\bar{s}(k+1)$ and  ${\cal R}(\bar{s}(k+1))$ for the knot-points by restricted sets ${\cal S}^{k+1,*}(d)$ and, for each $\bar{s}(k+1)$ with $\mid s_{k+1}\mid =0$, ${\cal R}^*(\bar{s}(k+1))$, respectively.  Such a restriction then defines a new total index set 
\[
{\cal R}^{k,*}(d)=\{(\bar{s}(k+1),u): \bar{s}(k+1)\in {\cal S}^{k+1,*}(d),\mid s_{k+1}\mid>0,u\in {\cal R}^*(\bar{s}(k+1))\}\cup{\cal R}_2^{k,*}(d),\]
where ${\cal R}_2^{k,*}(d)\equiv \{\bar{s}(k+1): \mid s_{k+1}\mid =0\}$.
 which  on its turn defines the submodel $D^{(k)}_{{\cal R}^{k,*}(d),M}([0,1]^d)$ of smoothness class $D^{(k)}_M([0,1]^d)$ consisting of the infinite  linear combinations of $\bar{\phi}_{\bar{s}(k+1)}\phi_u^k$ and $\bar{\phi}_{\bar{s}(k+1)}$ indexed by elements in ${\cal  R}^{k,*}(d)$ with $L_1$-norm representing $\pl Q\pl_{k,v}^*\leq M$. 

The subset ${\cal S}^{k+1,*}(d)$  can be chosen to exclude a collection of $\bar{s}(k+1)$-elements and, given a $\bar{s}(k+1)\in  {\cal S}^{k+1,*}(d)$, one can still restrict the set of non-zero coefficients $Q^{(k))}_{\bar{s}(k+1)}(du)$  to a subset ${\cal R}^*(\bar{s}(k+1))$ of $(0(s_{k+1}),1(s_{k+1})]$. Note that   ${\cal R}^*(\bar{s}(k+1))$ is the empty set when $\bar{s}(k+1)$ is not an element of ${\cal S}^{k+1,*}(d)$.

In this manner, by varying $k$ as well as the set of potential basis functions identified by ${\cal R}^{k,*}(d)$ one obtains a general class of models of form $D^{(k)}({\cal R}^{k,*}(d))$  that varies smoothness as well as support-restrictions on higher order derivatives. Recall from our zero-order spline section that restricting $s_1$ already implies that one works in a general additive model within the class of zero-order splines. Similarly, one can enforce higher order derivatives being zero on the left edges of $[0,1]^d$, resulting in restrictions on $s_j$, $j=2,\ldots,k+1$. Such submodels are discussed in some detail in  Appendix \ref{AppendixL}. 

{\bf Uniform approximation results generalized to submodels:}
Analogue to Section \ref{section4}, to obtain a corresponding finite dimensional approximation of this submodel $D^{(k)}({\cal R}^{k,*}(d))$,  one would select a set of ${\bf J}(\bar{s}(k+1))$ knot-points contained in ${\cal R}^*(\bar{s}(k+1))$ for all $\bar{s}(k+1)\in {\cal S}^{k+1,*}(d)$ with $s_{k+1}$ non-empty. One would naturally expect that our uniform approximation result  of Section \ref{section5} for $D^{(k)}({\cal R}^k(d))=D^{(k)}([0,1]^d)$ generalize to these smaller models, and  our MLE results w.r.t its oracle projection are ignorant of the precise basis/working model so they directly apply as well. Indeed, in Appendix \ref{AppendixM} we generalize our uniform approximation results for $k$-th order primitives  to general subsets ${\cal R}^*(\bar{s}(k+1))$ of ${\cal R}(\bar{s}(k+1))=(0(s_{k+1}),1(s_{k+1})]$ and thereby general index sets ${\cal R}^{k,*}(d)$ defined by ${\cal S}^{k+1,*}(d)$ and ${\cal R}^*(\bar{s}(k+1))$ across $\bar{s}(k+1)\in {\cal S}^{k+1,*}(d)$.


\section{Defining the finite dimensional $k$-th order spline approximation of a $k$-th order smooth function (and thereby candidate parametric working models for MLE).}\label{section4}


Suppose $Q\in D^{(k)}([0,1]^d)$ so that by Theorem \ref{theoremkthorderspline} we have
\begin{eqnarray}
Q&=&\Phi(\mu^{k}(\tilde{Q}^{(k)}_{\bar{s}(k+1)}):\bar{s}(k+1))\nonumber\\
&\equiv& \sum_{\bar{s}(k+1)\subset{\cal S}^{k+1}(d)} \bar{\phi}_{\bar{s}(k+1)}\mu^{k}(\tilde{Q}^{(k)}_{\bar{s}(k+1)})\label{repra} \\
&=&\sum_{\bar{s}(k+1),\mid s_{k+1}\mid\geq 1}\bar{\phi}_{\bar{s}(k+1)}\mu^k(\tilde{Q}^{(k)}_{\bar{s}(k+1)})+
\sum_{\bar{s}(k+1),\mid s_{k+1}\mid =0}\bar{\phi}_{\bar{s}(m+1)}Q^{(m)}_{\bar{s}(m)}(0(s_m))
\nonumber,
\end{eqnarray}
where the notation suppresses that $m=m(\bar{s}(k+1))$.
We note that
\[
\begin{array}{l}
\sum_{\bar{s}(k+1),\mid s_{k+1}\mid =0}\bar{\phi}_{\bar{s}(m+1)}Q^{(m)}_{\bar{s}(m)}(0(s_m))=
\sum_{j=1}^k \sum_{\bar{s}(k+1),m(\bar{s}(k+1))=j}\bar{\phi}_{\bar{s}(j+1)}Q^{(j)}_{\bar{s}(j)}(0(s_j))\\
=\sum_{j=1}^k\sum_{\bar{s}(j),\mid s_j\mid>0,\mid s_{j+1}\mid =0}\bar{\phi}_{\bar{s}(j+1)}Q^{(j)}_{\bar{s}(j)}(0(s_j)).
\end{array}
\]

For each $\bar{s}(k+1)$ with $s_{k+1}$ non-empty, we have that
$\tilde{Q}^{(k)}_{\bar{s}(k+1)}\in {\cal F}^{(0)}((0(s_{k+1}),1(s_{k+1})] )$, since it equals the $s_{k+1}$-section of $Q^{(k)}_{\bar{s}(k)}$ and can thus be represented as $\tilde{Q}^{(k)}_{\bar{s}(k+1)}(x(s_{k+1}))=\int_{(0(s_{k+1}),x(s_{k+1})]}  Q^{(k)}_{\bar{s}(k)}(du(s_{k+1}),0(s_k/s_{k+1}) )$.
Thus, each $\tilde{Q}^{(k)}_{\bar{s}(k+1)}$ with $s_{k+1}$ non-empty can be approximated by a linear combination of zero order splines $\phi_u^0$ with knot-points $u$ in  \[
 {\cal R}(\bar{s}(k+1),{\bf J}(\bar{s}(k+1)))= \{u\in (0(s_{k+1}),1(s_{k+1})]: u(s_{k+1})\in {\cal R}(\mid s_{k+1}\mid,{\bf J}(\bar{s}(k+1))\},\] 
implied by knot-point set ${\cal R}(\mid s_{k+1}\mid ,{\bf J}(\bar{s}(k+1))$ as defined in  corollary \ref{defknots}, using a $\bar{s}(k+1)$-specific size ${\bf J}(\bar{s}(k+1))$ and corresponding  variation independent coefficient vector
$\beta(\bar{s}(k+1),\cdot)$. By properties of Corollary \ref{defknots}, this then implies that this linear model has an  $O(r(d,J))$-$L^2$-approximation of $\tilde{Q}^k_{\bar{s}(k+1)}$.


For $\bar{s}(k+1)$ with $s_{k+1}$ empty set, we simply approximate $ \bar{\phi}_{\bar{s}(m+1)}Q^{(m)}_{\bar{s}(m)}(0(s_m))$ with $\bar{\phi}_{\bar{s}(m+1)}\beta((\bar{s}(k+1))$ with $\beta(\bar{s}(k+1))$ a constant that can be set equal to $Q^{(m)}_{\bar{s}(m)}(0(s_m))$ to get an exact equality (no approximation error), where $m=m(\bar{s}(k+1))$.
For $\bar{s}(k+1)$ with $s_{k+1}$ empty set, let ${\bf J}(\bar{s}(k+1))=1$ and we define ${\cal R}(\bar{s}(k+1),{\bf J}(\bar{s}(k+1)))=\{0(s_m)\}$ as the singleton with $m=m(\bar{s}(k+1))$.

Let ${\bf J}=({\bf J}(\bar{s}(k+1)):\bar{s}(k+1))$ be the combined vector of sizes ${\bf J}(\bar{s}(k+1))$ of these knot-point index  sets ${\cal R}(\bar{s}(k+1),{\bf J}(\bar{s}(k+1))$ across all $\bar{s}(k+1)\in {\cal S}^{k+1}(d)$.
Let $\beta_{\bf J}=(\beta_{{\bf J}}(\bar{s}(k+1),\cdot):\bar{s}(k+1))$ be the combined vector of coefficients obtained by stacking the $\beta_{{\bf J}}(\bar{s}(k+1),\cdot)$-coefficient-vector of size ${\bf J}(\bar{s}(k+1))$. 
We will use the following notation for the $\bar{s}(k+1)$-specific index set:
\[
{\cal R}(\bar{s}(k+1),{\bf J})\equiv {\cal R}(\bar{s}(k+1),{\bf J}(\bar{s}(k+1)))
.\]
Let \[
{\cal R}^k(d,{\bf J})\equiv \{(\bar{s}(k+1),u):u\in {\cal R}(\bar{s}(k+1),{\bf J}),\bar{s}(k+1)\subset{\cal S}^{k+1}(d)\}\]
 be the set of indices that identify both $\bar{s}(k+1)$ and the knot-point $u\in {\cal R}(\bar{s}(k+1),{\bf J})$.

Given $\bar{s}(k+1)$ with $s_{k+1}$ non-empty, let \[
\tilde{Q}^{(k)}_{\bar{s}(k+1),{\bf J},\beta_{\bf J}}=\sum_{u\in {\cal R}(\bar{s}(k+1),{\bf J}) }\beta_{\bf J}(\bar{s}(k+1),u)\phi_u^0\]
 be the finite dimensional zero-order spline approximation of $\tilde{Q}^{(k)}_{\bar{s}(k+1)}$.
Plugging this in $\mu^k(\tilde{Q}^{(k)}_{\bar{s}(k+1)})$ yields the corresponding approximation $\mu^{k}\left( \tilde{Q}^{(k)}_{\bar{s}(k+1),{\bf J},\beta_{\bf J}}\right)$.
By Lemma \ref{lemmahandyz1}, for $\bar{s}(k+1)$ with $s_{k+1}$ non-empty
we have
\[
\mu^{k}\left( \tilde{Q}^{(k)}_{\bar{s}(k+1),{\bf J},\beta_{\bf J}}\right)=
\sum_{u\in {\cal R}(\bar{s}(k+1),{\bf J})} \beta_{\bf J}(\bar{s}(k+1),u) \phi_u^k.\]
 


Plugging these finite dimensional approximations of $\mu^k(\tilde{Q}^{(k)}_{\bar{s}(k+1)})$ into our representation (\ref{repra}) of $Q$ gives the resulting $k$-th order spline approximation of $Q$: 
\begin{eqnarray*}
Q^k_{{\bf J},\beta_{\bf J} }
&=& \sum_{\bar{s}(k+1),\mid s_{k+1}\mid>0} \bar{\phi}_{\bar{s}(k+1)}\mu^{k}\left( \tilde{Q}^{(k)}_{\bar{s}(k+1),{\bf J},\beta_{\bf J}}\right) \\
&&+
\sum_{\bar{s}(k+1),\mid s_{k+1}\mid =0}
\bar{\phi}_{\bar{s}(k+1)}\beta(\bar{s}(k+1))\\
&=& \sum_{\bar{s}(k+1) ,\mid s_{k+1}\mid >0} \sum_{u\in {\cal R}(\bar{s}(k+1),{\bf J} ) }\beta_{{\bf J}}(\bar{s}(k+1),u)
\bar{\phi}_{\bar{s}(k+1)} \phi^k_u \\
&&+\sum_{\bar{s}(k+1),\mid s_{k+1}\mid =0}
\beta(\bar{s}(k+1)) \bar{\phi}_{\bar{s}(k+1)}
,
\end{eqnarray*}
where we are reminded that last sum equals $\sum_{\bar{s}(k+1),\mid s_{k+1}\mid =0}
\bar{\phi}_{\bar{s}(m(\bar{s}(k+1))+1)}\beta(\bar{s}(k+1))$.
Note also that this last sum is simply a finite (uniformly in ${\bf J}$) dimensional parametric form.

We note that $Q^k_{{\bf J},\beta_{\bf J} }$ is a linear combination over basis functions $\bar{\phi}_{\bar{s}(k+1)}\phi^k_{u}$ with $(\bar{s}(k+1),u)\in {\cal R}^k(d,{\bf J})$, where these basis functions reduce to $\bar{\phi}_{\bar{s}(k+1)}=\bar{\phi}_{\bar{s}(m(\bar{s}(k+1))+1)}$ for all $\bar{s}(k+1)$ with $s_{k+1}$ empty.
The total number of basis functions is given by $\bar{J}\equiv \sum_{\bar{s}(k+1)}{\bf J}(\bar{s}(k+1))$. 
\begin{definition}\label{definitionbasis}
Recall $\bar{\phi}_{\bar{s}(k+1)}=\bar{\phi}_{\bar{s}(m(\bar{s}(k+1))+1)}$ for all $\bar{s}(k+1)$ with $s_{k+1}$ empty.
Let \[
\phi_{\bar{s}(k+1),u}^k\equiv I(\mid s_{k+1}\mid >0)\bar{\phi}_{\bar{s}(k+1)}\phi^k_u+
I(\mid s_{k+1}\mid =0)\bar{\phi}_{\bar{s}(k+1)}.\]
\end{definition}
To emphasize that it is just a linear combination of $\bar{J}$ basis functions $\bar{\phi}_{\bar{s}(k+1),u}^k$ across all $(\bar{s}(k+1),u)\in {\cal R}^k(d,{\bf J})$, we can also write it as
\[
Q^k_{{\bf J},\beta_{\bf J}}=\sum_{ (\bar{s}(k+1),u)\in {\cal R}^k(d,{\bf J})}\beta_{\bf J}(\bar{s}(k+1),u) \phi^k_{\bar{s}(k+1),u} .\]

\section{Uniform approximation error of finite dimensional $k$-th order spline approximation of general $k$-th order smooth function}\label{section5}
We  assume that there exists a $\delta>0$ and $M<\infty$ so that for all $\bar{s}(k+1)$ with $s_{k+1}$ non-empty  $\delta J<{\bf J}(\bar{s}(k+1))<MJ$ across all $\bar{s}(k+1)$ so that all these knot-point sets ${\cal R}(\bar{s}(k+1),{\bf J})$ grow linearly in $J$. 
In addition, for all $\bar{s}(k+1)$ with $s_{k+1}$ empty set, let $\beta(\bar{s}(k+1))=Q^{(m)}_{\bar{s}(m)}(0(s_m)$ so that the corresponding parametric term does not contribute any approximation error. 
Assuming $Q\in D^{(k)}_M([0,1]^d)$ and $\beta(\bar{s}(k+1))$ for $s_{k+1}$ empty set defined as above, we have
\begin{eqnarray*}
(Q_{{\bf J},\beta_{\bf J}}^k-Q)
&=&\sum_{\bar{s}(k+1),\mid s_{k+1}\mid>0}\bar{\phi}_{\bar{s}(k+1)}
\left\{ \mu^{k}\left(\tilde{Q}^{(k)}_{\bar{s}(k+1),{\bf J},\beta_{\bf J} }\right)-
\mu^{k}(\tilde{Q}^{(k)}_{\bar{s}(k+1)})\right\}.
\end{eqnarray*}
Thus, it follows that
\[
\pl Q_{{\bf J},\beta_{\bf J}}^k-Q\pl_{\infty}
\leq  C(d,k) \max_{\bar{s}(k+1),\mid s_{k+1}\mid>0}\pl  \mu^{k}(\tilde{Q}^{(k)}_{\bar{s}(k+1),{\bf J},\beta_{\bf J}})-
\mu^{k}(\tilde{Q}^{(k)}_{\bar{s}(k+1)})\pl_{\infty},
\]
where $C(d,k)\equiv \sum_{\bar{s}(k+1),\mid s_{k+1}\mid>0}$. In Appendix \ref{AppendixD} we prove the following lemma.
\begin{lemma}

Let  $\tilde{Q}_{J,\beta}=\sum_{u\in {\cal R}(d,J)}\beta(u)\phi_u^0$ and ${\cal R}(d,J)$ satisfies (\ref{Rbm}).
Then, we have
\begin{equation}\label{todoa}
\sup_{\tilde{Q}\in {\cal F}^{(0)}_M((0,1]^d)}\inf_{\beta} \pl \mu^{k}(\tilde{Q}_{J,\beta})-\mu^{k}(\tilde{Q})\pl_{\infty}=O(C(M)r(d,J)^{k+1}).
\end{equation}
\end{lemma}
{\bf Proof: }
See Theorem \ref{mainkthordercdf} for the statement relying on ${\cal R}(d,J)$ satisfying both (\ref{Ra}) and (\ref{Rb}), while Theorem \ref{mainkthordercdfnew} establishes the stronger result (mostly relying on same proof) for ${\cal R}(d,J)$ only relying on (\ref{Rbm}).  $\Box$

This proves the following main  theorem.
\begin{theorem}\label{theoremuniformapprox}
Let $k\in \{1,2,\ldots,\}$ be given.
Recall the construction of set ${\cal R}^k(d,{\bf J})$ defining the approximation $Q^k_{{\bf J},\beta}=\sum_{(\bar{s}(k+1),u)\in {\cal R}^k(d,{\bf J})}\beta(u,\bar{s}(k+1))\phi_{\bar{s}(k+1),u}$ in terms of sets $({\cal R}(\bar{s}_{k+1},{\bf J}(\bar{s}(k+1)):\bar{s}(k+1))$ determined by  $({\cal R}(m,J):J)$ satisfying condition (\ref{Rb}) of Corollary \ref{defknots} for all $m=1,\ldots,d$. In addition, assume that for some $\delta>0$ and $M<\infty$,  $\delta J<{\bf J}(\bar{s}(k+1))<M J$ across all $\bar{s}(k+1)\in {\cal S}^{k+1}(d)$ with $\mid s_{k+1}\mid>0$.
We have
\[
\sup_{Q\in D^{(k)}_M([0,1]^d)}\inf_{\beta}\pl Q^k_{{\bf J},\beta}-Q\pl_{\infty}=O(C(M) r(d,J)^{k+1}).\]
\end{theorem}

\subsection{Global idea of overall proof  of uniform approximation error for $k$-th order primitive functions by $k$-th order splines.}
Detailed proofs are presented in Appendix \ref{AppendixD}.
Here we present the outline of the proof for the first order spline representation, presented by Theorem \ref{1storder}, while the general proof is a proof by induction: i.e., after having proved the result for $k=1$, assuming we have the result for $k$, prove it for $k+1$.
Therefore,  the proof for $k=1$ is the most essential part and is essentially repeated. 

Let $\tilde{Q}_{J,\beta}=\sum_{u\in {\cal R}(d,J)}\beta(u)\phi_u^0$ viewed as an approximation of $\tilde{Q}\in {\cal F}^{(0)}_M([0,1]^d)$, where ${\cal R}(d,J)$ satisfies the key $L^2$-approximation property (\ref{Rbm}) and is of size $J$, so that one can select $\beta$ so that $\pl \tilde{Q}_{J,\beta}-\tilde{Q}\pl_{\mu}=O(r(d,J))$.
Our goal is to show that $\inf_{\beta}\pl \mu(\tilde{Q}_{J,\beta})-\mu(\tilde{Q})\pl_{\infty}=O(r(d,J)^2)$. 
Define the criterion:
\[
R^1_J(\beta)\equiv  J^{-1}\sum_{v\in {\cal R}(d,J)}\left( \int_{(0,v]}\left(\sum_{u\in {\cal R}(d,J)}\beta(u)\phi_u^0-\tilde{Q}\right )d\mu \right)^2,\]
which thus equals $1/J \sum_{v\in {\cal R}(d,J)}(\mu(\tilde{Q}_{J,\beta})-\mu(\tilde{Q}))^2(v)$.
Let $\beta_J^*=\arg\min_{\beta}R^1_J(\beta)$ be the minimizer.
We claim that $\pl \mu(\tilde{Q}_{J,\beta_J^*})-\mu(\tilde{Q})\pl_{\infty}=O(r(d,J)^2)$, which would then complete the proof.
We have $H_v(\beta_J^*)\equiv \int_{(0,v]}\{ \tilde{Q}_{J,\beta_J^*(\tilde{Q})}-\tilde{Q}\} d\mu=0$
 for all $v\in {\cal R}(d,J)$, i.e., the least square minimizer sets the squared residuals
 $(\mu(\tilde{Q}_{J,\beta_J^*})-\mu(\tilde{Q}))(v)$ equal to zero for all $v\in {\cal R}(d,J)$.
Note that $R^1_J(\beta)=1/J \sum_{u\in {\cal R}(d,J)}\{H_v(\beta)\}^2$ so that the equations $H_v(\beta_J^*)=0$ are equations implied by the derivative equations $d/d\beta R^1_J(\beta)=0$ at $\beta=\beta_J^*$ solved by a minimum $\beta_J^*$ (details are presented in proof of Theorem \ref{1storder}).
Note that
\[
H_v(\beta_J^*)=\int \tilde{\phi}_v(y)(\tilde{Q}_{J,\beta_J^*}-\tilde{Q})  d\mu(y),\]
where $\tilde{\phi}_v(y)=I(y\leq v)$ is survivor analogue of the zero order spline $\phi_v(y)=I(y\geq v)$.
Since $H_v(\beta_J^*)=0$ for all $v\in {\cal R}(d,J)$, we have for all vectors $\alpha$
\[
0=\int \sum_{v\in {\cal R}(d,J)}\alpha(v)\tilde{\phi}_v(y)(\tilde{Q}_{J,\beta_J^*}-\tilde{Q}) d\mu(y).\]
Suppose we already showed that $\beta_J^*$ still satisfies $\pl \tilde{Q}_{J,\beta_J^*}-\tilde{Q}\pl_{\mu}=O(r(d,J))$, which is shown to be true in our proof of Theorem \ref{1storder} by viewing $\beta_J^*$ as the result of applying a universal steepest descent algorithm to an initial $\beta_J$ satisfying 
$\pl \tilde{Q}_{J,\beta_J}-\tilde{Q}\pl_{\mu}=O(r(d,J))$, and showing that the resulting update preserves the $L^2(\mu)$-rate of convergence (analogue to how the TMLE preserves the rate of convergence of an initial estimator).
 We can then carry out the  following proof. For any vector $\alpha_x$ we have
\begin{eqnarray*}
(\mu(\tilde{Q}_{J,\beta_J^*})-\mu(\tilde{Q}))(x)&=&\int_{(0,x]}(\tilde{Q}_{J,\beta_J^*}-\tilde{Q})) d\mu(y)\\
&=&\int \tilde{\phi}_x(y) (\tilde{Q}_{J,\beta_J^*}-\tilde{Q})(y)d\mu(y)\\
&=&\int (\tilde{\phi}_x(y)-\sum_{v\in {\cal R}(d,J)}\alpha_x(v)\tilde{\phi}_v(y) )    (\tilde{Q}_{J,\beta_J^*}-\tilde{Q})(y)       d\mu(y)\\
&\leq&  \pl \tilde{\phi}_x-\sum_{v\in {\cal R}(d,J)}\alpha_x(v)\tilde{\phi}_v\pl_{\mu} \pl \tilde{Q}_{J,\beta_J^*} -\tilde{Q}\pl_{\mu}.
\end{eqnarray*}
By the property (\ref{Rbm}) of ${\cal R}(d,J)$, we know that we can approximate the zero order splines $\phi_v$ by a linear combination  of the $J$ zero order splines in ${\cal R}(d,J)$ with an $L^2(\mu)$-norm error $\sim r(d,J)$,  uniformly in $v$. In this case we need this result for $\tilde{\phi}_v$. We could simply augment ${\cal R}(d,J)$ with the set of  $J$  knots that is tailored to approximate survivor functions instead of cumulative distribution functions to satisfy this property (\ref{Rb}) as well and Theorem \ref{1storder} relies on this. Subsequently, in Appendix \ref{AppendixD}  we generalize the proof to only relying on (\ref{Rbm})  (the result might not be surprising if one notes that $\tilde{\phi}_v$ can be written as a generalized difference of $\phi_u$ over corners of cube $(0,v]$). 
So then we have
\begin{eqnarray*}
(\mu(\tilde{Q}_{J,\beta_J^*})-\mu(\tilde{Q}))(x)&\leq&
\inf_{\alpha} \pl \tilde{\phi}_x-\sum_{v\in {\cal R}(d,J)}\alpha(v)\tilde{\phi}_v\pl_{\mu} \pl \tilde{Q}_{J,\beta_J^*}-\tilde{Q}\pl_{\mu}\\
&=& O(r(d,J)^{2}).
\end{eqnarray*}
This bound is uniformly in $x$ and uniformly in all $\tilde{Q}$ with $\pl\tilde{ f}\pl_{v}<1$. Thus, this proves $\pl \mu(\tilde{Q}_{J,\beta_J^*}) -\mu(\tilde{Q})\pl_{\infty}=O(r(d,J)^2)$. To conclude, we have constructed a $\mu_{J,\beta_J^*}(\tilde{Q})$ so that uniformly in $\tilde{Q}$ with $\pl \tilde{Q}\pl_{v}<1$:
1) $\pl \tilde{Q}_{J,\beta_J^*}-\tilde{Q} \pl_{\mu} =O(r(d,J))$; 2) $\pl \mu(\tilde{Q}_{J,\beta_J^*})-\mu(\tilde{Q})\pl_{\infty}=O(r(d,J)^{2})$. This then completes  the proof of the theorem.

\section{Uniform approximation error of  oracle MLE of $k$-th order spline working model w.r.t. true target function}\label{section6}
  
 Let ${\cal R}^k(d,{\bf J})$ be our index set  for the basis functions as defined in our construction of the finite dimensional $k$-th order spline approximation. 
Consider the corresponding working model $D^{(k)}({\cal R}(d,{\bf J}))=\{\sum_{j\in {\cal R}^k(d,{\bf J}) }\beta(j)\phi_j^k:\beta\}$, and $D_M^{(k)}({\cal R}^k(d,{\bf J}))$ is the subset that enforces the $L_1$-norm of $\beta$ to be bounded by $M$. Let $Q_{J,\beta}^k=\sum_{j\in {\cal R}^k(d,{\bf J})}\beta(j)\phi_j^k$. It was assumed that all components ${\bf J}(\bar{s}(k+1))$ are behaving as a constant times a $J$, and, for notational, convenience, we also express dependence on ${\bf J}$ only in terms of $J$. Theorem \ref{theoremuniformapprox} shows that for a given $Q\in D^{(k)}_{M}([0,1]^d)$, there exists a $\tilde{Q}_J\in D_M^{(k)}({\cal R}^k(d,{\bf J}))$ so that $\pl \tilde{Q}_J-Q\pl_{\infty}=O(C(M) r(d,J)^{k+1})$. 

Let $L(Q)$ be a loss function for $Q_0=\arg\min_fP_0L(Q)$ and $R_0(\beta)=P_0 L(Q_{J,\beta}^k)$. 
Let $Q_{0,J}^k=\arg\min_{Q\in D^{(k)}({\cal R}^k(d,{\bf J}))}P_0L(Q)$, and let $\beta_J$ be so that $Q_{0,J}^k=Q_{J,\beta_J}^k$. In particular, for the true $Q_0\in D_M^{(k)}([0,1]^d)$ we have that there exists a $\tilde{Q}_{0,J}\in D_M^{(k)}({\cal R}^k(d,{\bf J}))$ so that $\pl \tilde{Q}_{0,J}-Q_0\pl_{\infty}=O(C(M)r(d,J)^{k+1})$.  Let $\tilde{\beta}_J$ be so that $\tilde{Q}_{0,J}=Q^k_{J,\tilde{\beta}_J}$. By assumption,  we have $d_0(\tilde{Q}_{0,J},Q_0)=O(P_0 \{L(\tilde{Q}_{0,J})-L(Q_0)\}^2)$, so that it is a weak assumption on the loss function that the latter is $O(C^2(M)r(d,J)^{2(k+1)})$. Since $d_0(Q_{0,J}^k,Q_0)\leq d_0(\tilde{Q}_{0,J},Q_0)$, an immediate consequence of our uniform approximation result is given by: for $M$ chosen so that $\pl \beta_J\pl_1<M$
\[
d_0(Q_{0,J}^k,Q_0)=P_0L(Q_{0,J}^k)-P_0L(Q_0)=O(C(M)^2r(d,J)^{2(k+1)}).\]
The following theorem proves that this rate of convergence also applies to the supremum norm:  $\pl Q_{0,J}^k-Q_0\pl_{\infty}=O(C(M)r(d,J)^{k+1})$. 
The proof is given in Appendix \ref{AppendixI}. It is based on viewing $\beta_J$ as an update of $\tilde{\beta}_J$ along a universal steepest descent algorithm applied to $R_0(\beta)$, which will result in an update $\beta_J^*$ that solves the gradient $D^*_{\beta_J^*}=0$ of $R_0()$ at $\beta_J^*$, which, by assuming it has a  unique solution, implies that $\beta_J^*$ must equal $\beta_J$. We then show that this $P_0$-MLE update preserves the  sup-norm rate of convergence of the initial $Q_{J,\tilde{\beta}_J}^k$, due to $Q_{J,\tilde{\beta}_J}^k$ already within sup-norm distance $O(C(M)r(d,J)^{k+1})$ of optimizer $Q_0$.

In general this theorem shows that if a family of finite dimensional subspaces $D({\cal R}(J))$ indexed by $J$ of a function space $D$ has a certain sup-norm approximation of the true target function $Q_0\in D$, and the risk function $Q\rightarrow R_0(Q)$ of $Q_0=\arg\min_{Q\in D}R_0(Q)$ is well behaved, then the minimizer $Q_{0,J}=\arg\min_{Q\in D({\cal R}(J))}R_0(Q)$ of the risk function over this family will satisfy the same sup-norm approximation error (in rate). 

\begin{theorem}\label{theoremoraclemle}
Let $R_0(\beta)\equiv P_0L(Q_{J,\beta}^k)-P_0L(Q_0)$; $\beta_J=\arg\min_{\beta}R_0(\beta)$; 
$Q_{0,J}^k=Q_{J,\beta_J}^k$, and let $D^*_{\beta}=(\frac{d}{d\beta(u)}R_0(\beta):u\in {\cal R}^k(d,{\bf J}))$ be the gradient of $R_0(\beta)$ at $\beta$. Note  $D^*_{\beta_J}=0$. Let $\pl \beta\pl_J^2\equiv \sum_{j\in {\cal R}^k(d,{\bf J})}\beta(u)^2$.
Let $\tilde{\beta}_J$ be so that $\tilde{Q}_{0,J}=Q_{J,\tilde{\beta}_J}^k$ with $\pl \tilde{Q}_{0,J}-Q_0\pl_{\infty}=O(C(M)r(d,J)^{k+1})$.  Note that under a weak regularity condition,  $R_0(\tilde{\beta}_J)=O(C(M)^2 r(d,J)^{2(k+1)})$.

{\bf Assumptions:}
Assume that $D^*_{\beta}=0$ implies $\beta=\beta_J$, and  (e.g., using that $d_0(Q_{0,J},Q_0)=O(C(M)^2r(d,J)^{2(k+1)})$ and Cauchy-Schwarz inequality)
\begin{equation}\label{z22m}
\pl D^*_{\tilde{\beta}_J}\pl_J=O(J^{1/2} C(M)r(d,J)^{(k+1)}).
\end{equation}

Then, $\pl Q_{0,J}^k-Q_0\pl_{\infty}=O(C(M) r(d,J)^{k+1})$. 
This results holds uniformly in $P_0$ with $Q_0\in D^{(k)}_M([0,1]^d)$.
Thus, 
\begin{equation}\label{uniformoraclemle}
\sup_{\{P_0: Q(P_0)\in D^{(k)}_M([0,1]^d)\}}\pl Q_{J,\beta_J(P_0)}^k-Q(P_0)\pl_{\infty}=O(C(M)r(d,J)^{k+1}).\end{equation}

\end{theorem} 
\paragraph{Showing that assumption (\ref{z22m}) holds for squared error loss}
For example, suppose that $L(Q)=(Y-Q(X))^2$ so that $R_0(\beta)=P_0(Q_{J,\beta}-Q_{J,\beta_J})^2$.
Then, $R_0(\tilde{\beta}_J)=P_0(\tilde{Q}_{0,J}-Q_{0,J})^2=O(C(M)^2r(d,J)^{2(k+1)})$ due to $\pl \tilde{Q}_{0,J}-Q_{0,J}\pl_{P_0}\leq \pl \tilde{Q}_{0,J}-Q_0\pl_{P_0}+\pl Q_{0,J}-Q_0\pl_{P_0}$; $\pl \tilde{Q}_{0,J}-Q_0\pl_{P_0}=O(C(M)r(d,J)^{k+1})$ and $\pl Q_{0,J}-Q_0\pl_{P_0}=O(d_0^{1/2}(Q_{0,J},Q_0))=
O_P(C(M)r(d,J)^{k+1})$. 
Then, $\frac{d}{d\tilde{\beta}_J(u)}R_0(\tilde{\beta}_J)=2P_0 (Q_{J,\tilde{\beta}_J}-Q_{J,\beta_J})\frac{d}{d\tilde{\beta}_J(u)}Q_{J,\tilde{\beta}_J}$ so that it immediately also follows that $\max_{u\in {\cal R}(d,{\bf J})}\mid D^*_{\tilde{\beta}_J}(u)\mid =O(C(M)r(d,J)^{k+1})$. Thus, $\pl D^*_{\tilde{\beta}}\pl_J^2=O(J C(M)^2 r(d,J)^{2(k+1)})$ so that assumption (\ref{z22m}) indeed holds.

\section{Rates of convergence of sieve and HAL-MLEs towards oracle MLE w.r.t. loss-based dissimilarity}
\label{section7}
Analogue to Section 2, 
given our $k$-th order spline working model $Q_{{\bf J},\beta}^k=\sum_{j\in {\cal R}^k(d,{\bf J})}\beta_j\phi_j^k$ with $O(r(d,J)^{k+1})$-uniform approximation error for $D^{(k)}_M([0,1]^d)$, we can define  the $k$-th order spline sieve MLE and HAL-MLEs. Specifically, for a given ${\bf J}$, 
let  $Q_{{\bf J},\beta_n}^k$ be the MLE over the working model, while
$Q_{{\bf J},\beta_{n,C}}^k$ be the MLE over working model under constraint that $\pl \beta\pl_1\leq C$, while the relax lasso $Q_{{\bf J},\beta_{n,C}^r}^k$ refits the non-zero coefficients in $\beta_{n,C}$ with non-penalized MLE. For the sieve MLE we generally have to define a data adaptive selector of ${\bf J}_n$, while
for the two HAL-MLEs we might set ${\bf J}={\bf J}_{max,n}$, and define a data adaptive selector $C_n$ for $C$, although, one might also data adaptively select ${\bf J}_{max,n}$. 
Here we recommend to select ${\bf J}_{max,n}$ so that $O(r(d,J_{max,n})^{k+1})$ is balanced with the rate $(n/J_{max,n})^{-1/2}$, up till a $\log n$-factor allowing for some undersmoothing so that the bias is asymptotically negligible. 
Our results in this and next sections prove that the sieve MLE using this fixed model $D^{(k)}({\cal R}^k(d,{\bf J}_{max,n}))$ is already achieving the rate of convergence $n^{-k^*/(2k^*+1)}$ up till $\log n$-factors, both pointwise as well as w.r.t. square-root of loss-based dissimilarity.  Choosing this ${\bf J}_{max,n}$ large enough, the cross-validation can then focus on selecting ${\bf J}$ and $C$ of smaller models, knowing that it will do an asymptotically optimal job and thus will be at least as good as choosing (no regularization) $C=\infty$ and ${\bf J}_{max,n}$.

From now on our notation will often suppress the dependence on $k$, since it plays no role in the proof of asymptotic normality w.r.t. oracle MLE. 
For notational convenience, in this and following sections the realization of all  three estimators as functions of $P_n$ are denoted with $Q_n$, suppressing dependence on $k$, which can also be viewed as a random variable $Q_n=\hat{Q}(P_n)$, with $\hat{Q}$ a mapping from $P_n$ to function space $D^{(k)}([0,1]^d)$.  These estimators correspond with a  data dependent set ${\cal R}_n$ of indices of the non-zero coefficients. Let $D^{(k)}({\cal R}_n)=\{\sum_{j\in {\cal R}_n}\beta_j\phi_j: \beta\}$ be the corresponding data dependent working model, where we suppress dependence of $\phi_j^k$ on $k$. 
The $k$-th order spline sieve MLE and $k$-th order relax HAL-MLE solve exactly the scores $P_n \frac{d}{dQ_n}L(Q_n)(\phi_{j})=0$ for $j\in {\cal R}_n$ for this data dependent working model, while the HAL-MLE solves these score equations up till an understood approximation error $r_n(j_0)$ (see Lemma \ref{lemmascoreeqnhal} in the Appendix), which represents a second order remainder.  
Another difference is that the HAL-MLE is guaranteed to satisfy a uniformly bounded $k$-th order sectional variation norm, $\pl Q_n\pl_{v,k}^*\leq C_n$ (i.e., $L_1$-norm of $\beta$), while the other two estimators rely on the size $\sim J_n$ of ${\cal R}_n$ to be appropriately chosen so that (e.g.) the sectional variation norm is controlled as proven by Theorem \ref{theoremconvsectvarnorm}.

\subsection{Rate of convergence of HAL-MLE w.r.t loss based dissimilarity}
Our uniform approximation error $O(C(M)r(d,J)^{k+1})$ for $D^{(k)}_M([0,1]^d)$ by a $\sim J$-dimensional parametric model $D^{((k)}({\cal R}(d,{\bf J}))$ implies a uniform/sup-norm covering number, which then implies a bound on the uniform entropy integral \[
J_{\infty}(\delta,D^{(k)}_M([0,1]^d))\sim \delta^{(2k+1)/(2k+2)} (-\log \delta)^{(d_1+1)/2}.\]
This result is shown in Lemma \ref{lemmasupnormcoveringnumber}. 
Consider the $k$-th order spline HAL-MLE as an MLE over $D^{(k)}_M([0,1]^d)$. We can then carry out a general proof for its rate of convergence as in \citep{Bibaut&vanderLaan19},  just replacing their covering number for $D^{(0)}_M([0,1]^d)$  by this new one for $D^{(k)}_M([0,1]^d)$. This general proof only relies on an $L^2$-covering number which is implied by our stronger uniform covering number and gives the remarkable rate $n^{-k^*/(2k^*+1)}$ up till $\log n$-factors. These results  generalize to defining the HAL-MLE over a submodel $D^{(k)}({\cal R}(d,{\bf J}_{max,n}))$ of $D^{(k)}([0,1]^d)$ that has a uniform approximation error $\sim r(d,J_{max,n})^{k+1}$   of smaller order than the rate of convergence $(n/J_{max,n})^{-1/2}$ of the  oracle MLE of this initial model.

Recently, \citep{Kietal21} derived $L^2$-covering numbers for $L_1$-norm constrained linear combinations of higher order splines, which suffice for obtaining these rates of convergence for $D^{(k)}([0,1]^d)$, so that we could also use their results instead of our stronger uniform covering number results. 

\begin{theorem}\label{theoremlossbaseda}
Let  $Q_n=\arg\min_{Q\in D^{(k)}_{C_n}([0,1]^d)}P_n L(Q)$ be the $k$-th order spline HAL-MLE over  $D^{(k)}_{C_n}([0,1]^d)$.
Assume Assumption \ref{assumptionepsnet} of Lemma \ref{lemmasupnormcoveringnumber}. We have for some $m=m(d,k)<\infty$ that
\[
d_0(Q_n,Q_0)=O(n^{-2k^*/(2k^*+1)}(\log n)^m),\]
These results hold as long as $C_n=O((\log n)^a)$ for some $a<\infty$.
The same results apply to $Q_n=\arg\min_{Q\in D^{(k)}_{C_n}({\cal R}(d,{\bf J}_{0,n}))}P_0L(Q)$ when
$J_{0,n}$ is chosen so that $r(d,J_{0,n})^{k+1}/(n^{-k^*/(2k^*+1)}) =O(n^{-\delta})$ for some $\delta>0$.
\end{theorem}
{\bf Proof:} 
We can copy the proof for the rate of convergence for the zero-order HAL-MLE but now using the entropy integral $
J_{\infty}(\delta,D^{(k)}_M([0,1]^d))$.
In other words,  if $d_0(Q_n,Q_{0})\leq \delta^2$, then it follows that, ignoring $\log \delta$-factors 
\begin{eqnarray*}
d_0(Q_n,Q_{0})&\leq&-(P_n-P_0)(L(Q_n)-L(Q_{0}))\leq n^{-1/2}\sup_{ \pl g\pl_{P_0}\leq \delta, g\in D^{(k)}_C([0,1]^d)}\mid G_n(g)\mid \\
&=& O(n^{-1/2} \delta^{(2k+1)/(2k+2)}) ,\end{eqnarray*}
where we use the empirical process bound in terms of the entropy integral as stated in 
\citep{vanderVaart:Wellner11}.

We can start with $\delta_1=1$. Given  $\delta_k$ at step $k$,  using that the right-hand side gives a new bound for $d_0(Q_n,Q_0)$ in terms of $\delta_k$,  and that $\pl L(Q_n)-L(Q_0)\pl_{P_0}^2\sim d_0(Q_n,Q_0)$, we obtain a new bound  $\delta_{k+1}$. 
The resulting final $\delta$ (fixed point) satisfies $\delta= n^{-1/4}\delta^{(2k+1)/(4k+4)}$ resulting in  rate of convergence given by $\delta^2$. This yields $\delta^2=n^{-1/2}\delta^{(2k+1)/(2k+2)}$ and solving for $\delta=\delta_n$ gives the rate $d_0(Q_n,Q_{0,n})\sim \delta_n^2$ giving the stated rate $n^{-2k^*/(2k^*+1)}$ up till $\log n$-factors.  The other statement follows straightforwardly. 
$\Box$

\ \newline
$Q_n$ can also be viewed as an MLE over $D^{(k)}({\cal R}_n)$ and an estimator of its $P_0$-MLE $Q_{0,n}=\arg\min_{Q\in D^{(k)}({\cal R}_n)}P_0L(Q)$. Then the above proof, using the covering number of $D^{(k)}([0,1]^d)$ as a  conservative covering number for $D^{(k)}({\cal R}_n)$, gives the same rate for $d_0(Q_n,Q_{0,n})$.

\subsection{Rate of convergence in loss based dissimilarity for sieve and relax HAL-MLE, and HAL-MLE using fixed starting working model}
The above proof for HAL-MLE relies on knowing that $Q_n$ falls in the class $D^{(k)}_{C_n}([0,1]^d)$ whose covering number drives the rate of convergence. We will now establish this rate of convergence for the sieve-MLE and relax HAL-MLE without reference to our covering number result about $D^{(k)}_M([0,1]^d)$ but just utilizing our uniform approximation error result. 
 We first consider the  $J_{0,n}$-dimensional fixed parametric working model $D^{(k)}({\cal R}(d,{\bf J}_{0,n}))$, and establish the rate of convergence for the corresponding sieve MLE. We then embed the HAL-MLEs in this model by defining the HAL-MLE over this fixed working model, and use the cross-validation selector for the $L_1$-norm in the HAL-MLE and ${\bf J}$ for the sieve-MLE. These cross-validation selectors consider in particular the sieve-MLE for this fixed working model. Combined with the asymptotic equivalence of the cross-validation selector and the oracle selector, this then shows that the cross-validated HAL-MLEs and the cross-validated sieve MLEs achieve the same rate of convergence as this fixed sieve MLE, but,  in fact, they will converge as fast up till the constant as the oracle selected candidate estimator.  
\begin{theorem}\label{theoremlossbasedb}\ \nl
{\bf Sieve MLE for fixed sieve:}
Let ${\bf J}_{0,n}$ be a fixed sequence, and let  $Q_{{\bf J}_{0,n},n}=\arg\min_{Q\in D^{(k)}({\cal R}(d,{\bf J}_{0,n}))}P_n L(Q)$ be the sieve MLE over $D^{(k)}({\cal R}(d,{\bf J}_{0,n})$ estimating the oracle MLE $Q_{{\bf J}_{0,n},0}=\arg\min_{Q\in D^{(k)}({\cal R}(d,{\bf J}_{0,n}))}P_0L(Q)$.
Let $\beta_{{\bf J}_{0,n},n}$ and $\beta_{{\bf J}_{0,n},0}$ be  the coefficient vectors identifying $Q_{{\bf J}_{0,n},n}$ and $Q_{{\bf J}_{0,n},0}$. Let $B_{0,n}$ be a set of $\beta$-values so that  $\beta_{{\bf J}_{0,n},n},\beta_{{\bf J}_{0,n},0}\in B_{0,n}$ with probability tending to 1   and let 
${\cal F}_{0,n}=\{L(Q_{J_{0,n},\beta})-L(Q_{J_{0,n},\beta_{{\bf J}_{0,n},0}}):\beta\in B_{0,n}\}$. 
Assume that ${\bf J}_{0,n}$ is such that the entropy integral $J(\delta,{\cal F}_{{\bf J}_{0,n},n},L^2)$ for this $d_{0,n}\sim J_{0,n}$-dimensional set ${\cal F}_{0,n}$ can be conservatively bounded by $ d_{0,n}^{1/2} \delta\log \delta\sim J_{0,n}^{1/2}\delta \log \delta$. 
Then,
\[
d_{0}^{1/2}(Q_{{\bf J}_{0,n},n},Q_{{\bf J}_{0,n},0})=O_P ((J_{0,n}/n)^{1/2}\log n).\]
Moreover, suppose that ${\bf J}_{0,n}$ (say constant across components at value $J_{0,n}$)  is  chosen so that $(J_{0,n}/n)^{1/2}$ is balanced with $r(d,J_{0,n})^{k+1}$ up till $\log n$-factors. Then, $J_{0,n}\approx n^{1/(2k^*+1)}$ up till a  power of $\log n$-factor.  Then, for some $m<\infty$
\[
d_0(Q_{{\bf J}_{0,n},n},Q_0)=O_P(n^{-k^*/(2k^*+1)})(\log n)^m).\]

{\bf Cross-validated Sieve-MLE and  HAL-MLEs:}
Suppose that the HAL-MLE is an MLE over $D^{(k)}_C({\cal R}(d,{\bf J}_{0,n}))$ and $C$ is chosen with cross-validation for both the relax HAL-MLE and HAL-MLE, and that the sieve MLE selects ${\bf J}$ with a cross-validation selector that includes the candidate ${\bf J}_{0,n}$ (e.g., ${\bf J}\leq {\bf J}_{0,n}$).
Then, we have for all three estimators, for some $m<\infty$
\begin{equation}\label{l2ratehalmles}
d_0(Q_n,Q_0)=O(n^{-k^*/(2k^*+1)})(\log n)^m) .\end{equation}


\end{theorem}
{\bf Proof of Theorem \ref{theoremlossbasedb}:}
Let $Q_n=Q_{{\bf J}_{0,n},n}$ and $Q_{0,n}=Q_{{\bf J}_{0,n},0}$.
By assumption, the entropy integral $J(\delta,{\cal F}_{0,n},L^2)$ for a $d_{0,n}\sim J_{0,n}$-dimensional parametric model ${\cal F}_{0,n}$ can be conservatively bounded by $ d_{0,n}^{1/2} \delta\log \delta\sim J_{0,n}^{1/2}\delta \log \delta$. 
Recall that empirical process theory \citep{vanderVaart:Wellner11} teaches us (for the rates $\delta$ we apply it for) that $\sup_{Q\in {\cal F}_{0,n},\pl Q\pl_{P_0}<\delta}\mid G_n(Q)\mid =O_P(J(\delta,{\cal F}_{0,n},L^2))$, where $G_n(Q)=n^{1/2}(P_n-P_0)Q$ (ref).
Therefore, if $d_0(Q_n,Q_{0,n})\leq \delta^2$, then it follows that 
\[
\begin{array}{l}
d_0(Q_n,Q_{0,n})\leq -(P_n-P_0)(L(Q_n)-L(Q_{0,n}))\leq n^{-1/2}\sup_{ \pl g\pl_{P_0}\leq \delta}\mid G_n(g)\mid \\
=O_P(n^{-1/2} J_{0,n}^{1/2} \delta \log \delta ).\end{array}
\]
Recursively applying this result, starting with $\delta =1$, the resulting final $\delta$ (fixed point) satisfies $\delta= n^{-1/4}J_{0,n}^{1/4}\delta^{1/2} \log^{1/2}\delta $ and thus $\delta^{1/2}=(n/J_{0,n})^{-1/4}\log^{1/2}\delta$. 
So $\delta=(n/J_{0,n})^{-1/2}\log \delta$, and thus $\delta\sim (n/J_{0,n})^{-1/2}\log n$. The rate of convergence for $d_0(Q_n,Q_{0,n})$ is then $\delta^2$. This proves the first statement about the sieve MLE. 

The last statement about the HAL-MLEs is shown as follows. Due to asymptotic equivalence of cross-validation selector with oracle selector,  we know that the rate of convergence $d_0(Q_n,Q_0)$  is at least as fast as the one achieved with the working model $D^{(k)}({\cal R}(d,{\bf J}_{0,n}))$ which represents the unpenalized candidate for the HAL-MLEs and one of candidates for the sieve MLE. However, the latter  rate is  already achieving the optimal rate of convergence $\sim n^{-2k^*/(2k^*+1)}$ up till  $\log n$-factors. So that then proves that all these MLEs will converge at least at  the same rate w.r.t. loss based dissimilarity.
$\Box$

The condition on ${\bf J}_{0,n}$ essentially states that the $J_{0,n}$-dimensional class ${\cal F}_{0,n}$ consists of functions that are uniformly bounded from above by some $M<\infty$, since then the covering number would generally be polynomial. This is a significantly weaker condition than assuming that $\pl \beta_n\pl_1<M$ for some $M<\infty$ since that controls the $k$-th order sectional variation norm. Even for $k=1$, this condition will still allow $J_{0,n}$ to grow at a faster rate than the one that optimal balances the bias and standard error so that this is not a constraint obstructing the stated optimal rates of convergence for the sieve and relax HAL-MLE. In fact, we show in Appendix \ref{AppendixL} that for the desired/optimal rates of $J_{0,n}$ we have convergence of $Q_n-Q_{0,n}$ in sectional variation norm to zero, so that even the sectional variation norm will be uniformly bounded, let alone the supremum norm. 

\section{Pointwise asymptotic normality of the sieve and HAL-MLEs of regression function}\label{section8}
\subsection{Defining the estimation problem with respect to an oracle MLE w.r.t. fixed working model $D^{(k)}({\cal R}_{0,n})$.}
Recall ${\cal R}_n={\cal R}(P_n)$ is the set of non-zero coefficients in the fit $Q_n$ of $Q_0$. 
 Let ${\cal R}_{0,n}$ be a fixed sequence of sets of size $d_{0,n}\approx d_n$  approximating  ${\cal R}_n$ (i.e., ${\cal R}_{0,n}$ is independent of $P_n$) in the sense that the score equations $P_n S_{Q_n}(\phi_j)=0$ for $j\in {\cal R}_n$  solved by $Q_n$ imply $P_n S_{Q_n}(\phi_j)\approx 0$, $j\in {\cal R}_{0,n}$ are solved as well at a specified level of approximation. Let $J_{0,n}$ be the corresponding representative number of basis functions used to model each function in the additive model. 
 The choice of this set ${\cal R}_{0,n}$ will be discussed in a next subsection, where our recommended choice is ${\cal R}(P_n^{\#})$ for an independent sample $O_i^{\#}\sim_{iid} P_0$. Importantly, the choice should satisfy that $D^{(k)}({\cal R}_{0,n})$ uniformly approximates $Q_0$ and $Q_n$ at rate $O(C(M_n)r(d,J_{0,n})^{k+1})$. Thus, ${\cal R}_{0,n}$ could be one of our constructions ${\cal R}(d,{\bf J}_{0,n})$ targeting the whole or  submodel of $D^{(k)}([0,1]^d)$ that is known to contain $Q_0$.
 
 \begin{definition}
We say that $D^{(k)}({\cal R}_{0,n})$ satisfies the nonparametric uniform approximation condition if $\sup_{Q\in D^{(k)}_M([0,1]^d)}\inf_{g\in D^{(k)}({\cal R}_{0,n})}\pl g-Q\pl_{\infty}=O_P(C(M)r(d,J_{0,n})^{k+1})$.
We say that $D^{(k)}({\cal R}_{0,n})$ satisfies the  adaptive uniform approximation condition  if
$\sup_{Q\in \{Q_n,Q_0\}}\inf_{g\in D^{(k)}({\cal R}_{0,n})}\pl g-Q\pl_{\infty}=O_P(C(M_n)r(d,J_{0,n})^{k+1})$, where $M_n=\pl Q_n\pl_{v,k}^*$, the $L_1$-norm of its non-zero coefficients. 
The latter condition is weaker than the earlier, by allowing that ${\cal R}_{0,n}$ is tailored to the actual support of $Q_0$. 
\end{definition}

 This weakening of the uniform approximation assumption for $D^{(k)}({\cal R}_n)$ (and thereby $D^{(k)}({\cal R}_{0,n})$) so that it does not need to uniformly approximate $D^{(k)}([0,1]^d)$ is particularly important for the HAL-MLEs, since these are able to adapt to the support of $Q_0$, in particular, when using the cross-validation selector or a slightly undersmoothed version of it. However, it equally applies to the sieve MLE when selecting ${\bf J}$ with the cross-validation selector, thereby not necessarily enforcing that $\min_{\bar{s}(k+1)}{\bf J}(\bar{s}(k+1))/d_n>\delta$ for some $\delta>0$.

 Let  $O=(X,Y)$ and, let $\sigma^2_n$ be an estimator of $\sigma^2_{0,n}=E_0((Y-\tilde{Q}_{0,n})^2\mid X)$, where $\tilde{Q}_{0,n}=\arg\min_{Q\in D^{(k)}({\cal R}_{0,n})}P_0 (Y-Q(X))^2$.
 Let $L_{\sigma^2_n}(Q)(O)=\sigma^{-2}_n(Y-Q(X))^2$ be the weighted least squares loss function, so that $Q_0(X)=E_0(Y\mid X)$; $S_{Q_n,j}(X,Y)=\frac{d}{dQ_n}L(Q_n)(\phi_j)=\sigma^{-2}_n\phi_j(X)(Y-Q_n(X))$, $j\in {\cal R}_{0,n}$ and $j\in {\cal R}_n$.  
 We  denote $Q_{0,n}=Q_{0,{\cal R}_{0,n}}=\arg\min_{Q\in D^{(k)}({\cal R}_{0,n})}P_0 L_{\sigma^2_{0,n}}(Q)$ as the corresponding oracle (fixed) MLE, while  $Q_n$ is the empirical MLE over $D^{(k)}({\cal R}_n)$ or $D^{(k)}_{C_n}({\cal R}_n)$. 
 
 For a given $x_0$, we want to establish convergence in distribution of  $(n/d_{0,n})^{1/2}(Q_n(x_0)-Q_{0,n}(x_0))$ to a normal mean zero limit distribution. By simple stacking the linear expansion for $(Q_n-Q_{0,n})(x_0)$ of our pointwise convergence theorem below for a single $x_0$ across a given  set of points $(x_1,\ldots,x_k)$,  it will also follow  that
$((n/d_n)^{1/2}(Q_n(x_j)-Q_{0,n}(x_j)):j=1,\ldots,k)$ converges to a multivariate normal distribution with mean zero and certain covariance matrix. Then, by selecting ${\bf J}_{0,n}$ appropriately so that $r(d,J_{0,n})^{k+1}/(d_{0,n}/n)^{1/2}\rightarrow_p 0$ (i.e., bias negligible relative to standard error), given $\pl Q_{0,n}-Q_0\pl_{\infty}=O_P(r(d,J_{0,n})^{k+1})$, it follows that also $((n/d_{0,n})^{1/2}(Q_n-Q_0)(x_j):j=1,\ldots,k)$ converges to a multivariate normal mean zero limit distribution, which is our desired result.

We also cover logistic regression in which case $O=(X,Y)$, $Y\in \{0,1\}$, $L(Q)=-\log \{m_Q(X)^Y(1-m_Q(X))^{1-Y}\}$, where $m_Q(x)=1/(1+\exp(-Q(x)))$. 
 In the Appendix \ref{AppendixK} we also provide analogue theorems for  asymptotic normality for  unweighted least squares. 
 Our results will be generalized in the  Section \ref{section10} to arbitrary loss functions. Since the results  for likelihood behaving loss functions are more elegant, we emphasize our results   likelihood behaving loss functions, such as the weighted least squares loss function.



\subsection{Proof of asymptotic normality}

Let $dP^*_0(x)=\sigma^{-2}_{0,n}(x)dP_{X,0}(x)$, and let $L^2(P_0^*)$ be the corresponding Hilbert space of functions of $X$ with inner product $\langle f,g\rangle_{P_0^*}=P_0^* fg$. 
For simplicity, we assume that $\inf_{x\in [0,1]^d}\sigma^2_{0,n}(x)>\delta>0$ for some $\delta>0$, although for pointwise convergence at a particular $x_0$ one only needs to assume $\lim\inf \sigma^2_{0,n}(x_0)>0$.
We consider the weighted least squares loss function $L_{\sigma^2_n}(Q)(O)=\sigma^{-2}_n(X)(Y-Q(X))^2$. 
Recall  $Q_{0,n}=\arg\min_{Q\in D^{(k)}({\cal R}_{0,n})}P_0 L_{\sigma^2_{0,n}}(Q)$, 
the scores for $Q_n$ are given by $S_{Q_n,j}(O)=\phi_j(X)/\sigma^2_n(X) (Y-Q_n(X))$, and the scores for $Q_{0,n}$ are given by $S_{Q_{0,n},j}(O)=\phi_j(X)/\sigma^2_{0,n}(X)(Y-Q_{0,n}(X))$.

 Let $\{\phi_j^*:j\in {\cal R}_{0,n}\}$ be an orthonormal basis for $\{\phi_j:j\in {\cal R}_{0,n}\}$ in $L^2(P_0^*)$. This does not change the fixed working model $D^{(k)}({\cal R}_{0,n})$ and $Q_{0,n}$, since it only changes the parametrization of its elements. We use this orthonormal basis to linearize $(Q_n-Q_{0,n})(x)$  in a direct manner, without the need for an inverse of a matrix, which implicitly is carried out  by the orthonormalization of the basis, and, to establish an explicit elegant expression for its asymptotic variance.  
 
{\bf Solving the efficient influence curve equation of parameter $Q_{0,n}(x)$ in nonparametric model:}
 Since $\phi_j^*$ is a linear combination of $\{\phi_j: j\in {\cal R}_{0,n}\}$, and ${\cal R}_{0,n}$ is chosen so that score equations $r_n(j)=P_n \frac{d}{dQ_n}L_{\sigma^2_n}(Q_n)(\phi_j)\approx 0$, $j\in {\cal R}_{0,n}$,  the latter should also imply the  solving of score equations $r_n^*(j)\equiv P_n \frac{d}{dQ_n}L_{\sigma^2_n}(Q_n)(\phi_j^*)\approx 0$, $j\in {\cal R}_{0,n}$.
Define the term: 
\begin{equation}\label{tildernx}
\tilde{r}_n(x)\equiv \sum_{j\in {\cal R}_{0,n}}r_n^*(j)\phi_j^*(x).\end{equation}
As one increases $C_n$ in the HAL-MLEs, this term will reduce in size and can generally be controlled at the desired level. 
One can also represent  $\tilde{r}_n(x)$ as $P_n D_{{\cal R}_{0,n},x,\beta_n}$, where $D_{{\cal R}_{0,n},x,\beta}$ is the efficient influence curve of the statistical target parameter $P\rightarrow Q_{{\cal R}_{0,n},P}(x)=\sum_{j\in {\cal R}_{0,n}}\beta_j(P)\phi_j(x)$ with $\beta(P)=\arg\min_{\beta}P L(\sum_{j\in {\cal R}_{0,n}}\beta_j\phi_j)$ for the nonparametric model for $P_0$. So $\tilde{r}_n(x)\approx 0$ states that $Q_n$ needs to solve the efficient influence curve/score equation for nonparametric  statistical target parameter 
$Q_{{\cal R}_{0,n},P}(x)$ with $Q_{{\cal R}_{0,n},P}=\arg\min_{Q\in D^{(k)}({\cal R}_{0,n})}PL(Q)$.  

As starting point for our analysis of $(Q_n-Q_{0,n})(x)$, using that $P_0 \phi_j^*/\sigma^2_{0,n} (Y-Q_{0,n})=0$ and $P_n \sigma^{-2}_n \phi_j^*(Y-Q_n)=r_n^*(j)$,  we have
\[
P_0\phi_j^*/\sigma^2_n(Y-Q_n)-P_0\phi_j^*/\sigma^2_{0,n}(Y-Q_{0,n})=-(P_n-P_0)\phi_j^*/\sigma^2_n(Y-Q_n)+r_n^*(j).\]
We now write $P_0\phi_j^*/\sigma^2_n(Y-Q_n)=P_0 \phi_j^*/\sigma^2_{0,n}(Y-Q_n)+
P_0\phi_j^*(\sigma^{-2}_n-\sigma^{-2}_{0,n})(Y-Q_n)$ and put the second term to the other side of the equation. This yields:
\begin{eqnarray*}
P_0 \phi_j^*/\sigma^2_{0,n} (Y-{Q_n})-P_0\phi_j^*/\sigma^2_{0,n} (Y-{Q_{0,n}})&=&-(P_n-P_0)\phi_j^*/\sigma^2_n (Y-{Q_n})\\
&&+P_0 \phi_j^*(\sigma^{-2}_{0n}-\sigma^{-2}_{n})(Y-Q_n)+r_n^*(j).\end{eqnarray*}
Thus,
\[
P_0 \phi_j^*/\sigma^2_{0,n} ({Q_n}-{Q_{0,n}} )=(P_n-P_0)\phi_j^*/\sigma^2_n(Y-{Q_n})-P_0 \phi_j^*(\sigma^{-2}_{0n}-\sigma^{-2}_{n})(Y-Q_n)-r_n^*(j).\]
We view the second term $P_0 \phi_j^*(\sigma^{-2}_{0n}-\sigma^{-2}_{n})(Y-Q_n)$ as a second order remainder, and note that it equals
\begin{equation}\label{R1naphi}
R_{n,1}(\phi_j^*)\equiv -P_0 \phi_j^*(1/\sigma^2_n-1/\sigma^2_{0,n})(Q_n-Q_{0}).\end{equation}
Then, 
\[
P_0^* \phi_j^* (Q_n-Q_{0,n})=(P_n-P_0)\phi_j^*/\sigma^2_n(Y-Q_n)-r_n^*(j)+R_{n,1}(\phi_j^*).\]
Let $\tilde{Q}_n=\Pi_{J_{0,n}}(Q_n)\equiv \arg\min_{Q\in D^{(k)}({\cal R}_{0,n})}P_0^*(f-Q_n)^2$ be
the projection of  $Q_n$ on the fixed working model $D^{(k)}({\cal R}_{0,n})$ in $L^2(P_0^*)$.
More generally, we could have defined $\Pi_{J_{0,n}}(Q_n)\in D^{(k)}({\cal R}_{0,n})$ as an element for which $\pl Q_n-\tilde{Q}_n\pl_{\infty}=O_P(C(M_n)r(d,J_{0,n})^{k+1})$, where, by our adaptive uniform approximation assumption on $D^{(k)}({\cal R}_{0,n})$ such an element exists. With our projection definition we need to apply Theorem \ref{theoremoraclemle} to actually prove that $\pl \tilde{Q}_n-Q_n\pl_{\infty}=O_P(C(M_n)r(d,J_{0,n})^{k+1})$, which will be done below.


We decompose $Q_n-Q_{0,n}=Q_n-\tilde{Q}_n+\tilde{Q}_n-Q_{0,n}$. Let
$R_{n,2}(\phi_j^*)=-P_0^*\phi_j^*(Q_n-\tilde{Q}_n)$. 
So we have
\[
P_0^*\phi_j^*(\tilde{Q}_n-Q_{0,n})=(P_n-P_0)\phi_j^*/\sigma^2_n(Y-Q_n)-r_n^*(j)+(R_{n,1}+R_{n,2})(\phi_j^*).\]

Since $\{\phi_j^*:j\in {\cal R}_{0,n}\}$, is an orthonormal basis of the linear subspace $D^{(k)}({\cal R}_{0,n})$ of $L^2(P_0^*)$, and $\tilde{Q}_n-Q_{0,n}\in D^{(k)}({\cal R}_{0,n})$, we have
\begin{eqnarray*}
(\tilde{Q}_n-Q_{0,n})(x)&=& \sum_{j\in {\cal R}_{0,n}}\{P_0^*(\tilde{Q}_n-Q_n) \phi_j^* \}\phi_j^*(x)\\
&=&\sum_{j\in {\cal R}_{0,n}}(P_n-P_0)\phi_j^*/\sigma^2_n(Y-{Q_n}) \phi_j^*(x)-\sum_{j\in {\cal R}_{0,n}}r_n^*(j) \phi_j^*(x)\\
&& +\sum_{j\in {\cal R}_{0,n}}(R_{n,1}+R_{n,2})(\phi_j^*)\phi_j^*(x)\\
&=&\sum_{j\in {\cal R}_{0,n}}(P_n-P_0)\phi_j^*/\sigma^2_n(Y-{Q_n}) \phi_j^*(x)-
\tilde{r}_n(x)\\
&&+\sum_{j\in {\cal R}_{0,n}}(R_{n,1}+R_{n,2})(\phi_j^*)\phi_j^*(x).
\end{eqnarray*}
Let's first consider the term $R_{n,2}(x)\equiv \sum_{j\in {\cal R}_{0,n}}R_{n,2}(\phi_j^*)\phi_j^*(x)$.
Note that this term equals $\sum_{j\in {\cal R}_{0,n}}P_0^*\phi_j^*(Q_n-\tilde{Q}_n)\phi_j^*(x)$ and thus equals the projection of $(Q_n-\tilde{Q}_n)$ onto $D^{(k)}({\cal R}_{0,n})$. So, by definition of $\tilde{Q}_n$,  this equals 
$\Pi_{J_{0,n}}(Q_n)-\tilde{Q}_n=0$. So $R_{n,2}(x)=0$.

The remaining remainder term is thereby $R_{n,1}(x)\equiv  \sum_{j\in {\cal R}_{0,n}}R_{n,1}(\phi_j^*)\phi_j^*(x)$. 
We now note that
\begin{eqnarray}
R_{n,1}(x)&= &-\sum_{j\in {\cal R}_{0,n}}\{ P_0(\sigma^{-2}_n-\sigma^{-2}_{0,n})(Q_n-Q_{0,n}) \phi_j^*\}\phi_j^*(x)\label{R1nx}\\
&=&-\Pi_{J_{0,n}}( e_n)\nonumber,
\end{eqnarray}
where
\[
e_n\equiv (\sigma^{-2}_n-\sigma^{-2}_{0,n})(Q_n-Q_{0,n}).\]
Therefore, the $L^2(P_0^*)$-norm of $R_{n,1}$ is bounded by the $L^2(P_0^*)$-norm of $e_n$.
By Cauchy-Schwarz inequality this is bounded by constant times $\pl Q_n-Q_{0,n}\pl_{P_0^*}\pl \sigma^2_n-\sigma^2_{0,n}\pl_{P_0^*}$ so that this will converge at rate $(d_{0,n}/n)^{1/2} \pl \sigma^2_n-\sigma^2_{0,n}\pl_{P_0^*}$. Therefore, it follows that  $\pl R_{n,1}\pl_{P_0}^*=o_P((d_{0,n}/n)^{1/2})$ even when only assuming consistency of $\sigma^2_n$ as estimator of $\sigma^2_{0,n}$ in $L^2(P_0^*)$-norm (i.e., no need for an actual rate). Based on this $L^2(P_0^*)$-rate, we expect that also  $R_{n,1}(x)=o_P((d_{0,n}/n)^{1/2})$ only relies on consistency of $\sigma^2_n$. To formally understand $\pl R_{n,1}\pl_{\infty}$, we can directly and conservatively bound $R_{n,1}(x)$ by $\sim d_{0,n}\pl Q_n-Q_{0,n}\pl_{P_0}^*\pl \sigma^2_n-\sigma^2_{0,n}\pl_{P_0^*}$. This shows that   $\sup_x\mid R_{n,1}(x)\mid =o_P((n/d_{0,n})^{-1/2})$ if 
$d_{0,n}\pl \sigma^2_n-\sigma^2_{0,n}\pl_{P_0^*}=o_P(1)$. For example, if $Q_0,E_0(Y^2\mid X)$ are elements of $D^{(1)}_M([0,1]^d)$, then we know that the first order spline HAL-MLE $\sigma^2_n$ achieve $\pl \sigma^2_n-\sigma^2_{0,n}\pl_{P_0^*}=O_P(n^{-2/5})$ up till $\log n$-factors, in which case, this requirement is satisfied if $d_{0,n} n^{-2/5}=O_P(n^{-\delta})$ for some $\delta>0$. The latter  holds for the optimal rate $d_{0,n}\sim n^{1/5}$ (up till $\log n$-factors), as well as under high levels of undersmoothing (not needed).

Thus, we have \[
(\tilde{Q}_n-Q_{0,n})(x)= \sum_{j\in {\cal R}_{0,n}}(P_n-P_0)\phi_j^*/\sigma^2_n(Y-{Q_n}) \phi_j^*(x)-
\tilde{r}_n(x)+R_{n,1}(x).\]

We now write the left-hand side as $\tilde{Q}_n-Q_n+(Q_n-Q_{0,n})$. 
By definition of $\tilde{Q}_n$ we have $\pl \tilde{Q}_n-Q_n\pl_{\infty}=O_P(C(M_n)r(d,J_{0,n})^{k+1})$.
So we then have
\begin{eqnarray*}
(Q_n-Q_{0,n})(x)&=& \sum_{j\in {\cal R}_{0,n}}(P_n-P_0)\phi_j^*/\sigma^2_n(Y-{Q_n}) \phi_j^*(x)-
\tilde{r}_n(x)+R_{n,1}(x)\\
&&+O_P(C(M_n)r(d,J_{0,n})^{k+1}).\end{eqnarray*}

Let 
\begin{equation}\label{Enx}
-E_n(x)\equiv \sum_{j\in {\cal R}_{0,n}}(P_n-P_0)\phi_j^*\left\{1/\sigma^2_n(Y-{Q_n}) -1/\sigma^2_{0,n}(Y-Q_{0,n}) \right\} \phi_j^*(x).\end{equation}
Then, we have
\begin{eqnarray*}
(Q_n-Q_{0,n})(x)&=& (P_n-P_0) \sum_{j\in {\cal R}_{0,n}}\phi_j^*/\sigma^2_{0,n}(Y-{Q_{0,n}}) \phi_j^*(x)-E_n(x)\\
&&-
\tilde{r}_n(x)+R_{n,1}(x) +O_P(C(M_n)r(d,J_{0,n})^{k+1}).\end{eqnarray*}

Multiply both sides with $(n/d_{0,n})^{1/2}$. 
Then, the leading term is a sum $1/n^{1/2}\sum_i D_{Q_{0,n},x}(O_i)$ of independent mean zero random variables \begin{equation}\label{Df0nx}
D_{Q_{0,n},x}(O)\equiv d_{0,n}^{-1/2}\sum_{j\in {\cal R}_{0,n}}\phi_j^*/\sigma^2_{0,n}(Y-{Q_{0,n}}) \phi_j^*(x).\end{equation}
The variance  of $D_{Q_{0,n},x}(O)$ is given by
\[
\tilde{\sigma}^2_{0,n}(x)=\frac{1}{d_{0,n}}\sum_{j_1,j_2\in {\cal R}_{0,n}}P_0\{ \phi_{j_1}^*\phi_{j_2}^*/\sigma^2_{0,n}\}\phi_{j_1}^*(x)\phi_{j_2}^*(x).\]
 Due to 
the orthonormality of $\{\phi_j^*: j\in {\cal R}_{0,n}\}$ in $L^2(P_0^*)$  with $dP_0^*=\sigma^{-2}_{0,n}dP_0$, this variance simplifies to
\begin{equation}
\tilde{\sigma}^2_{0,n}(x)=
\frac{1}{d_{0,n}}\sum_{j\in {\cal R}_{0,n}}\{\phi_{j}^*(x)\}^2.\end{equation}
 Lemma \ref{lemmaexplicit} shows that $\tilde{\sigma}^2_{0,n}=O(1)$, and provides an explicit expression for this variance. In particular, it shows that if $D^{(k)}({\cal R}_{0,n})$ approximates $D^{(k)}([0,1]^d)$ w.r.t. a supremum norm around $x$ at rate $r(d,J_{0,n})^{k+1})$,  then this expression $\tilde{\sigma}^2_{0,n}(x)$ approximates $1/p_0^*(x)$. The latter then implies that $\tilde{\sigma}^2_{0,n}(x)\rightarrow_p \sigma^2_0(x)/p_0(x)$, where $\sigma^2_0(x)=E_0( (Y-Q_0)^2\mid X=x)$.
 If the model $D^{(k)}({\cal R}_{0,n})$ is aimed to approximate a real subset of $D^{(k)}([0,1]^d)$ (like one of our proposed submodels $D^{(k)}({\cal R}^{k,*}(d))$), thereby allowing for extrapolation (e.g., an additive model), then this asymptotic variance would be smaller so that $\sigma^2_0(x)/p_0(x)$ represents an upper bound for the actual asymptotic variance.

{\bf Proving that the projection $\tilde{Q}_n=\Pi_{J_{0,n}}(Q_n)$  approximates $Q_n$ uniformly at rate $O(C(M_n)r(d,J_{0,n})^{k+1})$:}
By assumption on ${\cal R}_{0,n}$ there exists a $\tilde{g}_n\in D^{(k)}({\cal R}_{0,n})$ so that $\pl Q_n-\tilde{g}_n\pl_{\infty}=O_P(C(M_n)r(d,J_{0,n})^{k+1})$,
 where $M_n$ is a bound on $\pl Q_n\pl_{v,k}^*$.
  Let $Q_{J_{0,n},\beta}=\sum_{j\in {\cal R}_{0,n}}\beta(j)\phi_j^*$ and $\tilde{\beta}_n$ identifies this $\tilde{g}_{n}=Q_{J_{0,n},\tilde{\beta}_n}$.
We have that $\tilde{Q}_n=\Pi_{J_{0,n}}(Q_n)=\arg\min_{g\in D^{(k)}({\cal R}_{0,n})}P_0^*(Q_n-g)^2$. We can apply Theorem \ref{theoremoraclemle} with $R_0(\beta)=P_0^*(Q_{J_{0,n},\beta}-Q_n)^2$; $\beta_n=\arg\min_{\beta}R_0(\beta)$; $\Pi_{J_{0,n}}(Q_{n})=Q_{J_{0,n},\beta_n}$, and $\tilde{g_n}=Q_{J_{0,n},\tilde{\beta}_n}$.  The conditions of Theorem \ref{theoremoraclemle} are  satisfied as demonstrated with the least squares loss function under Theorem \ref{theoremoraclemle}. Applying Theorem \ref{theoremoraclemle} proves that the uniform approximation error of $\tilde{g}_{n}$ w.r.t approximating $Q_n$  carries over to the minimizer $\tilde{Q}_n=Q_{J_{0,n},\beta_n}$ of the squared error risk $R_0(\beta)$. Specifically, it proves that if $\pl Q_n\pl_{v,k}^*<M$,  then $\pl \Pi_{J_{0,n}}(Q_n)-Q_n\pl_{\infty}=O_P(C(M)r(d,J_{0,n})^{k+1})$.



\subsection{Asymptotic normality theorem for weighted least squares regression}
Consider/recall the following definitions:
\begin{eqnarray}
\tilde{r}_n(x)&\equiv& \sum_{j\in {\cal R}_{0,n}}\{P_n \phi_j^*/\sigma^2_n(Y-Q_n)\}\phi_j^*(x)\label{tildernxf}\\
R_{n,1}(\phi_j^*)&\equiv& -P_0 \phi_j^*(1/\sigma^2_n-1/\sigma^2_{0,n})(Q_n-Q_{0})\nonumber\\
R_{n,1}(x)&\equiv&  \sum_{j\in {\cal R}_{0,n}}R_{n,1}(\phi_j^*)\phi_j^*(x)\nonumber\\
&=&-\Pi_{J_{0,n}}( e_n)\label{R1nxf}\\
e_n&\equiv& (\sigma^{-2}_n-\sigma^{-2}_{0,n})(Q_n-Q_{0,n})\nonumber \\
E_n(x)&\equiv& -\sum_{j\in {\cal R}_{0,n}}(P_n-P_0)\phi_j^*\left\{\sigma^{-2}_n(Y-{Q_n}) -\sigma^{-2}_{0,n}(Y-Q_{0,n}) \right\} \phi_j^*(x)\label{Enxf}\\
D_{Q_{0,n},x}&\equiv&
d_{0,n}^{-1/2}\sum_{j\in {\cal R}_{0,n}}\phi_j^*/\sigma^2_{0,n}(Y-{Q_{0,n}}) \phi_j^*(x).\label{Df0nxf}
\end{eqnarray}

Let ${\cal R}_{0,n}={\cal R}(P_n^{\#})$ the independent copy of ${\cal R}_n={\cal R}(P_n)$ with $P_n^{\#}$ empirical measure of independent (from $P_n$) i.i.d. sample from $P_0$.
\newline
{\bf Assumption A0:}
\begin{equation}\label{A0}
\sup_{Q\in \{Q_n,Q_0\}}\inf_{g\in D^{(k)}({\cal R}_{n})}\pl Q-g\pl_{\infty}=O_P(C(M_n)r(d,J_{0,n})^{k+1}),\end{equation}
 where $M_n=\pl Q_n\pl_{v,k}^*$, the $L_1$-norm of its non-zero coefficients. \newline
{\bf Assumption A1:}
\begin{equation}\label{undersmoothm}
\tilde{r}_n(x)=o_P((n/d_{0,n})^{-1/2}).
\end{equation}
{\bf Assumption A2:}
\begin{equation}\label{A2}
R_{n,1}(x)=o_P((n/d_{0,n})^{-1/2}).\end{equation}
{\bf Assumption A3:}\begin{equation}\label{keyconditiona1m}
E_n(x)=o_P((n/d_{0,n})^{-1/2}).\end{equation}
{\bf Assumption A4:}
\begin{equation}\label{neglbias}
C(M_n)r(d,J_{0,n})^{k+1}=o((n/d_{0,n})^{-1/2}).
\end{equation}
The following lemma provides a sufficient condition for Assumption A0. It states that one wants ${\cal R}_n$ to be such that $D^{(k)}({\cal R}_n)$ uniformly approximates a space $D^{(k)}({\cal R}_0)$ that includes $Q_0$, and with probability tending to 1 also includes $Q_n$.

\begin{lemma}\label{lemmasufficientA0}
{\bf Assumption A0*:}
Let ${\cal R}_0$ be a support set identifying a set of basis functions indexed by $(\bar{s}(k+1),u)$ in our $k$-th order spline representation  so that for a given $M<\infty$ $\{Q_n,Q_0\}\subset D^{(k)}({\cal R}_0)$ with probability tending to 1. 
Assume that ${\cal R}_n={\cal R}(P_n)$ satisfies $\max_{g\in D^{(k)}_M({\cal R}_0)} \min_{Q\in D^{(k)}({\cal R}_{n})}\pl g-Q\pl_{\infty}=O_P(C(M)r(d,J_{n})^{k+1})$. 

Then Assumption A0 holds.
\end{lemma}
\paragraph{Discussion of Assumption A0*:}
For example, for the sieve MLE, we can define ${\cal R}_n={\cal R}^k(d,{\bf J}_n)$ with ${\bf J}_n$ a cross-validation selector over all ${\bf J}$ satisfying $\min_{\bar{s}(k+1)}{\bf J}(\bar{s}(k+1))>\delta d_n$ for some $\delta>0$, which case Assumption A0* holds with ${\cal R}_0$ the maximal set ${\cal R}^k(d)$ so that $D^{(k)}({\cal R}_0)=D^{(k)}([0,1]^d)$. This trivially generalizes to submodels $D^{(k)}({\cal R}^{k,*}(d))$ known to hold by our statistical model in which case ${\cal R}_0$ just reduces to this set ${\cal R}^{k,*}(d)$. If ${\bf J}_n$ is a cross-validation selector that allows for certain $\bar{s}(k+1)$ (allowed by the model) ${\bf J}_n(\bar{s}(k+1))/d_n\rightarrow_p 0$, then Assumption A0* might still hold for a restricted ${\cal R}_0$ since, beyond capturing $Q_n$ with $D^{(k)}({\cal R}_0)$ with probability tending to 1,  we only need to approximate the true $Q_0$ and the cross-validation selector is asymptotically equivalent with the oracle selector. One natural candidate for ${\cal R}_0$ is actual support of the true $Q_0$ so that $Q_0\in D^{(k)}({\cal R}_0)$: this set includes all $\bar{s}(k+1)$ for which $Q^{(k)}_{0,\bar{s}(k+1)}\not =0$, and, for each such $\bar{s}(k+1)$ with $\mid s_{k+1}\mid >0$, it includes the $u$ for which  $Q^{(k)}_{0,\bar{s}(k+1)}(du)\not =0$, and one might always include the parametric form basis functions $\bar{\phi}_{\bar{s}(k+1)}$ with $\mid s_{k+1}\mid =0$.  With this definition we have $Q_0\in D^{(k)}({\cal R}_0)$. However, we also need that $Q_n\in D^{(k)}({\cal R}_0)$ with probability tending to 1. We know an oracle selector ${\bf J}^*_{0,n}\equiv \arg\min_{{\bf J}}\sum_v P_0 L(\hat{Q}_{{\bf J}}(P_{n,v}))$ of ${\bf J}$ would aim to select ${\bf J}$ so that it avoids any of the unnecessary basis functions outside support of $Q_0$, but nonetheless, it might select some with asymptotically negligible non-zero coefficients. This suggests that the choice ${\cal R}_0$ might have to include a small buffer of extra basis functions outside the support of $Q_0$. On the other hand, we could consider a ''truncated'' version $Q_{n,tr}\in D^{(k)}({\cal R}_0)$ of $Q_n$, defined by setting non-zero coefficients equal to zero if they concern basis functions outside the support of $Q_0$,  and assume that the difference between $Q_{n,tr}$ and $Q_n$ is negligible so that the analysis can focus on $Q_{n,tr}-Q_0$. This is beyond the scope of this article, but it indicates if $Q_0$ is highly restricted in its support, then there will be highly restricted ${\cal R}_0$, possibly as small as the actual support  of $Q_0$ defined above,  that will make assumption A0* true.

\begin{theorem}\label{thasnormal}\ \nl
Consider the $k$-th order sieve MLE, HAL-MLE or relax HAL-MLE $Q_n$.  Let $k^*=k+1$.
Let ${\cal R}_n$ be the set of $d_n$ non-zero coefficients in $Q_n=\sum_{j\in {\cal R}_n}\beta_n(j)\phi_j$. Assume $\inf_{x\in [0,1]^d}\sigma^{2}_{0,n}(x)>0$. 
Recall the definitions of $\tilde{r}_n(x)$ (\ref{tildernxf}), $R_{n,1}(x)$ (\ref{R1nxf}), $E_n(x)$ (\ref{Enxf}), $D_{Q_{0,n},x}(O)$ (\ref{Df0nxf}) above. 
Assume assumptions A0,A1,A2, A3 and A4.

Due to Assumption A0, we have the following expansion:
\begin{eqnarray}
(n/d_{0,n})^{1/2}(Q_n-Q_{0,n})(x)&=& n^{1/2}(P_n-P_0) d_{0,n}^{-1/2}\sum_{j\in {\cal R}_{0,n}}\phi_j^*/\sigma^2_{0,n}(Y-{Q_{0,n}}) \phi_j^*(x)\nonumber\\
&&-(n/d_{0,n})^{1/2}E_n(x) -
(n/d_{0,n})^{1/2}\tilde{r}_n(x)
\nonumber \\
&&\hspace*{-3cm}+(n/d_{0,n})^{1/2}R_{n,1}(x) +(n/d_{0,n})^{1/2}O_P(C(M_n)r(d,J_{0,n})^{k+1}),\label{keyexpansion}\end{eqnarray}
where $M_n=\pl Q_n\pl_{v,k}^*$, which equals the $L_1$-norm  of the vector of non-zero coefficients in $Q_n$.

Due to Assumptions A0,A1,A2,A3 and A4 we have
\begin{eqnarray}
\frac{(Q_n-Q_{0,n})(x)}{(n/d_{0,n})^{-1/2}}= n^{1/2}(P_n-P_0) d_{0,n}^{-1/2}\sum_{j\in {\cal R}_{0,n}}\phi_j^*/\sigma^2_{0,n}(Y-{Q_{0,n}}) \phi_j^*(x)+o_P(1). \label{keyexpansion} 
\end{eqnarray}

For a given $x$, the leading term is now  a sum of independent mean zero random variables $D_{Q_{0,n},x}(O)$ (\ref{Df0nx}) so that the central limit theorem can be applied. 
The variance is given by:
\begin{equation}\label{sigma2n}
\tilde{\sigma}^2_{0,n}(x)=
\frac{1}{d_{0,n}}\sum_{j\in {\cal R}_{0,n}}\{\phi_{j}^*(x)\}^2.\end{equation}

{\bf Conclusion:}
We have
\[
\tilde{\sigma}^{-1}_{0,n} (n/d_{0,n})^{1/2}({Q_n}-{Q_{0,n}})(x)\Rightarrow_d N(0,1),\]
and
\[
\tilde{\sigma}^{-1}_{0,n}(n/d_{0,n})^{1/2}(Q_n-Q_0)(x)\Rightarrow_d N(0,1).\]
Assuming $M_n=O(\log^{m}n)$ for some $m<\infty$, by choosing $J_{0,n}=n^{1/(2k^*+1)}\log^{m_1} n$ for some  $m_1<\infty$,  it follows that for  some  $m<\infty$, \[
\mid Q_n-Q_0\mid(x)=O_P(n^{-k^*/(2k^*+1)}\log^m n),\]
while Assumption A4 holds too. 
\end{theorem}
Assumption A4 is controlled by the user by selecting $d_n$ so that $C(M_n)r(d,d_n)^{k+1}=o((n/d_n)^{-1/2})$. We also established sufficient conditions for conditions A1,A2, A3.
In particular, sufficient conditions for A3 follow from Lemma \ref{lemmasuffEnx}, while a sufficient condition for A2 stated in next lemma was proved in our proof of Theorem \ref{thasnormal}. 
Sufficient conditions for A1 were established by Lemma  \ref{lemmatildernx}.
We summarize these sufficient conditions in the next lemma. 
\begin{lemma}
\begin{itemize}
\item
If $d_{0,n}\pl \sigma^2_n-\sigma^2_{0,n}\pl_{P_0^*}=o_P(1)$, then $R_{n,1}(x)=o_P((n/d_{0,n})^{-1/2})$ (i.e., A2 holds).
\item If $\sigma^2_n,Q_n,Q_{0,n},\sigma^2_0\in D^{(0)}_M([0,1]^d)$ and
$\pl \sigma^2_n-\sigma^2_0\pl_{P_0}=O_P(n^{-1/3}(\log n)^{d/2})$ (i.e., if $\sigma^2_n$ is an HAL-MLE of $\sigma^2_0$), and  $d_{0,n}=O_P(n^{1/3-\delta})$ for some $\delta>0$, then 
$E_n(x)=o_P((n/d_{0,n})^{-1/2})$ (i.e., A3 holds).
\item 
So using this HAL-MLE $\sigma^2_n$ guarantees that both $R_{1n}(x)$ and $E_n(x)$ are $o_P((n/d_{0,n})^{-1/2})$ if $d_{0,n}=O_P(n^{1/3-\delta})$ for some $\delta>0$. 
\item
Recall $S_j^*(Q_n)=\phi_j^*(Y-Q_n(X))$, $j\in {\cal R}_{0,n}$. Let $\tilde{S}_j(Q_n)$ be a linear combination of scores $S_j(Q_n)$ solved by $Q_n$ (i.e, $P_n S_j(Q_n)=0$) so that $P_n \tilde{S}_j(Q_n)=0$ approximating $S_j^*(Q_n)$ for $j\in {\cal R}_{0,n}$. For example, for relax HAL-MLE and sieve MLE, we have $S_j(Q_n)=S_{Q_n}(\phi_j)$, $j\in {\cal R}_n$.  Here $S_j^*(Q_n)-\tilde{S}_j(Q_n)=(\phi_j^*-\tilde{\phi}_j)(Y-Q_n(X))$, $j\in {\cal R}_{0,n}$. Lemma  \ref{lemmatildernx} shows  \begin{equation}\label{tildernxnice}
\tilde{r}_n(x)=\sum_{j\in {\cal R}_{0,n}}(P_n-\tilde{P}_n)\{S_j^*(Q_n)-\tilde{S}_j(Q_n)\} \phi_j^*(x).\end{equation}
Let 
$\delta_n\equiv  \pl d_{0,n}^{-1}\sum_{j\in {\cal R}_{0,n}}\{S_j(Q_n)-\tilde{S}_j(Q_n)\}\pl_{\mu}$.
If $d_{0,n}^{1/2}\delta_n^{(2k+1)/(2k+2)}=O_P(n^{-\delta})$ for some $\delta>0$, then (\ref{undersmoothm}) holds (i.e., A1 holds).
\item Also recall the sufficient condition for Assumption A0 as stated in Lemma \ref{lemmasufficientA0}.
\end{itemize}
\end{lemma}
So regarding assumption A1 we might have to obtain some rate for $\delta_n$. For the HAL-MLE this is more subtle, but for the relax HAL-MLE and sieve-MLE it is easily verified since we know $P_n S_{Q_n}(\phi_j)=0$ (exact) for all $j\in {\cal R}_n$, and we know the approximation error  of the linear span of these $\phi_j$'s.

{\bf Explicit asymptotic variance under nonparametric uniform approximation condition:}
Lemma \ref{lemmaexplicit} below shows that
 \[
 \tilde{\sigma}^2_{0,n}(x)=\{p_0^*(x)\}^{-1}  \Pi_{J_{0,n}}(K((\cdot-x)/h)(x)+o_P(1),\]
 where $\Pi_{J_{0,n}}$ denotes  the projection operator onto  space spanned by $\{\phi_j^*:j\in {\cal R}_{0,n}\}$, and  $K((\cdot-x)/h)$ is a $d$-variate kernel  satisfying $K(0)=1$, $\int K(x)d\mu(x)=1$, and bandwidth $h=(1/J_{0,n})^{1/d}$. Moreover, under the condition that $D^{(k)}_{M}({\cal R}_{0,n})$ satisfies the uniform approximation error $O(C(M(h))r(d,J)^{k+1})$ for function $K((\cdot-x)/h)$ with $h$ such that sectional variation norm of $K((\cdot- x)/h)$ equals $M(h)$, then $\Pi_{J_{0,n}}(K((\cdot-x)/h)(x)=1+o(1)$, so that
 this expression approximates $1/p_0^*(x)$. In particular, if $D^{(k)}_M({\cal R}_{0,n})$ uniformly approximates $D^{(k)}_M([0,1]^d)$ at rate $O(C(M)r(d,J_{0,n})^{k+1})$, then this holds. However, note that  this only requires that it behaves as a nonparametric model for local functions around $x$. This then implies that $\tilde{\sigma}^2_{0,n}(x)\rightarrow_p \sigma^2_0(x)/p_0(x)$, where $\sigma^2_0(x)=E_0( (Y-Q_0)^2\mid X=x)$.

{\bf Behavior of sectional variation norm of $Q_n$ under optimal tuning of $d_n$ up till $\log n$-factors:}
To use these sufficient conditions, it is helpful if the  sectional variation norms of $Q_n,Q_{0,n}$ are uniformly bounded. That is automatically guaranteed for the HAL-MLE. In addition, the  rate  $J_{0,n}$ stated in Theorem \ref{thasnormal}  is slow enough so that this  result follows  from Theorem \ref{theoremconvsectvarnorm}: in fact, we will have convergence in sectional variation norm of $Q_n-Q_{0,n}$. 

Since establishing a desired rate for $\delta_n$ is easily handled for the relax HAL-MLE and sieve MLE, we state here the resulting corollary for these two classes of estimators. 
\begin{corollary}\label{corthasnormal}
Consider the $k$-th order sieve MLE or relax HAL-MLE $Q_n$, $k\in \{1,2,\ldots\}$,
so that $P_n S_j(Q_n)=0$ with $S_j(Q_n)=\phi_j/\sigma^2_n (Y-Q_n)$ for all $j\in {\cal R}_n$. 
Assume that for some $M<\infty$, $Q_0,\sigma^2_0\in D^{(k)}_M([0,1]^d)$; 
and that $\sigma^2_n$ is a zero-order spline HAL-MLE so that
$\pl \sigma^2_n-\sigma^2_0\pl_{P_0}=O_P(n^{-1/3}(\log n)^{d/2})$.
Let $k^*=k+1$.
Assume that for some $\delta>0$ $\inf_{x\in [0,1]^d}\sigma^{2}_{0,n}(x)>\delta$ with probability tending to 1. 

Let ${\cal R}_0$ be a support set identifying a set of basis functions indexed by $(\bar{s}(k+1),u)$ in our $k$-th order spline representation  so that for a given $M<\infty$ $\{Q_n,Q_0\}\subset D^{(k)}({\cal R}_0)$ with probability tending to 1. 
Assume that ${\cal R}_n={\cal R}(P_n)$ satisfies $\max_{g\in D^{(k)}_M({\cal R}_0)} \min_{Q\in D^{(k)}({\cal R}_{n})}\pl g-Q\pl_{\infty}=O_P(C(M)r(d,J_{n})^{k+1})$. Then Assumption A0 holds.

Let $d_{0,n}=n^{1/(2k^*+1)}\log^{m_1} n$ for some $m_1$ with $m_1$ chosen so hat $r(d,J_{0,n})^{k+1}/(n/d_{0,n})^{-1/2}\rightarrow 0$. 
Then, 
\[
\tilde{\sigma}^{-1}_{0,n} (n/d_{0,n})^{1/2}({Q_n}-{Q_{0,n}})(x)\Rightarrow_d N(0,1),\]
and
\[
\tilde{\sigma}^{-1}_{0,n}(n/d_{0,n})^{1/2}(Q_n-Q_0)(x)\Rightarrow_d N(0,1).\]
Moreover,  for some $m<\infty$, 
\[
\mid Q_n-Q_0\mid(x)=O_P(n^{-k^*/(2k^*+1)}\log^m n).\]
These results also apply to the HAL-MLE if one verifies A1 with help of Lemma \ref{lemmatildernx}.\end{corollary}

 \subsection{Key  condition (\ref{undersmoothm})  that  sieve/HAL/relaxHAL-MLEs solve the efficient influence curve equation for a fixed working model $D^{(k)}({\cal R}_{0,n})$.}

The following lemma expresses $\tilde{r}_n(x)$ and clarifies why one generally expect (\ref{undersmoothm}) to hold when ${\cal R}_{0,n}$ approximates ${\cal R}_n$, and that, either way, undersmoothing will shrink it to the desired size. The lemma is presented in general so that it can also be  used in the next section for general loss functions. We refer to Lemma \ref{lemmascoreeqnhal} in Appendix \ref{AppendixJ} for  an explicit bound on the regular score equations $P_n S_j(Q_n)=P_n \frac{d}{dQ_n}L(Q_n)(\phi_j)$, $j\in {\cal R}_n$.
\begin{lemma}\label{lemmatildernx}
Let $S_j(Q)=\frac{d}{dQ}L(Q)(\phi_j)$, $j\in {\cal R}_n$, and let $S_j^*(Q)=\frac{d}{dQ}L(Q)(\phi_j^*)$, $j\in {\cal R}_{0,n}$.
Suppose $P_n S_j(Q_n)=0$ for $j\in {\cal R}_n$. Let $\tilde{P}_n\in {\cal M}$ be a probability distribution so that $Q(\tilde{P}_n)=Q_n$. Then, we have $\tilde{P}_n S_j(Q_n)=0$ for all $j\in {\cal R}_n$, and $\tilde{P}_n S_j^*(Q_n)=0$ for all $j\in {\cal R}_{0,n}$. For a given $S_j^*(Q_n)$ with $j\in {\cal R}_{0,n}$ let $\tilde{S}_j(Q_n)$ be a linear combination of $\{S_k(Q_n): k\in {\cal R}_n\}$ approximating $S_j^*(Q_n)$ in $L^2(P_0)$, $j\in {\cal R}_{0,n}$.
We have \[
\tilde{r}_n(x)=\sum_{j\in {\cal R}_{0,n}}(P_n-\tilde{P}_n)\{S_j^*(Q_n)-\tilde{S}_j(Q_n)\} \phi_j^*(x).\]
Thus, \[
\begin{array}{l}
(n/d_{0,n})^{1/2}\tilde{r}_n(x)=n^{1/2}(P_n-P_0)d_{0,n}^{-1/2}\sum_{j\in {\cal R}_{0,n}}\{S_j^*(Q_n)-\tilde{S}_j(Q_n)\}\phi_j^*(x)\\
\hfill -
n^{1/2}(\tilde{P}_n-P_0)d_{0,n}^{-1/2}\sum_{j\in {\cal R}_{0,n}}\{S_j^*(Q_n)-\tilde{S}_j(Q_n)\}\phi_j^*(x).
\end{array}
\]
The second term can be bounded by $\pl \tilde{p}_n-p_0\pl_{\mu}(n/d_0)^{1/2} d_{0,n}\pl d_{0,n}^{-1}\sum_{j\in {\cal R}_{0,n}}\{S_j(Q_n)-\tilde{S}_j(Q_n)\}\pl_{\mu}$, with $\mu$ Lebesgue measure, so that it suffices 
$\delta_n\equiv  \pl d_{0,n}^{-1}\sum_{j\in {\cal R}_{0,n}}\{S_j(Q_n)-\tilde{S}_j(Q_n)\}\pl_{\mu}=o_P(1)$.
The first term is an empirical process term and can be bounded by $d_{0,n}^{1/2}J(\delta_n,{\cal F}^k_n)$, where ${\cal F}^k_n\equiv \{1/d_{0,n}\sum_{j\in {\cal R}_{0,n}}\{S_j(Q)-\tilde{S}_j(Q)\}\phi_j^*(x): Q\in D^{(k)}_M([0,1]^d)\}$. Assume that  covering number of ${\cal F}^k_n$ is of same order as the covering number of $D^{(k)}_M([0,1]^d)$ (true for regression case), so that by  Lemma \ref{lemmasupnormcoveringnumber} this term is $O_P(d_{0,n}^{1/2} \delta_n^{(2k+1)/(2k+2)})$ up till $\log \delta$-factors. Then,  the first term is $o_P(1)$ if $d_{0,n}^{1/2}\delta_n^{(2k+1)/(2k+2)}=o_P(1)$.

If $Q_n$ is the HAL-MLE so that we only have $P_n S_j(Q_n)\approx 0$, $j\in {\cal R}_n$, the same bound applies with $\tilde{S}_j(Q_n)$ replaced by linear combination of scores  $S
^j(Q_n)$ for which $P_n S^j(Q_n)=0$ for $j\in {\cal R}_n$. 
These scores solved by the $L_1$-norm constrained MLE $Q_n$ correspond with paths $\{ (1+\delta h(j))\beta_n(j):j\in {\cal R}_n):\delta\}$ with  $h$ satisfying $\sum_{j\in {\cal R}_n}h(j)\mid \beta_n\mid (j)=0$. 
\end{lemma}
{\bf Proof of Lemma \ref{lemmatildernx}:}
We have $(\tilde{P}_n-P_n)\tilde{S}_j(Q_n)=0$. 
So, for each $j\in {\cal R}_{0,n}$,  we have $P_n S_j^*(Q_n)=(P_n-\tilde{P}_n)\{S_j^*(Q_n)-\tilde{S}_j(Q_n))\}$.
Therefore, 
\[
\tilde{r}_n(x)=\sum_{j\in {\cal R}_{0,n}}(P_n-\tilde{P}_n)\{S_j^*(Q_n)-\tilde{S}_j(Q_n)\} \phi_j^*(x).\]
We also notice that
\[
(P_n-\tilde{P}_n)\{S_j^*(Q_n)-\tilde{S}_j(Q_n)\}=(P_n-P_0)\{S_j^*(Q_n)-\tilde{S}_j(Q_n)\}-
(\tilde{P}_n-P_0)\{S_j^*(Q_n)-\tilde{S}_j(Q_n)\}.\]
$\Box$

The following lemma shows that condition $d_{0,n}^{1/2}\delta_n^{(2k+1)/(2k+2)}=o_P(1)$ holds if $D^{(k)}({\cal R}_{0,n})$ satisfies the nonparametric uniform approximation condition, but this latter is by no means necessary.

\begin{lemma}
We have for each $j\in {\cal R}_{0,n}$ $S_j^*(Q_n)-\tilde{S}_j(Q_n)=(\phi_j^*-\tilde{\phi}_j)/\sigma^2_n(Y-Q_n(X))$, where
$\tilde{\phi}_j$ is a linear combination of $\{\phi_j:j\in {\cal R}_n\}$, approximating $\phi_j^*$.
If $D^{(k)}_M({\cal R}_n)$ is such that it provides a $O(C(M)r(d,J_n)^{k+1})$-uniform approximation of
$D^{(k)}_M([0,1]^d)$ (which holds if ${\cal R}_n={\cal R}^k(d,{\bf J}_n)$), then, given that $\phi_j^*\in D^{(k)}_{M_n}({\cal R}_{0,n})$ for all $j\in {\cal R}_{0,n}$, and $J_n\sim J_{0,n}$, this then implies that also $\max_{j\in {\cal R}_{0,n}}\pl S_j^*(Q_n)-\tilde{S}_j(Q_n)\pl_{\infty}=O(C(M_n)r(d,J_{0,n})^{k+1})$. So then $\delta_n=O_P(C(M_n)r(d,J_{0,n})^{k+1})$, so that condition $d_{0,n}^{1/2}\delta_n^{(2k+1)/(2k+2)}=o_P(1)$ in previous lemma trivially holds.
\end{lemma}



\subsection{Understanding second order empirical process term $E_n(x)$.}

The following lemma avoids any smoothness beyond $D^{(0)}([0,1]^d)$, since it is simply not needed to control $E_n(x)$. In this manner, these lemmas also apply to the zero order spline sieve and HAL-MLEs.
\begin{lemma}\label{lemmasuffEnx}
Suppose that $\pl \sigma^2_n\pl_v^*=O_P(1)$ and $\sigma^2_0\in D^{(0)}_M([0,1]^d)$ for some $M<\infty$. 
Assume $\sigma^2_n$, based on an estimator of $E_0(Y\mid X)$ and $E_0(Y^2\mid X)$,  satisfies $\pl \sigma^2_n-\sigma^2_0\pl_{P_0}=O_P( n^{-1/3}(\log n)^{d/2})$, as would be the case when using a zero-order HAL-MLE and only assuming $E_0(Y^2\mid X)\in D^{(0)}_M([0,1]^d)$ for some $M<\infty$. Also assume $\pl Q_n-Q_0\pl_{P_0}=O_P(n^{-1/3}(\log n)^{d/2})$.
Then, $E_n(x)=O_P(d_{0,n} n^{-2/3})$ up till $\log n$-factors. Thus, $E_n(x)=o_P((d_{0,n}/n)^{1/2})$ if $d_{0,n}<n^{1/3-\delta}$ for some $\delta>0$. 
\end{lemma}
{\bf Proof:}
Note that $E_n(x)=d_{0,n} (P_n-P_0) g_n$, where 
\[
g_n=\frac{1}{d_{0,n}} \sum_{j\in {\cal R}_{0,n}}\phi_j^*\left\{1/\sigma^2_n(Y-{Q_n}) -1/\sigma^2_{0,n}(Y-Q_{0,n}) \right\} \phi_j^*(x).\]
 Let ${\cal F}_n\equiv \{ d_{0,n}^{-1}\sum_{j\in {\cal R}_{0,n}}\phi_j^*\{1/\sigma^2(Y-Q)-1/\sigma^2_{0,n}(Y-Q_{0,n}): f,\sigma^2\in D^{(0)}_M([0,1]^d)\}$.  We have that ${\cal F}_n$ is subset of $D^{(0)}_M([0,1]^{d+1})$ as functions of $(X,Y)$, where, for simplicity, we assume $Y\in [0,1]$.
We have $\pl g_n\pl_{P_0}=O_P(n^{-1/3}(\log n)^{d/2})$. Moreover, $g_n\in {\cal F}_n\subset D^{(0)}_M([0,1]^d)$ for some $M<\infty$, with probability tending to 1. We have $J(\delta,{\cal F}_n)\sim \delta^{1/2}$ up till $\log \delta$-factor. So this shows that $E_n(x)\sim d_{0,n}n^{-1/2} n^{-1/6}$ up till $\log n$-factors. $\Box$

The following lemma provides a seemingly weaker condition for controlling $E_n(x)=o_P((n/d_{0,n})^{-1/2})$, but either way, the above lemma suffices for our purpose. 
\begin{lemma}
Let $e_n=1/\sigma^2_n(Y-{Q_n}) -1/\sigma^2_{0,n}(Y-Q_{0,n})$.
For $(j_1,j_2)\in {\cal R}_{0,n}\times {\cal R}_{0,n}$, let \[
E_n(j_1,j_2)\equiv
P_0 \left\{ \phi_{j_1}^*\phi_{j_2}^*e_n^2 \right\} .\]
Let $\lambda_n$ be the maximal eigenvalue of matrix $E_n$. 
By Lemma \ref{lemmasuffkeyconditiona}, (\ref{keyconditiona1m}) holds if   $d_{0,n} \lambda_n=o_P(1)$, where we are reminded $d_{0,n}\sim J_{0,n}$. 
\end{lemma}

\subsection{Reasoning behind our concrete proposal for ${\cal R}_{0,n}$.} 
Firstly, ${\cal R}_n$ needs to be chosen  so that $D^{(k)}({\cal R}_n)$ uniformly approximates  $Q_0$ at rate  $r(d,J_{n})^{k+1}$ negligible relative to $(n/d_{0,n})^{-1/2}$. By selecting ${\cal R}_{0,n}$ as an approximation (or conservative approximation) of ${\cal R}_n$, this  uniform approximation property will then also be inherited by ${\cal R}_{0,n}$.
Secondly, ${\cal R}_{0,n}$ should be chosen so that $S_j^*(Q_n)$, $j\in {\cal R}_{0,n}$, is well approximated by scores $\{\tilde{S}_j(Q_n):j\in {\cal R}_{0,n}\}$ that are linear combinations of $\{S_j(Q_n): j\in {\cal R}_n\}$ satisfying $P_n S_j(Q_n)=0$, $j\in {\cal R}_n$, so that $P_n \tilde{S}_j(Q_n)=0$ for $j\in {\cal R}_{0,n}$. 

{\bf A particular proposal for ${\cal R}_{0,n}$:}
In order to make the score equations $r_n(j)$, $j\in {\cal R}_n$, solved by $Q_n$ easily span the score equations, $r_n^*(j)$, $j\in {\cal R}_{0,n}$, and thereby its linear combination $\tilde{r}_n(x)$,  a particular choice of interest is to select ${\cal R}_{0,n}$ as the realization of the same algorithm $P_n\rightarrow {\cal R}_n=\widehat{\cal R}(P_n)$ but applied to the empirical measure $P_n^{\#}$ of an independent sample $O_i^{\#}\sim_{iid} P_0$, so that ${\cal R}_{0,n}=\widehat{\cal R}(P_n^{\#})$ is independent of ${\cal R}_n$.
In that case, ${\cal R}_n$ has the same distribution as ${\cal R}_{0,n}$ so that the set of basis functions $\widehat{\cal R}(P_n^{\#})$ with non-zero coefficients  should be very comparable and live in the same overall  support set of the $k$-th order spline representation of the target function $Q_0$. Of course, choosing ${\cal R}_{0,n}$ as a shrunk version of ${\cal R}_n$ would satisfy the same. 

Moreover, due to the HAL-MLE  optimal rate of convergence  w.r.t. loss-based dissimilarity established in previous section, it is expected that the HAL-MLE fit applied to $P_n^{\#}$ generates sets ${\cal R}(d,{\bf J}_{0,n})$ satisfying the adaptive uniform approximation error $O(r(d,J_{0,n})^{k+1})$  (otherwise it could not achieve that rate of convergence to $Q_0$, see Appendix \ref{AppendixB}).  Therefore this proposal for ${\cal R}_{0,n}$ is certainly an appropriate choice for the sieve MLE and relax HAL-MLE since these exactly solve $r_n(j)=0$, $j\in {\cal R}_n$. If one uses HAL-MLE, one might note that ${\cal R}_{0,n}={\cal R}(d,{\bf J}_{0,n})$ for some ${\bf J}_{0,n}$ and replace ${\bf J}_{0,n}$ by $\delta {\bf J}_{0,n}$ for some shrinking factor $\delta\in (0,1)$. Alternatively, if HAL uses an undersmoother for $C_n>C_{n,cv}$, one might select ${\cal R}_{0,n}$ as applying HAL with the cross-validation selector $C_{n,cv}$ or slightly less aggressive undersmoother such as $C_n/\log n$  to an independent sample $P_n^{\#}$. 
In that manner, one takes into account that HAL uses $L_1$-regularization  so that it might only solve the desired efficient score equation $\tilde{r}_n(x)$ for a slightly smaller model $D^{(k)}({\cal R}_{0,n})$ than $D^{(k)}({\cal R}_n)$, although our Lemma \ref{lemmatildernx} suggests that this might not be needed at all.

{\bf Appropriate choice of ${\cal R}_{0,n}$ that allows conservative inference based on $D^{(k)}({\cal R}_n)$:}
To assist in our inference, we like to restrict to ${\cal R}_{0,n}$ with size $d_{0,n}\leq d_n$ with probability tending to 1,  or $d_{0,n}/d_n\rightarrow_p 1$, making the above choices ${\cal R}_{0,n}$ appropriate.  Although this is not needed for establishing  asymptotic normality of $(Q_n-Q_{0,{\cal R}_{0,n}})(x)$, it allows for conservative variance estimation based on working model $D^{(k)}({\cal R}_n)$ without having to know ${\cal R}_{0,n}$. 
Since ${\cal R}(P_n^{\#})=_d {\cal R}(P_n)$ it should follow that $d_{0,n}/d_n\rightarrow_p 1$, so that our proposed choice ${\cal R}_{0,n}={\cal R}(P_n^{\#})$ appears to also satisfy this desired property. 
To conclude, ${\cal R}_{0,n}={\cal R}(P_n^{\#})$ or a (proportionally) shrunk version of it is the most sensible appropriate choice for ${\cal R}_{0,n}$, and the uniform approximation property required on $D^{(k)}({\cal R}_{0,n})$ at most requires undersmoothing $Q_n$ relative to the cross-validation selector of its tuning parameters.

\subsection{Convergence in distribution of sieve and HAL-MLEs across multiple points} Consider Theorem \ref{thasnormal}. Consider two points $x,y$ and assume the conditions of this theorem for both points.  The asymptotic variance of $(n/d_{0,n})^{1/2}(Q_n-Q_{0,n})(x)$ is given by $\tilde{\sigma}^2_{0,n}(x)=d_{0,n}^{-1}\sum_{j\in {\cal R}_{0,n}}\{\phi_j^*(x)\}^2$, and similarly for $(n/d_{0,n})^{1/2}(Q_n-Q_{0,n})(y)$. 
By applying the bivariate central limit theorem to the bivariate mean zero empirical  mean of independent random variable approximation stated in the theorem, we have that $(n/d_{0,n})^{1/2}((Q_n-Q_{0,n})(x),(Q_n-Q_{0,n})(y))$ converges to a bivariate normal mean zero distribution with covariance matrix that has $\tilde{\sigma}^2_{0,n}(x)$ and $\tilde{\sigma}^2_{0,n}(y)$ on the diagonal, and asymptotic covariance $\rho_{0,n}(x,y)$   given by \[
\rho_{0,n}(x,y)= \frac{1}{d_{0,n}}\sum_{j\in {\cal R}_{0,n}}\phi_j^*(x)\phi_j^*(y).
\]
Lemma \ref{lemmacovariance} shows that
\[
\rho_{0,n}(x,y)= (p_0^*(x)p_0^*(y))^{-1} J_{0,n}^{-1}P_0^*\Pi_{J_{0,n}}(K_{J_{0,n},x})\Pi_{J_{0,n}}(K_{J_{0,n},y}),\]
where $K$ is any kernel so that $K(0)=1$; $\int K(x)d\mu(x)=1$, and $K_{J,x}(\cdot)=K((\cdot-x)/h)/h^d$ with $h=J^{-1/d}$, and $\Pi_{J_{0,n}}$ is the projection operator onto linear span of $\{\phi_j^*:j\in {\cal R}_{0,n}\}$ (i.e., onto $D^{(k)}({\cal R}_{0,n})$).
In particular, if  $D^{(k)}_M({\cal R}_{0,n})$ uniformly approximates $D^{(k)}([0,1]^d)$ at rate $O(C(M)r(d,J_{0,n})^{k+1})$ (we refer to this as the uniform nonparametric approximation condition on $D^{(k)}({\cal R}_{0,n})$), then $\Pi_J(K_{J,x})$ and $\Pi_J(K_{J,y})$ can be replaced by $K_{J,x}$ and $K_{J,y}$ so that it follows that $\rho_{0,n}(x,y)=o_P(1)$. Thus, if our working model is tailored to approximate the whole nonparametric model $D^{(k)}([0,1]^d)$ (as in our construction of ${\cal R}(d,{\bf J})$ in Section \ref{section4}), then the estimators $Q_n(x)$ and $Q_n(y)$ are asymptotically independent.

For simplicity we state the convergence in distribution at two points but this has a trivial generalization to a general finite set of points. 
\begin{corollary}
Consider Theorem \ref{thasnormal}. Consider two points $x,y$ and assume the conditions A0,A1,A2,A3,A4 hold at both points.   Then, both $(n/d_{0,n})^{1/2}((Q_n-Q_{0,n})(x),(Q_n-Q_{0,n})(y))$ and $(n/d_{0,n})^{1/2}((Q_n-Q_0)(x),(Q_n-Q_0)(y))$ converge jointly distribution to a bivariate normal mean zero with covariance the limit $\rho_0(x,y)$ of  
\begin{equation}\label{rho0n}
\rho_{0,n}(x,y)=\frac{1}{d_{0,n}}\sum_{j\in {\cal R}_{0,n}}\phi_j^*(x)\phi_j^*(y)
\end{equation}
 and variance elements the limits $\tilde{\sigma}^2_0(x)$ and $\tilde{\sigma}^2_0(y)$  of $\tilde{\sigma}^2_{0,n}(x)=d_{0,n}^{-1}\sum_{j\in {\cal R}_{0,n}}\{\phi_j^*(x)\}^2$ and $\tilde{\sigma}^2_{0,n}(y)$, respectively.
We have $\sigma^2_{0,n}(x)=p_0^*(x)^{-1} \Pi_{J_{0,n}}(K_{J_{0,n},x})(x)+o(1)$; $\sigma^2_{0,n}(y)=p_0^*(y)^{-1} \Pi_{J_{0,n}}(K_{J_{0,n},y})(y)+o(1)$; and \[
\rho_{0,n}(x,y)=(p_0^*(x)p_0^*(y))^{-1} J_{0,n}^{-1}P_0^*\Pi_{J_{0,n}}(K_{J_{0,n},x})\Pi_{J_{0,n}}(K_{J_{0,n},y}),\]
 where $K$ is any kernel so that $K(0)=1$; $\int K(x)d\mu(x)=1$, and $K_{J,x}(\cdot)=K((\cdot-x)/h)/h^d$ with $h=J^{-1/d}$, and $\Pi_{J_{0,n}}$ is the projection operator onto linear span of $\{\phi_j^*:j\in {\cal R}_{0,n}\}$ (i.e., onto $D^{(k)}({\cal R}_{0,n})$). If $D^{(k)}_M({\cal R}_{0,n})$ uniformly approximates $D^{(k)}([0,1]^d)$ at rate $O_P(C(M)r(d,J_{0,n})^{k+1})$, then  $\tilde{\sigma}^2_0(x)=\sigma^2_0(x)/p_0(x)$, $\tilde{\sigma}^2_0(y)=\sigma^2_0(y)/p_0(y)$ and $\rho_0(x,y)=0$. 
\end{corollary}

\subsection{Understanding asymptotic variance and covariance expressions (\ref{sigma2n}) and (\ref{rho0n})}

The next lemma establishes that $J^{-1}\sum_{j\in {\cal R}_{0,n}} \{\phi_j^*\}^2(x)=p_0^*(x)^{-1}\Pi_J(K((\cdot-x))/h))(x)+o(1)$, which is an expression that avoids having to know the orthonormal basis and could be easily estimated from the original basis $\phi_j$. In addition, if the nonparametric uniform approximation condition on $D^{(k)}({\cal R}_{0,n})$, then this equals $\{\sigma^{-2}_{0,n}(x)p_0(x)\}^{-1}+o(1)$,  providing a closed form expression for asymptotic variance of $Q_n(x)$. If $D^{(k)}({\cal R}_{0,n})$ approximates a real submodel $D^{(k)}({\cal R}_0)$ of $D^{(k)}([0,1]^d)$ containing $Q_0$, then this expression (\ref{p1})  still applies.

\begin{lemma}\label{lemmaexplicit}
Let $x\in [0,1]^d$ with $p_0^*(x)>0$.
Let $\{\phi_j^*:j=1,\ldots,J\}$ be an orthonormal basis in $L^2(P_0^*)$. 
Let $K:[-1,1]^d\rightarrow \openr_{\geq 0}$ be a kernel so that $K(0)=1$, $\int K(x)dx=1$. Let $K_{J,x}(z)=K((z-x)/h)/h^d$ with $h^d=1/J$. Let 
\[
r_J(j,x)\equiv \langle K_{J,x},\phi_j^*\rangle_{P_0^*}-\phi_j^*(x)p_0^*(x).\] 
{\bf Assumption:}
 Suppose that 
  \begin{equation}\label{lemmavar2}\frac{1}{J}\sum_{j=1}^Jr_J(j,x) \phi_j^*(x)=o_P(1).\end{equation} 
Then  \[ \frac{1}{J}\sum_j \langle K_{x,J},\phi_j^*\rangle_{P_0^*}\phi_j^*(x)=p_0^*(x) \frac{1}{J}\sum_j \{\phi_j^*(x)\}^2+o(1),\]
and thus
 \begin{equation}\label{p1}
 \frac{1}{J}\sum_{j=1}^J\{\phi_j^*(x)\}^2=\{p_0^*(x)\}^{-1}  \Pi_J(K((\cdot-x)/h) )(x)+o(1).\end{equation}
 If 
 \begin{equation}\label{kernelapprox}
 \pl K((\cdot-x)/h)-\Pi_J(K((\cdot-x)/h)\pl_{\infty}=o_P(1),\end{equation} then it follows that the right-hand side of (\ref{p1}) equals $p_0^*(x)+o_P(1)$.
 
 {\bf Particular assumption under which (\ref{kernelapprox}) holds:}
 If $K$ is $k$-times continuously differentiable and  $D^{(k)}({\cal R}_{0,n})$ satisfies the nonparametric uniform approximation condition, then $\pl K((\cdot-x)/h)-\Pi_J(K((\cdot-x)/h)\pl_{\infty}=O((\log J)^m/J)$ so that assumption (\ref{kernelapprox}) holds.
 
  \end{lemma}
  Regarding the last statement, we have $K_{J,x}=\Pi_J(K_{J,x})+E_{J,x}$, where 
$\Pi_J(K_{J,x})=\sum_j \langle K_{J,x},\phi_j^*\rangle_{P_0^*}\phi_j^*$ is the projection of $K_{J,x}$ onto the linear span of $\{\phi_j^*:j\}$ in $L^2(P_0^*)$,
and $E_{J,x}\equiv K_{J,x}-\Pi_J(K_{J,x})$ is the approximation error. 
If $K$ is $k$-times continuously differentiable, then it follows that
$K_{J,x}\in D^{(k)}_{(1/J)^k}([0,1]^d)$. 
If also $\sup_{Q\in D^{(k)}_M([0,1]^d)}\inf_{\beta}\pl Q-\sum_j\beta(j)\phi_j^*\pl_{\infty}=O(C(M) r(d,J)^{k+1})$, then it follows that $\pl E_{J,x}\pl_{\infty}=O( J^{k+1} r(d,J)^{k+1})$ which is $O((\log J)^m)$ for some $m<\infty$, which proves the last statement. 
The proof of other statements in this lemma is straightforward and  is omitted here.



The following lemma proves that $\rho_{0,n}(x,y)=o(1)$ under the nonparametric uniform approximation condition so that we can conclude that the standardized estimator at two different points converges to a bivariate normal with mean zero and diagonal covariance matrix, demonstrating that the two point estimators are asymptotically independent. 
If only the adaptive uniform approximation condition holds, then we still obtain an expression that avoids having to know the orthonormal basis, but in that case there is no reason to expect that $Q_n(x)$ and $Q_n(y)$ are asymptotically uncorrelated due to the limit model $D^{(k)}({\cal R}_0)$ allowing for extrapolation. 

\begin{lemma}\label{lemmacovariance}
Let $x,y\in [0,1]^d$ and  $x\not =y$ with $p_0^*(x)>0$ and $p_0^*(y)>0$. Consider setting of above lemma including definition of  $K_{J,x}$ and $r_J(j,x)$. 
We have $K_{J,x}=\Pi_J(K_{J,x}) +E_{J,x}$, where 
$\Pi_J(K_{J,x})= \sum_j \langle K_{J,x},\phi_j^*\rangle_{P_0^*}\phi_j^*$ is the projection of $K_{J,x}$ onto the linear span of $\{\phi_j^*:j\}$ in $L^2(P_0^*)$,
and $E_{J,x}=K_{J,x}-\Pi_J(K_{J,x})$ is the approximation error. Similarly, we define $E_{J,y}$. 

Assume 
 \begin{equation}\label{approxb}
 J^{-1}\left\{ \sum_j r_J(j,x)r_J(j,y)+r_J(j,x)\phi_j^*(y)p_0^*(y)+r_J(j,y)\phi_j^*(x)p_0^*(x) \right\}=o(1).
 \end{equation}
Then we have
 \[
 J^{-1}\sum_j \phi_j^*(x)\phi_j^*(y)=(p_0^*(x)p_0^*(y))^{-1} J^{-1}P_0^*\Pi_J(K_{J,x})\Pi_J(K_{J,y})=O(1).\]

Suppose we also have
\begin{equation}\label{approxa}
 J^{-1}\left\{ P_0^*(E_{J,x} E_{J,y})+P_0^*(E_{J,x}\Pi_J(K_{J,y}))+P_0^*(E_{J,y}\Pi_J(K_{J,x}))\right\}=o(1),\end{equation}
which holds, in particular, under the nonparametric uniform approximation condition on $D^{(k)}({\cal R}_{0,n})$.

 Then, 
\[J^{-1}\sum_j \phi_j^*(x)\phi_j^*(y)=o(1).\]
\end{lemma}
Condition (\ref{approxa}) holds if $K$ is $k$-times continuously differentiable, and is thus a trivial condition.

{\bf Proof of Lemma \ref{lemmacovariance}:}
We have that $K_{J,x}\approx \sum_j \langle K_{J,x},\phi_j^*\rangle_{P_0^*}\phi_j^*$, where the right-hand size equals the projection $\Pi_J(K_{J,x})$ of $K_{J,x}$ onto the linear span of $\{\phi_j^*:j\}$, and similarly for $K_{J,y}$. We defined $E_{J,x}=K_{J,x}-\Pi_J(K_{J,x})$, and similarly we defined $E_{J,y}$. Regarding this $E_{J,x}$ and $E_{J,y}$ we assumed (\ref{approxa}).
 We defined  $r_J(j,x)\equiv \langle K_{J,x},\phi_j^*\rangle_{P_0^*}-\phi_j^*(x)p_0^*(x)$, 
and regarding this approximation error $r_J(j,x)$ and $r_J(j,y)$ we assumed (\ref{approxb}).   
  Using that $\{\phi_j^*:j\}$ is an orthonormal basis in $L^2(P_0^*)$, it follows that 
  \[
  \begin{array}{l}
    J^{-1}P_0^*(\Pi_J(K_{J,x})\Pi_J(K_{J,y})\\
=J^{-1} P_0^*\left( \sum_{j_1,j_2} \langle K_{J,x},\phi_{j_1}^*\rangle_{P_0^*}\langle K_{J,y},\phi_{j_2}^*\rangle_{P_0^*}\phi_{j_1}^*\phi_{j_2}^* \right ) \\
=J^{-1}\sum_{j_1,j_2} \langle K_{J,x},\phi_{j_1}^*\rangle_{P_0^*}\langle K_{J,y},\phi_{j_2}^*\rangle_{P_0^*}
P_0^*( \phi_{j_1}^*\phi_{j_2}^*)\\
=J^{-1} \sum_j \langle K_{J,x},\phi_{j}^*\rangle_{P_0^*}\langle K_{J,y},\phi_{j}^*\rangle_{P_0^*}\\
=J^{-1} \sum_j \left\{ \phi_j^*(x)p_0^*(x) +r_J(j,x)\right\}\left\{ \phi_j^*(y)p_0^*(y)+r_J(j,y)\right\}\\
=J^{-1}p_0^*(x)p_0^*(y)\sum_j \phi_j^*(x)\phi_j^*(y), 
\end{array}
\]
by condition (\ref{approxb}). This proves the first statement of the lemma. 

Suppose now that also condition {\ref{approxa}) holds.   
Consider the quantity $P_0^* K_{J,x}K_{J,y}$. As $J$ converges to infinity,  $K_{J,x}K_{J,y}=0$ for large enough $J$. In particular $P_0^* K_{J,x}K_{J,y}=0$ for large enough $J$.
Thus, for large enough $J$, we have 
\[
\begin{array}{l}
0=J^{-1}P_0^* K_{J,x}K_{J,y}\\
=J^{-1}P_0^*(\Pi_J(K_{J,x} )+E_{J,x})(\Pi_J(K_{J,y})+E_{J,y})\\
=J^{-1}P_0^*(\Pi_J(K_{J,x})\Pi_J(K_{J,y})\\
+
J^{-1}\left\{ P_0^*(E_{J,x} E_{J,y})+P_0^*(E_{J,x}\Pi_J(K_{J,y}))+P_0^*(E_{J,y}\Pi_J(K_{J,x}))\right\}\\
=J^{-1}P_0^*(\Pi_J(K_{J,x})\Pi_J(K_{J,y})+o(1),
\end{array}
\]
by (\ref{approxa}).
The proof above showed that the latter expression equals $J^{-1}p_0^*(x)p_0^*(y)\sum_j \phi_j^*(x)\phi_j^*(y)+o(1)$, under assumption (\ref{approxb}). This proves that $J^{-1}\sum_j \phi_j^*(x)\phi_j^*(y)=o(1)$ for $x\not =y$, under (\ref{approxa}) and (\ref{approxb}). 
 This completes the proof of the lemma. $\Box$

\subsection{Inference}\label{subsectinference}
{\bf Confidence interval based on our formulas for the asymptotic variance:}
Theorem \ref{thasnormal} provides an asymptotic variance  for $(n/d_{0,n})^{1/2}(Q_n-Q_{0,n})(x)$ given by $\tilde{\sigma}^2_{0,n}(x)=
\sigma^2_{0,n}(x)/p_0(x)   \Pi_{J_{0,n}}(K((\cdot-x)/h_{0,n}) )(x)$ with $h_{0,n}=d_{0,n}^{-1/d}$. This can  be directly estimated by estimating $p_0(x)$, $\sigma^2_0(x)$, and computing the projection of a kernel onto the  linear working model $D^{(k)}({\cal R}_{n})$ (as an approximation of $D^{(k)}({\cal R}_{0,n})$, and evaluating it at $x$. If the working model satisfies the nonparametric uniform approximation condition, then we showed that  $\tilde{\sigma}^2_{0,n}(x)\rightarrow \sigma^2_0(x)/p_0(x)$, where $\sigma^2_0(x)$ is the conditional variance of $(Y-Q_0(x))$, given $X=x$, so that could be used as a conservative variance estimator. 
The variance of $Q_n(x)$ is then estimated as $(d_n/n)\sigma^2_n(x)/p_n(x)$, using that $d_n\approx d_{0,n}$. This method for constructing a pointwise confidence interval avoids any need for knowing ${\cal R}_{0,n}$. 
We can also use the expression $1/d_n\sum_{j\in {\cal R}_n}\{\phi_j^*\}^2(x)$, where $\{\phi_j^*:j\in {\cal R}_n\}$ is an orthonormal basis of $\{\phi_j: j\in {\cal R}_n\}$ in $L^2(\sigma^{-2}_ndP_n)$, as an approximation of $1/d_{0,n}\sum_{j\in {\cal R}_{0,n}}\{\phi_j^*\}^2(x)$ with $\phi_j^*$ an orthonormal basis of $D^{(k)}({\cal R}_{0,n})$ in $L^2(P_0^*)=L^2(\sigma^{-2}_{0,n}dP_0)$.  This then requires calculating an orthonormal basis from $\{\phi_j:j\in {\cal R}_n\}$, which is straightforward. 

{\bf Confidence interval based on $\delta$-method applied to working model $D^{(k)}({\cal R}_n)$:}
However, it might be desirable to obtain an asymptotic variance estimator that acts as if the working model is fixed, which should be more reflective of the actual variance of $Q_n(x)$ in finite samples.  
Again, due to $\lim\sup d_n/d_{0,n}\geq 1$, we should be able to estimate this working model based asymptotic variance with the working model is $D^{(k)}({\cal R}_n)$. 
For example, we could directly aim to estimate  $\tilde{\sigma}^2_{0,n}(x)=1/d_{0,n}\sum_{j\in {\cal R}_{0,n}}\{\phi_j^*\}^2(x)$ with $\phi_j^*$ an orthonormal basis of $D^{(k)}({\cal R}_{0,n})$ in $L^2(P_0^*)=L^2(\sigma^{-2}_{0,n}dP_0)$. 
The natural estimator of $\tilde{\sigma}^2_{0,n}(x)$ is given by $\tilde{\sigma}^2_n=1/d_n\sum_{j\in {\cal R}_n}\{\phi_j^*\}^2(x)$, where $\{\phi_j^*:j\in {\cal R}_n\}$ is an orthonormal basis of $\{\phi_j: j\in {\cal R}_n\}$ in $L^2(\sigma^{-2}_ndP_n)$. This then requires calculating an orthonormal basis from $\{\phi_j:j\in {\cal R}_n\}$, which is straightforward.

Alternatively and similarly, we can work with the original basis $\{\phi_j: j\in {\cal R}_n\}$. Under the conditions of Theorem \ref{thasnormal}, in terms of the original basis we have \[
Q_n(x)-Q_{0,n}(x)\approx P_n D_{{\cal R}_{n},x,\beta_{0,n}},\]
where 
\begin{equation}\label{deltamethodic}
D_{{\cal R}_{n},x,\beta_{0,n}}=\left\{ I(\beta_{0,n})^{-1}S_{{\cal R}_n,\beta_{0,n}}\right\}^{\top}{\bf \phi}(x),\end{equation}
$S_{{\cal R}_{n},\beta_{0,n}}=\left( \frac{d}{dQ_{0,n}}L(Q_{0,n})(\phi_j): j\in {\cal R}_n\right)^{\top}$ is the score for working model $D^{(k)}({\cal R}_{n})$, and $I(\beta_{0,n})$ is the information matrix defined by the minus second derivative of the empirical risk:
\[
I(\beta_{0,n})=-P_n \frac{d^2}{d\beta_{0,n}}L(\sum_{j\in {\cal R}_{n}}\beta(j)\phi_j).\]
That is, we can base statistical inference on the standard $\delta$-method influence curve for MLE $\Phi(\beta_n)\equiv \beta_n^{\top}{\bf \phi}(x)$ of $\Phi(\beta_{0,n})\equiv \beta_{0,n}^{\top}{\bf \phi}(x)$ defined by the parametric working model $D^{(k)}({\cal R}_{n})=\{\sum_{j\in {\cal R}_{n}} \beta(j)\phi_j:\beta\}$ assuming a nonparametric model for $P_0$. 
so that $Q_n(x)\sim N(Q_{0,n}(x),\sigma^2_{{\cal R}_{n}}(x))$ with
\[
\tilde{\sigma}^2_n(x)=P_n D_{{\cal R}_{n},x,\beta_n}^2/n.\]

Through simulations we have been able to demonstrated that indeed, under undersmoothing condition solving the efficient score equation  $P_n D_{{\cal R}_{n},\beta_n}=o_P(n^{-1/2})$, which corresponds with (\ref{tildernx}) for parameter $Q_{0,{\cal R}_n}(x)$,  this delta-method based variance estimator consistently estimates the asymptotic variance of the HAL-MLE.

Therefore, we conclude that for all three classes of estimators, assuming that $C_n$ appropriately undersmooths for the HAL-MLEs, we can use the variance estimator $\tilde{\sigma}^2_n$ corresponding with applying the $\delta$-method to obtain the influence curve of $\Phi(\beta_n)$ of  parameter $\Phi(\beta_{0,n})$ defined by the parametric working model $D^{(k)}({\cal R}_n)$ for $Q_0$ within the nonparametric model.

\subsection{Rate of convergence w.r.t. supremum norm and simultaneous confidence bands for regression function}\label{section11}

Our previous Theorem \ref{thasnormal} actually provided a linear approximation for $(Q_n-Q_{0,n})(x)$, which can then be used to determine the rate of convergence w.r..t supremum norm and obtain simultaneous confidence bands. Of course, we then need to make sure that all the conditions on the total remainder $R_n(x)$ of Theorem \ref{thasnormal} apply to each $x\in [0,1]^d$ and, in facts, hold uniformly in $x$. The linear approximation can then be viewed as an empirical process indexed by the $d$-dimensional class $\{D_{Q_{0,n},x}:x\}$ which is easily handled through empirical process theory, and only means that we lose a $\log n$-rate. 

Let Assumption A1* ,A2* and A3* be the uniform analogues of A1,A2 and A3 defined by $\pl  \tilde{r}_n\pl_{\infty} =o_P((n/d_{0,n})^{1/2})$; $\pl R_{n,1}\pl_{\infty}=o_P((n/d_{0,n})^{-1/2})$, and $\pl E_n\pl_{\infty}=o_P((n/d_{0,n})^{1/2})$, respectively.

\begin{theorem}\label{thunifregr}\ \nl
Consider the $k$-th order sieve MLE, HAL-MLE or relax HAL-MLE $Q_n$.  Let $k^*=k+1$.
Let ${\cal R}_n$ be the set of $d_n$ non-zero coefficients in $Q_n=\sum_{j\in {\cal R}_n}\beta_n(j)\phi_j$. Assume that for some $\delta>0$ $\inf_{x\in [0,1]^d}\sigma^{2}_{0,n}(x)>\delta$ with probability tending to 1. 
Recall the definitions of $\tilde{r}_n(x)$ (\ref{tildernx}), $R_{1n}(x)$ (\ref{R1nx}), $E_n(x)$ (\ref{Enx}), $D_{Q_{0,n},x}(O)$ (\ref{Df0nx}) above. Assume assumptions A0,A1*,A2*,A3*,A4. 

By  Assumption A0, we have the following expansion:
\begin{eqnarray}
(n/d_{0,n})^{1/2}(Q_n-Q_{0,n})(x)&=& n^{1/2}(P_n-P_0) d_{0,n}^{-1/2}\sum_{j\in {\cal R}_{0,n}}\phi_j^*/\sigma^2_{0,n}(Y-{Q_{0,n}}) \phi_j^*(x)\nonumber\\
&&-(n/d_{0,n})^{1/2}E_n(x) -
(n/d_{0,n})^{1/2}\tilde{r}_n(x)
\nonumber \\
&&\hspace*{-3cm}+(n/d_{0,n})^{1/2}R_{1n}(x) +(n/d_{0,n})^{1/2}O_P(C(M_n)r(d,J_{0,n})^{k+1}),\label{keyexpansion}\end{eqnarray}
where $M_n=\pl Q_n\pl_{v,k}^*$, which equals the $L_1$-norm  of the vector of non-zero coefficients in $Q_n$.

By Assumptions A0,A1*,A2*,A3*,A4, we have
\begin{eqnarray}
\frac{(Q_n-Q_{0,n})(x)}{(n/d_{0,n})^{-1/2}}= n^{1/2}(P_n-P_0) d_{0,n}^{-1/2}\sum_{j\in {\cal R}_{0,n}}\phi_j^*/\sigma^2_{0,n}(Y-{Q_{0,n}}) \phi_j^*(x)+o_P(1), \label{keyexpansion} 
\end{eqnarray}
where remainder is $o_P(1)$ uniformly in $x$.
For a given $x$, the leading term is now  a sum of independent mean zero random variables $D_{Q_{0,n},x}(O)$ (\ref{Df0nx}) so that the central limit theorem can be applied. 
The variance is given by:
\begin{equation}\label{sigma2na}
\tilde{\sigma}^2_{0,n}(x)=
\frac{1}{d_{0,n}}\sum_{j\in {\cal R}_{0,n}}\{\phi_{j}^*(x)\}^2.\end{equation}
Thus, $(n/d_{0,n})^{1/2}(Q_n-Q_{0,n})(x)/\tilde{\sigma}_{0,n}(x)\Rightarrow_d N(0,1)$ for all $x$.

We have 
\[
\sup_{x\in [0,1]^d}(n/d_{0,n})^{1/2}\mid {Q_n}-{Q_{0,n}}\mid (x)/\tilde{\sigma}_{0,n}(x)= o_P(\log n),
\]
and  
\[
\sup_{x\in [0,1]^d}(n/d_{0,n})^{1/2}\mid Q_n-Q_0\mid(x)/\tilde{\sigma}_{0,n}(x)=o_P(\log n).\]
Note that, by choosing $J_{0,n}=n^{1/(2k^*+1)}\log^{m_1}n$ for some specified $m_1$,  this implies that for some specified $m<\infty$, 
\[
\pl  Q_n-Q_0\pl_{\infty}=O_P(n^{-k^*/(2k^*+1)}\log^m n) .\]
\end{theorem}
{\bf Proof:}
We have that $(n/d_{0,n})^{1/2}(Q_n-Q_{0,n})(x)/\tilde{\sigma}_{0,n}(x)=n^{1/2}(P_n-P_0) D^*_{Q_{0,n},x}+R_n(x)$, where $\sup_x \mid R_n(x)\mid=o_P(1)$ and 
$D^*_{Q_{0,n},x}=d_{0,n}^{-1/2}\sum_{j\in {\cal R}_{0,n}}\phi_j^*(Y-Q_{0,n})\phi_j^*(x) /\tilde{\sigma}_{0,n}(x)$.
Note $D^*_{Q_{0,n},x}$ has mean zero and variance $1$.
$d_{0,n}^{-1/2}D^*_{Q_{0,n},x}$ is included in the fixed $d$-dimensional class
\[
{\cal F}_{n}=\left\{ \frac{1}{d_{0,n}}\sum_{j\in {\cal R}_{0,n}}\phi_j^*(Y-{Q_{0,n}} )\phi_j^*(x)/\tilde{\sigma}_{0,n}(x): x\right\}.\]
That is, these functions are indexed by $x\in [0,1]^d$.
Let $g_{x}$ be one element in ${\cal F}_n$. Then
\[
\sup_x P_0 g_{x}^2=\frac{1}{d_{0,n}},\]
due to our extra dividing by $d_{0,n}^{-1/2}$.

Thus, $\sup_{Q\in {\cal F}_n}P_0 Q^2\sim d_{0,n}^{-1}$.
Moreover, ${\cal F}_n$ has  universally  bounded sup-norm and its covering number $N(\epsilon,{\cal F}_n,L^2)=O(\epsilon^{-d})$.
Therefore, the entropy integral $J(\delta,{\cal F}_n)\equiv \int_0^{\delta}\sup_Q\sqrt{\log N(\epsilon,{\cal F}_n,L^2(Q))} d\epsilon$ of $ {\cal F}_n$ is bounded by $\sim -\int_0^{\delta} (\log \epsilon)^{1/2} d\epsilon$.  This can be {\em conservatively} bounded by \[
\sim -\int_0^{\delta}\log \epsilon d\epsilon = -\delta \log \delta +\delta\sim -\delta \log \delta.\]
So we can state that $J(\delta,{\cal F}_n)=o(\delta \log \delta)$.
We can can set $\delta=\delta_n=d_{0,n}^{-1/2}$.

Applying empirical process inequalities for $E \sup_{Q\in {\cal F}_n,\pl Q\pl_{P_0}\leq \delta_n}\mid n^{1/2}(P_n-P_0)Q\mid \sim J(\delta_n,{\cal F}_n,L^2)$ with $\delta_n=d_{0,n}^{-1/2}$ yields now:
\[
\begin{array}{l}
\sup_x (n/d_{0,n})^{1/2}\mid {Q_n}-{Q_{0,n}}\mid (x)/\tilde{\sigma}_{0,n}(x)\\
\approx 
d_{0,n}^{1/2}\sup_x\left | n^{1/2}(P_n-P_0) \frac{1}{d_{0,n}}\sum_{j\in {\cal R}_{0,n}}\phi_j^*(Y-{Q_{0,n}}) \phi_j^*(x)/\tilde{\sigma}_{0,n}(x) \right | \\
 \leq
d_{0,n}^{1/2} \sup_{g\in {\cal F}_n,\pl g\pl_{P_0}\leq \delta_n} \mid n^{1/2}(P_n-P_0)g\mid \\
\sim  d_{0,n}^{1/2}O_P( J(\delta_n,{\cal F}_n,L^2))\\
\sim  d_{0,n}^{1/2} o_P( \delta_n \log \delta_n)\\
\sim o_P(\log n) . \Box
\end{array}
\]

Our sufficient conditions for Assumptions A1,A2 and A3 are also sufficient for A1*,A2* and A3*. Therefore we can state the analogue of Corollary \ref{corthasnormal}.
\begin{corollary}\label{corthasnormalu}
Consider the $k$-th order sieve MLE or relax HAL-MLE $Q_n$, $k\in \{1,2,\ldots\}$,
so that $P_n S_j(Q_n)=0$ with $S_j(Q_n)=\phi_j/\sigma^2_n (Y-Q_n)$ for all $j\in {\cal R}_n$. 
Assume that for some $M<\infty$, $Q_0,\sigma^2_0\in D^{(k)}_M([0,1]^d)$; 
and that $\sigma^2_n$ is a zero-order spline HAL-MLE so that
$\pl \sigma^2_n-\sigma^2_0\pl_{P_0}=O_P(n^{-1/3}(\log n)^{d/2})$.
Let $k^*=k+1$.
Assume that for some $\delta>0$ $\inf_{x\in [0,1]^d}\sigma^{2}_{0,n}(x)>\delta$ with probability tending to 1. 

Let ${\cal R}_0$ be a support set identifying a set of basis functions indexed by $(\bar{s}(k+1),u)$ in our $k$-th order spline representation  so that for a given $M<\infty$ $\{Q_n,Q_0\}\subset D^{(k)}({\cal R}_0)$ with probability tending to 1. 
Assume that ${\cal R}_n={\cal R}(P_n)$ satisfies $\max_{g\in D^{(k)}_M({\cal R}_0)} \min_{Q\in D^{(k)}({\cal R}_{n})}\pl g-Q\pl_{\infty}=O_P(C(M)r(d,J_{n})^{k+1})$. Then Assumption A0 holds.

Let $d_{0,n}=n^{1/(2k^*+1)}\log^{m_1} n$ for some $m_1$ with $m_1$ chosen so hat $r(d,J_{0,n})^{k+1}/(n/d_{0,n})^{-1/2}\rightarrow 0$. 

We have 
\[
\sup_{x\in [0,1]^d}(n/d_{0,n})^{1/2}\mid {Q_n}-{Q_{0,n}}\mid (x)/\tilde{\sigma}_{0,n}(x)= o_P(\log n),
\]
and  
\[
\sup_{x\in [0,1]^d}(n/d_{0,n})^{1/2}\mid Q_n-Q_0\mid(x)/\tilde{\sigma}_{0,n}(x)=o_P(\log n).\]
Note that, by choosing $J_{0,n}=n^{1/(2k^*+1)}\log^{m_1}n$ for some specified $m_1$,  this implies that for some specified $m<\infty$, 
\[
\pl  Q_n-Q_0\pl_{\infty}=O_P(n^{-k^*/(2k^*+1)}\log^m n) .\]
The same results apply to HAL-MLE if one verifies Assumption A1* with help of Lemma \ref{lemmatildernx}.
\end{corollary}

\subsection{Simultaneous confidence band}

{\bf Scaling up a pointwise confidence band ignoring correlations:}
We showed that for each $x$:
\[
(n/d_n)^{1/2}\tilde{\sigma}^{-1}_{0,n}({Q_n}-{Q_{0,n}})(x)\rightarrow_d N(0,1).\]
Therefore a pointwise $0.95$-confidence interval for $Q_{0,n}(x)$  is given by ${Q_n}(x)\pm 1.96 (n/d_n)^{-1/2}\tilde{\sigma}_{0,n}(x)$. By our result that the supremum norm converges at a rate faster than $\log n (n/d_n)^{-1/2}$, we can now make it an asymptotically valid $0.95$-simultaneous confidence band by multiplying it by $\log d_n\sim \log n$:
\[
{Q_n}(x)\pm 1.96 \log d_n (n/d_n)^{-1/2}\tilde{\sigma}_{0,n}(x).\]
Moreover, since $\pl Q_{0,n}-Q_{0}\pl_{\infty}=O(r(d,J_n)^{k+1})$, by choosing $J_n$ so that $(n/J_n)^{1/2}(\log n) r(d,J_n)^{k+1}=o(1)$, then this  is also a simultaneous confidence band for the true $Q_{0}$.

{\bf Constructing directly a simultaneous confidence band based on estimating correlation structure:}
Let $D_{n,P_0,x}\equiv 
d_n^{-1/2}\sum_{j\in {\cal R}_n}\phi_u^*(Y-{Q_{0,n}}) \phi_j^*(x)$ so that $({Q_n}-{Q_{0,n}})(x)\approx 1/n \sum_i D_{n,P_0,x}(O_i)$.
 As discussed in Subsection \ref{subsectinference},  we can replace this influence curve by the delta-method based influence curve (\ref{deltamethodic}) in terms of the original basis.
Let $\sigma_n(x)=\tilde{\sigma}_n^2(x)/n$ be the corresponding estimator of $\tilde{\sigma}^2_{0,n}/n$, where $\sigma_{0,n}(x)=\tilde{\sigma}^2_{0,n}(x)=P_0D_{n,P_0,x}^2$.  Note that this variance estimator  behaves as $d_n/n$.
We could  construct a simultaneous confidence band by estimating the correlation matrix of 
the standardized estimator $(({Q_n}-{Q_{0,n}})(x_m)/\sigma_n(x_m): m=1,\ldots,M)$ ( scaled to converge pointwise to $N(0,1)$) based on the empirical correlation matrix $\rho_n$ of the vector influence curve  $(D^*_{n,P_0,x_m}:m=1,\ldots,M)$ across the user supplied set of grid of time-points $x_1,\ldots,x_M$. This empirical correlation matrix $\rho_n$, for a given set of grid points, would converge to a diagonal covariance matrix as sample size increases in the case that the nonparametric uniform approximation condition on $D^{(k)}({\cal R}_n)$ holds. But even in that case $\rho_n$ will  demonstrate real correlations across the grid-points not far apart. We could now use as working model that the vector of standardized estimators is multivariate normal with mean zero and covariance matrix equal to this correlation matrix $\rho_n$. This then implies a simultaneous confidence band for ${Q_{0,n}}$ of form ${Q_n(x_m)}\pm q_n(0.95) \sigma_n(x_m)$, where $q_n(0.95)$ is the 0.95-quantile of a large sample  of $\max_m \mid Z(m)\mid$ with $Z\sim N(0,\rho_n)$.  We are claiming that this will still be a valid simultaneous confidence band for $Q_{0,n}$ since the correlation matrix will approximate independence across points far enough apart, so that $q_n(0.95)$ will be of the order of $\log n$, thereby behaving as the above proven simultaneous confidence band.  We do not have a formal proof of this and simulations can shed more light on this conjecture. 

\subsection{Proof of asymptotic normality of sieve and HAL-MLEs of logistic regression}

Our proof of Theorem \ref{thasnormal} for least squares regression above mostly generalizes to logistic regression. For completeness, this analogue proof is given here. 
Let   $m_Q(x)=1/(1+\exp(-Q(x)))$; $Q_0(X)=\mbox{Logit}P_0(Y=1\mid X)$; $m_{Q_0}(x)=P_0(Y=1\mid X=x)$; $L(Q)(X,Y)=-Y\log m_Q(X)-(1-Y)\log(1-m_Q(X))$; $\frac{d}{dQ}L(Q)(\phi_j)=\phi_j(X)(Y-m_Q(X))$; $Q_n=\arg\min_{Q\in D^{(k)}({\cal R}_{0,n})}P_n L(Q)$; $Q_{0,n}=\arg\min_{Q\in D^{(k)}({\cal R}_{0,n})}P_0L(Q)$.  
Let $\{\phi_j:j\in {\cal R}_{0,n}\}$ be the original basis of $D^{(k)}({\cal R}_{0,n})$ and let $\{\phi_j^*: j\in {\cal R}_{0,n}\}$ be an orthonormal basis of $D^{(k)}({\cal R}_{0,n})$
in $L^2(P_0^*)$, where $dP_0^*=\sigma^2_{0,n}dP_0$, and $\sigma^2_{0,n}=m_{Q_{0,n}}(1-m_{Q_{0,n}})$. Let
\begin{equation}\label{tildernxl}
\tilde{r}_n(x)=\sum_{j\in {\cal R}_{0,n}}P_n \phi_j^* (Y-Q_n) \phi_j^*(x).
\end{equation}

As starting point for our analysis of $(Q_n-Q_{0,n})(x)$, using that $P_0 \phi_j^* (Y-m_{Q_{0,n}})=0$ and $r_n^*(j)=P_n \phi_j^*(Y-m_{Q_n})$, $j\in {\cal R}_{0,n}$,  we have
\[
P_0 \phi_j^*(Y-m_{Q_n})-P_0\phi_j^*(Y-m_{Q_{0,n}} )=-(P_n-P_0)\phi_j^*(Y-m_{{Q_n}})+r_n^*(j),\]
so that
\[
P_0 \phi_j^*(m_{Q_n}-m_{Q_{0,n} } )=(P_n-P_0)\phi_j^*(Y-m_{Q_n})-r_n^*(j).\]
We  linearize $m_{Q_n}-m_{Q_{0,n}}$ in $(Q_n-Q_{0,n})$: 
A tailor expansion of $m$ at $Q_{0,n}(x)$ shows that  $m(Q_n(x))-m(Q_{0,n}(x))=(m_{Q_{0,n}}(1-m_{Q_{0,n}})) (Q_n-Q_{0,n})(x)$ plus a second order term $R_{1n,a}(x)$ in $(Q_n-Q_{0,n})(x)$ that is $O((Q_n-Q_{0,n})^2(x))$.
Thus, $P_0\phi_j^*(m_{Q_n}-m_{Q_{0,n}})=P_0 \phi_j^* m_{Q_{0,n}}(1-m_{Q_{0,n}}) (Q_n-Q_{0,n})-R_{1,n}(\phi_j^*)$, where $R_{1,n}(\phi_j^*)=-P_0\phi_j^* R_{1n,a}$.
Thus $P_0\phi_j^*(m_{Q_n}-m_{Q_{0,n}})=P_0^* \phi_j^* (Q_n-Q_{0,n})-R_{1,n}(\phi_j^*)$.
We can write $R_{1,n}(\phi_j^*)=-P_0^*\phi_j^* \tilde{R}_{1n,a}$, where $\tilde{R}_{1n,a}=R_{1n,a}/\sigma^2_{0,n}$.

So at this point we have
\[
P_0^* \phi_j^*(Q_n-Q_{0,n})=(P_n-P_0)\phi_j^*(Y-m_{Q_n}) -r_n^*(j)+R_{1,n}(\phi_j^*).\]
From now on the analysis is identical to our analysis for mean regression. 
For sake of self-contained proof, we still proceed here.

Let $\tilde{Q}_n=\Pi_{J_{0,n}}(Q_n)$, where $\Pi_{J_{0,n}}(Q_n)=\arg\min_{Q\in D^{(k)}({\cal R}_{0,n})}P_0^*(f-Q_n)^2$ is the projection of $Q_n$ on the fixed working model $D^{(k)}({\cal R}_{0,n})$ in $L^2(P_0^*)$. 
We decompose $Q_n-Q_{0,n}=Q_n-\tilde{Q}_n+\tilde{Q}_n-Q_{0,n}$. Let
$R_{2,n}(\phi_j^*)=-P_0^*\phi_j^*(Q_n-\tilde{Q}_n)$. 
So we have
\[
P_0^*\phi_j^*(\tilde{Q}_n-Q_{0,n})=(P_n-P_0)\phi_j^*/\sigma^2_n(Y-Q_n)-r_n^*(j)+(R_{1,n}+R_{2,n})(\phi_j^*).\]

Since $\{\phi_j^*:j\in {\cal R}_{0,n}\}$, is an orthonormal basis of the linear subspace $D^{(k)}({\cal R}_{0,n})$ of $L^2(P_0^*)$, and $\tilde{Q}_n-Q_{0,n}\in D^{(k)}({\cal R}_{0,n})$, we have
\begin{eqnarray*}
(\tilde{Q}_n-Q_{0,n})(x)&=& \sum_{j\in {\cal R}_{0,n}}\{P_0^*(\tilde{Q}_n-Q_n) \phi_j^* \}\phi_j^*(x)\\
&=&\sum_{j\in {\cal R}_{0,n}}(P_n-P_0)\phi_j^*/\sigma^2_n(Y-{Q_n}) \phi_j^*(x)-\sum_{j\in {\cal R}_{0,n}}r_n^*(j) \phi_j^*(x)\\
&& +\sum_{j\in {\cal R}_{0,n}}(R_{1,n}+R_{2,n})(\phi_j^*)\phi_j^*(x)\\
&=&\sum_{j\in {\cal R}_{0,n}}(P_n-P_0)\phi_j^*/\sigma^2_n(Y-{Q_n}) \phi_j^*(x)-
\tilde{r}_n(x)\\
&&+\sum_{j\in {\cal R}_{0,n}}(R_{1,n}+R_{2,n})(\phi_j^*)\phi_j^*(x).
\end{eqnarray*}
Let's consider the term $R_{2,n}(x)\equiv \sum_{j\in {\cal R}_{0,n}}R_{2,n}(\phi_j^*)\phi_j^*(x)$.
Note that this term equals $\sum_{j\in {\cal R}_{0,n}}P_0^*\phi_j^*(Q_n-\tilde{Q}_n)\phi_j^*(x)$ and thus equals the projection of $(Q_n-\tilde{Q}_n)$ onto $D^{(k)}({\cal R}_{0,n})$. So, by definition of $\tilde{Q}_n$,  this equals 
$\Pi_{J_{0,n}}(Q_n)-\tilde{Q}_n=0$. So $R_{2,n}(x)=0$.

Let 
\begin{equation}\label{R1nxl}
R_{1,n}(x)\equiv\sum_{j\in {\cal R}_{0,n}}R_{1,n}(\phi_j^*)\phi_j^*(x).\end{equation}
We note that
\[
R_{1,n}(x)=-\Pi_{J_{0,n}}(\tilde{R}_{1n,a}).\]

Thus, we have \[
(\tilde{Q}_n-Q_{0,n})(x)= \sum_{j\in {\cal R}_{0,n}}(P_n-P_0)\phi_j^*(Y-{Q_n}) \phi_j^*(x)-
\tilde{r}_n(x)+R_{1,n}(x).\]

We now write the left-hand side as $\tilde{Q}_n-Q_n+(Q_n-Q_{0,n})$. We already showed that $\pl \tilde{Q}_n-Q_n\pl_{\infty}=O_P(C(M_n)r(d,J_{0,n})^{k+1})$.
So we then have
\begin{eqnarray*}
(Q_n-Q_{0,n})(x)&=& \sum_{j\in {\cal R}_{0,n}}(P_n-P_0)\phi_j^*(Y-{Q_n}) \phi_j^*(x)-
\tilde{r}_n(x)+R_{1,n}(x)\\
&&+O_P(C(M_n)r(d,J_{0,n})^{k+1}).\end{eqnarray*}

Let 
\begin{equation}\label{Enxl}
E_n(x)\equiv \sum_{j\in {\cal R}_{0,n}}(P_n-P_0)\phi_j^* (Q_n-Q_{0,n})   \phi_j^*(x).\end{equation}
Then, we have
\begin{eqnarray*}
(n/d_{0,n})^{1/2}(Q_n-Q_{0,n})(x)&=& n^{1/2}(P_n-P_0) d_{0,n}^{-1/2}\sum_{j\in {\cal R}_{0,n}}\phi_j^*(Y-{Q_{0,n}}) \phi_j^*(x)\\
&&\hspace*{-5cm}-
-(n/d_{0,n})^{1/2}\tilde{r}_n(x)+(n/d_{0,n})^{1/2}R_{1,n}(x) +(n/d_{0,n})^{1/2}O_P(C(M_n)r(d,J_{0,n})^{k+1})\\
&&-(n/d_{0,n})^{1/2}E_n(x).\end{eqnarray*}

The leading term is a sum of independent mean zero random variables \begin{equation}\label{Df0nxl}
D_{Q_{0,n},x}(O)\equiv d_{0,n}^{-1/2}\sum_{j\in {\cal R}_{0,n}}\phi_j^*(Y-{Q_{0,n}}) \phi_j^*(x).\end{equation}
The variance  of $D_{Q_{0,n},x}(O)$ is given by
\[
\sigma^2_{0,n}(x)=\frac{1}{d_{0,n}}\sum_{j_1,j_2\in {\cal R}_{0,n}}P_0\{ \phi_{j_1}^*\phi_{j_2}^*\sigma^2_{0,n}\}\phi_{j_1}^*(x)\phi_{j_2}^*(x).\]
 Due to 
the orthonormality of $\{\phi_j^*: j\in {\cal R}_{0,n}\}$ in $L^2(P_0^*)$  with $dP_0^*=\sigma^{2}_{0,n}dP_0$, this variance simplifies to
\begin{equation}\label{sigma2nb}
\tilde{\sigma}^2_{0,n}(x)=
\frac{1}{d_{0,n}}\sum_{j\in {\cal R}_{0,n}}\{\phi_{j}^*(x)\}^2.\end{equation}
Lemma \ref{lemmaexplicit} shows that under the nonparametric uniform approximation condition on $D^{(k)}({\cal R}_{0,n})$, this expression approximates $1/p_0^*(x)=\sigma^{-2}_0(x)p_0(x)^{-1}=\{p_0(x)Q_0(x)(1-Q_0(x))\}^{-1}$. If the model $D^{(k)}({\cal R}_{0,n})$ is aimed to approximate a real subset $D^{(k)}({\cal R}_0)$ of $D^{(k)}([0,1]^d)$, thereby allowing for extrapolation towards local behavior at $x$ (e.g., an additive model), then this asymptotic variance would be smaller and the same lemma provides an expression for that as well.

\subsection{Asymptotic normality theorem for logistic regression}
Consider/recall the following definitions:
\begin{eqnarray}
\tilde{r}_n(x)&\equiv& \sum_{j\in {\cal R}_{0,n}}\{P_n \phi_j^*(Y-Q_n)\}\phi_j^*(x)\label{tildernxlf}\\
R_{n,1}(\phi_j^*)&\equiv& -P_0 \phi_j^*\tilde{R}_{1n,a}\nonumber\\
\tilde{R}_{1n,a}&\equiv& R_{1n,a}/\sigma^2_{0,n}\nonumber\\
R_{1n,a}&\equiv& m(Q_n)(x)-m(Q_{0,n})(x)-(m_{Q_{0,n}}(1-m_{Q_{0,n}}))(x) (Q_n-Q_{0,n})(x)\nonumber \\
R_{n,1}(x)&\equiv&  \sum_{j\in {\cal R}_{0,n}}R_{n,1}(\phi_j^*)\phi_j^*(x)\label{R1nxlf}\\
E_n(x)&\equiv& \sum_{j\in {\cal R}_{0,n}}(P_n-P_0)\phi_j^* (Q_n-Q_{0,n})   \phi_j^*(x) \label{Enxlf}\\
D_{Q_{0,n},x}&\equiv& d_{0,n}^{-1/2}\sum_{j\in {\cal R}_{0,n}}\phi_j^*(Y-{Q_{0,n}}) \phi_j^*(x).\label{Df0nxlf}
\end{eqnarray}
Let ${\cal R}_{0,n}={\cal R}(P_n^{\#})$ the independent copy of ${\cal R}_n={\cal R}(P_n)$ with $P_n^{\#}$ empirical measure of independent (from $P_n$) i.i.d. sample from $P_0$.\newline
{\bf Assumption B0:}
\begin{equation}\label{B0}
\sup_{Q\in \{Q_n,Q_0\}}\inf_{g\in D^{(k)}({\cal R}_{n})}\pl Q-g\pl_{\infty}=O_P(C(M_n)r(d,J_{0,n})^{k+1}),\end{equation}
 where $M_n=\pl Q_n\pl_{v,k}^*$, the $L_1$-norm of its non-zero coefficients. \newline
{\bf Assumption B1:}
\begin{equation}\label{undersmoothml}
\tilde{r}_n(x)=o_P((n/d_{0,n})^{-1/2}).
\end{equation}
{\bf Assumption B2:}
\begin{equation}
\label{B2}
R_{n,1}(x)=o_P((n/d_{0,n})^{-1/2}).
\end{equation}
{\bf Assumption B3:}\begin{equation}\label{keyconditiona1ml}
E_n(x)=o_P((n/d_{0,n})^{-1/2}).\end{equation}
{\bf Assumption B4:}
\begin{equation}\label{neglbiasl}
C(M_n)r(d,J_{0,n})^{k+1}=o((n/d_{0,n})^{-1/2}).
\end{equation}

This proves the following theorem. 
\begin{theorem}\label{thasnormal1}\ \nl
Consider the $k$-th order sieve MLE, HAL-MLE or relax HAL-MLE $Q_n$ of $Q_0=\mbox{Logit}P_0(Y=1\mid X)$.  Let $k^*=k+1$.
Let ${\cal R}_n$ be the set of $d_n$ non-zero coefficients in $Q_n=\sum_{j\in {\cal R}_n}\beta_n(j)\phi_j$. Assume $\inf_{x\in [0,1]^d}m_{Q_0(x)}>0$. 
Recall the definitions of $\tilde{r}_n(x)$ (\ref{tildernxlf}), $R_{1,n}(x)$ (\ref{R1nxlf}), $E_n(x)$ (\ref{Enxlf}), $D_{Q_{0,n},x}(O)$ (\ref{Df0nxlf}) above. 
Assume Assumptions B0,B1,B2,B3 and B4.

Due to Assumption B0, we have the following expansion:
\begin{eqnarray}
(n/d_{0,n})^{1/2}(Q_n-Q_{0,n})(x)&=& n^{1/2}(P_n-P_0) d_{0,n}^{-1/2}\sum_{j\in {\cal R}_{0,n}}\phi_j^*(Y-{Q_{0,n}}) \phi_j^*(x)\nonumber\\
&&-(n/d_{0,n})^{1/2}E_n(x) -
(n/d_{0,n})^{1/2}\tilde{r}_n(x)
\nonumber \\
&&\hspace*{-3cm}+(n/d_{0,n})^{1/2}R_{1,n}(x) +(n/d_{0,n})^{1/2}O_P(C(M_n)r(d,J_{0,n})^{k+1}),\label{keyexpansion}\end{eqnarray}
where $M_n=\pl Q_n\pl_{v,k}^*$, which equals the $L_1$-norm  of the vector of non-zero coefficients in $Q_n$.


By Assumption B1-B4, we have
\begin{eqnarray}
\frac{(Q_n-Q_{0,n})(x)}{(n/d_{0,n})^{-1/2}}= n^{1/2}(P_n-P_0) d_{0,n}^{-1/2}\sum_{j\in {\cal R}_{0,n}}\phi_j^*(Y-{Q_{0,n}}) \phi_j^*(x)+o_P(1). \label{keyexpansion} 
\end{eqnarray}

For a given $x$, the leading term is now  a sum of independent mean zero random variables $D_{Q_{0,n},x}(O)$ (\ref{Df0nx}) so that the central limit theorem can be applied. 
The variance is given by:
\begin{equation}\label{sigma2nc}
\tilde{\sigma}^2_{0,n}(x)=
\frac{1}{d_{0,n}}\sum_{j\in {\cal R}_{0,n}}\{\phi_{j}^*(x)\}^2.\end{equation}



{\bf Conclusion:}
We have
\[
\tilde{\sigma}^{-1}_{0,n} (n/d_{0,n})^{1/2}({Q_n}-{Q_{0,n}})(x)\Rightarrow_d N(0,1),\]
and
\[
\tilde{\sigma}^{-1}_{0,n}(n/d_{0,n})^{1/2}(Q_n-Q_0)(x)\Rightarrow_d N(0,1).\]
Assuming $M_n=O(\log^{m}n)$ for some $m<\infty$, by choosing $J_{0,n}=n^{1/(2k^*+1)}\log^{m_1} n$ for some  $m_1<\infty$,  it follows that for  some  $m<\infty$, \[
\mid Q_n-Q_0\mid(x)=O_P(n^{-k^*/(2k^*+1)}\log^m n).\]

By the $\delta$-method, these results imply the analogue results for $(m_{Q_n}-m_{Q_{0,n}})(x)$ and $(m_{Q_n}-m_{Q_0})(x)$, but with  $D_{Q_{0,n},x}$ replaced by $m_{Q_0}(1-m_{Q_0})(x) D_{Q_{0,n},x}$. Thus, in this case the asymptotic variance of $(m_{Q_n}-m_{Q_{0,n}})(x)$ is given by  \[
 \tilde{\sigma}^2_{0,n}(x)=\frac{\sigma^2_0(x)}{p_0(x)} \Pi_{J_{0,n}}(K((\cdot-x)/h)(x)\]
 with $\sigma^2_0=m_{Q_0}(1-m_{Q_0})$. 
\end{theorem}
{\bf Sufficient condition for B2 and B3: }
Analogue to the least squares regression case, we have that $R_{1n}(x)$ and $E_n(x)$ are $o_P((n/d_{0,n})^{-1/2})$ if $d_{0,n}=O_P(n^{1/3-\delta})$ for some $\delta>0$. 
\nl
{\bf Sufficient condition for B1:}
Recall $S_j^*(Q_n)=\phi_j^*(Y-m_{Q_n}(X))$, $j\in {\cal R}_{0,n}$. Let $\tilde{S}_j(Q_n)$ be a linear combination of scores solved by $Q_n$ so that $P_n \tilde{S}_j(Q_n)=0$ approximating $S_j^*(Q_n)$ for $j\in {\cal R}_{0,n}$.  Here $S_j^*(Q_n)-\tilde{S}_j(Q_n)=(\phi_j^*-\tilde{\phi}_j)(Y-m_{Q_n}(X))$ for a $\tilde{\phi}_j$, $j\in {\cal R}_{0,n}$. Lemma  \ref{lemmatildernx} shows  \[
\tilde{r}_n(x)=\sum_{j\in {\cal R}_{0,n}}(P_n-\tilde{P}_n)\{S_j^*(Q_n)-\tilde{S}_j(Q_n)\} \phi_j^*(x).\]
Let 
$\delta_n\equiv  \pl d_{0,n}^{-1}\sum_{j\in {\cal R}_{0,n}}\{S_j(Q_n)-\tilde{S}_j(Q_n)\}\pl_{\mu}$.
If $d_{0,n}^{1/2}\delta_n^{(2k+1)/(2k+2)}=O_P(n^{-\delta})$ for some $\delta>0$, then Assumption B1 (\ref{undersmoothml})  holds.


{\bf Explicit asymptotic variance under nonparametric uniform approximation assumption:}
Lemma \ref{lemmaexplicit} shows that 
 \[
 \tilde{\sigma}^2_{0,n}(x)=\{p_0^*(x)\}^{-1}  \Pi_{J_{0,n}}(K((\cdot-x)/h)(x)+o_P(1),\]
 where $\Pi_{J_{0,n}}$ denotes  the projection operator onto  space spanned by $\{\phi_j^*:j\in {\cal R}_{0,n}\}$, and  $K((\cdot-x)/h)$ is a $d$-variate kernel  satisfying $K(0)=1$, $\int K(x)d\mu(x)=1$, and bandwidth $h=(1/J_{0,n})^{1/d}$. 
 Under the nonparametric uniform approximation condition on $D^{(k)}_{M}({\cal R}_{0,n})$ we have that $\tilde{\sigma}^2_{0,n}(x)= 1/p_0^*(x)+o_P(1)=\{Q_0(1-Q_0)p_0(x)\}^{-1}+o_P(1)$.

\begin{corollary}\label{corthasnormall}
Consider the $k$-th order sieve MLE or relax HAL-MLE $Q_n$, $k\in \{1,2,\ldots\}$,
so that $P_n S_j(Q_n)=0$ with $S_j(Q_n)=\phi_j (Y-m_{Q_n})$ for all $j\in {\cal R}_n$. 
Assume that for some $M<\infty$, $Q_0\in D^{(k)}_M([0,1]^d)$ and assume $\inf_{x\in [0,1]^d}m_{Q_0}(x)>0$.
Let $k^*=k+1$.

Let ${\cal R}_0$ be a support set identifying a set of basis functions indexed by $(\bar{s}(k+1),u)$ in our $k$-th order spline representation  so that for a given $M<\infty$ $\{Q_n,Q_0\}\subset D^{(k)}({\cal R}_0)$ with probability tending to 1. 
Assume that ${\cal R}_n={\cal R}(P_n)$ satisfies $\max_{g\in D^{(k)}_M({\cal R}_0)} \min_{Q\in D^{(k)}({\cal R}_{n})}\pl g-Q\pl_{\infty}=O_P(C(M)r(d,J_{n})^{k+1})$. Then Assumption B0 holds.

Let $d_{0,n}=n^{1/(2k^*+1)}\log^{m_1} n$ for some $m_1$ with $m_1$ chosen so hat $r(d,J_{0,n})^{k+1}/(n/d_{0,n})^{-1/2}\rightarrow 0$. 
Then, 
\[
\tilde{\sigma}^{-1}_{0,n} (n/d_{0,n})^{1/2}({Q_n}-{Q_{0,n}})(x)\Rightarrow_d N(0,1),\]
and
\[
\tilde{\sigma}^{-1}_{0,n}(n/d_{0,n})^{1/2}(Q_n-Q_0)(x)\Rightarrow_d N(0,1).\]
Moreover,  for some $m<\infty$, 
\[
\mid Q_n-Q_0\mid(x)=O_P(n^{-k^*/(2k^*+1)}\log^m n).\]
These results also apply to the HAL-MLE if one verifies B1 with help of Lemma \ref{lemmatildernx}.\end{corollary}

\section{Asymptotic linearity (and super-efficiency) of plug-in HAL-MLE, relax-HAL-MLE, and sieve-MLE of  pathwise differentiable features of regression function $Q_0$.}\label{section9}
The purpose of this section is the analysis of plug-in estimator $\Phi(Q_n)-\Phi(Q_{0,n})$ and $\Phi(Q_n)-\Phi(Q_0)$ for a pathwise differentiable target parameter $P\rightarrow \Phi(Q(P))$. 
Consider the inverse weighted least squares regression setting of previous section. 
In this section we can replace ${\cal R}_{0,n}$ by ${\cal R}_n$ since there will not be an issue to deal with the data dependence of ${\cal R}_{n}$. As a consequence, we then have that $Q_n-Q_{0,n}\in D^{(k)}({\cal R}_{0,n})$, thereby also avoiding having to deal with the projection operator $\Pi_{J_{0,n}}$ in our proof of Theorem \ref{thasnormal}.

Due to our ${\cal R}_{0,n}$ is set equal to ${\cal R}_n$ in this proof, our starting equality is 
\begin{eqnarray*}
(Q_n-Q_{0,n})(x)&=&\sum_{j\in {\cal R}_{0,n}}(P_n-P_0)\phi_j^*\sigma^{-2}_n(Y-{Q_n}) \phi_j^*(x)\\
&&-\sum_{j\in {\cal R}_{0,n}}\{P_n \sigma^{-2}_n \phi_j^*(Y-Q_n)\} \phi^*_j(x)+R_{1,n}(x),
\end{eqnarray*}
where $R_{1,n}(x)$ is defined by (\ref{R1nx}), an easy second order remainder, essentially only relying on consistency of $\sigma^2_n(x)$ to $\sigma^2_0(x)$.

{\bf Analyzing linear functions $\Phi$:}
Firstly, consider the case that $\Phi(Q)=\int h_0(x) Q(x) dP_0(x)$ for some specified $h_0$ possibly depending on $P_0$. Applying this linear operator $\Phi$ on both sides implies.
\[
\begin{array}{l}
\Phi(Q_n)-\Phi(Q_{0,n})=(P_n-P_0) \left\{\sum_{j\in {\cal R}_{0,n}}P_0 (h_0\phi_j^*) \phi_j^*(X)\right\}
\sigma^{-2}_n(Y-Q_n)\\
+P_n \left\{ \sum_{j\in {\cal R}_{0,n}}P_0( h_0\phi_j^*)\phi_j^*(X)\right\} \sigma^{-2}_n(Y-Q_n)+
P_0 h_0 R_{1,n}.
\end{array}
\]
We can write $P_0 h_0\phi_j^*=P_0^* \{(h_0\sigma^2_{0,n})\phi_j^*\}$. 
Let $h^*_{0,n}=h_0\sigma^2_{0,n}$ and let $h_0^*=h_0\sigma^2_0$ be its limit.
Therefore, we can also write this as:
\[
\begin{array}{l}
\Phi(Q_n)-\Phi(Q_{0,n})=(P_n-P_0) \left\{\sum_{j\in {\cal R}_{0,n}}P_0^* (h_{0,n}^*\phi_j^*) \phi_j^*(X)\right\}\sigma^{-2}_n(Y-Q_n)\\
+P_n \left\{ \sum_{j\in {\cal R}_{0,n}}P_0^*( h_{0,n}^*\phi_j^*)\phi_j^*(X)\right\}\sigma^{-2}_n (Y-Q_n)-P_0 h_0 R_{1,n}.
\end{array}
\]
Due to $\phi_j^*$ being an orthonormal basis in $L^2(P_0^*)$, note that
\[
h_{0,n}(X)\equiv \sum_{j\in {\cal R}_{0,n}}P_0^* (h_{0,n}^*\phi_j^*) \phi_j^*(X)\]
is just the projection of $h_{0,n}^*$ onto the linear space spanned by $\{\phi_j^*:j\in {\cal R}_{0,n}\}\subset L^2(P_0^*)$.
Let $\tilde{h}_0$ represent an $L^2(P_0^*)$-limit of $h_{0,n}$. Note that $\tilde{h}_0=h_0^*=\sigma^2_0h_0$ in the case that ${\cal R}_{0,n}$ satisfies the nonparametric uniform approximation assumption but could be the projection onto a smaller space $D^{(k)}({\cal R}_0)$ if it only satisfies the adaptive uniform approximation assumption. 
So we have
\begin{eqnarray*}
\Phi(Q_n)-\Phi(Q_{0,n})&=&(P_n-P_0) h_{0,n}\sigma^{-2}_n(Y-Q_n)+P_n h_{0,n}\sigma^{-2}_n(Y-Q_n)\\
&&-P_0 h_0 R_{1,n}.\end{eqnarray*}
Thus, efficient influence curve condition for $D^{(k)}({\cal R}_{0,n})$  is that $P_n h_{0,n}\sigma^{-2}_n(Y-Q_n)=o_P(n^{-1/2})$, which holds at $0$ for the sieve and relax HAL-MLE, and represents a weak undersmoothing condition  for the HAL-MLE. 
We assume the second order remainder condition
$(P_n-P_0)\{h_{0,n}\sigma^{-2}_n(Y-Q_n)-\tilde{h}_0\sigma^{-2}_0(Y-Q_0)\} =o_P(n^{-1/2})$,
which holds if $Q_n,Q_0,\sigma^2_n,\sigma^2_0\in D^{(0)}_M([0,1]^d)$ and $d_0(Q_n,Q_0)\rightarrow_p 0$, and $\pl \sigma^2_n-\sigma^2_0\pl_{P_0}\rightarrow_p 0$.
 Thus, this is  a very weak condition given the known rates of convergence of $Q_n,\sigma^2_n$, and even convergence in sectional variation norm (Appendix \ref{AppendixM}).  In addition, we assume $P_0 h_0 R_{1,n}=o_P(n^{-1/2})$. If  sup-norm of $\sigma^2_{-2n},\sigma^2_{0,n}, Q_n,Q_{0,n}$ are uniformly bounded, then  $P_0 h_0R_{1,n}$ can be bounded by the  $L^1(P_0)$-norm $\pl (\sigma^{-2}_n-\sigma^{-2}_{0,n})(Q_n-Q_{0,n})\pl_{P_0,1}$.
 So is suffices the assume the latter is $o_P(n^{-1/2})$.   
This then proves $\Phi(Q_n)-\Phi(Q_{0,n})=P_n \tilde{h}_0\sigma^{-2}_0(Y-Q_0)+o_P(n^{-1/2})$.

It then remains to analyze $\Phi(Q_{0,n})-\Phi(Q_0)=\int h_0 (Q_{0,n}-Q_0) dP_0$. However, 
\[
P_0 h_0(Q_{0,n}-Q_0)=-P_0 h_0(Y-Q_{0,n})=-P_0^* h_0\sigma^2_{0,n}(Y-Q_{0,n}),\]
and we know that $P_0^* \phi_j^*(Y-Q_{0,n})=0$ for all $j\in {\cal R}_{0,n}$ (i.e., the score equations of the MLE $Q_{0,n}$ are given by $P_0 \phi_j/\sigma^2_{0,n}(Y-Q_{0,n})$, and $\phi_j^*$ is a linear combination of $\{\phi_j: j\in {\cal R}_{0,n}\}$). Recall $h_{0,n}=\Pi(h_0\sigma^2_{0,n}\mid D^{(k)}({\cal R}_{0,n}))$ is the projection of $h_0\sigma^2_{0,n}$ onto the linear span of $\phi_j^*$, $j\in {\cal R}_{0,n}$, in $L^2(P_0^*)$. Then, we have $P_0^* h_{0,n} (Y-Q_{0,n})=0$. 
Thus, 
\[
P_0 h_0(Q_{0,n}-Q_0)=P_0^*({h}_{0,n}-h_0\sigma^2_{0,n})(Q_0-Q_{0,n}).\]
So, if $P_0^*(h_{0,n}-h_0\sigma^2_{0,n})(Q_0-Q_{0,n})=o_P(n^{-1/2})$, then $\Phi(Q_{0,n})-\Phi(Q_0)=o_P(n^{-1/2})$.
We could simply bound this term in $L^2(P_0^*)$-norms of $h_{0,n}-h_0\sigma^2_{0,n}$ and $Q_{0,n}-Q_0$, but that is still stronger than needed, as we show now. 

Suppose that $Q_0,Q_{0,n}\in D^{(k)}({\cal R}_0)$ for a support set ${\cal R}_0$ containing ${\cal R}_{0,n}$ with probability tending to 1. Note that this support set just has to be rich enough so that it spans $Q_0$. Then, $P_0^*({h}_{0,n}-h_0\sigma^2_{0,n})(Q_0-Q_{0,n})=
P_0^*\{\Pi({h}_{0,n}-h_0\sigma^2_{0,n}\mid D^{(k)}({\cal R}_0))\} (Q_0-Q_{0,n})$, where the projection operator is in $L^2(P_0^*)$. Note also that $\Pi({h}_{0,n}\mid D^{(k)}({\cal R}_0))={h}_{0,n}$. Let $\tilde{h}_0=\Pi(h_0\sigma^2_{0,n}\mid D^{(k)}({\cal R}_0)$. So we need
\[
P_0 ({h}_{0,n}-\tilde{h}_0)(Q_0-Q_{0,n})=o_P(n^{-1/2}).\]
We could assume $\pl {h}_{0,n}-\tilde{h}_0\pl_{P_0}\pl Q_{0,n}-Q_0\pl_{P_0}=o_P(n^{-1/2})$. We have $\pl Q_{0,n}-Q_0\pl_{\infty}=O(r(d,J_{0,n})^{k+1})$ under the adaptive uniform approximation assumption on $D^{(k)}({\cal R}_{0,n})$. If $h_0,\sigma^2_0\in D^{(k)}_M([0,1]^d)$ for some $M<\infty$, so that we have $\tilde{h}_0\in D^{(k)}_M({\cal R}_0)$, then it follows that we also have $\pl \tilde{h}_{0,n}-\tilde{h}_0\pl_{\infty}=O_P(r(d,J_{0,n})^{k+1})$. 
This proves the following theorem. 
\begin{theorem}\label{theoremsmoothregr}
Let ${\cal R}_{0,n}={\cal R}_n$ and consider setting of Theorem \ref{thasnormal}. Let $k\in \{1,\ldots\}$.
We have 
\begin{eqnarray*}
(Q_n-Q_{0,n})(x)&=&\sum_{j\in {\cal R}_{0,n}}(P_n-P_0)\phi_j^*\sigma^{-2}_n(Y-{Q_n}) \phi_j^*(x)\\
&&-
\sum_{j\in {\cal R}_{0,n}}\{P_n \phi_j^*\sigma^{-2}_n(Y-Q_n)\} \phi^*_{j,0}(x)+R_{1,n}(x).
\end{eqnarray*} 
Consider $\Phi(Q)=P_0 h_0Q$ for a given $h_0$.

{\bf Definitions:} Let \[
h_{0,n}(X)\equiv \sum_{j\in {\cal R}_{0,n}}P_0^* (h_{0,n}^*\phi_j^*) \phi_j^*(X),\]
where $h_{0,n}^*=h_0\sigma^2_{0,n}$ and note that $h_{0,n}=\Pi(h_{0,n}^*\mid D^{(k)}({\cal R}_{0,n})$ in $L^2(P_0^*)$. Let $h_0^*=h_0\sigma^2_0$ be the $L^2(P_0^*)$ limit of $h_{0,n}^*$.
Let $\tilde{h}_0$ represent an $L^2(P_0^*)$-limit of $h_{0,n}$. 
We have
\begin{eqnarray*}
\Phi(Q_n)-\Phi(Q_{0,n})&=&(P_n-P_0) h_{0,n}\sigma^{-2}_n(Y-Q_n)-P_n h_{0,n}\sigma^{-2}_n(Y-Q_n)\\
&&+P_0 h_0R_{1,n}.
\end{eqnarray*}
{\bf Assumptions:}
$\inf_x \sigma^2_0(x)>\delta>0$ for some $\delta>0$;  the adaptive uniform approximation condition on $D^{(k)}({\cal R}_{0,n})$ w.r.t. limit model $D^{(k)}({\cal R}_0)$, where ${\cal R}_0$ is such that $Q_n,Q_{0,n},Q_0\in D^{(k)}({\cal R}_0)$ with probability tending to 1; $h_0,\sigma^2_0\in D^{(k)}_M([0,1]^d)$ for some $M<\infty$; score equation $P_n h_{0,n}\sigma^{-2}_n(Y-Q_n)=o_P(n^{-1/2})$; 
$Q_n,Q_{0,n},Q_0,\sigma^2_n,\sigma^2_{0,n}\in D^{(0)}_M([0,1]^d)$; $\pl Q_n-Q_{0,n}\pl_{P_0}=o_P(1)$; $\pl \sigma^2_n-\sigma^2_0\pl_{P_0}=o_P(1)$; $\pl (\sigma^{-2}_n-\sigma^{-2}_{0,n})(Q_n-Q_{0,n})\pl_{P_0,1}=o_P(n^{-1/2})$; $r(d,J_{0,n})^{2(k+1)}=o_P(n^{-1/2})$.

{\bf Conclusion:}
Then, $\Phi(Q_n)-\Phi(Q_{0,n})=P_n \tilde{h}_0\sigma^{-2}_0(Y-Q_0)+o_P(n^{-1/2})$.
We also have $\Phi(Q_{0,n})-\Phi(Q_0)=O_P(r(d,J_{0,n})^{2(k+1)})$. 
Therefore, $\Phi(Q_n)$ is an asymptotically linear estimator of $\Phi(Q_0)$ at $P_0$ with influence curve 
$D_{P_0}(O)=\tilde{h}_0\sigma^{-2}_0(X)(Y-Q_0(X))$.
\end{theorem}
The conditions allow for a large degree of undersmoothing, as long as the sectional variation norm $\pl Q_n\pl_{v,0}^*=O_P(1)$ is controlled and the product of the rates of convergence of $Q_n-Q_0$  and $\sigma^2_n-\sigma^2_0$ in $L^2(P_0)$-norm  is $o_P(n^{-1/2})$. In addition, if one chooses the $d_{0,n}$ to optimize the rate of convergence w.r.t. $Q_0$, then the conditions of this theorem hold as well. 

So this proves asymptotic linearity of linear functions $\Phi$ of $Q_n$ as estimator of $\Phi(Q_0)$, under very mild conditions. 
We note that if $\tilde{h}_0=h_0\sigma^2_0$, which is the case if $D^{(k)}({\cal R}_{0,n})$ satisfies the nonparametric uniform approximation condition, then, $D_{P_0}$ is the efficient influence curve of $\Phi(Q_0)$ for the nonparametric model given by $h_0 (Y-Q_0)$.
However,  if ${\cal R}_{0,n}$ grows to a support set ${\cal R}_0$ so that $D^{(k)}({\cal R}_0)$ is a real subset of $D^{(k)}([0,1]^d)$, then ${h}_{0,n}$ would converge to the projection of $h_0\sigma^2_0$ onto this smaller linear space $D^{(k)}({\cal R}_0)$ so that $\tilde{h}_0\sigma^{-2}_0$ will be some projection of $h_0$ onto a smaller model $D^{(k)}({\cal R}_0)$. 
This is precisely the reason for $D_{P_0}$ not necessarily be equal to the efficient influence curve, but instead would equal a super-efficient influence curve, the efficient influence curve for a smaller model than ${\cal M}$ implied by assuming $Q_0\in D^{(k)}({\cal R}_0)$. 

The above theorem also provides the basis for asymptotic linearity of general differentiable functionals $\Phi(Q_n)$, shown as follows. Assume 
\[
\Phi(Q_n)-\Phi(Q_0)=\dot{\Phi}_{Q_0}(Q_n-Q_0)+R_{\Phi}(Q_n,Q_0),\]
where $\dot{\Phi}_{Q_0}:L^2(P_0)\rightarrow\openr$ is a bounded linear mapping, so that by the Riesz-representation theorem we have
\[
\dot{\Phi}_{Q_0}(h)=P_0 D_{\Phi,Q_0} h\] for a gradient $D_{\Phi,Q_0}\in L^2(P_0)$. 
So then we have $\Phi(Q_n)-\Phi(Q_0)=P_0 D_{\Phi,Q_0} (Q_n-Q_0)+o_P(n^{-1/2})$, and 
$P_0 D_{\Phi,Q_0}(Q_n-Q_0)$ is now analyzed by the previous theorem with $h_0=D_{\Phi,Q_0}$. This proves the following more general theorem.

\begin{theorem}
Consider a $\Phi$ so that 
\[
\Phi(Q_n)-\Phi(Q_0)=\dot{\Phi}_{Q_0}(Q_n-Q_0)+R_{\Phi}(Q_n,Q_0),\]
where $\dot{\Phi}_{Q_0}:L^2(P_0)\rightarrow\openr$ is a bounded linear mapping, so that by the Riesz-representation theorem we have
\[
\dot{\Phi}_{Q_0}(h)=P_0 D_{\Phi,Q_0} h\] for a gradient $D_{\Phi,Q_0}\in L^2(P_0)$. 
For notational convenience, let $h_0=D_{\Phi,Q_0}$. 
Assume $R_{\Phi}(Q_n,Q_0)=o_P(n^{-1/2})$ and the conditions of above Theorem \ref{theoremsmoothregr} with $h_0=D_{\Phi,Q_0}$ and corresponding definitions of $h_{0,n}$ and its limit $\tilde{h}_0)$. 

Then, $\Phi(Q_n)$ is an asymptotically linear estimator of $\Phi(Q_0)$ at $P_0$ with influence curve 
$D_{P_0}(O)=\tilde{h}_0\sigma^{-2}_0(X)(Y-Q_0(X))$.
\end{theorem}

{\bf Super efficiency:}
Analogue to linear $\Phi$-case, it follows that under the uniform nonparametric approximation condition on $D^{(k)}({\cal R}_{0,n})$, $D_{P_0}$ is the efficient influence curve for $\Phi(Q_0)$ for the nonparametric model. On the other hand, under the weaker adaptive uniform approximation condition on $D^{(k)}({\cal R}_{0,n})$, the influence curve $D_{P_0}$ will generally be a  super-efficient influence curve with a variance smaller than the efficient influence curve.

\section{Asymptotic normality of $k$-th order sieve and HAL-MLEs for general loss functions}\label{section10}


Recall our fixed working model $D^{(k)}({\cal R}_{0,n})$ and the oracle $Q_{0,n}=\arg\min_{Q\in D^{(k)}({\cal R}_{0,n})}P_0L(Q)$. We are also given one of our sieve or HAL-MLEs $Q_n\in D^{(k)}({\cal R}_n)$, where ${\cal R}_{0,n}$ is an independent random set approximating ${\cal R}_n$ as discussed earlier, so that $Q_n$ approximately solves the score equations for the fixed working model $D^{(k)}({\cal R}_{0,n})$.  
We assume  ${\cal R}_{0,n}$ satisfies the adaptive uniform approximation error condition  so that $\sup_{Q\in \{Q_0,Q_n\}}\inf_{\tilde{Q}\in D^{(k)}({\cal R}_{0,n})}\pl \tilde{Q}-Q\pl_{\infty}=O(C(M_n)r(d,J_{0,n})^{k+1})$ with $M_n=\pl Q\pl_{v,k}^*$. The latter is implied by the nonparametric uniform approximation error in which the supremum is over $Q\in D^{(k)}([0,1]^d)$.  Let $S_{Q}(\phi)= \left . \frac{d}{d\delta} L(f+\delta \phi)\right |_{\delta =0}$. 
We will have a special focus on the case that $L(Q)=-\log Q$ is a log-likelihood loss or at least has the behavior of the log-likelihood loss in the sense that
\begin{equation}\label{likelihoodloss}
-\frac{d}{d\delta_0} P S_{Q(P)+\delta_0 h}(\phi)=P S_{Q(P)}(\phi)S_{Q(P)}(h), \end{equation}
which is well known to be true for the log-likelihood loss. However, we present the results also for general loss functions that might not satisfy this property. 

\subsection{Proving pointwise asymptotic normality of sieve and HAL-MLEs}

{\bf Starting equation of analysis:}
Since $Q_{0,n}=\arg\min_{Q\in D^{(k)}({\cal R}_{0,n})}P_0L(Q)$ we have $P_0 S_{Q_{0,n}}(\phi_j) =0$ for all $j\in {\cal R}_{0,n}$. Let $r_n(\phi_j)\equiv P_n S_{Q_n}(\phi_j) $ be the score of $\beta(j)$ in model $D^{(k)}({\cal R}_{0,n})$ at $Q_n$, $j\in {\cal R}_{0,n}$. Our starting equation is given by:
\[
P_0 \{ S_{Q_n}(\phi_j)-S_{Q_{0,n}}(\phi_j)\}=-(P_n-P_0)S_{Q_n}(\phi_j) +r_n(\phi_j).\]
{\bf Tailor expansion of left-hand side to obtain second equality:}
We apply a first order Tailor expansion to the left-hand side.
Specifically, let $R_{1,n}(\phi_j)$ be defined by
 \begin{equation}\label{R1n}
P_0\{ S_{Q_n}(\phi_j)-S_{Q_{0,n}}(\phi_j)\}=P_0 dS_{Q_{0,n}}(\phi_j) (Q_n-Q_{0,n})+R_{1,n}(\phi_j),\end{equation}
where \[
dS_{Q_{0,n}}(\phi_j)(Q_n-Q_{0,n})=\left . \frac{d}{d\delta}S_{Q_{0,n}+\delta (Q_n-Q_{0,n})}(\phi_j)
\right |_{\delta =0}\]
is the  directional derivative of $Q\rightarrow S_{Q}(\phi_j)$ at $Q=Q_{0,n}$ in the direction $Q_n-Q_{0,n}$. 
Given this definition of the second order remainder $R_{1,n}(\phi_j)$, we have our next equality: 
\[P_0 dS_{Q_{0,n}}(\phi_j)(Q_n-Q_{0,n})=-(P_n-P_0) S_{Q_n}(\phi_j)+r_n(\phi_j)-R_{1,n}(\phi_j).\]

{\bf Approximating $Q_n$ by an element in $D^{(k)}({\cal R}_{0,n})$ so that the left-hand side derivative becomes a mapping on $D^{(k)}({\cal R}_{0,n})$:}
Let $\Pi_{J_{0,n}}(Q_n)$ be an element in $D^{(k)}({\cal R}_{0,n})$ so that $\pl \Pi_{J_{0,n}}(Q_n)-Q_n\pl_{\infty}=O(C(M_n) r(d,J_{0,n})^{k+1})$, where $M_n$ is the $L_1$-norm of the vector of non-zero coefficients in $Q_n$. By our assumption on ${\cal R}_{0,n}$ this element exists. We already proved that we can define $\Pi_{J_{0,n}}(Q_n)$ equal to projection of $Q_n$ onto $D^{(k)}({\cal R}_{0,n})$ in $L^2(P_0)$. 

We can write $Q_n=\Pi_{J_{0,n}}(Q_n)+(Q_n-\Pi_{J_{0,n}}(Q_n))$.
Let 
\begin{equation}\label{R2n}
R_{2,n}(\phi_j) \equiv P_0dS_{Q_{0,n}}(\phi_j)(Q_n-\Pi_{J_{0,n}}(Q_n)).\end{equation}
Then, we have
\[-P_0 dS_{Q_{0,n}}(\phi_j)(\Pi_{J_{0,n}}(Q_n)-Q_{0,n})=(P_n-P_0) S_{Q_n}(\phi_j)-r_n(\phi_j)+R_{1,n}(\phi_j) +R_{2,n}(\phi_j) .\]

{\bf Riesz representation theorem applied to left-hand linear real valued operator on $D^{(k)}({\cal R}_{0,n})$:}
For a given $j$, the left-hand side is a linear real valued mapping on the linear space $D^{(k)}({\cal R}_{0,n})$ evaluated at $\Pi_{J_{0,n}}(Q_n)-Q_{0,n}\in D^{(k)}({\cal R}_{0,n})$. 
Suppose that we put a particular inner product $\langle f,g\rangle_{J_{0,n}}$ on this finite $d_{0,n}$ dimensional linear space $D^{(k)}({\cal R}_{0,n})$, where our chosen definition will be presented and motivated later below. At this point we keep it general.
Given a definition of the inner product, by the Riesz-representation theorem, there exists a $D_{Q_{0,n}}(\phi_j)\in D^{(k)}({\cal R}_{0,n})$ so that the left-hand size equals $\langle D_{Q_{0,n}}(\phi_j),\Pi_{J_{0,n}}(Q_n)-Q_{0,n}\rangle_{J_{0,n}}$. 
So we now have
\[
\langle D_{Q_{0,n}}(\phi_j),\Pi_{J_{0,n}}(Q_n)-Q_{0,n}\rangle_{J_{0,n}}=
(P_n-P_0) S_{Q_n}(\phi_j)-r_n(\phi_j)+R_{1,n}(\phi_j) +R_{2,n}(\phi_j) .\]
Here we note that $D_{Q_{0,n}}(\phi_j)$ is a $j$-specific linear combination of $\{\phi_j:j\in {\cal R}_{0,n}\}$ and this linear combination is linear in $\phi_j$. 
We will allow to use an approximate inner product representation to naturally handle the log-likelihood loss function case covered in detail in the next subsection: $D_{Q_{0,n}}(\phi_j)\in D^{(k)}({\cal R}_{0,n})$ is defined so that
\begin{equation}\label{R3n}
-R_{3,n}(\phi_j)\equiv -P_0 dS_{Q_{0,n}}(\phi_j)(\Pi_{J_{0,n}}(Q_n)-Q_{0,n})-\langle D_{Q_{0,n}}(\phi_j), \Pi_{J_{0,n}}(Q_n)-Q_{0,n}\rangle_{J_{0,n}}\end{equation}
will end up contributing a second order (negligible) remainder.
Then,
\[
\langle D_{Q_{0,n}}(\phi_j),\Pi_{J_{0,n}}(Q_n)-Q_{0,n}\rangle_{J_{0,n}}=
(P_n-P_0) S_{Q_n}(\phi_j)-r_n(\phi_j)+\sum_{m=1}^3R_{m,n}(\phi_j) .\]
{\bf Particularly convenient and recommended inner product:}
We suggest that in general a good choice is
\begin{equation}\label{recommendinnerproduct}
\langle \phi,g\rangle_{J_{0,n}}\equiv -P_0 \frac{d}{dQ_{0,n}}S_{Q_{0,n}}(\phi) g,\end{equation}
or an approximation thereoff. In that case, $D_{Q_{0,n}}(\phi)=\phi$, and in the case of the log-likelihood loss creating an orthonormal basis $\{\phi_j^*:j\in {\cal R}_{0,n}\}$ from $\{\phi_j:j\in {\cal R}_{0,n}\}$ w.r.t. this inner product corresponds with creating an orthonormal basis of scores $\{ S_{Q_{0,n}}(\phi_{j}^*): j\in{\cal R}_{0,n}\}$ in $L^2(P_0)$, thereby creating a nice asymptotic variance expression. Even for non-likelihood loss functions, it simplifies matters to choose (\ref{recommendinnerproduct}).

{\bf Orthonormalizing the gradient basis functions $D_{Q_{0,n}}(\phi_j)$:}
Let $\tilde{\phi}_j=\sum_{k\in {\cal R}_{0,n}} \alpha_{j,0}(k)\phi_k$ be chosen so that $\{\phi_j^*\equiv D_{Q_{0,n}}(\tilde{\phi}_j):j\in {\cal R}_{0,n}\}$ is an orthonormal basis  of $D^{(k)}({\cal R}_{0,n})$ w.r.t. inner product $\langle g_1,g_2\rangle_{J_{0,n}}$, so that $\langle \phi_j^*,\phi_j^*\rangle_{J_{0,n}}=1$ and $\langle \phi_{j_1}^*,\phi_{j_2}^*\rangle_{J_{0,n}}=0$ for all $j_1\not =j_2$. 
If we select (\ref{recommendinnerproduct}) as inner product, then $\tilde{\phi}_j=\phi_j^*$.

Due to the linearity in $\phi_j$ of the last expression,  we have
\begin{equation}\label{tempidb}
\langle \phi_j^*,\Pi_{J_{0,n}}(Q_n)-Q_{0,n}\rangle_{J_{0,n}}=(P_n-P_0) S_{Q_n}(\tilde{\phi}_j)-r_n(\tilde{\phi}_j)+\sum_{k=1}^3R_{k,n}(\tilde{\phi}_j).\end{equation}
If we select (\ref{recommendinnerproduct}) as inner product, then this becomes
\begin{equation}\label{idbnice}
\langle \phi_j^*,\Pi_{J_{0,n}}(Q_n)-Q_{0,n}\rangle_{J_{0,n}}=(P_n-P_0) S_{Q_n}({\phi}_{j}^*)-r_n(\phi_j^*)+\sum_{k=1}^3R_{k,n}({\phi}_{j}^*).\end{equation}

{\bf Representing $\Pi_{J_{0,n}}(Q_n)-Q_{0,n}$ in terms of orthonormal basis:}
We now return to the general case, not just the log-likelihood loss case above or  the case obtained by selecting (\ref{recommendinnerproduct}).
Since $\{\phi_j^*:j\}$ is an orthonormal basis of $(D^{(k)}({\cal R}_{0,n}),\langle \cdot,\cdot \rangle_{J_{0,n}})$ and $\Pi_{J_{0,n}}(Q_n)-Q_{0,n}\in D^{(k)}({\cal R}_{0,n})$, we have 
\[
\Pi_{J_{0,n}}(Q_n)-Q_{0,n}=\sum_{j\in {\cal R}_{0,n}}\langle \Pi_{J_{0,n}}(Q_n)-Q_{0,n},\phi_j^*\rangle_{J_{0,n}}\phi_j^*.\]
We can now substitute  expression (\ref{tempidb}) for $\langle \Pi_{J_{0,n}}(Q_n)-Q_{0,n},\phi_j^*\rangle_{J_{0,n}}$ which yields:
\begin{eqnarray}
(\Pi_{J_{0,n}}(Q_n)-Q_{0,n})(x)&=&\sum_{j\in {\cal R}_{0,n}} (P_n-P_0) S_{Q_n}(\tilde{\phi}_j
)
\phi_j^*(x)\nonumber \\
&& -\sum_{j\in {\cal R}_{0,n}} r_n(\tilde{\phi}_j)\phi_j^*(x)\nonumber \\
&&+\sum_{m=1}^3\sum_{j\in {\cal R}_{0,n}}R_{m,n}(\tilde{\phi}_j)\phi_j^*(x)\nonumber
 \end{eqnarray}
 Let $\tilde{r}_n(x)\equiv \sum_{j\in {\cal R}_{0,n}} r_n(\tilde{\phi}_j)\phi_j^*(x)$, which is a score equation. 
  Let \[
  R_{m,n}(x)\equiv \sum_{j\in {\cal R}_{0,n}}R_{m,n}(\tilde{\phi}_j)\phi_j^*(x),\]
   $m=1,2,3$, and
\begin{equation}\label{En}
E_n(x)\equiv \sum_{j\in {\cal R}_{0,n}} (P_n-P_0) (S_{Q_n}-S_{Q_{0,n}}) (\tilde{\phi}_j)
\phi_j^*(x).\end{equation}

Regarding bounding $E_n(x)$ we can easily generalize Lemma \ref{lemmasuffEnx} to the following lemma. 
\begin{lemma}\label{lemmasuffEnxgen}

Assume $\pl Q_n-Q_0\pl_{P_0}=O_P(n^{-1/3}(\log n)^{d/2})$.
Assume that this implies
\[
\pl \frac{1}{d_{0,n}} \sum_{j\in {\cal R}_{0,n}}(S_{Q_n}-S_{Q_{0,n}})(\tilde{\phi}_j)\phi_j^*(x)\pl_{P_0}=O_P(n^{-1/3}(\log n)^{d/2}).\]

Then, $E_n(x)=O_P(d_{0,n} n^{-2/3})$ up till $\log n$-factors. Thus, $E_n(x)=o_P((d_{0,n}/n)^{1/2})$ if $d_{0,n}<n^{1/3-\delta}$ for some $\delta>0$. 
\end{lemma}
{\bf Proof:}
Note that $E_n(x)=d_{0,n} (P_n-P_0) g_n$, where 
\[
g_n=\frac{1}{d_{0,n}} \sum_{j\in {\cal R}_{0,n}}(S_{Q_n}-S_{Q_{0,n}})(\tilde{\phi}_j)\phi_j^*(x).\]
 Let ${\cal F}_n\equiv \{ d_{0,n}^{-1}\sum_{j\in {\cal R}_{0,n}}(S_{Q}-S_{Q_{0,n}})(\tilde{\phi}_j)\phi_j^*(x): Q\in  D^{(0)}_M([0,1]^d)\}$. 
 Assume $J(\delta,{\cal F}_n,L^2)\sim \delta^{1/2}$ up till $\log \delta$-factor. For example, ${\cal F}_n\subset D^{(0)}_{M_1}([a,b]^m)$ for some $M_1<\infty$ and $m$-dimensional cube $[a,b]^m$.  
We have $\pl g_n\pl_{P_0}=O_P(n^{-1/3}(\log n)^{d/2})$. Moreover, $g_n\in {\cal F}_n\subset D^{(0)}_M([0,1]^d)$ for some $M<\infty$, with probability tending to 1. Using the bound for $J(\delta,{\cal F}_n,L^2)\sim \delta^{1/2}$ up till $\log \delta$-factor with $\delta=\delta_n=n^{-1/3}(\log n)^{d/2}$ gives
 $E_n(x)\sim d_{0,n}n^{-1/2} n^{-1/6}$ up till $\log n$-factors. $\Box$

We also note that by the adaptive uniform approximation condition on $D^{(k)}({\cal R}_{0,n})$ \begin{eqnarray*}
(Q_n-Q_{0,n})(x)&=&(\Pi_{J_{0,n}}(Q_n)-Q_{0,n})(x)+(Q_n-\Pi_{J_{0,n}}(Q_n))(x)\\
&=&(\Pi_{J_{0,n}}(Q_n)-Q_{0,n})(x)+O_P(C(M_n)r(d,J_{0,n})^{k+1}).\end{eqnarray*}
{\bf Resulting expansion:}
So we can conclude the  following equation:
\begin{eqnarray*}
(Q_n-Q_{0,n})(x)&=&\sum_{j\in {\cal R}_{0,n}} (P_n-P_0) S_{Q_{0,n}}(\tilde{\phi}_j)
\phi_j^*(x)\\
&& -\tilde{r}_n(x)+\sum_{m=1}^3R_{m,n}(x) \\
&&-E_n(x)+O_P(C(M_n) r(d,J_{0,n})^{k+1}). \end{eqnarray*}
If we would choose (\ref{recommendinnerproduct}), then we can replace $\tilde{\phi}_j=\phi_j^*$.

{\bf Defining the remainder:}
 Let \begin{eqnarray}
 \bar{R}_n(x)=\sum_{m=1}^3\sum_{j\in {\cal R}_{0,n}} R_{m,n}(\tilde{\phi}_j)\phi_j^*(x) ,\label{barRn}
\end{eqnarray}
 and the total remainder
\begin{equation}\label{Rn}
R_n(x)\equiv- \tilde{r}_n(x)+\bar{R}_n(x)-E_n(x).\end{equation}
 
{\bf Asymptotically linear approximation with specified influence curve:}
Thus, we have shown  the following asymptotic linearity expansion:
\begin{eqnarray*}
(Q_n-Q_{0,n})(x)&=&\sum_{j\in {\cal R}_{0,n}} (P_n-P_0) S_{Q_{0,n}}\left(\tilde{\phi}_j\right)
\phi_j^*(x)\\
&& +R_n(x)+O_P(C(M_n) r(d,J_{0,n})^{k+1}).
\end{eqnarray*}
Our theorem assumes $R_n(x)=o_P((d_{0,n}/n)^{1/2})$. We also assume that $J_{0,n}$ is such that $C(M_n) r(d,J_{0,n})^{k+1}=o_P ( (d_{0,n}/n)^{1/2})$. 

Under these assumptions \[
(Q_n-Q_{0,n})(x)=P_n D_{Q_{0,n},x}+o_P((d_{0,n}/n)^{1/2}),\] where 
the influence curve is given by 
\begin{equation}\label{Df0n}
D_{Q_{0,n},x}\equiv \sum_{j\in {\cal R}_{0,n}}  S_{Q_{0,n}}\left(\tilde{\phi}_j\right)
\phi_j^*(x).\end{equation}

{\bf Expression for asymptotic variance:}
Let \begin{equation}\label{Sigma0}
\Sigma_{0,n}(j_1,j_2)\equiv P_0 S_{Q_{0,n}}\left(\tilde{\phi}_{j_1,0}\right)S_{Q_{0,n}}\left(\tilde{\phi}_{j_2,0}\right).\end{equation}
Then,
 $\sqrt{n/d_{0,n}}\sum_{j\in {\cal R}_{0,n}} (P_n-P_0) S_{Q_{0,n}}\left(\tilde{\phi}_j\right)
\phi_j^*(x)$ is a sum of independent mean zero random variables with variance given by 
\begin{eqnarray*}
\tilde{\sigma}^2_{0,n}(x)&=&\frac{1}{d_{0,n}}\sum_{j_1,j_2}P_0 S_{Q_{0,n}}\left(\tilde{\phi}_j\right)S_{Q_{0,n}}\left(\tilde{\phi}_j\right)\phi_{j_1}^*(x)\phi_{j_2}^*(x)\\
&=&\frac{1}{d_{0,n}}\sum_{j_1,j_2}\Sigma_{0,n}(j_1,j_2)\phi_{j_1}^*(x)\phi_{j_2}^*(x).
\end{eqnarray*}

{\bf Understanding asymptotic variance if  covariance matrix $\Sigma_{0,n}$ is non-diagonal:}
Let's now provide some additional insight about $\tilde{\sigma}^2_{0,n}(x)$.
We note that the asymptotic variance $\tilde{\sigma}^2_{0,n}(x)$ can be represented as
\[
\tilde{\sigma}^2_{0,n}(x)=d_n^{-1}\phi^{*,\top}(x)\Sigma_{0,n}\phi^*(x)=O(1).\]
 
 The following lemma bounds $\tilde{\sigma}^2_{0,n}(x)$ by the maximal eigenvalue of $\Sigma_{0,n}$ times $1/d_{0,n}\sum_{j\in {\cal R}_{0,n}}\phi^{*,2}_j(x)$.

 \begin{lemma}\label{lemmaboundedvar}
Consider an expression $\tilde{\sigma}^2_{0,n}(x)\equiv 1/d_n\sum_{j_1,j_2\in {\cal R}_{0,n}}\Sigma_{0,n}(j_1,j_2)\phi_{j_1}^*(x)\phi_{j_2}^*(x)$, where $\Sigma_{0,n}$ symmetric positive definite covariance matrix with maximal eigenvalue $\lambda_n$.
Note that $\tilde{\sigma}^2_{0,n}(x)=d_n^{-1}\phi^{*,\top}(x)\Sigma_0\phi^*(x)$, where $\phi^*(x)=(\phi_j^*(x):j\in {\cal R}_{0,n})$. 
Then, $\tilde{\sigma}^2_{0,n}(x)=O\left(\lambda_n d_n^{-1}\sum_{j\in {\cal R}_{0,n}}\phi^{*2}_j(x)\right)$. 
\end{lemma}
{\bf Proof:}
 Let $E_n^{\top}D_nE_n$ be the eigenvalue decomposition of $\Sigma_{0,n}$, where $D_n$ is the diagonal matrix with eigenvalues and $E_n$ is a the matrix with columns being the eigenvectors. Then, we
have 
\begin{eqnarray*}
\tilde{\sigma}^2_{0,n}(x)&=&d_n^{-1}\phi^{*,\top}(x)E_n^{\top} D_nE_n\phi^*(x)\\
&=& d_n^{-1} (E_n\phi^*)^{\top}D_n(E_n\phi^*)\\
&=& d_n^{-1}(D_n^{1/2}E_n\phi^*)^{\top}(D_n^{1/2}E_n\phi^*)\\
&=& d_n^{-1}\pl D_n^{1/2}E_n\phi^*\pl^2\\
&\leq &\max_j\mid \lambda_j\mid d_n^{-1}\pl E_n\phi^*\pl^2\\
&=&\lambda_n d_n^{-1} \pl \phi^*\pl^2\\
&=&\lambda_n d_n^{-1}\sum_{j\in {\cal R}_{0,n}}\{\phi^*_j(x)\}^2. \Box
\end{eqnarray*}

\subsection{Special focus on log-likelihood behaving loss function:}
{\bf For log-likelihood behaving loss functions (\ref{recommendinnerproduct}) corresponds with orthogonalizing scores:}
When selecting (\ref{recommendinnerproduct}) as inner product we obtain the nice expression (\ref{idbnice}). If the loss function is log-likelihood behaving, then the additional benefit is that this choice corresponds with asymptotically orthogonalizing the scores. Moreover, in this case, we could actually select the inner product equal to the covariance of the scores to obtain exact orthogonalization, at cost of introducing a remainder $R_{3n}(\phi_j)$. This remainder is well understood and creates a negligible contribution, as shown by next lemma. 
\begin{lemma}\label{lemmaR3nphi}
Suppose that $L(Q)=-\log p_Q$ is a log-likelihood loss function, with densities defined w.r.t. a dominating measure $\mu$,  for some link function $Q\rightarrow p_Q$, where $\{p_Q: Q\in D^{(k)}({\cal R}_{0,n})\}$ is a working model for ${\cal M}=\{p_Q: Q\in D^{(k)}([0,1]^d)\}$ or ${\cal M}=\{p_Q : Q\in D^{(k)}({\cal R})\}$ for some reduced support ${\cal R}$ in our $k$-th order spline representation with $D^{(k)}({\cal R})\subset D^{(k)}([0,1]^d)$. 
Define
\[
-R_{3n}(\phi_j)\equiv -P_0 dS_{Q_{0,n}}(\phi_j)(\Pi_{J_{0,n}}(Q_n)-Q_{0,n})-
P_0 S_{Q_{0,n}}(\phi_j)S_{Q_{0,n}}(\Pi_{J_{0,n}}(Q_n)-Q_{0,n}).
\]
{\bf Assumptions:}
$\pl p_{Q_{0,n}}-p_0\pl_{\mu}=O(r(d,J_{0,n})^{k+1})$ (e.g., implied by $Q_{0,n}-Q_0\pl_{\infty}=O(r(d,J_{0,n})^{k+1})$); 
$dS_{Q_{0,n}}(\phi_j)()$ and $S_{Q_{0,n}}()$ are uniformly bounded linear operators (also uniform in $n$) on $L^2(\mu)$. 

Then,
\[
R_{3,n}(\phi_j)= O(r(d,J_{0,n})^{k+1}) \pl \Pi_{J_{0,n}}(Q_n)-Q_{0,n}\pl_{\mu}.\]
More generally, for $g_n\in D^{(k)}({\cal R}_{0,n}\cup {\cal R}_n)$
\[
-P_0 dS_{Q_{0,n}}(\phi_j)(g_n)-
P_0 S_{Q_{0,n}}(\phi_j)S_{Q_{0,n}}(g_n)= O(r(d,J_{0,n})^{k+1}) \pl g_n \pl_{\mu}.\]


\end{lemma}
{\bf Proof: }
Suppose $g_n\in D^{(k)}({\cal R}_{0,n}\cup {\cal R}_n)$.
By Cauchy-Schwarz inequality, we have
\[
\begin{array}{l}
(P_0-P_{Q_{0,n}})dS_{Q_{0,n}}(\phi_j)(g_n)\\
\leq
\pl p_0-p_{Q_{0,n}}\pl_{\mu}\pl dS_{Q_{0,n}}(\phi_j)(g_n)\pl_{\mu}\\
=
O(r(d,J_{0,n})^{k+1}) \pl dS_{Q_{0,n}}(\phi_j)(g_n)\pl_{\mu}\\
=O_P(r(d,J_{0,n})^{k+1})\pl g_n\pl_{\mu},\end{array}
\]
by assumption that the operator $dS_{Q_{0,n}}(\phi_j)$ is a uniformly bounded linear operator on $L^2(\mu)$. 

Since for a correctly specified parametric model ${\cal M}_n\equiv \{p_Q: Q\in D^{(k)}({\cal R}_{0,n}\cup{\cal R}_n)\}$ with $Q_{0,n}\in D^{(k)}({\cal R}_{0,n}\cup{\cal R}_n)$ so that $P_{Q_{0,n}}\in {\cal M}_n$, the minus second derivative of the log likelihood at $p_{Q_{0,n}}$ equals the covariance of the scores, we have the following identity: for $g_1,g_2\in D^{(k)}({\cal R}_{0,n}\cup{\cal R}_n)$, 
\begin{equation}\label{keyidentitya}
-P_{Q_{0,n}}dS_{Q_{0,n}}(g_1)(g_2)=P_{Q_{0,n}}S_{Q_{0,n}}(g_1)S_{Q_{0,n}}(g_2).
\end{equation}
Thus,
\[
\begin{array}{l}
P_0 dS_{Q_{0,n}}(\phi_j)(g_n)\\
=
(P_0-P_{Q_{0,n}})dS_{Q_{0,n}}(\phi_j)(g_n)+
P_{Q_{0,n}}dS_{Q_{0,n}}(\phi_j)(g_n)\\
= O(r(d,J_{0,n})^{k+1}) \pl dS_{Q_{0,n}}(\phi_j)(g_n)\pl_{\mu}\\
\hfill -P_{Q_{0,n}} S_{Q_{0,n}}(\phi_j)S_{Q_{0,n}}(g_n)\\
=O(r(d,J_{0,n})^{k+1}) \pl dS_{Q_{0,n}}(\phi_j)(g_n)\pl_{\mu}\\
\hfill -(P_{Q_{0,n}}-P_0) S_{Q_{0,n}}(\phi_j)S_{Q_{0,n}}(g_n)\\
\hfill -P_0 S_{Q_{0,n}}(\phi_j)S_{Q_{0,n}}(g_n)\\
= -P_0 S_{Q_{0,n}}(\phi_j)S_{Q_{0,n}}(g_n)\\
\hfill +
O(r(d,J_{0,n})^{k+1}) \pl dS_{Q_{0,n}}(\phi_j)(g_n)\pl_{\mu}\\
+O(r(d,J_{0,n})^{k+1})\pl S_{Q_{0,n}}(g_n)\pl_{\mu}\\
=-P_0 S_{Q_{0,n}}(\phi_j)S_{Q_{0,n}}(g_n)+O(r(d,J_{0,n})^{k+1}\pl g_n\pl_{\mu}).
\end{array}
\]
This completes the proof of the lemma. $\Box$

{\bf Analogue of (\ref{tempidb}) for log-likelihood behaving loss functions:}
So for such log-likelihood behaving loss functions we have
\begin{equation}\label{tempida}
P_0 S_{Q_{0,n}}(\phi_j)S_{Q_{0,n}}(\Pi_{J_{0,n}}(Q_n)-Q_{0,n})=(P_n-P_0) S_{Q_n}(\phi_j)-r_n(\phi_j)+\sum_{m=1}^3 R_{m,n}(\phi_j) ,\end{equation}
where the left-hand side equals an inner product $\langle \phi_j,\Pi_{J_{0,n}}(Q_n)-Q_{0,n}\rangle_{J_{0,n}}$ on $D^{(k)}({\cal R}_{0,n})$ given by 
\begin{equation}\label{scoreinner}
\langle g_1,g_2\rangle_{J_{0,n}}\equiv P_0 S_{Q_{0,n}}(g_1)S_{Q_{0,n}}(g_2).\end{equation}
In this case we relied on an approximate inner product representation: 
\[-P_0 dS_{Q_{0,n}}(\phi_j)(\Pi_{J_{0,n}}(Q_n)-Q_{0,n})=
\langle \phi_j,\Pi_{J_{0,n}}(Q_n)-Q_{0,n}\rangle_{J_{0,n}}-R_{3,n}(\phi_j).\]
So we obtain the expression (\ref{idbnice})
\begin{equation}\label{tempidbloglik}
\langle \phi_j^*,\Pi_{J_{0,n}}(Q_n)-Q_{0,n}\rangle_{J_{0,n}}=(P_n-P_0) S_{Q_n}({\phi}_{j}^*)-r_n({\phi}_{j}^*)+\sum_{m=1}^3R_{m,n}({\phi}_{j}^*),\end{equation}
but with the additional benefit that we now know that $P_0 S_{Q_{0,n}}(\phi_{j_1}^*)S_{Q_{0,n}}(\phi_{j_2}^*)=0$ for $j_1,j_2\in {\cal R}_{0,n}$ with $j_1\not =j_2$.

{\bf Understanding $R_{m,n}(x)$ for $m=2,3$ and the log-likelihood loss:}
The following lemma shows that for log-likelihood loss we have that $R_{2,n}(x)$ and $R_{3,n}(x)$ are generally $o_P((n/d_{0,n})^{-1/2})$ under a weak regularity condition.
Recall $R_{2,n}(x)=\sum_{j\in {\cal R}_{0,n}}P_0 dS_{Q_{0,n}}(\phi_j^*)(Q_n-\tilde{Q}_n) \phi_j^*(x)$, 
where $\tilde{Q}_n=\Pi_{J_{0,n}}(Q_n)$, and recall $R_{2,n}(\phi_j^*)=P_0 dS_{Q_{0,n}}(\phi_j^*)(Q_n-\tilde{Q}_n)$. 
By Lemma \ref{lemmaR3nphi} we have 
$R_{2,n}(\phi_j^*)=-P_0 S_{Q_{0,n}}({\phi}_j^*)S_{Q_{0,n}}(Q_n-\tilde{Q}_n)+
O_P( r(d,J_{0,n})^{2(k+1)})$, using that $\pl Q_n-\tilde{Q}_n\pl_{\mu}=O(C(M_n)r(d,J_{0,n})^{k+1})$.
Therefore, 
\begin{eqnarray*}
-R_{2,n}(x)&=&\sum_{j\in {\cal R}_{0,n}}P_0 S_{Q_{0,n}}(\tilde{\phi}_j) S_{Q_{0,n}}(Q_n-\tilde{Q}_n) \phi_j^*(x)\\
&&+O_P\left( \sum_{j\in {\cal R}_{0,n}}\phi_j^*(x) r(d,J_{0,n})^{2(k+1)}\right).\end{eqnarray*}
The leading term equals $\sum_{j\in {\cal R}_{0,n}}\langle Q_n-\tilde{Q}_n, \phi_j^*\rangle_{J_{0,n}} \phi_j^*$ which is the projection of $Q_n-\Pi_{J_{0,n}}(Q_n)$ onto $D^{(k)}({\cal R}_{0,n})$ and thus equals zero. 
So, we conclude
\[
R_{2,n}(x)=O_P\left( r(d,J_{0,n})^{2(k+1)}\right)\sum_{j\in {\cal R}_{0,n}}\phi_j^*(x) 
.\]
By Lemma \ref{lemmaR3nphi} we also have
\[
R_{3,n}(x)=O_P\left( r(d,J_{0,n})^{2(k+1)}\right) \sum_{j\in {\cal R}_{0,n}}\phi_j^*(x). \]
This proves the following lemma for the log-likelihood loss.


\begin{lemma} \label{lemmaR2nR3nloglik}
Suppose that $L(Q)=-\log p_Q$ is a log-likelihood loss function, with densities defined w.r.t. a dominating measure $\mu$,  for some link function $Q\rightarrow p_Q$, where $\{p_Q: Q\in D^{(k)}({\cal R}_{0,n})\}$ is a working model for ${\cal M}=\{p_Q: Q\in D^{(k)}([0,1]^d)\}$ or ${\cal M}=\{p_Q : Q\in D^{(k)}({\cal R})\}$ for some reduced support ${\cal R}$ in our $k$-th order spline representation with $D^{(k)}({\cal R})\subset D^{(k)}([0,1]^d)$. 
Define
\[
-R_{3n}(\phi_j)\equiv -P_0 dS_{Q_{0,n}}(\phi_j)(\Pi_{J_{0,n}}(Q_n)-Q_{0,n})-
P_0 S_{Q_{0,n}}(\phi_j)S_{Q_{0,n}}(\Pi_{J_{0,n}}(Q_n)-Q_{0,n}).
\]
{\bf Assumptions:}
$\pl p_{Q_{0,n}}-p_0\pl_{\mu}=O(r(d,J_{0,n})^{k+1})$ (e.g., implied by $Q_{0,n}-Q_0\pl_{\infty}=O(r(d,J_{0,n})^{k+1})$); 
$dS_{Q_{0,n}}(\phi_j)()$ and $S_{Q_{0,n}}()$ are uniformly bounded linear operators (also uniform in $n$) on $L^2(\mu)$.


Let $\Pi_{J_{0,n}}(Q_n)$ be the projection of $Q_n$ onto $D^{(k)}({\cal R}_{0,n})$ in $L^2(P_0)$. 

Then,
\[
R_{2,n}(x)=O_P\left( r(d,J_{0,n})^{2(k+1)}\right)\sum_{j\in {\cal R}_{0,n}}\phi_j^*(x) 
,\]
and 
\[
R_{3,n}(x)=O_P\left( r(d,J_{0,n})^{2(k+1)}\right) \sum_{j\in {\cal R}_{0,n}}\phi_j^*(x). \]
\end{lemma}
If $J_{0,n}$ is chosen so that $r(d,J_{0,n})^{k+1}=o_P((d_{0,n}/n)^{1/2})$, then it follows that
$R_{3,n}(x)=O_P( d_{0,n}^2n^{-1})$, which is $o_P((d_{0,n}/n)^{1/2})$ if $d_{0,n}=O(n^{1/3-\delta})$ for some $\delta>0$.

{\bf Particularly nice expression of variance for log-likelihood behaving loss function:}
In the case that $D_{Q_{0,n}}(\phi_j)=\phi_j$ and $\langle g_1,g_2\rangle_{J_{0,n}}=P_0 S_{Q_{0,n}}(g_1)S_{g_{0,n}}(g_2)$ we have that the influence curve $D_{Q_{0,n},x}$ is a sum over $j$ of uncorrelated random variables. 
So in that case we have
\[
\tilde{\sigma}^2_{0,n}(x)=\frac{1}{d_{0,n}}\sum_{j\in {\cal R}_{0,n}} \{\phi_j^*(x)\}^2.\]

\subsection{Asymptotic normality theorem}
We state the theorem for the most important case in which we select (\ref{recommendinnerproduct}) as inner product $\langle \cdot,\cdot\rangle_{J_{0,n}}$, or an approximation thereof as with the log-likelihood loss so that $\phi_j^*=\tilde{\phi}_j$. 
So this still covers general loss functions. 

Recall that ${\cal R}_{0,n}={\cal R}(P_n^{\#})$ is an independent copy of ${\cal R}_n={\cal R}(P_n)$ with $P_n^{\#}$ empirical measure of independent (from $P_n$) i.i.d. sample from $P_0$.
In addition, recall that $\{\phi_j^*:j\in {\cal R}_{0,n}\}$ is  an orthonormal basis of $\{\phi_j: j\in {\cal R}_{0,n}\}$ w.r.t. the inner product $\langle \cdot,\cdot\rangle_{J_{0,n}}$ (\ref{recommendinnerproduct}).

Consider/recall the following definitions:
\begin{eqnarray}
\tilde{r}_n(x)&\equiv& \sum_{j\in {\cal R}_{0,n}} \{P_n S_{Q_n}(\phi_j^*)\}\phi_j^*(x)\label{tildernxgf}
\\
R_{1,n}(\phi)&\equiv& P_0\{ S_{Q_n}(\phi)-S_{Q_{0,n}}(\phi)\}-P_0 dS_{Q_{0,n}}(\phi) (Q_n-Q_{0,n})\label{R1ngf}\\
R_{2,n}(\phi)& \equiv& P_0dS_{Q_{0,n}}(\phi)(Q_n-\Pi_{J_{0,n}}(Q_n))\label{R2ngf}\\
-R_{3,n}(\phi)&\equiv& -P_0 dS_{Q_{0,n}}(\phi)(\Pi_{J_{0,n}}(Q_n)-Q_{0,n})-\langle \phi, \Pi_{J_{0,n}}(Q_n)-Q_{0,n}\rangle_{J_{0,n}}\label{R3ngf}\\
E_n(x)&\equiv& \sum_{j\in {\cal R}_{0,n}} (P_n-P_0) (S_{Q_n}-S_{Q_{0,n}}) ({\phi}_j^*)
\phi_j^*(x)\label{Engf}\\
 \bar{R}_n(x)&\equiv&\sum_{m=1}^3\sum_{j\in {\cal R}_{0,n}} R_{m,n}(\phi_j^*)\phi_j^*(x) \label{barRngf}
\\
R_n(x)&\equiv&- \tilde{r}_n(x)+\bar{R}_n(x)-E_n(x)\label{Rngf}\\
D_{Q_{0,n},x}&\equiv& \sum_{j\in {\cal R}_{0,n}}  S_{Q_{0,n}}\left({\phi}_j^*\right)
\phi_j^*(x)\label{Df0ngf}\\
\Sigma_0(j_1,j_2)&\equiv& P_0 S_{Q_{0,n}}\left({\phi}_{j_1}^*\right)S_{Q_{0,n}}\left({\phi}_{j_2}^*\right).\label{Sigma0gf}
\end{eqnarray}
\newline
{\bf Assumption G0:}
\begin{equation}\label{G0}
\sup_{Q\in \{Q_n,Q_0\}}\inf_{g\in D^{(k)}({\cal R}_{n})}\pl Q-g\pl_{\infty}=O_P(C(M_n)r(d,J_{0,n})^{k+1}),\end{equation}
 where $M_n=\pl Q_n\pl_{v,k}^*$, the $L_1$-norm of its non-zero coefficients. \newline
{\bf Assumption G1:}
\begin{equation}\label{G1}
\tilde{r}_n(x)=o_P((n/d_{0,n})^{-1/2}).
\end{equation}
{\bf Assumption G2:}
\begin{equation}
\label{G2}
\bar{R}_n(x)=o_P((n/d_{0,n})^{-1/2}).
\end{equation}
{\bf Assumption G3:}\begin{equation}\label{G3}
E_n(x)=o_P((n/d_{0,n})^{-1/2}).\end{equation}
{\bf Assumption G4:}
\begin{equation}\label{G4}
C(M_n)r(d,J_{0,n})^{k+1}=o((n/d_{0,n})^{-1/2}).
\end{equation}
We have proven the following theorem. 


\begin{theorem}\label{theoremasnormgen}
Consider the $k$-th order sieve MLE, HAL-MLE or relax HAL-MLE $Q_n$.  Let $k^*=k+1$.
Let ${\cal R}_n$ be the set of $d_n$ non-zero coefficients in $Q_n=\sum_{j\in {\cal R}_n}\beta_n(j)\phi_j$.
Consider the definitions above and assume  G0-G4.


We have
\begin{eqnarray*}
(Q_n-Q_{0,n})(x)&=&\sum_{j\in {\cal R}_{0,n}} (P_n-P_0) D_{Q_{0,n},x}\\
&& +R_n(x)+O_P(C(M_n) r(d,J_{0,n})^{k+1}).
\end{eqnarray*}

Then, 
\begin{eqnarray*}
\frac{(Q_n-Q_{0,n})(x)}{(d_{0,n}/n)^{1/2}}&=& P_n d_{0,n}^{-1/2}D_{Q_{0,n},x}+o_P(1).
\end{eqnarray*}
For a given $x$, the leading term is now  a sum of independent mean zero random variables so that the central limit theorem can be applied. 
The variance is given by:
\[
\tilde{\sigma}^2_{0,n}(x)=\frac{1}{d_{0,n}}\sum_{j_1,j_2\in {\cal R}_{0,n}}\Sigma_{0,n}(j_1,j_2)\phi_{j_1}^*(x)\phi_{j_2}^*(x).\]
In the case of log-likelihood behaving loss in the sense that the inner product representation (\ref{Df0n}) holds  with$\langle g_1,g_2\rangle_{J_{0,n}}=P_0 S_{Q_{0,n}}(g_1)S_{Q_{0,n}}(g_2)$, we have
\begin{equation}\label{sigma2nd}
\tilde{\sigma}^2_{0,n}(x)=\frac{1}{d_{0,n}}\sum_{j\in {\cal R}_{0,n}}\{\phi_j^*(x)\}^2.\end{equation}

{\bf Conclusion:}
Under these assumptions, we have
\[
\tilde{\sigma}^{-1}_{0,n} (n/d_{0,n})^{1/2}({Q_n}-{Q_{0,n}})(x)\Rightarrow_d N(0,1),\]
and
\[
\tilde{\sigma}^{-1}_{0,n}(n/d_{0,n})^{1/2}(Q_n-Q_0)(x)\Rightarrow_d N(0,1).\]
Assuming $M_n=O(\log^{m}n)$ for some $m<\infty$, by choosing $J_{0,n}=n^{1/(2k^*+1)}\log^{m_1} n$ for some  $m_1<\infty$,  it follows that for  some  $m<\infty$, \[
\mid Q_n-Q_0\mid(x)=O_P(n^{-k^*/(2k^*+1)}\log^m n).\]
\end{theorem}
\paragraph{Discussion of assumptions:}
Assumption G4 is controlled by user and therefore holds by choosing optimal rate up till $\log n$-factor.
Our lemmas also establish sufficient conditions for G1-G3.\newline
{\bf Assumption G3:}
Lemma \ref{lemmasuffEnxgen} shows that under weak regularity condition we have $E_n(x)=o_P((d_{0,n}/n)^{1/2})$ if $d_{0,n}<n^{1/3-\delta}$ for some $\delta>0$. \newline
{\bf Assumption G2:}
We also note that $R_{1,n}(x)=\sum_{j\in {\cal R}_{0,n}}R_{1,n}(\phi_j^*)\phi_j^*(x)$ can be represented as $\Pi_{J_{0,n}}(H(Q_n,Q_{0,n})\mid D^{(k)}({\cal R}_{0,n}))$ for a function $H(Q_n,Q_{0,n})$ involving a square difference  of $Q_n-Q_{0,n}$ and will therefore typically be  $o_P(\mid Q_n-Q_{0,n}\mid (x))$.
In addition,  $R_{2,n}(x)$ will be $O_P(r(d,J_{0,n})^{k+1})$ under weak regularity conditions, in particular, we have $\pl R_{2,n}\pl_{J_{0,n}}^2=O_P( r(d,J_{0,n})^{2(k+1)})$. Finally,  $R_{3,n}(x)=\sum_{j\in {\cal R}_{0,n}}R_{3,n}({\phi}_j^*) \phi_j^*(x)$ where $R_{3,n}(\phi_j^*)$ will be zero or a second order remainder. 
In particular, if $L(Q)=-\log p_Q$ is the log-likelihood loss, then one can apply Lemma \ref{lemmaR2nR3nloglik} to obtain, under weak regularity conditions,
\[
R_{m,n}(x)=O_P\left( \sum_{j\in {\cal R}_{0,n}}\phi_j^*(x) r(d,J_{0,n})^{2(k+1)}\right),\mbox{ m=2,3}.\]
{\bf Assumption G1:}
Lemma  \ref{lemmatildernx} shows  \[
\tilde{r}_n(x)=\sum_{j\in {\cal R}_{0,n}}(P_n-\tilde{P}_n)\{S_j^*(Q_n)-\tilde{S}_j(Q_n)\} \phi_j^*(x).\]
Let  ${\cal F}^k_n\equiv \{1/d_{0,n}\sum_{j\in {\cal R}_{0,n}}\{S_j(Q)-\tilde{S}_j(Q)\}\phi_j^*(x): Q\in D^{(k)}_M([0,1]^d)\}$. Suppose the covering number of ${\cal F}^k_n$  is of same order as the covering number of $D^{(k)}_M([0,1]^d)$ given by  Lemma \ref{lemmasupnormcoveringnumber}.
 Let $\delta_n\equiv  \pl d_{0,n}^{-1}\sum_{j\in {\cal R}_{0,n}}\{S_j(Q_n)-\tilde{S}_j(Q_n)\}\pl_{\mu}$ with $\mu$ Lebesgue measure.
Then, by Lemma \ref{lemmatildernx}, if  $d_{0,n}^{1/2}\delta_n^{(2k+1)/(2k+2)}=O_P(n^{-\delta})$ for some $\delta>0$, then (\ref{G1}) holds.



\subsection{Uniform convergence rate for sieve and HAL-MLEs}
Let G1*, G2* and G3* be the sup-norm analogues of G1,G2 and G3:\newline
{\bf Assumption G1*:}
\begin{equation}\label{G1*}
\pl \tilde{r}_n\pl_{\infty}=o_P((n/d_{0,n})^{-1/2}).
\end{equation}
{\bf Assumption G2*:}
\begin{equation}
\label{G2*}
\pl \bar{R}_n\pl_{\infty}=o_P((n/d_{0,n})^{-1/2}).
\end{equation}
{\bf Assumption G3*:}\begin{equation}\label{G3*}
\pl E_n\pl_{\infty}=o_P((n/d_{0,n})^{-1/2}).\end{equation}
In additon, consider
\newline
{\bf Assumption G4*:}
\begin{equation}\label{G4*}
\log n(n/d_{0,n})^{1/2}C(M_n) r(d,J_{0,n})^{k+1}=o_P(1).\end{equation}

\begin{theorem}\label{thunifregr}\ \nl
Consider the $k$-th order sieve MLE, HAL-MLE or relax HAL-MLE $Q_n$.  Let $k^*=k+1$.
Consider the setting described in Theorem \ref{theoremasnormgen} and Assume G0,G1*,G2*,G3*,G4*.
By Theorem \ref{theoremasnormgen} we have the following expansion:
\begin{eqnarray}
(n/d_{0,n})^{1/2}(Q_n-Q_{0,n})(x)&=& n^{1/2}(P_n-P_0) d_{0,n}^{-1/2}\sum_{j\in {\cal R}_{0,n}}S_{Q_{0,n}}({\phi}^*_j) \phi_j^*(x) \nonumber\\
&&+(n/d_{0,n})^{1/2}E_n(x) -
(n/d_{0,n})^{1/2}\tilde{r}_n(x)
\nonumber \\
&&\hspace*{-3cm}-(n/d_{0,n})^{1/2}\bar{R}_{n}(x) +(n/d_{0,n})^{1/2}O_P(C(M_n)r(d,J_{0,n})^{k+1}),\label{keyexpansion}\end{eqnarray}
where $M_n=\pl Q_n\pl_{v,k}^*$, which equals the $L_1$-norm  of the vector of non-zero coefficients in $Q_n$.

{\bf Regularity condition:}
Let \[
{\cal F}_{n}=\left\{ \frac{1}{d_{0,n}}\sum_{j\in {\cal R}_{0,n}}S_{Q_{0,n}}({\phi}_j^*))\phi_j^*(x)/\tilde{\sigma}_{0,n}(x): x\right\}.\]
That is, these functions are indexed by $x\in [0,1]^d$.
Assume $\sup_{Q\in {\cal F}_n}\sup_o\mid Q(o)\mid =O(1)$ and that the covering number $N(\epsilon,{\cal F}_n,L^2)=O(\epsilon^{-d})$.

Then, 
\begin{eqnarray}
\frac{(Q_n-Q_{0,n})(x)}{(n/d_{0,n})^{-1/2}}= n^{1/2}(P_n-P_0) d_{0,n}^{-1/2}\sum_{j\in {\cal R}_{0,n}}
S_{Q_{0,n}}({\phi}_j^*) \phi_j^*(x)+o_P(1). \label{keyexpansiongen} 
\end{eqnarray}
For a given $x$, the leading term is now  a sum of independent mean zero random variables $D_{Q_{0,n},x}(O)$ (\ref{Df0nx}) so that the central limit theorem can be applied. 
The variance is given by:
\begin{equation}\label{sigma2ne}
\tilde{\sigma}^2_{0,n}(x)=
\frac{1}{d_{0,n}}\sum_{j_1,j_2\in {\cal R}_{0,n}}\Sigma_{0,n}(j_1,j_2)\phi_{j_1}^*(x)\phi_{j_2}^*(x),\end{equation}
which simplifies under the log-likelihood loss so that $\Sigma_{0,n}$ is the identity matrix. 
Thus, $(n/d_{0,n})^{1/2}(Q_n-Q_{0,n})(x)/\tilde{\sigma}_{0,n}(x)\Rightarrow_d N(0,1)$ for all $x$.

We have 
\[
\sup_{x\in [0,1]^d}(n/d_{0,n})^{1/2}\mid {Q_n}-{Q_{0,n}}\mid (x)/\tilde{\sigma}_{0,n}(x)= o_P(\log n),
\]
and  
\[
\sup_{x\in [0,1]^d}(n/d_{0,n})^{1/2}\mid Q_n-Q_0\mid(x)/\tilde{\sigma}_{0,n}(x)=o_P(\log n).\]
Note that, by choosing $J_{0,n}=n^{1/(2k^*+1)}\log^{m_1}n$ for some specified $m_1$,  this implies that for some specified $m<\infty$, 
\[
\pl  Q_n-Q_0\pl_{\infty}=O_P(n^{-k^*/(2k^*+1)}\log^m n) .\]
\end{theorem}
{\bf Proof:}
We have that $(n/d_{0,n})^{1/2}(Q_n-Q_{0,n})(x)/\tilde{\sigma}_{0,n}(x)=n^{1/2}(P_n-P_0) D^*_{Q_{0,n},x}+R_n(x)$, where $\sup_x \mid R_n(x)\mid=o_P(1)$ and 
$D^*_{Q_{0,n},x}=d_{0,n}^{-1/2}\sum_{j\in {\cal R}_{0,n}}S_{Q_{0,n}}({\phi}_j^*)\phi_j^*(x) /\tilde{\sigma}_{0,n}(x)$.
Note $D^*_{Q_{0,n},x}$ has mean zero and variance $1$.
$d_{0,n}^{-1/2}D^*_{Q_{0,n},x}$ is included in the fixed $d$-dimensional class
\[
{\cal F}_{n}=\left\{ \frac{1}{d_{0,n}}\sum_{j\in {\cal R}_{0,n}}S_{Q_{0,n}}({\phi}_j^*))\phi_j^*(x)/\tilde{\sigma}_{0,n}(x): x\right\}.\]
That is, these functions are indexed by $x\in [0,1]^d$.
Let $g_{x}$ be one element in ${\cal F}_n$. Then
\[
\sup_x P_0 g_{x}^2=\frac{1}{d_{0,n}},\]
due to our extra dividing by $d_{0,n}^{-1/2}$.

Thus, $\sup_{Q\in {\cal F}_n}P_0 Q^2\sim d_{0,n}^{-1}$.
Moreover, by assumption, ${\cal F}_n$ has  universally  bounded sup-norm and its covering number $N(\epsilon,{\cal F}_n,L^2)=O(\epsilon^{-d})$.
Therefore, the entropy integral $J(\delta,{\cal F}_n)\equiv \int_0^{\delta}\sup_Q\sqrt{\log N(\epsilon,{\cal F}_n,L^2(Q))} d\epsilon$ of $ {\cal F}_n$ is bounded by $\sim -\int_0^{\delta} (\log \epsilon)^{1/2} d\epsilon$.  This can be {\em conservatively} bounded by \[
\sim -\int_0^{\delta}\log \epsilon d\epsilon = -\delta \log \delta +\delta\sim -\delta \log \delta.\]
So we can state that $J(\delta,{\cal F}_n)=o(\delta \log \delta)$.
We can can set $\delta=\delta_n=d_{0,n}^{-1/2}$.

Applying empirical process inequalities for $E \sup_{Q\in {\cal F}_n,\pl Q\pl_{P_0}\leq \delta_n}\mid n^{1/2}(P_n-P_0)Q\mid \sim J(\delta_n,{\cal F}_n,L^2)$ with $\delta_n=d_{0,n}^{-1/2}$ yields now:
\[
\begin{array}{l}
\sup_x (n/d_{0,n})^{1/2}\mid {Q_n}-{Q_{0,n}}\mid (x)/\tilde{\sigma}_{0,n}(x)\\
\approx 
d_{0,n}^{1/2}\sup_x\left | n^{1/2}(P_n-P_0) \frac{1}{d_{0,n}}\sum_{j\in {\cal R}_{0,n}}S_{Q_{0,n}}(\tilde{\phi}_j) \phi_j^*(x)/\tilde{\sigma}_{0,n}(x) \right | \\
 \leq
d_{0,n}^{1/2} \sup_{g\in {\cal F}_n,\pl g\pl_{P_0}\leq \delta_n} \mid n^{1/2}(P_n-P_0)g\mid \\
\sim  d_{0,n}^{1/2}O_P( J(\delta_n,{\cal F}_n,L^2))\\
\sim  d_{0,n}^{1/2} o_P( \delta_n \log \delta_n)\\
\sim o_P(\log n) . \Box
\end{array}
\]

In Theorem \ref{theoremasnormgen} we stated sufficient conditions for G1,G2,G3 that are also sufficient for the uniform versions G1*,G2*,G3*.
In Appendix \ref{AppendixK1} we prove asymptotic normality of arbitrary functions of $Q_n$ building on our asymptotic linearity expansion for $Q_n-Q_{0,n}$. 
\section{Asymptotic linearity (and super-efficiency) of plug-in HAL-MLE, relax-HAL-MLE, and sieve-MLE of  pathwise differentiable features of $Q_0$}\label{section11}

\subsection{Delta-method approach, handling general loss functions.}
Let $\Phi(Q_0)$ be a pathwise differentiable function of $P_0$ on the model ${\cal M}$ for $P_0$. 
When analyzing  a pathwise differentiable parameter $\Phi(Q_n)-\Phi(Q_{0,n})$, using a $\delta$-method approximation $d\Phi_{Q_{0,n}}(Q_n-Q_{0,n})$,  it will not be problematic if the influence curve of $(Q_n-Q_{0,n})$ involves a sum over an informative set ${\cal R}_n$ (which we had to avoid for analyzing $(Q_n-Q_{0,n})(x)$ itself), as long as we can argue that this influence curve lives in a Donsker class and converges to a fixed function in $L^2(P_0)$, as one expects to be true. 
As a consequence, we can carry out the same analysis as in previous sections to obtain an expansion for $(Q_n-Q_{0,n})(x)$,  but in this case just setting ${\cal R}_{0,n}={\cal R}_n$.
Imitating our analysis, our expansion for $(Q_n-Q_{0,n})(x)$ when setting ${\cal R}_{0,n}={\cal R}_n$ does not have to deal with the term $\Pi_{J_{0,n}}(Q_n)-Q_n$ anymore, since it equals zero now, while this term  contributed a remainder term $O_P(C(M_n)r(d,J_{0,n})^{k+1})$ and a similar $R_{2,n}(x)$ , which represents a problematic bias term in an analysis of $\Phi(Q_n)-\Phi(Q_{0,n})$ in which bias has to converge to zero faster than $n^{-1/2}$.
Thus, we obtain the same expansion but without $R_{2,n}(x)$ and without $O_P(C(M_n)r(d,J_{0,n})^{k+1})$. The next lemma makes this point. 
We present the result for the case that we select for inner product (\ref{recommendinnerproduct}) or, in the case of log-likelihood behaving loss functions,  its covariance of score (\ref{scoreinner}) approximation.

\begin{lemma}
Suppose that we set ${\cal R}_{0,n}={\cal R}_n$, which is now allowed to be informative. 
Define
\begin{eqnarray*}
R_{1n}({\phi}_j^*)&\equiv& P_0\{ S_{Q_n}({\phi}_j^*)-S_{Q_{0,n}}({\phi}_j^*)\}-P_0 dS_{Q_{0,n}}(\phi_j^*) (Q_n-Q_{0,n})\\
-R_{3,n}(\phi_j^*)&\equiv& -P_0 dS_{Q_{0,n}}(\phi_j^*)(Q_n-Q_{0,n})-\langle \phi_j^*, Q_n-Q_{0,n}\rangle_{J_{0,n}}.
\end{eqnarray*}
Define $R_{m,n}(x)=\sum_{j\in {\cal R}_{0,n}}R_{k,n}(\phi_j^*)\phi_j^*(x)$, $m=1,3$.
Let $\bar{R}_n(x)=R_{1,n}(x)+R_{3,n}(x)$.
Recall definitions of $E_n(x)=\sum_{j\in {\cal R}_{0,n}}E_n(\phi_j^*)\phi_j^*(x)$ and $\tilde{r}_n(x)\equiv \sum_{j\in {\cal R}_{0,n}}P_n S_{Q_n}(\phi_j^*)\phi_j^*(x)$, where $E_n(\phi_j^*)=(P_n-P_0)(S_{Q_n}(\phi_j^*)-S_{Q_{0,n}}(\phi_j^*))$. 
Let $R_n(x)=\bar{R}_n(x)-E_n(x)-\tilde{r}_n(x)$. 
We have
\begin{eqnarray*}
(Q_n-Q_{0,n})(x)&=&\sum_{j\in {\cal R}_{0,n}} (P_n-P_0) S_{Q_{0,n}}\left(\phi_j^*\right)
\phi_j^*(x) +R_n(x)
\end{eqnarray*}
Thus the approximate (due to ${\cal R}_{0,n}$ being informative) influence curve is given by 
\[
D_{Q_{0,n},x}\equiv \sum_{j\in {\cal R}_{0,n}}  S_{Q_{0,n}}\left(\phi_j^*\right)
\phi_j^*(x),\]
\end{lemma} 
In this subsection we do not assume that the loss function $L$ generates scores in the tangent space of the model ${\cal M}$ for $P_0$ and thereby behaves as a log-likelihood loss, but that case is treated separately in the next subsection. On the other hand, we still know that the scores $S_{Q_{0,n}}(\phi_j^*)$ are elements in $L^2_0(P_0)$ and thus represents scores in a nonparametric model. 
Instead we utilize the expansion we have for $Q_n(x)-Q_{0,n}(x)=P_n D_{Q_{0,n},x}+R_n(x)$, and use the delta-method to establish asymptotic linearity for $\Phi(Q_n)-\Phi(Q_{0,n})$. 
The analysis of $\Phi(Q_{0,n})-\Phi(Q_0)$ is carried out separately.

{\bf Asymptotic linearity of $\Phi(Q_n)-\Phi(Q_{0,n})$ with (data dependent) influence curve $D_{\Phi,Q_{0,n}}$:}
We can restrict $\Phi: D^{(k)}({\cal R}_{0,n})\rightarrow\openr$  and $d\Phi_{Q_{0,n}}: D^{(k)}({\cal R}_{0,n})\rightarrow\openr$, since we apply it to $Q_n,Q_{0,n}$ and $Q_n-Q_{0,n}$  which are elements of the linear space $D^{(k)}({\cal R}_{0,n})$. 
Given we have $(Q_n-Q_{0,n})(x)=P_n D_{Q_{0,n},x}+R_n(x)$, and $\Phi$ is differentiable we obtain \[
\Phi(Q_n)-\Phi(Q_{0,n})=P_n d\Phi_{Q_{0,n}}(D_{Q_{0,n},\cdot}) +d\Phi_{Q_{0,n}}(R_n)+R_{\Phi}(Q_n,Q_{0,n}).\]  By the fact that (conservatively) $d_0^{1/2}(Q_n,Q_{0,n})=o_P(n^{-1/4})$, it is reasonable to assume that $R_{\Phi}(Q_n,Q_{0,n})=o_P(n^{-1/2})$. In fact, we expect this remainder to behave as $O_P(n^{-2k^*/(2k^*+1)})$ up till $\log n$-factors, if we select the tuning parameters so that our optimal rate of convergence (in sup-norm) is achieved.

Using the smoothness of the linear operator $d\Phi_{Q_{0,n}}$, one then needs to show $d\Phi_{Q_{0,n}}(R_n)=o_P(n^{-1/2})$. This requires working out the formulas for $d\Phi_{Q_{0,n}}(\tilde{r}_n)$; $d\Phi_{Q_{0,n}}(E_n)$; $d\Phi_{Q_{0,n}}(R_{m,n})$, $m=1,3$, and showing that each one is indeed $o_P(n^{-1/2})$, which generally speaking will be relatively straightforward. In particular, $R_{1,n}(x)$, $R_{3,n}(x)$ will generally be $O(d_0(Q_n,Q_{0,n}))$ and thus be $o_P(n^{-1/2})$. 
Moreover, $d\Phi_{Q_{0,n}}(\tilde{r}_n)=\sum_{j\in {\cal R}_{0,n}}P_n S_{Q_n}(\phi_j^*) d\Phi_{Q_{0,n}}(\phi_j^*)$, which equals zero for the relax and sieve HAL-MLE and is well controlled by the HAL-MLE as well (possibly some undersmoothing needed). 
Let's now consider $d\Phi_{Q_{0,n}}(E_n)$.
We can view $d\Phi_{Q_{0,n}}$ as a bounded linear operator on $D^{(k)}({\cal R}_{0,n})$ with inner product $\langle Q_1,Q_2\rangle_{J_{0,n}}$ that was used in our representation of $(Q_n-Q_{0,n})$ as a linear combination of $\phi_j^*$. Thus, $\phi_j^*$ is an orthonormal basis w.r.t. this inner product. By the Riesz-representation theorem we have 
$d\Phi_{Q_{0,n}}(h)=\langle \tilde{D}_{\Phi,Q_{0,n}},h\rangle_{J_{0,n}}$. By linearity of $\phi\rightarrow S_{Q_{0,n}}(\phi)$ it follows that the empirical process term can be written as:
\[
d\Phi_{Q_{0,n}}(E_n)=(P_n-P_0)(S_{Q_n}-S_{Q_{0,n}})(g_n),
\]
where $g_n\equiv \sum_{j\in {\cal R}_{0,n}}\langle \tilde{D}_{\Phi,Q_{0,n}},\phi_j^*\rangle_{J_{0,n}}\phi_j^*$.
Note that $g_n=\Pi_{J_{0,n}}(\tilde{D}_{\Phi,Q_{0,n}}\mid D^{(k)}({\cal R}_{0,n}))$, but that equals $\tilde{D}_{\Phi,Q_{0,n}}$ itself. So we have
\[
d\Phi_{Q_{0,n}}(E_n)=(P_n-P_0)(S_{Q_n}-S_{Q_{0,n}})(\tilde{D}_{\Phi,Q_{0,n}}).\]
This can then be analyzed with standard empirical process theory, using $Q_n,Q_{0,n},\tilde{D}_{\Phi,Q_{0,n}}\in D^{(0)}_M([0,1]^d)$ for some $M<\infty$, and that $d_0(Q_n,Q_{0,n})=O_P(d_n/n)$. 

 One has then shown that $\Phi(Q_n)$ is asymptotically linear with (data dependent) influence curve $D_{\Phi,Q_{0,n}}=d\Phi_{Q_{0,n}}(D_{Q_{0,n},\cdot})$, where we are reminded that $D_{Q_{0,n},x}= \sum_{j\in {\cal R}_{0,n}} S_{Q_{0,n}}\left(\phi_j^*\right)
\phi_j^*(x)$.

{\bf Asymptotic convergence of $D_{\Phi,Q_{0,n}}$ to a fixed function $D_{\Phi,Q_0}$:}
We now want to show that $D_{\Phi,Q_{0,n}}$ converges in $L^2(P_0)$ to a fixed function $D_{\Phi,Q_0}\in L^2_0(P_0)$, which then represents the influence curve of $\Phi(Q_n)-\Phi(Q_{0,n})$. We could simply assume this, but it is helpful to obtain some insight about why one now expects that this influence curve has a stable variance, contrary to the influence curve $D_{Q_{0,n},x}$ of $(Q_n-Q_{0,n})(x)$, which needed to be divided by $d_{0,n}^{1/2}$ in order to keep it bounded in $L^2(P_0)$. 


As above for $d\Phi_{Q_{0,n}}(E_n)$, we have $D_{\Phi,Q_{0,n}}=S_{Q_{0,n}}(g_n)$ while $g_n=\tilde{D}_{\Phi,Q_{0,n}})$, so that
\[
D_{\Phi,Q_{0,n}}=S_{Q_{0,n}}(\tilde{D}_{\Phi,Q_{0,n}}).\]
So we need that $P_0(\tilde{D}_{\Phi,Q_{0,n}}-\tilde{D}_{\Phi,Q_0})^2\rightarrow_p 0$ for a limit $\tilde{D}_{\Phi,Q_0}$ and that $S_{Q_{0,n}}(\tilde{D}_{\Phi,Q_{0,n}})$ falls in  a Donsker class with probability tending to 1.

Above we proved the following lemma. 
\begin{lemma}\label{lemmagn}
Suppose that $d\Phi_{Q_{0,n}}:D^{(k)}({\cal R}_{0,n})\rightarrow \openr$ is represented as 
\[
d\Phi_{Q_{0,n}}(h)=\langle \tilde{D}_{\Phi,Q_{0,n}},h\rangle_{J_{0,n}}\]
with $\tilde{D}_{\Phi,Q_{0,n}}\in D^{(k)}({\cal R}_{0,n})$.
We have \[
D_{\Phi,Q_{0,n}}=S_{Q_{0,n}}\tilde{D}_{\Phi,Q_{0,n}}.\]
If for a limit $D_{\Phi,Q_0}$ we have $P_0\{S_{Q_{0,n}}\tilde{D}_{\Phi,Q_{0,n}})-D_{\Phi,Q_0}\}^2\rightarrow_p 0$ and $S_{Q_{0,n}}\tilde{D}_{\Phi,Q_{0,n}}$ falls in a Donsker class with probability tending to 1, then 
$(P_n-P_0)D_{\Phi,Q_{0,n}}=(P_n-P_0)D_{\Phi,Q_0}+o_P(n^{-1/2})$.
\end{lemma}

We also have $\sigma^2_{\Phi,0,n}=P_0 D_{\Phi,Q_{0,n}}^2=O(1)$ and converges to a fixed $\sigma_{\Phi,0}$, so that the CLT establishes $n^{1/2}(\Phi(Q_n)-\Phi(Q_{0,n}))\Rightarrow_d N(0,\sigma^2_0)$. More importantly,  $\Phi(Q_n)$ is asymptotically linear with fixed influence curve $D_{\Phi,Q_0}$ being the $L^2(P_0)$-limit of $d\Phi_{Q_{0,n}}(D_{Q_{0,n},\cdot})$.

If $Q_n$ acts as a sieve MLE for a sequence of working models $D^{(k)}({\cal R}_{0,n})$ approximating an oracle model $D^{(k)}({\cal R}_0)$ (i.e., ${\cal R}_0$ depending on $Q_0$ and satisfying that $Q_0\in D^{(k)}({\cal R}_0)$) strictly smaller than the assumed model for $Q_0$ implied by the statistical model ${\cal M}$, such as $D^{(k)}([0,1]^d)$, then $Q_n$ will typically be super-efficient.  

Finally, one needs to establish that $\Phi(Q_{0,n})-\Phi(Q_0)=o_P(n^{-1/2})$. 
In the theorem below we simply assume the convergence of $D_{\Phi,Q_{0,n}}$ while the above can be used to actually prove this in a particular example. 

\begin{theorem}
Consider the setting of Theorem \ref{theoremasnormgen} so that
 $(Q_n-Q_{0,n})(x)=P_n D_{Q_{0,n},x}+R_n(x)$, where 
 $R_n(x)=\sum_{j\in {\cal R}_{0,n} }R_n({\phi}_j^*)\phi_j^*(x)$, where
 $R_n({\phi}_j^*)=-r_n({\phi}_j^*)-E_n({\phi}_j^*)+
 (R_{1,n}+R_{3,n})({\phi}_j^*)$. 
 Assume $\Phi$ is differentiable at $Q_{0,n}$ so that
 $\Phi(Q_n)-\Phi(Q_{0,n})=P_n d\Phi_{Q_{0,n}}(D_{Q_{0,n},\cdot}) +d\Phi_{Q_{0,n}}(R_n)+R_{\Phi}(Q_n,Q_{0,n})$.  
 Let $D_{\Phi,Q_{0,n}}=d\Phi_{Q_{0,n}}(D_{Q_{0,n},\cdot})$.
 Define $\tilde{D}_{\Phi,Q_{0,n}}\in D^{(k)}({\cal R}_{0,n})$ be the canonical gradient of $d\Phi_{Q_{0,n}}(h)=\langle \tilde{D}_{\Phi,Q_{0,n}},h\rangle_{J_{0,n}}$. 
 We have \[
 D_{\Phi,Q_{0,n}}=S_{Q_{0,n}}\tilde{D}_{\Phi,Q_{0,n}}.\]
  
 {\bf Assumptions:}
 $d\Phi_{Q_{0,n}}(R_n)=o_P(n^{-1/2})$; $R_{\Phi}(Q_n,Q_{0,n})=o_P(n^{-1/2})$;
for a limit $D_{\Phi,Q_0}$ we have $P_0 \{D_{\Phi,Q_{0,n}}-D_{\Phi,Q_0}\}^2\rightarrow_p 0$, where one can use Lemma \ref{lemmagn}; $D_{\Phi,Q_{0,n}}$ is an element of a fixed $P_0$-Donsker class with probability tending to 1. 
 
 Then, $\Phi(Q_n)-\Phi(Q_{0,n})=P_n D_{\Phi,Q_0}+o_P(n^{-1/2})$.
 
 If also $\Phi(Q_{0,n})-\Phi(Q_0)=o_P(n^{-1/2})$ (by applying Lemma \ref{lemmaoraclebias} below), then $\Phi(Q_n)$ is asymptotically linear at $P_0$ with influence curve $D_{\Phi,Q_0,g_0}$.
 \end{theorem}

 {\bf Analysis of $\Phi(Q_{0,n})-\Phi(Q_0)$:}
 The following lemma provides the conditions under which $\Phi(Q_{0,n})-\Phi(Q_0)=o_P(n^{-1/2})$.
 \begin{lemma} \label{lemmaoraclebias}
 Let ${\cal M}_0\subset {\cal M}$ be an oracle model depending on $D^{(k)}({\cal R}_0)$ so that $p_0\in {\cal M}_0$. For example, in the case that ${\cal M}=\{p_Q: Q\in D^{(k)}([0,1]^d)\}$ we would have ${\cal M}_0=\{p_Q: Q\in D^{(k)}({\cal R}_0)\}$ with $Q_0\in D^{(k)}({\cal R}_0)$, where ${\cal R}_{0,n}\subset {\cal R}_0$ (with probability tending to 1). 
 
 We can then consider the oracle parameter 
$\Phi_0:{\cal M}_0\rightarrow\openr$. The fact that $\Phi:{\cal M}\rightarrow \openr$ is pathwise differentiable implies that $\Phi_0:{\cal M}_0\rightarrow\openr$ (just $\Phi$ but on smaller set of probability distributions)  is pathwise differentiable but possibly with a smaller tangent space $T_{{\cal M}_0}(P)$ instead of the model tangent space $T_{{\cal M}}(P)$ of model ${\cal M}$.
 Let $D_{\Phi_0,P}$ be a gradient of $\Phi_0$ at $P$, not necessarily the canonical gradient $D_{\Phi_0,P}^*\in T_{{\cal M}_0}(P)$. We might also denote this gradient with $D_{{\cal R}_0,P}$.  Suppose that $D_{\Phi_0,P}$ depends on $P$ through $Q(P)$ and $G(P)$. Then we can write $D_{\Phi_0,Q,G}$. 
 Let $D_{\Phi_{0,n},Q_{0,n},G_0}$ be a  gradient of $\Phi_{0,n}:{\cal M}\rightarrow\openr$ defined by $\Phi_{0,n}(Q_0)=\Phi(Q_{0,n})$.      
 Let 
 $R_{\Phi_0,0}(f,Q_0)\equiv \Phi_0(Q)-\Phi_0(Q_0)+P_0 D_{\Phi_0,Q,G_0}$ be the exact second order remainder for this parameter $\Phi_0$. 
 Then, we have
 \[
 \Phi(Q_{0,n})-\Phi(Q_0)=-P_0 D_{\Phi_0,Q_{0,n},G_0}+R_{\Phi_0,0}(Q_{0,n},Q_0).\] 
 
 {\bf Assumptions:}
 $R_{\Phi_0,0}(Q_{0,n},Q_0)=o_P(n^{-1/2})$, where one can use that $\pl Q_{0,n}-Q_0\pl_{\infty}=O(r(d,J_{0,n})^{k+1})$; $P_0 D_{\Phi_0,Q_{0,n},G_0}=o_P(n^{-1/2})$. 
 
 Then, 
 \[
 \Phi(Q_{0,n})-\Phi(Q_0)=o_P(n^{-1/2}).\]
 
{\bf Sufficient condition for $P_0 D_{\Phi_0,Q_{0,n},G_0}=o_P(n^{-1/2})$:}
Assume $P_0 D_{\Phi_{0,n},Q_{0,n},G_0}=0$ or $o_P(n^{-1/2})$. 
  Let $P_{Q_{0,n},G_0}$ represent an approximation of $P_0$ compatible with $Q_{0,n}$ and $G_0$ in the sense that $Q_{{\cal R}_{0,n}}(P_{Q_{0,n},G_0})=Q_{0,n}$ with  $Q_{{\cal R}_{0,n}}(P)=\arg\min_{Q\in D^{(k)}({\cal R}_{0,n})}P L(Q)$ and $G(P_{Q_{0,n},G_0})=G_0$. 
 Then,
 \begin{eqnarray*}
 P_0 D_{\Phi_0,Q_{0,n},G_0} &=&(P_0-P_{Q_{0,n},G_0})\{D_{\Phi_0,Q_{0,n},G_0}-D_{\Phi_{0,n},Q_{0,n},G_0}\}.
 \end{eqnarray*}
 Thus, if $p_0=dP_0/d\mu>\delta>0$ for some $\delta>0$, then we have
 \[
 P_0 D_{\Phi_0,Q_{0,n},G_0}  =O\left( 
  \pl p_{Q_{0,n},G_0}-p_0\pl_{\mu}\pl D_{\Phi_0,Q_{0,n},G_0}-D_{\Phi_{0,n},Q_{0,n},G_0}\pl_{\mu}\right).\]
\end{lemma}
The first $L^2(\mu)$-norm can be bounded in terms of $\pl Q_{0,n}-Q_0\pl_{\infty}$, while the second concerns the approximation error of the linear span of $\{S_j(Q_{0,n}):j\in {\cal R}_{0,n}\}$ for approximating an element in the span of $\{S_j(Q_{0,n}):j\in {\cal R}_0\}$. One expects both to be bounded by $r(d,J_{0,n})^{k+1}$ so that this term will behave as $r(d,J_{0,n})^{2(k+1)}$.

We have that $Q_{0,n}$ is an oracle MLE for $D^{(k)}({\cal R}_{0,n})$ with the latter working model growing towards the oracle model $D^{(k)}({\cal R}_0)$.  The linear span of score equations $P_0 S_{j}(Q_{0,n})=0$, $j\in {\cal R}_{0,n}$, solved by $Q_{0,n}$ approximate  the closure of the linear span of $S_j(Q_{0,n})$, $j\in {\cal R}_0$. 
We need that this linear span approximates a gradient of $D_{\Phi_{0,n},Q_{0,n},G_0}$ of $\Phi_{0,n}$ which itself then approximates a gradient $D_{\Phi_0,Q_{0,n},G_0}$ of $\Phi_0$.
For the log-likelihood loss,  one can set the latter gradients equal to the canonical gradients in the tangent space $T_{{\cal M}_0}(P_{Q_{0,n}})$. 

\subsection{Log-likelihood behaving loss functions}

In this subsection we carry out the typical analysis of MLE/TMLE based on solving the efficient influence curve equation, which thereby relies on the score equations solved by $Q_n$ to approximate the efficient influence curve equation. Thus, the loss function has to generate log-likelihood scores. 

Recall $Q_n$ is an MLE for working model $D^{(k)}({\cal R}_{0,n})$.
Define the parameter $P_0\rightarrow \Phi(Q_{0,n})$ on the actual statistical model ${\cal M}$ for $P_0$, which could be the nonparametric model or a model of form $\{p_Q: Q\in D^{(k)}([0,1]^d)\}$.  That is, this parameter maps $P_0$ into $Q_{0,n}=\arg\min_{Q\in D^{(k)}({\cal R}_{0,n})}P_0L(Q)$ and then plugs this in the parameter mapping $\Phi(Q_{0,n})$. 
Let's denote this parameter with $\Phi_{0,n}:{\cal M}\rightarrow\openr$ so that $\Phi_{0,n}(P)=\Phi(Q_{{\cal R}_{0,n}}(P))$ with $Q_{{\cal R}_{0,n}}(P)=\arg\min_{Q\in D^{(k)}({\cal R}_{0,n})}P_0L(Q)$. Let $D_{{\cal R}_{0,n},P}$ be the canonical gradient of $\Phi_{0,n}$, and assume that the generalized $L(\cdot)$-scores solved by $Q_{0,n}$ imply that $P_0 D_{{\cal R}_{0,n},P_0}=0$ or $o_P(n^{-1/2})$.
Suppose that $D_{{\cal R}_{0,n},P_0}$ depends on $P_0$ through $Q_0$  and 
a nuisance parameter $G_0$. Then, we can write $D_{{\cal R}_{0,n},Q_0,G_0}$ while we have $P_0 D_{{\cal R}_{0,n},Q_{0,n},G_0}=0$. 
Moreover, due to $Q_n$ solving scores, we assume that $P_n D_{{\cal R}_{0,n},Q_n,G_0}=o_P(n^{-1/2})$ as well. 

We also consider the oracle parameter $\Phi_0$ defined on the oracle model ${\cal M}_0$ that is defined by $\Phi_{0}:{\cal M}\rightarrow\openr$ so that $\Phi_{0}(P)=\Phi(Q_{{\cal R}_{0}}(P))$ with $Q_{{\cal R}_{0}}(P)=\arg\min_{Q\in D^{(k)}({\cal R}_{0})}P_0L(Q)$. By assumption $Q_0\in D^{(k)}({\cal R}_0)$ so that $\Phi_0(P_0)=\Phi(Q_0)$. Let $D_{{\cal R}_{0},P}$ be the canonical gradient of $\Phi_{0}$ at $P$.

As an example consider the case that ${\cal M}=\{p_Q: Q\in D^{(k)}([0,1]^d)\}$. 
For the given working model $D^{(k)}({\cal R}_{0,n})$, the efficient influence curve of $\beta_{0,n}$ is given by $I_{n,\beta_{0,n}}^{-1}S_{\beta_{0,n}}$, where
$S_{\beta_{0,n}}=\frac{d}{d\beta_{0,n}}\log p_{\sum_{j\in {\cal R}_{0,n}}\beta_{0,n}(j)\phi_j}$ is the score vector, and $I_{n,\beta_{0,n}}=-P_0\frac{d^2}{d\beta_{0,n}}\log p_{\sum_{j\in {\cal R}_{0,n}}\beta_{0,n}(j)\phi_j}$ is the information matrix. Let $g_{n,\beta_{0,n}}=(d/d\beta_{0,n}(j) \Phi(\sum_{l\in {\cal R}_{0,n}}\beta_{0,n}(l)\phi_l):j)$ be the gradient of the target parameter as a function of $\beta_{0,n}$. Then the efficient influence curve is given by $D_{{\cal R}_{0,n},Q_{0,n}}=g_{n,\beta_{0,n}}^{\top} I_{n,\beta_{0,n}}^{-1}S_{\beta_{0,n}}$, which thus only depends on $Q_{0,n}$ or equivalently $\beta_{0,n}$.
Moreover, this efficient influence curve is a linear combination of scores $S_{Q_{0,n}}(\phi_j)$ so that indeed $P_0 D_{{\cal R}_{0,n},Q_{0,n}}=0$. Moreover, $P_n D_{{\cal R}_{0,n},Q_n}=0$ for the sieve and relax HAL-MLE and will be $o_P(n^{-1/2})$ for the HAL-MLE possibly requiring a little undersmoothing. 

Let $R_{\Phi_{0,n}}(Q_n,Q_{0,n})\equiv \Phi(Q_n)-\Phi(Q_{0,n})+P_0 D_{{\cal R}_{0,n},Q_{0,n},G_0}$ be the exact second order remainder for the parameter $\Phi_{0,n}$. This is generally a second order term so that we can assume $R_{\Phi_{0,n}}(Q_n ,Q_{0,n})=o_P(n^{-1/2})$, given we already have a rate of convergence $d_0(Q_n,Q_{0,n})$ faster than $n^{-1/2}$. 
Then, we have
\[
\Phi(Q_n)-\Phi(Q_{0,n})=-P_0 D_{{\cal R}_{0,n},Q_n,G_0}+R_{\Phi_{0,n}}(Q_n,Q_{0,n}).\]
Moreover, combined with $P_n D_{{\cal R}_{0,n},Q_n,G_0}=o_P(n^{-1/2})$ this yields
\[
\Phi(Q_n)-\Phi(Q_{0,n})=(P_n-P_0) D_{{\cal R}_{0,n},Q_n,G_0}+o_P(n^{-1/2}).\]
Assume $(P_n-P_0)\{D_{{\cal R}_{0,n},Q_n,G_0}-D_{{\cal R}_{0,n},Q_{0,n},G_0}\}=o_P(n^{-1/2})$. Then, we have shown the desired asymptotic linearity w.r.t. the data adaptive target parameter (i.e., ${\cal R}_{0,n}={\cal R}_n$ depends on the data)
\[
\Phi(Q_n)-\Phi(Q_{0,n})=(P_n-P_0)D_{{\cal R}_{0,n},Q_{0,n},G_0}+o_P(n^{-1/2}).\]
Finally, assume that $P_0 \{D_{{\cal R}_{0,n},Q_{0,n},G_0}-D_{{\cal R}_0,Q_0,G_0}\}^2\rightarrow_p 0$ for a fixed $D_{{\cal R}_0,Q_0,G_0}$, presumably the canonical gradient of oracle statistical target parameter $\Phi_0$ at $P_0$ as defined above, and that $D_{{\cal R}_{0,n},Q_{0,n},G_0}$ is an element of  a $P_0$-Donsker class with probability tending to 1. 
So then we have shown
\[
\Phi(Q_n)-\Phi(Q_{0,n})=P_n D_{{\cal R}_0,Q_0,G_0}+o_P(n^{-1/2}).\]
If $D^{(k)}({\cal R}_0)$ is a stronger model for $Q_0$ than the one assumed by the actual statistical model ${\cal M}$, then $D_{{\cal R}_0,Q_0,G_0}$ represents an efficient influence curve for a smaller statistical model ${\cal M}_0$, so that in that case $\Phi(Q_n)$ is super-efficient. 

 We proved the following theorem.
 \begin{theorem}\ \nl
  {\bf Definitions:}
Define the parameter $P_0\rightarrow \Phi(Q_{0,n})$ on the actual statistical model ${\cal M}$ for $P_0$, which could be the nonparametric model or a model of form $\{p_Q: Q\in D^{(k)}([0,1]^d)\}$.  That is, this parameter maps $P_0$ into $Q_{0,n}=\arg\min_{Q\in D^{(k)}({\cal R}_{0,n})}P_0L(Q)$ and then plugs this in the parameter mapping $\Phi(Q_{0,n})$. 
Let's denote this parameter with $\Phi_{0,n}:{\cal M}\rightarrow\openr$ so that $\Phi_{0,n}(P)=\Phi(Q_{{\cal R}_{0,n}}(P))$ with $Q_{{\cal R}_{0,n}}(P)=\arg\min_{Q\in D^{(k)}({\cal R}_{0,n})}P_0L(Q)$. Let $D_{{\cal R}_{0,n},P}$ be the canonical gradient of $\Phi_{{\cal R}_{0,n}}$, and we assume below that the generalized $L(\cdot)$-scores solved by $Q_{0,n}$ imply that $P_0 D_{{\cal R}_{0,n},P_0}=o_P(n^{-1/2})$.
Suppose that $D_{{\cal R}_{0,n},P_0}$ depends on $P_0$ through $Q_0$  and 
a nuisance parameter $G_0$. Then, we can denote $D_{{\cal R}_{0,n},P_0}$ with  $D_{{\cal R}_{0,n},Q_0,G_0}$. 
Let $R_{\Phi_{0,n}}(Q_n,Q_{0,n})\equiv \Phi(Q_n)-\Phi(Q_{0,n})+P_0 D_{{\cal R}_{0,n},Q_{0,n},G_0}$ be the exact second order remainder for target parameter $\Phi_{0,n}:{\cal M}\rightarrow\openr$.
\newline
{\bf Assumptions:}  $P_n D_{{\cal R}_{0,n},Q_{0,n},G_0}=o_P(n^{-1/2})$; $P_n D_{{\cal R}_{0,n},Q_n,G_0}=o_P(n^{-1/2})$;
$R_{\Phi_{0,n}}(Q_n ,Q_{0,n})=o_P(n^{-1/2})$; $(P_n-P_0)\{D_{{\cal R}_{0,n},Q_n,G_0}-D_{{\cal R}_{0,n},Q_{0,n},G_0}\}=o_P(n^{-1/2})$; $P_0 \{D_{{\cal R}_{0,n},Q_{0,n},G_0}-D_{{\cal R}_0,Q_0,G_0}\}^2\rightarrow_p 0$ for a fixed $D_{{\cal R}_0,Q_0,G_0}$; and that $D_{{\cal R}_{0,n},Q_{0,n},G_0}$ is an element of  a $P_0$-Donsker class with probability tending to 1. 

Then,
\[
\Phi(Q_n)-\Phi(Q_{0,n})=P_n D_{{\cal R}_0,Q_0,G_0}+o_P(n^{-1/2}).\]
If, in addition, $\Phi(Q_{0,n})-\Phi(Q_0)=o_P(n^{-1/2})$ (by applying Lemma \ref{lemmaoraclebias}), then
\[
\Phi(Q_n)-\Phi(Q_0)=P_n D_{{\cal R}_0,Q_0,G_0}+o_P(n^{-1/2}).\]
\end{theorem}

\section{Summary and concluding remarks. }\label{section12}
In this article we proposed and analyzed three general $k$-th order spline sieve/relax lasso/lasso minimum loss estimators of  a functional parameter. Due to the adaptivity of the $k$-th order spline HAL-MLE, and the computational advantage of only having to tune the $L_1$-norm, we recommend the HAL-MLEs in practice, but as we highlighted the specification of index sets for the sieve-MLE provide  good initial large working models ${\cal D}^{(k)}({\cal R}(d,{\bf J}_{max}))$ for the lasso based estimators. 

{\bf Remarkable asymptotic rates of convergence for these $k$-th order spline sieve and HAL-MLEs:}
In this article we showed that $\sim J$-dimensional $k$-th order spline working models $D^{(k)}({\cal R}(d,{\bf J}))$ provide very strong $O(C(M)r(d,J)^{k+1})\approx O(1/J^{k+1})$ supremum  norm approximations of  the $k$-th order smoothness class $D^{(k)}_M([0,1]^d)$ for any $M<\infty$, where $M$ represents a bound on the $k$-th order sectional variation norm. In combination with the  MLE rate of pointwise convergence $(n/d)^{-1/2}$ for $Q_n-Q_{0,n}$, this then resulted in supremum norm convergence of these sieve and HAL-MLEs at the remarkable rate $n^{-k^*/(2k^*+1)}$ up till power of $\log n$-factors, where $k^*=k+1$. That is, these estimators of $d$-variate functions achieve the same (up till power of $\log n$-factor)  well known minimax rates of convergence for estimation of a univariate function that is $k$-times continuously differentiable, where $d$ only enters through the power of $\log n$-factor.  In addition we even showed convergence in sectional variation norm at rate $n^{-k/(2k+1)}$ up till power of $\log n$-factors, showing that these $k$-th order spline sieve and HAL-MLEs can effectively estimate derivatives of the target function as well. 

Obviously, our ''curse of dimensionality free''-rates of convergence results imply that the curse of dimensionality enters through the constant. Nonetheless, the asymptotic superiority relative to other algorithms  suggests strong practical performance relative to other machine learning approaches. In addition, the constant is mostly driven by the actual $k$-th order sectional variation norm $M_0$ of the true target function $Q_0$, so that the performance is highly adaptive to the actual complexity of the true function. 

{\bf Asymptotic normality of sieve and HAL-MLEs at same rate, when undersmoothing the $L_1$-norm:}
Moreover, we show that if the $L_1$-norm is chosen large enough so that
$r(d, d_n)^{k+1}$, $d_n$ size of working model $D^{(k)}({\cal R}_n)$,  divided by $(n/d_n)^{-1/2})$ approximates zero at a (power of) $\log n$-rate, then these estimators are asymptotically normally distributed pointwise, at the same rates as stated above. 
As a consequence, these estimators result in asymptotically valid confidence intervals for $Q_0(x)$ or smooth functions thereof. 

{\bf Statistical inference for non-pathwise differentiable statistical target parameters:}
That is, our results for these estimators establish statistical inference for non-pathwise differentiable parameters of the true distribution $P_0\in {\cal M}$, uniform rates of convergence, and corresponding simultaneous confidence bands for functional parameters. Indeed in some prior simulation studies not reported here we observe for reasonable sample sizes valid finite sample coverage for estimation of functionals  up till dimension 3, such as conditional treatment effects adjusting for some key covariates.

{\bf Cross-validation to select choice of smoothness, submodel, initial screening: }
These asymptotic results were proved assuming that the true function $Q_0$ is an element of a given $k$-th order smoothness class $D^{(k)}_M([0,1]^d)$ or submodel thereof. In practice, one prefers to not have to a priori commit to these choices. 
 
We characterized $D^{(k)}([0,1]^d)$ as an infinite linear combinations of basis functions indexed by a set ${\cal R}^k(d)$. This allowed us to define general submodels by replacing ${\cal R}^k(d)$ by a subset ${\cal R}^{k,*}(d)\subset {\cal R}^k(d)$, resulting in submodels $D^{(k)}( {\cal R}^{k,*}(d))$ of $D^{(k)}([0,1]^d)$. Beyond the $L_1$-norm $C$, the $k$-th order spline relax HAL-MLE (as do the others) still relies on the choice of smoothness $k$, the choice of submodel $D^{(k)}_{{\cal R}^{k,*}(d)}([0,1]^d)\subset D^{(k)}([0,1]^d)$ of the $k$-th order smoothness class, and the $J_{max}$ to define the initial index set ${\cal R}^k(d,{\bf J}_{max})$ with ${\bf J}_{max}=(J_{max},\ldots,J_{max})$. Regarding the latter, we already proposed such a set  in our prior work, specifically, for any (by model) allowed $\bar{s}(k+1)$ with $\mid s_{k+1}\mid>0$, use the knot-point set ${\cal R}(\bar{s}(k+1),n)=\{X_i(s_{k+1}):i=1,\ldots,n\}$, in which case ${\bf J}_{max}=(n,\ldots,n)$, across all $\bar{s}((k+1)$ our model allows. Of course,  we can then define subsets of the latter to generate smaller choices ${\cal R}(d,{\bf J}_{max})$ with $J_{max}<n$.  Suppose that the choice of subsets ${\cal R}^{k,*}_m(d)$ are indexed by another integer $m$.  We could select smoothness index $k$, the submodel index $m$, and prescreening index $J_{max}$ with the cross-validation selector, while using an undersmoother $C_n$ for $C$ across values larger or equal than its cross-validation selector. In that manner, for each $C$,  we implemented a discrete super-learner $Q_{n,C}$, where the candidate estimators $\hat{Q}_{k,m,J_{max},C}(P_n)$ are indexed by smoothness, submodel, and initial (non-informative) screening, beyond $L_1$-norm $C$. 

{\bf Adaptivity to smoothness/submodel:}
By asymptotic equivalence of the  $C$-specific cross-validation selector with the oracle selector, this  $C$-specific discrete super-learner will asymptotically select the best possible choice, so that the estimator will achieve a rate of convergence as if the true smoothness and  correct submodel would have been known a priori.

{\bf Undersmooth selector of $L_1$-norm $C$:}
We propose plateau (Lepski or variance plateaus) selectors for $C_n$, which can be applied to the collection of estimators $\{Q_{n,C}(x):C\}$ with corresponding variance estimator $\{\sigma^2_{n,C}(x):C\}$ (e.g., Chapter 25 of \citep{vanderLaan&Rose18}). For a given $C$, one selects a size $d_n(C)$ for the working model, which then provides our variance estimators  $\sigma^2_{n,C}(x)$ for $Q_n(x)$ according to our asymptotic normality theorems. 
Such selectors could be pointwise for a given $x$, but  could also be tailored for a global set of points. In past work we have seen good performance of the plateaus selectors and we are currently practically investigating global undersmoothers so that one compatible HAL-fit can be used for all points $x$. 

{\bf Asymptotic normality should be preserved under adaptive selection of tuning parameters:}
Our asymptotic normality theorems should still apply as long as the cross-validation selector $(k_n,m_n)$ converges to a fixed choice $(k_0,m_0)$ with probability tending to 1.

{\bf Generalization to non-linear risk functions:}
In this article we obtained  rates of  uniform convergence and asymptotic normality of sieve and HAL-MLEs for general {\em linear} risk functions $R_P(Q)\equiv P L(Q)$ implied by a loss function $L(Q)(O)$ where the target function $Q_0=Q(P_0)=\arg\min_{Q}R_{P_0}(Q)$ is defined as minimizer of the true risk function. In many applications of interest the risk depends on $P$ in a non-linear fashion. For example, the risk could be defined as the expectation of an inverse of probability censoring weighted loss function so that it also requires estimation of the censoring distribution (e.g., \citep{vanderLaan&Dudoit03}). In particular, the conditional treatment effect can be represented as a minimizer of a general risk function, and similarly, for the optimal individualized treatment rule \citep{vanderLaan:Luedtke15,luedtke2016super}.

Given that 1)  our uniform approximation results for $D^{(k)}([0,1]^d)$ have nothing to do with the choice of risk function, and 2) that our proof of the Theorem in Section \ref{section7}, to transfer the approximation error of a working model $D^{(k)}({\cal R}(d,{\bf J}))$ w.r.t. $D^{(k)}_M([0,1]^d)$ to the approximation error of the loss-based projection  $Q_{0,n}$ of $Q_0$ onto this finite dimensional working model w.r.t $Q_0$ itself, is general, our rate and asymptotic normality results should generalize to non-linear risk functions as well. In that case,  one estimates $R_{P_0}(Q)$ with an estimator $R_n(Q)$ that is asymptotically linear uniformly in $Q\in D^{(0)}([0,1]^d)$, and $Q_n$ is defined as a sieve or HAL-MLE of $Q\rightarrow R_n(Q)$ over finite dimensional $k$-th order spline working models $D^{(k)}({\cal R}(d,{\bf J}_{0,n,max}))$. Such global estimators of $(R_n(Q):f)$ could be obtained  by estimating $R_{P_0}(Q)$ by plugging in an  undersmoothed HAL-MLE $\tilde{P}_n$ and estimating the risk with $R_n(Q)=R_{\tilde{P}_n}(Q)$. 
In future work we will address these generalizations of $k$-th order HAL-MLE to essentially estimate  and provide inference for any target function of interest, not just once that can be defined as minimizers of linear risk functions. 

\bibliographystyle{plainnat}
\bibliography{TMLE1_TLB}

\begin{thebibliography}{26}
\providecommand{\natexlab}[1]{#1}
\providecommand{\url}[1]{\texttt{#1}}
\expandafter\ifx\csname urlstyle\endcsname\relax
  \providecommand{\doi}[1]{doi: #1}\else
  \providecommand{\doi}{doi: \begingroup \urlstyle{rm}\Url}\fi

\bibitem[Benkeser and van~der Laan(2016)]{Benkeser&vanderLaan16}
D.~Benkeser and M.J. van~der Laan.
\newblock The highly adaptive lasso estimator.
\newblock \emph{Proceedings of the IEEE Conference on Data Science and Advanced
  Analytics}, 2016.
\newblock To appear.

\bibitem[Bibaut and van~der Laan(2019)]{Bibaut&vanderLaan19}
A.~Bibaut and M.J. van~der Laan.
\newblock Fast rates for empirical risk minimization over cadlag functions with
  bounded sectional variation norm.
\newblock Technical report, Division of Biostatistics, University of
  California, Berkeley, 2019.

\bibitem[Bickel et~al.(1997)Bickel, Klaassen, Ritov, and Wellner]{Bickeletal97}
P.J. Bickel, C.A.J. Klaassen, Y.~Ritov, and J.~Wellner.
\newblock \emph{Efficient and adaptive estimation for semiparametric models}.
\newblock Springer, Berlin Heidelberg New York, 1997.

\bibitem[Fang et~al.(2019)Fang, Guntuboyina, and Sen]{Fangetal19}
B.~Fang, A.~Guntuboyina, and B.~Sen.
\newblock Multivariate extensions of isotonic regression and total variation
  denoising via entire monotonicity and hardy-krause variation.
\newblock \emph{https://arxiv.org/abs/1903.01395}, 2019.
\newblock Technical report.

\bibitem[Gao(2013)]{Gao13}
F.~Gao.
\newblock \emph{Bracketing entropy of high dimensional distributions}.
\newblock Birkhauser, In High Dimensional Probability VI, Volume 66, 2013.

\bibitem[Gill et~al.(1995)Gill, van~der Laan, and
  Wellner]{Gill&vanderLaan&Wellner95}
R.D. Gill, M.J. van~der Laan, and J.A. Wellner.
\newblock Inefficient estimators of the bivariate survival function for three
  models.
\newblock \emph{Annales de l'Institut Henri Poincar{\'e}}, 31:\penalty0
  545--597, 1995.

\bibitem[Ki et~al.(2021)Ki, Fang, and Guntuboyina]{Kietal21}
Dohyeong Ki, Billy Fang, and Adityanand Guntuboyina.
\newblock Mars via lasso, 2021.
\newblock URL \url{https://arxiv.org/abs/2111.11694}.

\bibitem[Luedtke and van~der Laan(2016)]{luedtke2016super}
A.R. Luedtke and M.J. van~der Laan.
\newblock Super-learning of an optimal dynamic treatment rule.
\newblock \emph{International Journal of Biostatistics}, 12\penalty0
  (1):\penalty0 305--332, 2016.

\bibitem[Neuhaus(1971)]{Neuhaus71}
G.~Neuhaus.
\newblock On weak convergence of stochastic processes with multidimensional
  time parameter.
\newblock \emph{Annals of Statistics}, 42:\penalty0 1285--1295, 1971.

\bibitem[Newey(2014)]{Newey94}
W.~Newey.
\newblock The asymptotic variance of semiparametric estimators.
\newblock \emph{Econometrica}, 62\penalty0 (6):\penalty0 1349--1382, 2014.

\bibitem[Qiu et~al.(2021)Qiu, Luedtke, and Carone]{Qiu&Luedtke&Carone21}
H.~Qiu, A.~Luedtke, and M.~Carone.
\newblock Universal sieve-based strategies for efficient estimation using
  machine learning tools.
\newblock \emph{Bernoulli}, 27\penalty0 (4):\penalty0 2300--2336, 2021.

\bibitem[Shen(1997)]{Shen97}
X.~Shen.
\newblock On methods of sieves and penalization.
\newblock \emph{Annals of Statitics}, 252\penalty0 (6):\penalty0 2555--2591,
  1997.

\bibitem[Shen(2007)]{Shen07}
X.~Shen.
\newblock Large sample sieve estimation of semiparametric models.
\newblock \emph{Chapter in Handbook of Econometrics}, 76\penalty0
  (00):\penalty0 0000, 2007.

\bibitem[van~der Laan et~al.(2019)van~der Laan, Benkeser, and
  Cai]{vanderLaan&Benkeser&Cai19}
Mark~J. van~der Laan, David Benkeser, and Weixin Cai.
\newblock Efficient estimation of pathwise differentiable target parameters
  with the undersmoothed highly adaptive lasso, 2019.

\bibitem[van~der Laan(2015)]{vanderLaan15}
M.J. van~der Laan.
\newblock A generally efficient targeted minimum loss-based estimator.
\newblock Technical Report 300, UC Berkeley, 2015.
\newblock http://biostats.bepress.com/ucbbiostat/paper343.

\bibitem[van~der Laan(2017)]{vanderLaan17}
M.J. van~der Laan.
\newblock A generally efficient targeted minimum loss estimator based on the
  highly adaptive lasso.
\newblock \emph{International Journal of Biostatistics}, 2017.

\bibitem[van~der Laan and Dudoit(2003)]{vanderLaan&Dudoit03}
M.J. van~der Laan and S.~Dudoit.
\newblock Unified cross-validation methodology for selection among estimators
  and a general cross-validated adaptive epsilon-net estimator: finite sample
  oracle inequalities and examples.
\newblock Technical Report 130, Division of Biostatistics, University of
  California, Berkeley, 2003.

\bibitem[van~der Laan and Luedtke(2015)]{vanderLaan:Luedtke15}
M.J. van~der Laan and A.~Luedtke.
\newblock Targeted learning of the mean outcome under an optimal dynamic
  treatment rule.
\newblock \emph{Journal of Causal Inference}, 3\penalty0 (1):\penalty0 61--95,
  2015.

\bibitem[van~der Laan and Robins(2003)]{vanderLaan&Robins03}
M.J. van~der Laan and J.M. Robins.
\newblock \emph{Unified methods for censored longitudinal data and causality}.
\newblock Springer, Berlin Heidelberg New York, 2003.

\bibitem[van~der Laan and Rose(2011)]{vanderLaan&Rose11}
M.J. van~der Laan and S.~Rose.
\newblock \emph{Targeted Learning: Causal Inference for Observational and
  Experimental Data}.
\newblock Springer, Berlin Heidelberg New York, 2011.

\bibitem[van~der Laan and Rose(2018)]{vanderLaan&Rose18}
M.J. van~der Laan and S.~Rose.
\newblock \emph{Targeted Learning in Data Science: Causal Inference for Complex
  Longitudinal Studies}.
\newblock Springer, Berlin Heidelberg New York, 2018.

\bibitem[van~der Laan and Rubin(2006)]{vanderLaan&Rubin06}
M.J. van~der Laan and Daniel~B. Rubin.
\newblock Targeted maximum likelihood learning.
\newblock \emph{Int J Biostat}, 2\penalty0 (1):\penalty0 Article 11, 2006.

\bibitem[van~der Laan et~al.(2006)van~der Laan, Dudoit, and van~der
  Vaart]{vanderLaan&Dudoit&vanderVaart06}
M.J. van~der Laan, S.~Dudoit, and A.W. van~der Vaart.
\newblock The cross-validated adaptive epsilon-net estimator.
\newblock \emph{Stat Decis}, 24\penalty0 (3):\penalty0 373--395, 2006.

\bibitem[van~der Laan et~al.(2021)van~der Laan, Wang, and van~der
  Laan]{vanderLaan&Wang21}
M.J. van~der Laan, Z.~Wang, and L.W.P. van~der Laan.
\newblock Higher order targeted maximum likelihood estimation, 2021.

\bibitem[van~der Vaart and Wellner(2011)]{vanderVaart:Wellner11}
A.W. van~der Vaart and J.A. Wellner.
\newblock A local maximal inequality under uniform entropy.
\newblock \emph{Electronic Journal of Statistics}, 5:\penalty0 192--203, 2011.
\newblock ISSN: 1935-7524, DOI: 10.1214/11-EJS605.

\bibitem[van~der Vaart et~al.(2006)van~der Vaart, Dudoit, and van~der
  Laan]{vanderVaart&Dudoit&vanderLaan06}
A.W. van~der Vaart, S.~Dudoit, and M.J. van~der Laan.
\newblock Oracle inequalities for multi-fold cross-validation.
\newblock \emph{Stat Decis}, 24\penalty0 (3):\penalty0 351--371, 2006.

\end{thebibliography}

\section{Notation Index}\label{section13}
\begin{description}
\item[\rm{cadlag function:} ]{multivariate real valued function that is right-continuous with left-hand limits}
\item[\rm CDF:]{cumulative distribution function}
\item[\rm HAL:]{Highly adaptive lasso}
\item[\rm MLE:]{Maximum likelihood estimation/estimator, Minimum loss estimation/estimator}
\item[$r_1(\epsilon)\sim r_2(\epsilon)$]{ short-notation for $r_1(\epsilon)=O(r_2(\epsilon))$ so that $r_2(\epsilon)$ represents an upper bound for the behavior of the expression $r_1(\epsilon)$ as $\epsilon\rightarrow 0$. Similarly, $r_1(n)\sim r_2(n)$ denotes $r_1(n)=O(r_2(n))$}
\item[$O$:]{Euclidean valued observed random variable}
\item[${\cal O}$:]{known support of $O$}
\item[\rm i.i.d:]{Independent and identically distributed}
\item[$P_0$:]{True data-generating distribution of $O$} 
\item[$O_1,\ldots,O_n\sim_{iid}P_0$:]{Sample of $n$ i.i.d. copies of random variable $O\sim P_0$}
\item[${\cal M}$:]{Statistical model, the set of possible probability distributions for $P_0$}
\item[$P_0 \in {\cal M}$:]{$P_0$ is known to be an element of the statistical model ${\cal M}$}
\item[$Q:{\cal M}\rightarrow {\cal Q}$:]{functional parameter of data distribution of interest}
\item[${\cal Q}$:]{parameter space for $Q$. So $Q_0=Q(P_0)\in {\cal Q}$. In our case, ${\cal Q}\subset D^{(k)}([0,1]^d)$ for some $k$}
\item[$L:{\cal Q}\times{\cal O}\rightarrow\openr$:]{loss function $(Q,o)\rightarrow L(Q,o)$ for $Q(P)$ so that $Q(P)=\arg\min_{Q\in {\cal Q}}PL(Q)$. We use interchangeably $L(Q)(o)$ and $L(Q,o)$ for the loss function at candidate $Q$ and observation $o$}
\item[$P$:]{Possible data-generating distribution in ${\cal M}$} 
\item[$P_n$:]{Empirical probability distribution; places probability $1/n$ on each observed $O_i\ldots, O_n\sim_{iid} P_0$}
\item[${\cal F}$:]{typically denotes a class of multivariate real valued functions of $O$}
\item[$\ell^{\infty}({\cal F})$:]{Banach space of real valued functions $G:{\cal F}\rightarrow\openr$ that map a function $f\in {\cal F}$ into a real number, endowed with the supremum norm $\pl G\pl_{\cal F}\equiv \sup_{Q\in {\cal F}}\mid G(Q)\mid$}
\item[$(G_n(f):f\in {\cal F})$:]{denotes empirical process $G_n=n^{1/2}(P_n-P)\in \ell^{\infty}({\cal F})$. So it uses notation $Pf\equiv \int f(o)dP(o)$ for a given measure $P$ and function $f$}
\item[${\cal M}_{np}$:]{Nonparametric statistical model that contains $P_n$ with probability 1}
\item[$\hat{Q}:{\cal M}_{np}\rightarrow {\cal Q}$:]{Estimator of $Q_0$ that maps $P_n$ into a function $\hat{Q}(P_n)$ in the parameter space}
\item[$Q_n=\hat{Q}(P_n)$:]{Realization of estimator when applied to $P_n$, typically still viewed as random element}
\item[$\mu$:]{ generally speaking, $\mu$ represents the Lebesgue measure on $(0,1]^d$, if not stated otherwise. For a subset $s\subset\{1,\ldots,d\}$, $\mu_s$ denotes the Lebesgue measure restricted to $(0(s),1(s)]=\{(x(j):j\in s): x\in (0,1]^d\}$. If we use notation $d\mu(x(s))$, it is meant to be $d\mu_s(x(s)=\prod_{j\in s}dx(j)$}
\item[$df/d\mu$:]{Radon-Nikodym derivative of measure implied by cadlag function $f$ w.r.t. Lebesgue measure $\mu$, so that $\int_A df=\int_A df/d\mu d\mu$ for all measurable subset $A\subset[0,1]^d$}
\item[$p=dP/d\mu$:]{Density of $P$ w.r.t. dominating measure $\mu$}
\item[$p_0$:]{True density of data-generating distribution $P_0$ w.r.t. appropriate dominating measure}
\item[$p$:]{Possible density in $\{dP/d\mu: P\in {\cal M}\}$ of data-generating distribution $P_0$ w.r.t. appropriate dominating measure}
\item[$L^2(P)$, $L^2(\mu)$:]{Hilbert space of square integrable functions with inner product $\langle h_1,h_2\rangle_P\equiv Ph_1h_2$, and similarly for $L^2(\mu)$}
\item[$L^2_0(P)$:]{Subspace of $L^2(P)$ of functions with mean zero w.r.t. $P$}
\item[$\phi_u^0$:]{zero order spline basis function $\phi_u^0(x)=I(x\geq u)$ with $u,x\in [0,1]^d$. Natural basis functions for modeling CDFs $\int_{(0,x]}dQ(u)=\int_u\phi_u^0(x)dQ(u)$. This is a tensor product of univariate zero-order spline basis functions}
\item[$\tilde{\phi}_u^0$:]{ $\tilde{\phi}_u^0(x)=I(x< u)$ with $\{u,x\}\subset [0,1]^d$. Natural basis functions for modeling survivor functions $\int_{(x,1]}dQ(u)=\int_u \tilde{\phi}_u^0(x)dQ(u)$}
\item[\mbox{$D^{(0)}([0,1 ]^d)$}:]{Space of real valued cadlag functions on $[0,1]^d$}
\item[$\pl Q\pl_v$]:{ variation norm of cadlag function $Q$ defined as $\pl Q\pl_v\equiv \int_{(0,1]^d} \mid dQ(x)\mid <\infty$. $dQ(u)$ can be interpret as the measure $Q$ assigns to infinitesimal cube $(u,u+du]$, which is defined as its generalized difference across the corners of $(u,u+du]$}
\item[$dQ$:]{short-hand notation for the measure on $(0,1]^d$ generated by  cadlag function $Q$ with finite variation norm $\pl Q\pl_v\equiv \int_{(0,1]^d} \mid dQ(x)\mid <\infty$. } 
\item[$x(s), x_s$:]{ for  a vector $x\in [0,1]^d$  and subset $s\subset\{1,\ldots,d\}$,  we define $x(s)=(x(j):  j\in s)$ as the subvector of $x$ with indices in subset $s\subset\{1,\ldots,d\}$. We alternate $x(s)$ and $x_s$ notation}
\item[$x(-s),x_{-s}$:]{ for a vector $x\in [0,1]^d$, we define $x(-s)=(x(j):j\not \in s)$ as the subvector of $x$ with indices in the complement of $s$}
\item[\mbox{$[0,1]^d=\cup_{s\subset\{1,\ldots,d\} }E_s$}:]{partitioning of the $d$-dimensional cube $[0,1]^d$ in zero-edges  $E_s=\{0(-s)\}\times (0(s),1(s)]$ across $2^d$-subsets $s$ with $E_{\emptyset}=\{0\}$}
\item[$Q_s:E_s\rightarrow\openr$:]{$s$-specific section of $Q:[0,1]^d\rightarrow\openr$ defined by $Q_s(0_{-s},x_{s})=Q(0_{-s},x_s)$ on $E_s$}
\item[$dQ_s$:]{short-hand notation for the measure generated by $\mid s\mid$-dimensional cadlag function $Q_s$ on $E_s$}
\item[$Q(x)=\int_{[0,1]^d}\phi_u^0(x)dQ(u)=\int_{[0,x]}dQ(u)$:]{representation of cadlag function w.r.t. its measure induced on $[0,1]^d$, when defining $Q(0_{-s}-,x_s)=0$ on the left of the $0$-edge $E_s$ across all subsets $s\subset\{1,\ldots,d\}$. Equivalently, the measure implied by $Q$ on each $E_s$ is the one generated by its section $Q_s(x_s)=Q(0(-s),x_s)$ restricted to $E_s$, so that $dQ=\sum_s I_{E_s} dQ_s$. This can be written as $Q(x)=Q(0)+\sum_{s\subset\{1,\ldots,d\}}\int_{(0(s),x(s)]}dQ_s(u)$, where $Q_s(u)=Q(u(s),0(-s))$ is the $s$-specific section of $Q$}
\item[$\pl Q\pl_v^*$:]{sectional variation norm $\pl Q\pl_v^*=\int_{[0,1]^d}\mid dQ(u)\mid$, which can be represented as $\pl Q\pl_v^*=\mid Q(0)\mid+\sum_{s\subset\{1,\ldots,d\}} \pl Q_s\pl_v$ with $\pl Q_s\pl_v=\int_{(0(s),1(s)]}\mid dQ_s(u)\mid$ the variation norm of section $Q_s$}
\item[\mbox{$D^{(0)}_M([0,1]^d)$}:]{set of cadlag functions with sectional variation norm bounded by $M$}
\item[${\cal R}^0(d)$:]{ complete set $\{(s,u): s\subset\{1,\ldots,d\},u\in(0(s),1(s)]\}$ indicating zero order spline basis functions $\phi^0_{s,u}(x)\equiv I(x(s)\geq u)$ indexed by subset $s\subset\{1,\ldots,d\}$ and knot-point $u\in (0(s),1(s)]$}
\item[\mbox{$D^{(0)}({\cal R}^{0,*}(d))$}:]{ for a given subset ${\cal R}^{0,*}(d)\subset {\cal R}^0(d)$,  this is defined as the space of functions that are in closure of linear span of $\{\phi_{s,u}: (s,u)\in {\cal R}^{0,*}(d)\}$. We have $D^{(0)}({\cal R}^0(d))=D^{(0)}([0,1]^d)$, by the representation of the class of cadlag functions as linear combinations of zero order splines}
\item[$N(\epsilon,{\cal F},d)$:]{Number of spheres of size $\epsilon$ w.r.t. metric $d$ that are needed to cover  a set  of functions ${\cal F}$}
\item[$J(\delta,{\cal F},\pl\cdot\pl)$:]{ this is the entropy integral for class ${\cal F}$ w.r.t. metric $d$ defined as $\int_0^{\delta}\sqrt{\log N(\epsilon,{\cal F},\pl\cdot\pl)} d\epsilon$. $J_{\infty}(\delta,{\cal F})$ is used when the norm is the supremum norm}
\item[$r(d,J,M)$:]{optimally $r(d,J,M)=\arg\min_{\{u_1,\ldots,u_J\}}\sup_{Q\in D^{(0)}_M([0,1]^d)}\inf_{\beta}\pl \sum_{j=1}^J\beta(j)\phi_{u_j}^0-Q\pl_{L^2(\mu)}$, with $\mu$ Lebesgue measure.  We establish the  upper bound $r(d,J,M)\sim MJ^{-1}(\log J)^{2(d-1)}$ based on the known upper bound on the  known covering number $N(\epsilon,D^{(0)}_1([0,1]^d),L^2)$. All our uniform approximation results apply with the optimal $r(d,J,M)$ definition above}
\item[$C(M)$, $r(d,J)$:]{ $r(d,J,M)=O( C(M)r(d,J))$, where $r(d,J)=1/J(\log J)^{2(d-1)}$ and $C(M)=M$. If the literature obtains improvements, then one can replace $C(M)$ and $r(d,J)$ by these, but it will at most change the power of $\log J$}
\item[\mbox{${\cal F}^{(0)}((0,1]^d)$}:]{ set of  differences of two cumulative distribution functions $Q$ on $(0,1]^d$ with finite variation norm $\pl Q\pl_v<\infty$. A cumulative distribution function $Q$ is a function with $Q(du)\geq 0$ for all $u\in (0,1]^d$ with $Q(1)<\infty$}
\item[\mbox{${\cal F}^{(0)}_M((0,1]^d)$}:]{ difference of two cumulative distribution functions with variation norm bounded by $M$}
\item[\mbox{${\cal S}^{(0)}_M((0,1]^d)$}:]{ set of differences of survivor functions $Q$ with variation norm $\pl Q\pl_v<M$. A survivor function is a function $Q$  with $Q(du)\leq 0$ for all $u$ and $Q(0)<\infty$}
\item[\mbox{${\cal F}^{(0)}_M({\cal R})$}:]{ for subset ${\cal R}\subset (0,1]^d$ this is the linear space of functions $\{x\rightarrow \int_{u\in {\cal R}}\phi_u^0(x)d\tilde{Q}(u):\pl \tilde{Q}\pl_v<M\}\subset {\cal F}^{(0)}_M((0,1]^d)$}
\item[${\cal R}(d,J)$:]{any set of $J$ knot-points in $(0,1]^d$ chosen so that $D^{(0)}({\cal R}(d,J))$ approximates  ${\cal F}^{(0)}_M((0,1]^d)$ w.r.t. $L^2(\mu)$-norm with rate $O(r(d,J,M))$, $\mu$ Lebesgue measure}
\item[${\cal R}(s,J)$, $s\subset\{1,\ldots,d\}$:]{same as ${\cal R}(\mid s\mid,J)$, but restricted to $[0(s),1(s)]$}
\item[${\cal R}^0(d,{\bf J})$:]{ Finite set of indices $\{(s,u): s\subset\{1,\ldots,d\}, u\in {\cal R}(s,{\bf J}(s))\}$ for zero-order spline basis functions $\phi_{s,u}(x)=I(x(s)\geq u)$, where for each subset $s$, the set ${\cal R}(s,{\bf J}(s))$ contains ${\bf J}(s)$ knot-points. It is generally assumed that each number ${\bf J}(s)$ of knot-points is proportional to the average $J=1/(\sum_{s\subset\{1,\ldots,d\}}) \sum_{s\subset\{1,\ldots,d\}}{\bf J}(s)$. The latter allows us to obtain rates of approximation for a function $Q$ in terms of $J$}
\item[\mbox{$D^{(0)}({\cal R}^0(d,{\bf J}))$}:]{space of linear combinations of $\{\phi_j:j\in {\cal R}^0(d,{\bf J})\}$.  It approximates $D^{(0)}_M([0,1]^d)$ w.r.t. $L^2(\mu)$-norm at rate $O(C(M)r(d,J))$}
\item[\mbox{$D^{(0)}_M({\cal R}^{0,*}(d))$}:]{subset of functions $Q\in D^{(0)}({\cal R}^{0,*}(d))$ for which sectional variation norm $\pl Q\pl_v^*=\int_{{\cal R}^{0,*}(d)}\phi_u^0 \mid dQ(u)\mid $ is bounded by $M$}
\item[$\phi^k_u$:]{$k$-th order spline basis functions identified by knot-point $u\in [0,1]^d$. It equals a tensor product over all $j\in \{1,\ldots,d\}$ with $u(j)>0$ of univariate $k$-th order spline basis function at knot-point $u(j)$}
\item[$\mu^k(f)$:]{ recursively defined by $\mu(f)(x)=\mu^1(f)\equiv \int_{(0,x]}f(y)d\mu(y)$; $\mu^{j+1}(f)(x)=\int_{(0,x]}\mu^j(f)(y)d\mu(y)$, $j=1,2,\ldots$. If $f$ only depends on $x$ through $x(s)$, then $\mu^k(f)=\mu^k_s(f)$. We have $\phi_u^k=\mu^k(\phi_u^0)$}
\item[\mbox{${\cal F}^{(k)}_M((0,1]^d)$}:]{ linear space $\{\mu^k(f): f\in {\cal F}^{(0)}_M((0,1]^d)\}$, which equals $\{x\rightarrow\int \phi_u^k(x)df(u): f\in {\cal F}^{(0)}_M((0,1]^d)\}$. Similarly, we define ${\cal F}^{(k)}_M({\cal R})=\{x\rightarrow\int_{\cal R}\phi_u^k(x)df(u): f\in {\cal F}^{(0)}_M({\cal R})\}$ for a subset ${\cal R}\subset(0,1]^d$}
\item[$\bar{s}(k+1)$:]{ a vector $(s_1,\ldots,s_{k+1})$ of nested subsets $s_1\supset \ldots\supset s_{k+1}$ of $\{1,\ldots,d\}$}
\item[${\cal S}^{k+1}(d)$:]{set of all possible $\bar{s}(k+1)=(s_1,\ldots,s_{k+1})$}
\item[${\cal S}^{k+1,*}(d)$:]{subset of ${\cal S}^{k+1}(d)$}
\item[$m(\bar{s}(k+1)$:]{the smallest integer $m$ in $\{1,\ldots,k+1\}$ for which $s_{m+1}$ is empty set}
\item[$Q^{(k)}_{\bar{s}(k+1)}$:]{ sequentially defined $k$-th order Radon-Nikodym derivative w.r.t. Lebesgue measure, by first defining $Q^{(1)}_{s_1}=dQ_{s_1}/d\mu_{s_1}$; constructing its section $Q^{(1)}_{s_1,s_2}$, then taking its derivative $Q^{(2)}_{s_1,s_2}=\frac{dQ^{(1)}_{s_1,s_2}}{d\mu_{s_2}}$, constructing its section $Q^{(2)}_{s_1,s_2,s_3}$, etc, till at step $k$ we have $Q^{(k)}_{\bar{s}(k+1)}$. The moment the $j+1$-th defined section is the evaluation at $0$ due to $s_{j+1}$ being the empty set, then $Q^{(j)}_{\bar{s}(j+1)}=Q^{(j)}_{\bar{s}(j)}(0(s_j))$, in which case $Q^{(k)}_{\bar{s}(k+1)}$ is the constant $Q^{(m)}_{\bar{s}(m)}(0(s_m))$ with $m=m(\bar{s}(k+1))$}
\item[$\tilde{Q}^{(k)}_{\bar{s}(k+1)}(x(s_{k+1}))$:]{ this is defined as $\tilde{Q}^{(k)}_{\bar{s}(k+1)}(x(s_{k+1}))=\int_{(0(s_{k+1}),x(s_{k+1})]}  Q^{(k)}_{\bar{s}(k)}(du(s_{k+1}),0(s_k/s_{k+1}) )$}
\item[$\pl Q\pl_{v,k}^*$:]{$k$-th order sectional variation norm defined as $\sum_{\bar{s}(k+1)}\pl Q_{\bar{s}(k+1)}^{(k)}\pl_v$, sum of variation norms across all $\bar{s}(k+1)$}
\item[\mbox{$D^{(k)}_M([0,1]^d)$}:]{functions in $D^{(k)}([0,1]^d)$ with $k$-th order sectional variation norm bounded by $M$}
\item[$\bar{\phi}_{\bar{s}(k+1)}$:]{ defined as $\prod_{j=1,\mid s_j\mid>0}\phi_0^j(x(s_j/s_{j+1}))$. Equivalently, this equals $\prod_{j=1}^{m(\bar{s}(k+1))-1}\phi_0^j(x(s_j/s_{j+1}))$. Note $\phi_0^j(x(s_j/s_{j+1}))$ denotes the $j$-th order spline on $(0(s_j/s_{j+1}),1(s_j/s_{j+1})]$ with knot-point $u=0$}
\item[${\cal R}^k(d)$:]{ defined as union of ${\cal R}^k_1(d)\equiv \{(\bar{s}(k+1),u): \bar{s}(k+1)\in {\cal S}^{k+1}(d),\mid s_{k+1}>0,u\in (0(s_{k+1}),1(s_{k+1})]\}$
and ${\cal R}^k_2(d)\equiv \{\bar{s}(k+1): \mid s_{k+1}\mid =0\}$}
\item[${\cal R}^{k,*}(d)$:]{subset of ${\cal R}^{k}(d)$ obtained by restricting in the definition of ${\cal R}^k_1(d)$  the set ${\cal S}^{k+1}(d)$  to a subset ${\cal S}^{k+1,*}(d)$, and, for each $\bar{s}(k+1)\in {\cal S}^{k+1,*}(d)$ with $\mid s_{k+1}\mid >0$, restricting knot-point set $(0(s_{k+1}),1(s_{k+1})]$ to a subset of this cube. This results in a subset ${\cal R}^{k,*}_1(d)\subset {\cal R}^k_1(d)$. In addition, it might also restrict ${\cal R}^k_2(d)$ to ${\cal R}^{k,*}_2(d)=\{\bar{s}(k+1)\subset{\cal S}^{k+1,*}(d):\mid s_{k+1}\mid =0\}$, but since this represents only a finite dimensional parametric form one might always include ${\cal R}^k_2(d)$}
\item[$\phi_{\bar{s}(k+1),u}$:]{ for $u\in (0(s_{k+1}),1(s_{k+1})]$, we have
$\phi_{\bar{s}(k+1),u}(x)=I(\mid s_{k+1}\mid >0)\bar{\phi}_{\bar{s}(k+1)}\phi_u^k+I(\mid s_{k+1}\mid =0)\bar{\phi}_{\bar{s}(k+1)}$. So if $s_{k+1}$ is empty, then $\phi_{\bar{s}(k+1),u}$ does not depend on $u$ and equals $\bar{\phi}_{\bar{s}(k+1)}$}
\item[\mbox{$D^{(k)}({\cal R}^{k,*}(d))$}:]{closure of linear combinations of $\{\phi_{\bar{s}(k+1),u}:(\bar{s}(k+1),u)\in {\cal R}^{k,*}(d)\}$. Our $k$-th order spline representation theorem states that $D^{(k)}_M([0,1]^d)=D^{(k)}_M({\cal R}^{k}(d))$, and   by choosing subsets ${\cal R}^{k,*}(d)\subset{\cal R}^k(d)$, we obtain subspaces of $D^{(k)}([0,1]^d)$}
\item[${\cal R}^k(d,{\bf J})$:]{ finite subset of ${\cal R}^k(d)$. It is defined as union of \[
\{(\bar{s}(k+1),u): \bar{s}(k+1)\in {\cal S}^{k+1}(d),\mid s_{k+1}\mid >0,u\in {\cal R}(s_{k+1},{\bf J}(\bar{s}(k+1))\}\]
and $\{\bar{s}(k+1)\in {\cal S}^{k+1}(d):\mid s_{k+1}\mid =0\}$.
So this index set identifies ${\bf J}(\bar{s}(k+1))$ $k$-th order spline basis functions $\phi_u^k(x(s_{k+1}))$ when $\mid s_{k+1}\mid >0$, and it identifies the finite set of $\bar{\phi}_{\bar{s}(k+1)}$ for each $\bar{s}(k+1)$ with $\mid s_{k+1}\mid =0$. Generally speaking ${\bf J}(\bar{s}(k+1))$ is proportional to $J$ across all $\bar{s}(k+1)$}
\item[${\cal R}^{k,*}(d,{\bf J})$:]{ subset of ${\cal R}^k(d,{\bf J})$ by intersecting it with ${\cal R}^{k,*}(d)$ so that it now yields the desired uniform approximation of $D^{(k)}({\cal R}^{k,*}(d))$}
\item[\mbox{$D^{(k)}({\cal R}^{k}(d,{\bf J}))$}:]{ linear combinations of $\{\phi_j: j\in {\cal R}^k(d,{\bf J})\}$. It yields a uniform approximation of any function $Q\in D^{(k)}_M([0,1]^d)$ at rate $
O(C(M)r(d,J)^{k+1})$, uniformly in all such $Q$}
\item[$Q_{{\bf J},\beta}$, $Q_{J,\beta}$:]{ $Q_{{\bf J},\beta}=\sum_{j\in {\cal R}^k(d,{\bf J})}\beta(j)\phi_j$, and $Q_{J,\beta}$ is short-hand notation for this assuming all non-zero components of ${\bf J}$ are proportional to $J$. This notation could also be applied to ${\cal R}^{k,*}(d,{\bf J})$}
\item[\mbox{$D^{(k)}({\cal R}^{k,*}(d,{\bf J}))$}:]{  linear combinations of $\{\phi_j: j\in {\cal R}^{k,*}(d,{\bf J})\}$. By our uniform approximation theorem, it approximates $ D^{(k)}_M({\cal R}^{k,*}(d))$ w.r.t. supremum norm at rate $
O(C(M)r(d,J)^{k+1})$}
\item[${\cal R}_n$:]{ data adaptive set indices $(\bar{s}(k+1),u)$ (or just $\bar{s}(k+1)$) identifying the basis functions of type $\phi_{\bar{s}(k+1),u}$  with non-zero coefficients in the $k$-th order spline sieve or HAL-MLEs $Q_n$ so that $Q_n=\sum_{j\in {\cal R}_n}\beta_n(j)\phi_j$. Sometimes also denoted with ${\cal R}_n^k$, if helpful. Note that ${\cal R}_n={\cal R}(P_n)$ identifies the whole estimator $\hat{Q}(P_n)$}
\item[${\cal R}_{0,n}$:]{ generally a fixed set defined as ${\cal R}(P_n^{\#})$ for an independent empirical measure  $P_n^{\#}$ of sample of $n$ i.i.d. copies of $O^{\#}\sim P_0$ independent of $P_n$, viewed as an independent approximation of ${\cal R}_n={\cal R}(P_n)$. Sometimes also denoted with ${\cal R}^k_{0,n}$}
\item[$D^{(k)}({\cal R}_n)$:]{ data adaptively determined working model consisting of linear combinations of $\{\phi_j: j\in {\cal R}_n\}$, while $D^{(k)}_M({\cal R}_n)$ enforces $L_1$-norm of coefficients to be bounded by $M$}
\item[$D^{(k)}({\cal R}_{0,n})$:]{ fixed working model  $\{\sum_{j\in {\cal R}_{0,n}}\beta(j)\phi_j:\beta\}$ as defined above, while $D^{(k)}_M({\cal R}_{0,n})$ enforces $\pl \beta\pl_1\leq M$}
\item[$Q_{{\cal R}^k,0}$:]{ oracle MLE over working model $D^{(k)}({\cal R}^k)$ at true data distribution $P_0$ defined as 
$Q_{{\cal R}^k,0}=\arg\min_{Q\in D^{(k)}({\cal R}^k)}P_0L(Q)$}
\item[$Q_{{\cal R}^k,P}$:]{ oracle MLE over working model $D^{(k)}({\cal R}^k)$ at $P$ defined as 
$Q_{{\cal R}^k,P}=\arg\min_{Q\in D^{(k)}({\cal R}^k)}PL(Q)$}
\item[$Q_{0,n}$:]{ $=\arg\min_{Q\in D^{(k)}({\cal R}_{0,n})}P_0L(Q)$, short-hand notation for $Q_{{\cal R}_{0,n},0}=Q_{{\cal R}_{0,n},P_0}$}
\item[$S_Q(\phi)$:]{ score $S_Q(\phi)=\frac{d}{dQ}L(Q)(\phi)$ for coefficient in front of basis function $\phi$}
\item[$S_j(Q)$:]{ short notation for $S_Q(\phi_j)=\frac{d}{dQ}L(Q)(\phi_j)$, $j\in {\cal R}_{0,n}$ or $j\in {\cal R}_n$}
\item[$\langle h_1,h_2\rangle_{J_{0,n}}$:]{ inner product chosen so that minus second derivative of expectation of loss (e.g, log-likelihood loss) is (approximately with remainder $R_{3n}(\phi)$ or exactly) represented by this inner product: 
\[
-P_0 \frac{d}{dQ_{0,n}}S_{Q_{0,n}}(\phi)(Q_n-Q_{0,n})=\langle \phi,Q_n-Q_{0,n}\rangle_{J_{0,n}}+R_{3,n}(\phi)\] }
\item[$\phi_j^*$, $j\in {\cal R}_{0,n}$:]{ orthonormal basis for linear span $D^{(k)}({\cal R}_{0,n})$ of $\{\phi_j:j\in {\cal R}_{0,n}\}$ w.r.t. inner product $\langle h_1,h_2\rangle_{J_{0,n}}$ on $D^{(k)}({\cal R}_{0,n}))$. We still index these by $j\in {\cal R}_{0,n}$, even though each one is a linear combination of $\{\phi_j:j\in {\cal R}_{0,n}\}$ (but in Gramm-Schmidt we could use an ordering and then  indeed create sequentially a $\phi_j^*$ for each $j\in {\cal R}_{0,n}$}
\item[$\Pi_{J_{0,n}}(Q)$:]{ projection of $Q$ onto $D^{(k)}({\cal R}_{0,n})$ w.r.t. inner product $\langle\cdot,\cdot\rangle_{J_{0,n}}$. So $\Pi_{J_{0,n}}(Q)=\sum_{j\in {\cal R}_{0,n}}\langle f,\phi_j^*\rangle_{J_{0,n}}\phi_j^*$}
\item[$S_j^*(Q)$:]{ short notation for $\frac{d}{dQ}L(Q)(\phi_j^*)$}
\item[$r_n(j)$:]{ empirical mean of score $P_n S_{Q_n}(\phi_j)$ for non-orthogonalized working model $D^{(k)}({\cal R}_{0,n})$}
\item[$r_n^*(j)$:]{ empirical mean of score  $P_n S_{Q_n}(\phi_j^*)$ for orthogonalized working model}
\item[$\tilde{r}_n(x)$:]{ $=\sum_{j\in {\cal R}_{0,n}}r_n^*(j)\phi_j^*(x)$, which equals the
empirical mean of efficient influence curve at estimator $Q_n$ for parameter $P\rightarrow Q_{{\cal R}_{0,n},P}(x)$ on nonparametric model ${\cal M}$}
\item[$D_{Q_{0,n},x}$:]{ $=\sum_{j\in {\cal R}_{0,n}}S_{Q_{0,n}}(\phi_j^*)\phi_j^*(x)$. It represents the influence curve of $(Q_n-Q_{0,n})(x)\approx P_n D_{Q_{0,n},x}$, and equals the efficient influence curve of $P\rightarrow Q_{{\cal R}_{0,n},P}(x)$ on nonparametric model ${\cal M}$}
\item[$d_n$:]{ number of elements in ${\cal R}_{n}$, i.e.,  dimension of model $D^{(k)}({\cal R}_n)$}
\item[$d_{0,n}$:]{ dimension of model $D^{(k)}({\cal R}_{0,n})$}
\item[$\tilde{\sigma}^2_{0,n}(x)$:]{ variance under $P_0$ of $D_{Q_{0,n},x}$, representing approximate  asymptotic variance of $(n/d_{0,n})^{1/2}(Q_n-Q_{0,n})(x)$}
\item[$\tilde{\sigma}^2_n(x)$:]{ estimator of $\tilde{\sigma}^2_{0,n}(x)$}
\item[$k^*=k+1$:]{ useful so that rate of convergence has form $n^{-k^*/(2k^*+1)} $ up till $\log n$-factors}
\item[$\tilde{Q}_n$:]{ $\tilde{Q}_n=\Pi_{J_{0,n}}(Q_n\mid D^{(k)}({\cal R}_{0,n}))$ is projection of $Q_n$ onto the fixed working model $D^{(k)}({\cal R}_{0,n})$ satisfying that $\pl \tilde{Q}_n-Q_n\pl_{\infty}=O(C(M_n)r(d,J_{0.n})^{k+1})$ with $M_n$ the $L_1$-norm of coefficient vector of $Q_n$}
\item[$\tilde{P}_n$:]{ a $P$ in model ${\cal M}$ that is compatible with $Q_n=Q(\tilde{P}_n)$}
\item[$\tilde{S}_j(Q_n)$:]{  a score of form $\tilde{S}_j(Q_n)=\frac{d}{dQ_n}L(Q_n)(\tilde{\phi}_j)$ solved by $Q_n$ so that $P_n \tilde{S}_j(Q_n)=0$ and that approximates $S_j(Q_n)$, $j\in {\cal R}_{0,n}$}
\item[$R_n(x)=\bar{R}_n(x)-\tilde{r}_n(x)-E_n(x)$:]{ total exact remainder so that $(Q_n-Q_{0,n})(x)=P_n D_{Q_{0,n},x}+R_n(x)$}
\item[$\bar{R}_n(x)$:]{ sum  of second order remainders, $\sum_{m=1}^3\sum_{j\in {\cal R}_{0,n}} R_{m,n}(\phi_j^*)\phi_j^*(x)$, with $R_{m,n}(\phi)$ linear in $\phi$ and, typically,  $\sum_jR_{m,n}(\phi_j^*)\phi_j^*$ represents a projection of a function (second order difference) onto $D^{(k)}({\cal R}_{0,n})$}
\item[$E_n(x)$:]{ empirical process second order remainder term $(P_n-P_0)(D_{Q_n,x}-D_{Q_{0,n},x})$}
\item[$Q_{{\bf J},P}$:]{ short-hand for oracle MLE $Q_{{\cal R}^k(d,{\bf J}),P}$ defined above for the case that the finite set of basis functions is given by  ${\cal R}^k(d,{\bf J})$}
\item[${\bf J}_{0,n}$:]{ fixed selector of ${\bf J}$ in ${\cal R}^k(d,{\bf J})$, such as the actual data adaptive selector ${\bf J}_n$ but applied to sample from $P_0$ independent of $P_n$}
\item[\rm $k$-th order spline HAL-MLE]{ generally defined as $\hat{Q}_C(P_n)=\arg\min_{Q\in D^{(k)}_C({\cal R}(d,{\bf J}_{max,n}) ) }P_n L(Q)$ indexed by tuning parameter the $L_1$-norm of the vector of coefficients, or equivalently, the $k$-th order sectional variation norm, enforced to be smaller or equal than $C$, for a given (large enough) initial working model. The initial working model could be chosen so that the unpenalized MLE for this working model already yields the optimal rate of convergence  $n^{-k*/(2k^*+1)}$ up till $\log n$-factors, while yielding negligible bias. $C$ is the tuning parameter that will be data adaptively selected (e.g., cross-validation selector(}
\item[{\rm $k$-th order spline HAL-MLE, extra tuned}:]{ defined as $\hat{Q}_{C,{\bf J}_{max}}(P_n)$ identical to our definition of $\hat{Q}_C(P_n)$ above but now viewing the size ${\bf J}_{max,n}$ of the initial working model $D^{(k)}({\cal R}(d,{\bf J}_{max,n}))$ as
a tuning parameter.  The joint tuning parameter is data adaptively selected (e.g., cross-validation selector)}
\item[{\rm $k$-th order spline sieve MLE}:]{ generally defined as $\hat{Q}_{{\bf J}}(P_n)=\arg\min_{Q\in D^{(k)}({\cal R}(d,{\bf J})}P_n L(Q)$, where ${\bf J}$ is the vector-valued tuning parameter, chosen with a data adaptive selector ${\bf J}_n$}
\item[$S_Q(\phi)$]{ $S_Q(\phi)=\frac{d}{dQ}L(Q)(\phi)$ is directional derivative of loss at $Q$ in direction $\phi$. When $Q=Q_{J,\beta}=\sum_{j\in {\cal R}}\beta(j)\phi_j$, then $S_Q(\phi_j)$ is the score of $\beta(j)$, $j\in {\cal R}$.
We might also index it by $k$, $S_Q^k(\phi)$, if not clear from context that we work with working model $D^{(k)}({\cal R}^k(d,{\bf J}))$ viewed as approximations of $D^{(k)}([0,1]^d)$ or $D^{(k)}({\cal R}^{k,*}(d))$ according to our $k$-th order spline representation theorem and uniform approximation theorem}

\item[{\rm $k$-th order spline relax HAL-MLE}:]{ denoted with $\hat{Q}_{C,r}(P_n)$ defined by   first computing  the HAL-MLE $\hat{Q}_C(P_n)$, then determining the subset ${\cal R}_n$ of basis functions with non-zero coefficient in its fit, and then refit it without $L_1$-penalization so that $\hat{Q}_{C,r}(P_n)=\arg\min_{Q\in D^{(k)}({\cal R}_n)}P_n L(Q)$ equals the regular (unpenalized) MLE over the data adaptive working model $D^{(k)}({\cal R}_n)$}
\item[$T_P({\cal M})$:]{Tangent space at $P$ in model ${\cal M}$, defined w.r.t. a class of paths $\{P^h_{\epsilon}:\epsilon\in (-\delta,\delta)\}\subset{\cal M}$  through $P\in {\cal M}$ at $\epsilon =0$ with score $S_h$. Tangent space is defined as closure of the linear span of all scores $S_h$ across class of all paths. We have $T_P({\cal M})\subset L^2_0(P)$ subspace of Hilbert space $L^2_0(P)$ with inner product $\langle h_1,h_2\rangle_P=Ph_1h_2$}
\item[$d\Psi_P:T_P({\cal M})\rightarrow\openr$:]{pathwise derivative of $\Psi:{\cal M}\rightarrow\openr$ along a specified class of paths through $P$ defined by $d\Psi_P(h)=\frac{d}{d\epsilon_0}\Psi(P_{\epsilon_0,h})$. It will be assumed that this pathwise derivative is bounded so that it can be represented by the inner product in $T_P({\cal M})$: $d\Psi_P(h)=\langle D^*_P,h\rangle_{P}$ for a $D^*_P\in T_P({\cal M})$}
\item[$D^*_P$:]{canonical gradient of $\Psi$, element of tangent space $T_P({\cal M})$, so that $d\Psi_P(h)=\langle D^*_P,h\rangle_{T_P({\cal M}}$. If we want to emphasize it is the canonical gradient of a particular function $\Psi$, then we also use $D^*_{\Psi,P}$}
\item[$D_P$:]{ gradient at $P$, i.e., any function that equals $D^*_P+S$ with $S$ an element of orthogonal complement of tangent space $T_P({\cal M})\subset L^2_0(P)$ (i.e., $S\perp T_P({\cal M})$)}
\end{description}

\appendix 
\section*{Appendix}\label{Appendix}
{\bf Appendix \ref{AppendixA}} establishes the key corollary \ref{defknots} for the existence of sets $({\cal R}(m,J):m,J)$ satisfying the $O(r(m,J))$ $L^2$-approximation error of $m$-variate cumulative distribution and survivor functions with $r(m,J)=J^{-1}\log^{2(d-1)}J$. 
 {\bf Appendix \ref{AppendixB}} provides a practical method  that uses  zero order HAL-MLE to fit a uniform CDF in order to   determine knot-point sets $({\cal R}(m,J):J)$ satisfying the desired $L^2$-approximation error $O(r(m,J))$ for approximating $m$-dimensional cumulative distribution functions, at least up till a $\log J$-factor.
  {\bf Appendix \ref{AppendixC}} presents the proof of the $k$-th order spline representation of $D^{(k)}([0,1]^d)$. {\bf Appendix \ref{AppendixD}} presents the general proof that linear combinations of $J$ $k$-th order splines uniformly approximates a $d$-dimensional $k$-th order primitive as $O(r(d,J)^{k+1})$, firstly focusses on making assuming that ${\cal R}(m,J)$ provides the desired $O(r(m,J))$ approximation error for both cumulative and survivor functions and subsequently showing that the latter requirement can be dropped. 
In {\bf Appendix \ref{AppendixE}} we discuss a variety of submodels of $D^{(k)}([0,1]^d)$ of interest  by restricting the space of basis functions. 
In {\bf Appendix \ref{AppendixF}} we generalize our uniform approximation results for finite dimensional $k$-th order spline working models  for the $k$-th order smoothness class to subsets of our $k$-th order smoothness class, thereby making our $k$-th order spline HAL-MLE and sieve MLE results applicable to a large variety of submodels that both restrict smoothness as well as the set of potential $k$-th order spline basis functions. 
In {\bf Appendix \ref{AppendixG}} we show that our uniform approximation results imply corresponding sup-norm covering numbers for the class of $k$-th order primitives and, as a consequence, for our $k$-th order smoothness class $D^{(k)}([0,1]^d)$.  
In {\bf Appendix \ref{AppendixI}} we provide our proof of Theorem \ref{theoremoraclemle}  showing that the oracle MLE $Q_{0,{\bf J}}$ implied by the $k$-th order spline working model $D^{(k)}({\cal R}(d,{\bf J}))$ achieves the same sup-norm approximation $O(r(d,J)^{k+1})$ of $Q_0$ as the best possible uniform approximation, thereby making our uniform approximation error results immediately apply to the oracle MLE/projection $Q_{0,{\bf J} }$. In {\bf Appendix \ref{AppendixJ}} we prove at which rate (and how to achieve this rate) an undersmoothed HAL-MLE solves the regular score equations (i.e., derivatives along paths not restricting the $L_1$-norm) for its non-zero coefficients, the ones that are solved exactly by the sieve MLE and relax HAL-MLE.  In {\bf Appendix \ref{AppendixK}} we present  the proof for asymptotic normality and uniform convergence of the $k$-th order spline sieve and HAL-MLEs  for the regression function when using {\em unweighted} least squares, while our results in \ref{section8} focussed on weighted least squares regression.  In {\bf Appendix \ref{AppendixK1}} we establish asymptotic normality of arbitrary functions of $Q_n$ (e.g., derivatives) building on our second order expansion $(Q_n-Q_{0,n})(x)\approx P_n D_{Q_{0,n},x}$ established in Section \ref{section10}.
In {\bf Appendix \ref{AppendixL}} we establish the approximation error of $D^{(k)}_M([0,1]^d)$ with our working models $D^{(k)}({\cal R}(d,{\bf J}) ) $ w.r.t. the sectional variation norm,  a much stronger norm than the supremum norm. It shows that this approximation error in sectional variation norm of the oracle  MLE  $Q_{0,{\bf J}}^k$ over  $D^{(k)}({\cal R}(d,{\bf J}))$ has the same rate $r(d,J)^k$ as the sup-norm rate of the approximation error for  the  oracle  MLE $Q_{0,{\bf J}}^{k-1}$ over $D^{(k-1)}_M([0,1]^d)$. 
In addition, we establish that the sectional variation norm of  $Q_{n,{\bf J}}^k-Q_{0,{\bf J}}^k$, $Q_{n,{\bf J}}^k$ being MLE over $D^{(k)}({\cal R}(d,{\bf J}))$,  converges at rate $(n/J)^{-1/2}$, showing that the sectional variation norm $\pl Q_{n,{\bf J}}^k-Q_0\pl_v^*$ converges at rate $(n/J)^{-1/2}+r(d,J)^k$, which for optimally tuned ${\bf J}$ gives the rate of convergence $n^{-k/(2k+1)}$ up till $\log n$-factors. 
These results also  generalize to the HAL-MLEs. For completeness, in {\bf Appendix \ref{AppendixM}} we also study the sectional variation norm behavior for the sieve zero-order spline-MLE, showing under what conditions it remains uniformly bounded. In {\bf Appendix \ref{AppendixN}} we provide special results for the regression case when the covariates are discrete valued.

 \section{Linear combinations of $J$ zero order splines that approximate differences of cumulative distribution and survivor functions with $L^2$-norm error $O(r(d,J))$.}
 \label{AppendixA}
 The following definitions and results for the zero-order spline approximations provide a basis for our approximation results for linear combinations of $\leq k$-th order splines of a $k$-th order smooth function in $D^{(k)}([0,1]^d)$.  It suffices to focus on approximating cumulative distribution functions and cumulative survivor functions, even though these results imply the same approximation errors for the general cadlag functions of bounded sectional variation, given that the latter are linear combinations of differences of cumulative distribution functions.  
 \subsection{Defining set of differences of cumulative distribution and survivor functions and $J$-dimensional zero order spline approximations}
 Firstly, we define the set of differences of cumulative distribution functions and survivor functions.
\begin{definition} {\bf (Differences of Cumulative Distribution Functions)}
Let $\mu$ denote the Lebesgue measure on $[0,1]^d$.
Let ${\cal F}^{(0)}((0,1]^d)=\{x\rightarrow\int_{(0,x]}dQ(u):\pl Q\pl_v<\infty\}$, which equals the set of differences of cumulative ''distribution'' functions  with finite  sup-norm (i.e its value at $1$ is finite).
Let ${\cal F}^{(0)}_M((0,1]^d)=\{x\rightarrow\int_{(0,x]}dQ(u):\pl Q\pl_v\leq M\}$.
We also define ${\cal S}^{(0)}((0,1]^d)=\{x\rightarrow \int_{(x,1]}dQ(u):\pl Q\pl_v<\infty\}$, which equals the  set of differences of ''survivor'' functions with finite sup-norm. Let ${\cal S}^{(0)}_M((0,1]^d)=\{x\rightarrow\int_{(x,1]}dQ(u):\pl Q\pl_v\leq M\}$.
\end{definition}

\begin{definition}{\bf (Separate knot-point sets for cumulative and survivor functions)} 
We note that $\int_{(0,x]}dQ(u)=\int\phi_u^0(x)dQ(u)$ is a linear combination of zero-order splines $x\rightarrow \phi_u^0(x)$ indexed by knot-points $u\in (0,1]^d$. Similarly,  $\int_{(x,1]}dQ(u)=\int\tilde{\phi}_u^0(x)dQ(u)$ is a  linear combination of 
$x\rightarrow\tilde{\phi}_u^0(x)\equiv I(x< u)$ indexed by knot-points $u\in (0,1]^d$. Thus, ${\cal F}^{(0)}((0,1]^d)$ can be represented as the space of linear combinations of $\phi_u^0$, $u\in (0,1]^d$, while ${\cal S}^{(0)}((0,1]^d)$ represents the space of linear combinations of $\tilde{\phi}_u^0$, $u\in (0,1]^d$. Therefore, we will denote knot-point sets for approximating ${\cal F}^{(0)}((0,1]^d)$ and ${\cal S}^{(0)}((0,1]^d)$ with ${\cal R}_{\phi^0}(d,J)$ and ${\cal R}_{\tilde{\phi}^0}(d,J)$, respectively, and their corresponding finite dimensional spaces with ${\cal F}^{(0)}({\cal R}_{\phi^0}(d,J))=\{\sum_{j\in {\cal R}_{\phi^0}(d,J)}\beta(j)\phi_j^0:\beta\}\subset {\cal F}^{(0)}((0,1]^d)$ and ${\cal S}^{(0)}({\cal R}_{\tilde{\phi}^0}(d,J))=\{\sum_{j\in {\cal R}_{\tilde{\phi}^0}(d,J)}\beta(j)\tilde{\phi}_j^0:\beta\}\subset{\cal S}^{(0)}((0,1]^d)$. Here $J$ denotes the size of the set. 
\end{definition}

\subsection{Existence of knot-point sets of size $J$ giving $O(r(d,J))$ $L^2$-approximation error}
We define the following rate $r(d,J,M)$ which represents the approximation error in $L^2(\mu)$ norm of ${\cal F}^{(0)}_M((0,1]^d)$ for an optimally chosen linear combination of $J$ zero-order splines $\phi_{u_j}^0$, $j=1,\ldots,J$. Let's first define this quantity formally.
\begin{definition}
Let 
\begin{equation}\label{defrdJM}
r(d,J,M)\equiv \arg\min_{u_1,\ldots,u_J\in (0,1]^d}\sup_{Q\in {\cal F}^{(0)}_M((0,1]^d)}\inf_{\beta}
\pl \sum_{j=1}^J\beta(j)\phi_{u_j}^0-Q\pl_{L^2(\mu)}.
\end{equation}
We will represent $r(d,J,M)=C(M)r(d,J)$. Below we prove that $r(d,J,M)=O(MJ^{-1}\log^{2(d-1)}J)$, so that we can choose as upper bounds $C(M)=M$ and $r(d,J)=J^{-1}\log^{2(d-1)}J$. In this article, the reader can replace $r(d,J)$ and $C(M)$ by these upper bounds $J^{-1}\log^{2(d-1)}J$ and $M$, respectively.
\end{definition}
It appears logical that this rate $r(d,J,M)$ should  be determined by the covering number of ${\cal F}^{(0)}((0,1]^d)$. 
We have that the $L^r$-covering number of ${\cal F}^{(0)}_M((0,1]^d)$ behaves (i.e., upper bounded by) as   $\log N(\epsilon,{\cal F}^{(0)}_1((0,1]^d) ,L^r)\sim C(M)\epsilon^{-1}(-\log \epsilon)^{2(d-1)}$  \citep{Bibaut&vanderLaan19}, where $C(M)\sim M (\log M)^{2d-2}$. For notational convenience, let $N_1(\epsilon)=N(\epsilon,{\cal F}^{(0)}_M((0,1]^d)),L^r)$.
This can also be written as $\log N_1(\epsilon)\sim (M /\epsilon)(\log (M/\epsilon))^{2(d-1)}$.
For our purpose, we can set $r=2$.
Based on this covering number we can show that that $r(d,J,M)$ behaves as $1/J$ up till $\log J$-factor. We can show this as follows. The $\epsilon$-net achieving this $\epsilon$-covering is represented by linear combination of $J(\epsilon)$ zero-order splines (i.e., discrete cumulative distribution functions with support the knot-point set) and discretizing each coefficient so that neighboring values for each coefficients are $O(\epsilon)$ apart. However, given the $L_1$-norm constraint $M$,  the $\beta$-vector values in the $\epsilon$-grid  that  make the $L_1$-norm bigger than $M$ should be removed, since these would be describe  functions outside ${\cal F}^{(0)}_M((0,1]^d)$. 
If we ignore this necessary removal of $\beta$-values, we would have $N_1(\epsilon)\sim (1/\epsilon)^{J(\epsilon})$. However, let's first determine how this order of $N_1(\epsilon)$ changes due to removing $\beta$-values that have $L_1$-norm larger than $M$.  Let $J=J(\epsilon)$.
For this proof, without loss of generality, we can simplify matters by restricting to actual cumulative distribution functions so that we can assume $\beta(j)\geq 0$ with $\sum_j\beta(j)\leq M$.
Let the knot-points $u_j$ be ordered so that we can represent them as $\mu_1,\ldots,\mu_J$ and corresponding $\beta(1),\ldots,\beta(J)$.
For $\beta(1)$ we have $M/\epsilon$ possible grid-values. Then, for $\beta(2)$ we have $M/\epsilon -\beta(1)$ possible values, due to $\sum_{j=1}^2\beta(j)\leq M$. And so on. So we can count all possible values in this $\epsilon$-grid with the following sum:
\[
\sum_{\beta(1)=1}^{M/\epsilon}\sum_{\beta(2)=(1-\beta(1))}^{M/\epsilon}\ldots\sum_{\beta(J)=1}^{M/\epsilon-\sum_{j=1}^{J-1}\beta(j)}.\]
To understand the behavior of this count, wen can represent the sums as integrals so that we have to solve 
\[
\int_{1}^{M/\epsilon}dx_1\int_{1-x_1}^{M/\epsilon} dx_2\ldots\int_{1-\sum_{j=1}^{J-1}x_j}^{M/\epsilon} dx_J.\]
This integral can be explicitly solved starting with the integral over $x_J$ and working backwards ending with integral over $x_1$. It follows that it equals
\[
\frac{(M/\epsilon)^J}{J!}.\]
So, $\log N_1(\epsilon)\sim J\log  (M/\epsilon)-\log J!$. We now note that
\[
\log J!=\sum_{j=1}^J \log j.\]
Again, this can be represented by an integral $\int_1^J \log x dx$ which is explicitly solved by $\left . (x\log x-x)\right |_1^J=J\log J-J$.
So we have $\log N_1(\epsilon)\sim J \log (M/\epsilon)-J\log J+J=J(1+\log M/\epsilon -\log J)$, where we can truncate the factor from below by $1$.
This has to be set equal to $(M/\epsilon)(\log( M/\epsilon))^{2(d-1)}$.
From this it follows that $J(\epsilon)=(M/\epsilon)(\log M/\epsilon)^k$ for some integer $k$.
Let's plug-in this expression for $J(\epsilon)$ and solve for $k$. For this expression $\log M/\epsilon -\log J=0$ due to cancellation of $\log M/\epsilon$,  so that $\log N_1(\epsilon)\sim J(\epsilon)$.
This yields the equation:
\[
(M/\epsilon)\log^k(M/\epsilon)=(M/\epsilon)\log^{2(d-1)} (M/\epsilon),\]
so that $k=2(d-1)$.
So we conclude that $J(\epsilon)=(M/\epsilon)(\log (M/\epsilon))^{2(d-1)}$.
This represents the number of zero-order spline basis functions required to obtain an $L^r$-norm  $\epsilon$-approximation of the set ${\cal F}^{(0)}_M((0,1]^d)$. Let ${\cal R}_{\phi}(d,\epsilon)$ represent a set of knot-points giving this $L^r$-$\epsilon$ approximation. 
For a given $J$ we can solve $J=J(\epsilon)$ in $\epsilon$ and we denote that solution with $\epsilon(J)$, which is the $L^r$-approximation error of the $J$-dimensional linear combination of zero-order splines $\phi_u$ with $u\in {\cal R}_{\phi}(d,\epsilon)$. We find
$M/\epsilon_J\sim J\log^{-2(d-1)}J$, and thus $\epsilon_J=MJ^{-1}\log^{2(d-1)}J$.
So with $J$ zero-order splines we obtain an $L^r$-approximation error of $MJ^{-1}\log^{2(d-1)}J$. Let ${\cal R}_{\phi}(d,J)\equiv {\cal R}_{\phi}(d,\epsilon(J))$. Then, \[
\sup_{Q\in {\cal F}^{(0)}_M((0,1]^d)}\inf_{Q_J\in {\cal F}^{(0)}({\cal R}_{\phi}(d,J))}\pl Q_J-Q\pl_{L^r(\mu)}=O(MJ^{-1}\log^{2(d-1)}J).\] 
So we have proved the following lemma.
\begin{lemma}
We have $r(d,J,M)=O(MJ^{-1}\log^{2(d-1)}J)=O(C(M)r(d,J))$ with $C(M)=M$ and $r(d,J)=J^{-1}\log^{2(d-1)}J$.
\end{lemma}

\begin{lemma}{\bf (Existence of knot-point set ${\cal R}(d,J)$  of size $J$ giving $O(Mr(d,J))$-approximation error)}
Consider ${\cal F}^{(0)}_{M}((0,1]^d)$ for a given $M<\infty$. Let $r(d,J,M)$ be defined as above. 
We have that there exist a set ${\cal R}_{\phi^0}(d,J)$ of $J$ knot-points so that 
\[
\sup_{Q\in {\cal F}^{(0)}_M((0,1]^d)}\inf_{Q_J\in {\cal F}^{(0)}({\cal R}_{\phi^0}(d,J) )}\pl Q_J-Q\pl_{L^2(\mu)}=O(r(d,J,M)).\]
We obtain the analogue result for ${\cal S}^0_M(0,1]^d$: There exists a possibly different set of knot-points ${\cal R}_{\tilde{\phi}^0}(d,J)$ but for basis functions $\tilde{\phi}_u^0(x)=I(x< u)$ so that 
\[
\sup_{Q\in {\cal S}^{(0)}_M((0,1]^d)}\inf_{Q_J\in {\cal S}^{(0)}({\cal R}_{\tilde{\phi}^0}(d,J))}\pl Q_J-Q\pl_{L^2(\mu)}=O(r(d,J,M)).\]
\end{lemma}

\begin{definition}
In this article, for a given integer $d$ and  subset $s\subset\{1,\ldots,d\}$,  let $\mu_s:[0,1]^{\mid s\mid}\rightarrow \openr$ be the lebesgue measure $\mu_s(dx(s))=\prod_{j\in s}dx_j$. We can also think of $\mu_s$ as the $\mid s\mid$-dimensional uniform cumulative distribution function with $\mu_s(x(s))=\prod_{j\in s}x_j$ on $\{x(j): j\in s,x\in (0,1]^d\}\subset (0,1]^{\mid s\mid}$.
We will also write $\mu(dx(s))$ since $x(s)$ already clarifies its restriction to $(x(j):j\in s)$. Thus, $\int Q(x(s))d\mu(x(s))=\int\ldots\int Q(x_j:j\in s )\prod_{j\in s}dx_j$ is a standard Riemann integral over all components $x_j$, $j\in s$.
\end{definition}

We have the following corollary.
\begin{corollary}\label{defknots} {\bf (Existence of knot-point set ${\cal R}(d,J)$ approximating both cumulative and survivor functions)}
Let $\pl Q\pl_{\mu}=\left( \int_{(0,1]^d} Q^2 d\mu\right)^{1/2}$ be the $L^2(\mu)$-norm.
There exists a set of $J$ knot-points ${\cal R}(d,J)$ so that  both
\begin{eqnarray}
\sup_x \inf_{\alpha}\pl \tilde{\phi}_x-\sum_{v\in {\cal R}(d,J)}\alpha(v)\tilde{\phi}_v\pl_{\mu}&=&O(C(M)r(d,J))\label{Ra}\\
\sup_{Q\in {\cal F}^0_M((0,1]^d)}\inf_{\alpha}\pl \sum_{v\in {\cal R}(d,J)}\alpha(v)\phi_v^0-Q\pl_{\mu}&=&O(C(M) r(d,J))\label{Rb} .
\end{eqnarray}
Specifically, one can define ${\cal R}(d,J)={\cal R}_{\phi}(d,J/2)\cup{\cal R}_{\tilde{\phi}}(d,J/2)$ as the union of the two sets of $J/2$  knot-points used to approximate ${\cal F}^{(0)}_M(0,1]^d$ and ${\cal S}^{(0)}_M(0,1]^d$  as defined in previous lemma. If one only needs (\ref{Rb}) then one defines ${\cal R}(d,J)={\cal R}_{\phi}(d,J)$.
\end{corollary}

\begin{definition}
In the remaining of this article ${\cal R}(d,J)$ represents a set of $J$ knot-points that satisfies (\ref{Rb}), and for one of the approximation error results  Theorem \ref{mainkthordercdf} we also need ${\cal R}(d,J)$ to satisfy  (\ref{Ra}). 
\end{definition}

The proof of Theorem \ref{mainkthordercdf} for our uniform approximation result appears more direct for ${\cal R}(d,J)$ satisfying both (\ref{Ra}) and (\ref{Rb}). However, with a trick, we are also able to proof the desired approximation error result for ${\cal R}(d,J)$ only satisfying (\ref{Rb}) (see Theorem \ref{mainkthordercdfnew}).
Our proof appears to suggest that making ${\cal R}(d,J)$ satisfy both (\ref{Ra}) and (\ref{Rb}) might help the constant in front of rate by a potentially large factor $2^d$, so it remains to evaluate what the practical impact is of choosing ${\cal R}(d,J)$ to satisfy both (\ref{Ra}) and (\ref{Rb}) w.r.t. behavior of our proposed sieve and HAL-MLEs.


\subsection{Practical method for determining knot-point sets of size $J$ giving desired $O(r(d,J))$ $L^2$-approximation error}
Corollary \ref{defknots} states the existence of ${\cal R}(d,J)$ of size $J$ satisfying (\ref{Ra}) and (\ref{Rb}), but it does not provide a practical tool for determining such knot-point sets.
For practical implementation of generating ${\cal R}_{\phi}(d,J)$ we suggest that one  may find such a set of $J$ knot-points by minimizing ${\cal R}(J)\rightarrow \inf_{\alpha}\int_{(0,1]^d}(\sum_{u\in {\cal R}(J)}\alpha(u)\phi_u^0(x)-\mu(x))^2d\mu(x)$ over all sets ${\cal R}(J)$ of $J$ knot-points, where $\mu(x)=\prod_{j=1}^d x(j)$ is the uniform cumulative distribution function. Similarly, we could determine  ${\cal R}_{\tilde{\phi}}(d,J)$ by minimizing ${\cal R}(J)\rightarrow \inf_{\alpha}\int_{(0,1]^d}(\sum_{u\in {\cal R}(J)}\alpha(u)\tilde{\phi}_u(x)-\int_{(x,1]}d\mu(y) )^2d\mu(x)$ over all sets ${\cal R}(J)$ of $J$ knot-points. That is, these sets of knot-points can be determined by optimizing the approximation error for approximating the uniform cumulative distribution function and uniform survivor function, respectively.  It is beyond the scope of the current article to implement such optimization algorithms. Instead, in  Appendix \ref{AppendixB} below we  suggest a practical method utilizing the lasso (HAL) for determining such sets $({\cal R}_{\phi}(d,J):J)$ and ${\cal R}_{\tilde{\phi}}(d,J):J)$ and prove that they indeed yield the $O(r(d,J))$ approximation error uniformly among all continuous elements ${\cal F}^{(0)}_M((0,1]^d)$, while we conjecture that it will extend to the whole ${\cal F}^{(0)}_M((0,1]^d)$.

\section{Determining knot-point sets with desired approximation error for approximating cumulative distribution functions with zero-order HAL}\label{AppendixB}

Our proposed $k$-th order spline sieve MLE relies on specifying knot-point sets $({\cal R}(\bar{s}(k+1),J):J)$ which can be defined as $\{x: x(s_{k+1})\in {\cal R}_{\phi}(\mid s_{k+1}\mid,J),x(-s_{k+1})=0\}$, where ${\cal R}_{\phi}(m=\mid s_{k+1}\mid,J)$  is the set of knot-points under which linear combinations of zero-order splines achieve an approximation error $O(r(m,J))$ of an $m$ dimensional cumulative distribution function.  
We first argue that this comes down to  fitting  uniform cumulative distribution functions, which then suggests a practical method for generating such sets ${\cal R}_{\phi}(m,J)$ and ${\cal R}_{\tilde{\phi}}(m,J)$ across $(m,J)$-values.

 \subsection*{Providing some insight why approximating a uniform cumulative distribution function implies same approximation error for class of cumulative distribution functions}
Let ${\cal R}_n$ be the knot-point set of size $J_n$ giving this approximation error $r(d,J_n)$ of a uniform cumulative distribution function $\mu$. 
Consider now a function $Q\in {\cal F}^{(0)}_M((0,1]^d)$. We have $Q=Q_d+Q_c$, where $Q_c$ is absolute continuous w.r.t. Lebesgue measure $\mu$ and $Q_d$ is discrete. 
So $Q_c(x)=\int_{(0,x]} dQ_c(y)/d\mu(y) d\mu(y)$. Let $Q_{c,J_n}=\int_{(0,x]} dQ_c/d\mu d\mu_{J_n}$, where we know $\pl \mu_{J_n}-\mu\pl_{\mu}=O(r(d,J_n))$.
 Integration by parts then proves that if $Q\in {\cal F}^{(1)}_M((0,1]^d)$, then $\pl Q_{c,J_n}-Q_c\pl_{\mu}=O(\pl \mu_{J_n}-\mu\pl_{\mu})$. This then shows that this set ${\cal R}_n$ works uniformly in all absolute continuous cumulative distribution functions. 
 This does not show yet that it also provides the desired approximation error uniformly in all discrete distributions. Since  $D^{(0)}({\cal R}_n)$ is itself a space of discrete distributions  one might argue that it should not have a harder time to approximate discrete distributions than continuous distributions, and therefore that this set will also approximate the  discrete cumulative distributions at the same rate.

\subsection{Utilizing HAL to generate knot-point sets ${\cal R}(m,J)$  across dimensions $m$ and sizes $J$ for approximating uniform cumulative distribution function at rate $r(m,J)$}
Let the dimension $m$ be given. Our goal is to fit a uniform cumulative distribution function $Q_K(x)=K\prod_{j=1}^{m} x(j)$ for integers $K=1,2,\ldots$, where $\pl Q_K\pl_v=K$. Alternatively, we can define $Q_K=Q_{1K}-Q_{2K}$ to be a difference of two of such uniform cumulative distribution function with $\pl Q_K\pl_v=2K$. First define a set of $n$ potential knot-points in $(0,1]^m$, and denote its values with $u_1,\ldots,u_n$. In our zero-order HAL-implementations for regression problems we have proposed the support of the empirical measure of the covariate distribution. In the next subsection we prove formally that this initial set of zero order spline basis functions indeed suffices. 


Secondly, we draw $Z_i\equiv Q_K(u_i)+e_i$ where $e_i\sim N(0,\sigma^2)$ is independent  mean zero noise, $i=1,\ldots,n$.  For example, we can select $\sigma^2=1$. We then use lasso regression to  fit $Q$ by minimizing the empirical mean of least squares over the $n$ i.i.d. observations $(Z_i,u_i)$, $i=1,\ldots,n$
\[
\beta\rightarrow \sum_{i=1}^n \left\{Z_i-\sum_{j=1}^n \beta(u_j)\phi_{u_j}^0(u_i)\right\}^2\]
under the $L_1$-constraint $\pl \beta\pl_1=\lambda$ and selecting $\lambda$ with the cross-validation selector $\lambda_N$. Let ${\cal R}_1(m,J_n(K))$ be the knot-points it selects, where $J_n(K)$ is the number of knot-points with non-zero coefficient.  The cross-validation selector $\lambda_N$ is asymptotically equivalent with the oracle selector and we know that the rate of of convergence of the resulting fit achieves the rate of convergence $n^{-1/3}(\log n)^{d/2}$ w.r.t. $Q_K$. The latter teaches us that the lasso is  able to select an effective set of $J$ knot-points to achieve the desired approximation  of $Q_K$ w.r.t. bias. 
We run this cross-validated lasso regression for a sequence of increasing  levels of variation norm $K$ with corresponding sets ${\cal R}_1(K,J_n(K)))$ of size $J_n(K)$, $K=1,\ldots,M$ for some large enough $M$.
By the known $L^2$-rate $n^{-1/3}(\log n)^{m/2}$ of convergence of zero-order HAL when using the cross-validation selector $J_n(K)$,  we know that $J_n(K)=C_n(K)n^{1/3}$ up till power of $\log n$-factor. Thus the size of these sets is $J_n(K)\sim n^{1/3}$ while its approximation error behaves as $n^{-1/3}\sim 1/J_n(K)$ up till $\log n$-factors, which thus shows that these sets of size $J$ give an approximation error $r(m,J)$, possibly up till a power of $\log J$.

The above method  essentially tabulate sets ${\cal R}_{\phi}(m,J(K))$ across $m$ and $J$. The same method can be applied to ${\cal R}_{\tilde{\phi}}(m,J)$ but now with $Q_K=K \prod_{j=1}^m (1-x_j)$ being a uniform survivor function or a difference of two uniform survivor functions.

If we want to tabulate ${\cal R}(s,J)$ for sets $s\subset\{1,\ldots,d\}$,
given the subset $s\subset\{1,\ldots,d\}$ of size $m$, we  can define the potential set $n$ of knot-points as $\{X_i(s): i=1,\ldots,n\}$. The advantage of the latter is that the choice of knot-points adapts to actual observed values that the empirical risk will depend on. However, in that case one would have to run our proposed procedure for generating these knot-point sets for each subset of components  $\{1,\ldots,d\}$ of size $m$, instead of using same support points across all subsets $s$ of size $m$ (but applied to the corresponding components).

 \subsection{A sufficient starting knot-point set for  zero-order HAL and proving HAL gives knot-point set with desired approximation error}
 In our work on zero-order HAL we have proposed an initial set of knot-points ${\cal R}^0_n(s,d)=\{X_i(s): i=1,\ldots,n\}$ for fitting  a cumulative distribution function on $(0(s),1(s)]\subset (0,1]^d$ for a given subset $s\subset \{1,\ldots,d\}$. In particular, this yields a starting set of zero-order spline basis functions for the zero-order HAL-MLE given by ${\cal R}^0_n(d)=\{(s,u): s\subset \{1,\ldots,d\},u\in {\cal R}^0_n(s,d)\}$. In this subsection we want to prove that indeed this initial set of knot-points ${\cal R}^0_n(s,d)$ suffices for obtaining the desired rate $n^{-1/3}(\log n)^{d/2}$ for fitting a difference of $s$-specific cumulative distribution functions, and thereby that ${\cal R}^0_n(d)$ suffices for fitting a cadlag function with bounded sectional variation norm. Let $P_{0,s}$ be the probability distribution of $X_i(s)$, $i=1,\ldots,n$.
 We rely in this proof on $Q_0$ be such that $dQ_{0,s}/dP_{0,s}$ exists for all subsets $s$. 
 

 The following lemma proves that indeed running HAL with this initial set of knot-points still achieves the desired rate $d_0(Q_n,Q_0)=O_P(n^{-2/3}(\log n)^d)$.
 In particular, it proves that applying HAL with initial knot-point set ${\cal R}^0_n=\{X_i: i=1,\ldots,n\}$ for fitting a $Q_0\in {\cal F}^{(0)}_M((0,1]^d)$ also yields this desired rate. 
 This holds uniformly in all $Q_0$ that satisfy that $dQ_{0,s}/dP_{0,s}$ exists and is uniformly bounded. 
  
  \begin{lemma}
 Let ${\cal R}^0_n(s,d)=\{X_i(s): i=1,\ldots,n\}$  for a given subset $s\subset \{1,\ldots,d\}$,
 where $\{X_i(s):i=1,\ldots,n\}$ is an i.i.d. sample of $n$ observations from $P_{0,s}$. Let ${\cal R}^0_n(d)=\cup_{s\subset\{1,\ldots,d\}}{\cal R}^0_n(s,d)$.
 Let ${D}^{(0)}({\cal R}^0_n(d))$ be the space of linear combinations of $\phi_u^0$ with $u\in {\cal R}^0_n(d)$ and $D_M^{(0)}({\cal R}^0_n)=\{Q\in D({\cal R}^0_n):\pl Q\pl_v^*\leq M\}$. Let $Q_0\in D_M^{(0)}([0,1]^d)$. 
 Let $Q_{n,C}=\arg\min_{Q\in D({\cal R}^0_n),\pl Q\pl_v^*\leq C}P_n L(Q)$ as estimator of $Q_{0,C}=\arg\min_{Q\in D_C^{(0)}([0,1]^d)}P_0L(Q)$.
 Assume that for each subset $s$ $Q_{0,s}\ll P_{0,s}$ and that $dQ_{0,s}/dP_{0,s}$  is uniformly bounded by some $M<\infty$. 
  We have $d_0(Q_{n,C},Q_{0,C})=O_P(n^{-2/3}\log^d n)$. In particular, if $C\geq M$, then $d_0(Q_{n,C},Q_0)=O_P(n^{-2/3}\log^d n)$. In addition, if $C_n$ is the cross-validation selector, and $Q_n=Q_{n,C_n}$, then $d_0(Q_n,Q_0)=O_P(n^{-2/3}\log^d n)$.

 In particular, this applies to the case that $Q_0\in {\cal F}^{(0)}_M((0,1]^d)$ is a difference of $d$-dimensional cumulative distribution functions and ${\cal R}^0_n(d)={\cal R}^0_n(\{1,\ldots,d\},d)=\{X_i:i=1,\ldots,n\}$, $X_i\sim_{iid} P_{0,\{1,\ldots,d\}}$, showing that this initial set ${\cal R}^0_n(d)$ suffices for approximating any multivariate cdf. 
 \end{lemma}
 {\bf Proof:}
  Let $\tilde{Q}_{0,s}=\int_{(0(s),x(s)]}dQ_{0,s}$ with $Q_{0,s}(x(s))=Q(0(-s),x(s))$, so that $Q_0=\sum_{s\subset\{1,\ldots,d\}} \tilde{Q}_{0,s}$.
 Let $P_{n,s}$ be the empirical measure of $\{X_i(s):i=1,\ldots,n\}$ and let $P_{0,s}$ be the true distribution of $X(s)$. We have $\tilde{Q}_{0,s}=\int_{(0(s),x(s)]} dQ_{0,s}/dP_{0,s} dP_{0,s}$, assuming that $Q_{0,s}$ is absolutely continuous w.r.t. $P_{0,s}$ and we can define $\tilde{Q}_{0,s,n}=\int_{(0(s),x(s)]} dQ_{0,s}/dP_{0,s} dP_{n,s}\in D^{(0)}({\cal R}^0_n(d,s))$. Let $Q_{0,n}^*=\sum_s \tilde{Q}_{0,s,n}\in D^{(0)}({\cal R}^0_n(d))$. We note that 
 \[
 \tilde{Q}_{0,s,n}-\tilde{Q}_{0,s}=\int_{(0(s),x(s)]} dQ_{0,s}/dP_{0,s} d(P_{n,s}-P_{0,s}).\]
 Assume that $\pl dQ_{0,s}/dP_{0,s}\pl_{\infty}<M$ for some $M<\infty$. Then, the latter is an empirical process indexed by the class of indicator functions $X(s)\rightarrow I(X(s)\leq x(s))$ indexed by $x(s)$  multiplied with a fixed bounded function. This is known to be a Donsker class, so that we know that $\sup_x \mid \tilde{Q}_{0,s,n}-\tilde{Q}_{0,s}\mid (x)=O_P(n^{-1/2})$. 
 This implies that we also have $\pl Q_{0,n}^*-Q_0\pl_{\infty}=O_P(n^{-1/2})$.  
 Therefore, it follows, in particular, that $d_0(Q_{0,n}^*,Q_0)=O_P(n^{-1})$.  Appendix \ref{AppendixI} (i.e., Theorem \ref{theoremoraclemle}) that this also implies $d_0(Q_{0,n},Q_0)=O_P(n^{-1})$  for the loss-based projection $Q_{0,n}=\arg\min_{Q\in D^{(0)}({\cal R}_n^0(d))}P_0 L(Q)$ as well. An MLE analysis of $d_0(Q_n,Q_{0,n})$ shows, analogue to the proof in \citep{Bibaut&vanderLaan19} of zero-order HAL-MLE for $d_0(Q_n,Q_0)=O_P(n^{-2/3}\log^d n)$ with $Q_n=\arg\min_{Q\in D^{(0)}_M([0,1]^d)}P_n L(Q)$, shows that $d_0(Q_n,Q_{0,n})=O_P(n^{-2/3}\log^d n)$. 
 $\Box$

 We now use this result to show that the knot-point set ${\cal R}_n$ for the non-zero coefficients selected by HAL, when starting with ${\cal R}_n^0$, and assuming that $Q_0$ is chosen difficult enough so that $d_0(Q_n,Q_{{\cal R}_n,0})\sim (J_n/n)$ up till $\log n$-factor, then $d_0(Q_{{\cal R}_n,0},Q_0)\sim J_n^{-1}$ up till $\log n$-factor as well, 
 where $J_n$ is the size of ${\cal R}_n$. In other words, indeed HAL is able to obtain sets ${\cal R}_n$ of size $J_n$ that satisfy the desired approximation error $r(d,J_n)$, at least up till a $\log n$-factor. 
 
 \begin{lemma} Recall $r(d,J)=J^{-1}(\log J)^{2(d-1)}$. Assume that for each subset $s$ $Q_{0,s}\ll P_{0,s}$ and that $dQ_{0,s}/dP_{0,s}$  is uniformly bounded by some $M<\infty$. 
 Let ${\cal R}_n$ be the set of knot-points with non-zero coefficients in $Q_n\arg\min_{Q\in D^{(0)}_{C_n}({\cal R}_n^0)}P_n L(Q)$, so that $Q_n=\sum_{u\in {\cal R}_n}\beta_n(u)\phi_u^0$, and let $Q_{0,{\cal R}_n}=\arg\min_{Q\in D^{(0)}({\cal R}_n)}P_0L(Q)$.
 Let $J_n$ be the size of ${\cal R}_n$. 
 We have $d_0(Q_n,Q_{0,{\cal R}_n})=O_P((J_n/n)\log^2 n)$. 
 Assume that it actually behaves at this rate so that $d_0(Q_n,Q_{0,{\cal R}_n})\sim (J_n/n)\log ^m n$ for some $m$, by, for example, choosing a difficult $Q_0$ such as the uniform CDF.
 
 Then, for some $m$ \[
 d_0(Q_{0,{\cal R}_n},Q_0)=O_P(J_n^{-2}\log^m J_n)),\] where $J_n$ is size of ${\cal R}_n$.
 This result can also be applied to special case that $Q_0\in {\cal F}^{(0)}((0,1]^d)$ and $Q_n$ is zero-order  HAL-MLE over ${\cal F}^{(0)}({\cal R}_n^0)$ (i.e constraining $L_1$-norm and selecting with cross-validation). 
 \end{lemma}
 {\bf Proof:} The result $d_0(Q_n,Q_{0,{\cal R}_n})=O_P((J_n/n)\log^2 n)$ follows from same proof as provided in  Theorem \ref{theoremlossbasedb}.
 We now actually assume $d_0(Q_n,Q_{0,n})\sim (J_n/n)\log^m n$ for some $m$.
 Since  $d_0(Q_n,Q_0)=d_0(Q_n,Q_{0,{\cal R}_n})+d_0(Q_{0,{\cal R}_n},Q_0)$ and $d_0(Q_n,Q_0)=O_P(n^{-2/3}\log^d n)$ it follows that then $(J_n/n)\log ^m n=O_P(n^{-2/3}\log^d n)$, and thereby that $J_n=O_P(n^{1/3}\log^{d-m}n)$. 
 It also follows that $d_0(Q_{0,{\cal R}_n},Q_0)=O_P(n^{-2/3}\log^d n)$ as well. 
 Suppose that $d_0(Q_{0,{\cal R}_n},Q_0)$ has a worse  rate than $n^{-2/3}$, say by some $n^{\delta}$ for some $\delta>0$. Plugging in $J_n=n^{1/3}\log^{d-m}n$  shows then that $d_0(Q_{0,{\cal R}_n},Q_0)^{1/2}$ has a worse rate than $n^{-1/3}$, while it should also be $O(n^{-1/3}\log^{d/2}n)$, which is a contradiction. This shows that $d_0(Q_{0,{\cal R}_n},Q_0)=O(J_n^{-1}\log^m J_n)$ for some $m$.
$  \Box$
\section{ Proof of Theorem \ref{theoremkthorderspline}}
\label{AppendixC}
In our proof we first ignore the issue that, if $s_{m+1}\subset s_m$ is the empty set (for first time, so that $m=m(\bar{s}(k+1))$) then  the $s_{m+1}$-section $Q^{(m)}_{\bar{s}(m+1)}$ of $Q^{(m)}_{\bar{s}(m)}$ equals a constant $Q^{(m)}_{\bar{s}(m+1)}=Q^{(m)}_{\bar{s}(m)}(0(s_m))$, and, as a consequence, the resulting term does not generate more terms anymore (integrating the constant just means the constant comes in front and the remaining factor corresponds with $\bar{\phi}_{\bar{s}(k+1)}$. We discuss this issue after the proof below and it will correspond with making the convention (\ref{conventiona}) stating that the term $\bar{\phi}_{\bar{s}(k+1)}\mu^k(\tilde{Q}^{(k)}_{\bar{s}(k+1)})$ in our result actually reduces to $\bar{\phi}_{\bar{s}(k+1)}Q^{(m)}_{\bar{s}(m)}(0(s_m))$, while
$\bar{\phi}_{\bar{s}(k+1)}=\prod_{j=1}^{\min(k,m(\bar{s}(k+1))}\phi_0^j(x(s_j/s_{j+1}))$.

For $k=0$, we already had the zero-order spline representation:
\[
Q(x)=\sum_{s_1\subset{\cal S}^{1}(d)}\int \phi^0_{u(s_1)}(x(s_1)) Q^{(0)}_{s_1}(du(s_1)).\]
Recall that for empty subset $s_1$, $Q^{(0)}_{s_1}$ is a point-mass at $0$ so that the integral reduces to $Q(0)$.
For notational convenience, in this proof we use $du(s)$ for $d\mu_{s}(u(s))=\prod_{j\in s}du(j)$.
Using this zero-order spline representation as start, for $k=1$, using $Q^{(0)}_{s_1}(du(s_1))=Q^{(1)}_{s_1}(u(s_1))du(s_1)$, and applying the zero-order spline representation to $Q^{(1)}_{s_1}(u(s_1))$ provides the next formula for $k=1$:
\begin{eqnarray*}
Q(x)&=&\sum_{s_1\subset{\cal S}^{1}(d)}\int \phi^0_{u(s_1)}(x(s_1)) Q^{(1)}_{s_1}(u(s_1))du(s_1)\\
&=&\sum_{\bar{s}(2)\subset{\cal S}^2(d)}\int_{(0(s_1),x(s_1)]}\int_{(0(s_2),u(s_2)]} Q^{(1)}_{\bar{s}(2)}(du_1(s_2))   du(s_1)\\
&=&\sum_{\bar{s}(2)\subset{\cal S}^2(d)}\left\{\int_{0(s_1/s_2)}^{x(s_1/s_2)} du(s_1/s_2)\right\}
\int_{0(s_2)}^{x(s_2)}\left\{\int_{[u_1(s_2),x(s_2)]}du(s_2)\right\} Q^{(1)}_{\bar{s}(2)}(du_1(s_2))\\
&=&\sum_{\bar{s}(2)\subset{\cal S}^2(d)}\phi^1_0(x(s_1/s_2))
\int_{(0(s_2),x(s_2)]}\phi^1_{u_1(s_2)}(x(s_2)) Q^{(1)}_{\bar{s}(2)}(du_1(s_2)).
\end{eqnarray*}
Note that, for $s_2$ the empty set, when $s_1$ is not the empty set, we would have \[
\begin{array}{l}
\int_{(0(s_1),x(s_1)]}\int_{(0(s_2),u(s_2)]} Q^{(1)}_{\bar{s}(2)}(du_1(s_2))   du(s_1)=
Q^{(1)}_{s_1}(0(s_1))
\int_{(0(s_1),x(s_1)]}  du(s_1)\\
=Q^{(1)}_{s_1}(0(s_1))\phi_0^1(x(s_1)).
\end{array}
\]
This proves the $k$-th order spline representation of $Q$ for $k=1$.
This method can be iterated to obtain the  $k=2$-representation: 
\[
\begin{array}{l}
Q(x)= \sum_{\bar{s}(2)}\phi^1_0
\int_{(0(s_2),x(s_2)]}\phi^1_{u_1(s_2)}(x(s_2)) Q^{(2)}_{\bar{s}(2)}(u_1(s_2))du_1(s_2)\\
=\sum_{\bar{s}(2)}\phi^1_0
\int_{(0(s_2),x(s_2)]}\phi^1_{u_1(s_2)}(x(s_2))
 \sum_{s_3}\int_{(0(s_3),u_1(s_3)]} Q^{(2)}_{\bar{s}(3)}(du_2(s_3)) du_1(s_2)\\
=\sum_{\bar{s}(3)}\phi^1_0 \int_{0(s_2)}^{x(s_2)}\phi^1_{u_1(s_2)}(x(s_2))\int_{0(s_3)}^{u_1(s_3)} Q^{(2)}_{\bar{s}(3)}(du_2(s_3)) du_1(s_2)\\
=\sum_{\bar{s}(3)}\phi^1_0\\
\int_{0(s_2/s_3)}^{x(s_2/s_3)}\int_{0(s_3)}^{x(s_3)} \phi^1_{u_1(s_2/s_3)}\phi^1_{u_1(s_3)} \int_{(0(s_3),u_1(s_3)]} Q^{(2)}_{\bar{s}(3)}(du_2(s_3)) du_1(s_2/s_3)du_1(s_3)\\
=\sum_{\bar{s}(3)}\phi^1_0\\
\int_{0(s_2/s_3)}^{x(s_2/s_3)}\phi^1_{u_1(s_2/s_3)}du_1(s_2/s_3)
\int_{0(s_3)}^{x(s_3)}  \int_{0(s_3)}^{u_1(s_3)}
\phi^1_{u_1(s_3)} Q^{(2)}_{\bar{s}(3)}(du_2(s_3)) du_1(s_3)\\
=\sum_{\bar{s}(3)}\phi^1_0\phi^2_0\\
 \int_{(0(s_3),x(s_3)]}  \int_{(0(s_3),u_1(s_3)]}
\phi^1_{u_1(s_3)}(x(s_3)) Q^{(2)}_{\bar{s}(3)}(du_2(s_3)) du_1(s_3)\\
=\sum_{\bar{s}(3)}\phi^1_0\phi^2_0
\int_{(0(s_3),x(s_3)]}  \int_{[u_2(s_3),x(s_3)]}\phi^1_{u_1(s_3)}du_1(s_3)       
              Q^{(2)}_{\bar{s}(3)}(du_2(s_3)) \\
=\sum_{\bar{s}(3)}\phi^1_0\phi^2_0
\int_{(0(s_3),x(s_3)]}  \phi^2_{u_2(s_3)}(x(s_3))    Q^{(2)}_{\bar{s}(3)}(du_2(s_3)) .
\end{array}
\]
This proves the $k$-th order spline representation of $Q$ for $k=2$.
This proof can be iterated establishing the first representation of the theorem. 

Regarding the alternative representation of $Q(x)$ in terms of $k$-th order primitives, we derive it as follows.
Starting with the $k=0$ representation, 
\[Q(x)=\sum_{s_1\subset{\cal S}^{1}(d)}\int \phi^0_{u(s_1)}(x(s_1)) Q^{(0)}_{s_1}(du(s_1)),\]
for $k=1$, we have
\[
\begin{array}{l}
Q(x)=\sum_{s_1}\int \phi^0_{u(s_1)}(x(s_1)) Q^{(1)}_{s_1}(u(s_1)) d\mu(u(s_1))\\
=\sum_{\bar{s}(2)}\int \phi^0_{u(s_1)}(x(s_1)) \int_{(0(s_2),u(s_2)]}Q^{(1)}_{s_1,s_2}(du_1(s_2)) d\mu(u(s_1))\\
=\sum_{\bar{s}(2)}\int_{(0(s_1/s_2),x(s_1/s_2)]}\phi^0_{u(s_1/s_2)}(x(s_1/s_2))d\mu(u(s_1/s_2))
\\
\int_{(0(s_2),x(s_2)]}\int_{(0(s_2),u(s_2)]}\phi^0_{u(s_2)}(x(s_2))  Q^{(1)}_{s_1,s_2}(du_1(s_2)) d\mu(s_2)\\
=\sum_{\bar{s}(2)}\phi^1_0(x(s_1/s_2))\int_{(0(s_2),x(s_2)]}\int_{(0(s_2),u(s_2)]} \phi^0_{u(s_2)}(x(s_2)) Q^{(1)}_{\bar{s}(2)}(du_1(s_2)) d\mu(u(s_2))\\
=\sum_{\bar{s}(2)}\phi^1_0(x(s_1/s_2))\int_{(0(s_2),x(s_2)]}\int_{(0(s_2),u(s_2)]}Q^{(1)}_{\bar{s}(2)}(du_1(s_2)) d\mu(u(s_2))\\
=\sum_{\bar{s}(2)}\phi^1_0(x(s_1/s_2))\int_{(0(s_2),x(s_2)]}\tilde{Q}^{(1)}_{\bar{s}(2)}(u(s_2))d\mu(u(s_2))\\
=\sum_{\bar{s}(2)}\phi^1_0(x(s_1/s_2))\mu(\tilde{Q}^{(1)}_{\bar{s}(2)})(x(s_2)).
\end{array}
\]
This proves the representation in terms of first order primitives for $k=1$.
For $k=2$, we proceed:
\[
\begin{array}{l}
Q(x)= \sum_{\bar{s}(2)}\phi^1_0(x(s_1/s_2))
\int_{(0(s_2),x(s_2)]}\int_{(0(s_2),u(s_2)]} Q^{(2)}_{\bar{s}(2)}(u_1(s_2)) d\mu(u_1(s_2))d\mu(u(s_2))\\
=\sum_{\bar{s}(3)}\phi^1_0
\int_{(0(s_2),x(s_2)]}\int_{(0(s_2),u(s_2)]}\int_{(0(s_3),u_1(s_3)]} Q^{(2)}_{\bar{s}(3)}(du_2(s_3)) d\mu(u_1(s_2))d\mu(u(s_2))\\
=\sum_{\bar{s}(3)}\phi^1_0
\int_{0(s_2/s_3)}^{x(s_2/s_3)}\int_{0(s_2/s_3)}^{u(s_2/s_3)}d\mu(u_1(s_2/s_3))d\mu(u(s_2/s_3)) \\
\int_{0(s_3)}^{x(s_3)}\int_{0(s_3)}^{u(s_3)}\int_{0(s_3)}^{u_1(s_3)} Q^{(2)}_{\bar{s}(3)}(du_2(s_3)) d\mu(u_1(s_3))d\mu(u(s_3))\\
=\sum_{\bar{s}(3)}\phi^1_0\phi^2_0(x(s_2/s_3))
\int_{0(s_3)}^{x(s_3)}\int_{0(s_3)}^{u(s_3)}\int_{0(s_3)}^{u_1(s_3)} Q^{(2)}_{\bar{s}(3)}(du_2(s_3)) d\mu(u_1(s_3))d\mu(u(s_3))\\
=\sum_{\bar{s}(3)}\phi^1_0\phi^2_0(x(s_2/s_3))
\int_{0(s_3)}^{x(s_3)}\int_{0(s_3)}^{u(s_3)}\tilde{Q}^{(2)}_{\bar{s}(3)}(u_1(s_3))  d\mu(u_1(s_3))d\mu(u(s_3))\\
=\sum_{\bar{s}(3)}\phi^1_0\phi^2_0(x(s_2/s_3))
\mu^2(\tilde{Q}^{(2)}_{\bar{s}(3)})(x(s_3)).
\end{array}
\]
This proves the representation of $Q$ in terms of second  order primitives for $k=2$. 
Again, this proof can be iterated giving the general presentation in terms of $k$-th order primitives, as stated in the theorem. 

Let's now discuss in more detail what happens when $s_{m+1}$ is for the first time the empty set.
The $\bar{s}(m)$-specific term to which we apply the generation of sections defined by $s_{m+1}$ is given by 
$\bar{\phi}_{\bar{s}(m)}\mu^m( Q^{(m)}_{\bar{s}(m)} )$: in fact, we first had a term $\bar{\phi}_{\bar{s}(m)}\mu^{m-1}(\tilde{Q}^{(m-1)}_{\bar{s}(m)})$, where \[
\tilde{Q}^{m-1}_{\bar{s}(m)}))(u_{s_m})=
\int_{(0(s_m),u(s_m)]} Q^{(m-1)}_{\bar{s}(m)}(du(s_m))
=\int_{(0(s_m),u(s_m)]} Q^{(m)}_{\bar{s}(m)}(u(s_m))d\mu(u(s_m)),\]
thereby becoming $\mu^m(Q^{(m)}_{\bar{s}(m)})$. 
At that point we replace $Q^{(m)}_{\bar{s}(m)} =\sum_{s_{m+1}\subset s_m}Q^{(m)}_{\bar{s}(m),s_{m+1}}$ by the sum of sections $Q^{(m)}_{\bar{s}(m+1)}$ defined by setting  the coordinates in $s_m/s_{m+1}$ in $Q^{(m)}_{\bar{s}(m)}$ equal to zero.  If $s_{m+1}$ is empty set, we generate the term  
$\bar{\phi}_{\bar{s}(m)}\mu^m (Q^{(m)}_{\bar{s}(m)}(0(s_m))  )$, which equals
$\bar{\phi}_{\bar{s}(m)}Q^{(m)}_{\bar{s}(m)}(0(s_m)) \mu^m(1)$, while $\mu^m(1)(x(s_m))=\phi^m_0(x(s_m))$.
Note $\bar{\phi}_{\bar{s}(m)}=\prod_{j=1}^{m-1}\phi^j_0(x(s_j/s_{j+1})$ so that multiplying the latter with $\phi^m_{0(s_m)}(x_{s_m})$ yields
$\bar{\phi}_{\bar{s}(m+1)}$.
So the term becomes $\bar{\phi}_{\bar{s}(m+1)} Q^{(m)}_{\bar{s}(m)}(0(s_m))$, which is precisely our convention (\ref{conventiona}).
$\Box$

 \section{Approximating $k$-th order primitives with linear combinations of $k$-th order splines}\label{AppendixD}

\subsection{Defining linear combinations of $k$-th order splines approximating a $k$-th order primitive function}

Recall that, for a given set ${\cal R}(d,J)$ as defined in Corollary \ref{defknots}, we define $\tilde{Q}_{J,\beta}=\sum_{u\in {\cal R}(d,J)}\beta(u)\phi_u^0$ as the discrete distribution with support ${\cal R}(d,J)$ and pointmasses $\beta(u)$, $u\in {\cal R}(d,J)$. 
By plugging $\tilde{Q}_{J,\beta}$ for $\tilde{Q}$ in the definition of $\mu^{k}(\tilde{Q})$ with $\tilde{Q}\in {\cal F}^{(0)}_M((0,1]^d)$,  using $\mu^k(\phi_u^0)=\phi_u^k$, we obtain linear approximations 
$\mu^k(\tilde{Q}_{J,\beta})$ of $\mu^k(\tilde{Q})$, where
\[
\mu^k(\tilde{Q}_{J,\beta})=\sum_{u\in {\cal R}(d,J)}\beta(u)\phi_u^k(x) .\]


By our definition of ${\cal R}(d,J)$ and Corollary \ref{defknots} we already established the following lemma.
\begin{lemma}
Let $\pl Q\pl_{\mu}\equiv\sqrt{\int Q^2d\mu}$.
We have
 \[
 \sup_{\tilde{Q}\in {\cal F}^{(0)}_M((0,1]^d)}\inf_{\beta}\pl \tilde{Q}_{J,\beta}-\tilde{Q}\pl_{\mu}=O(C(M) r(d,J)).\]
 \end{lemma}
  
  In the following subsection Theorem \ref{mainkthordercdf} establishes the approximation error for $\mu^{k}(\tilde{Q}_{J,\beta})$ of $\mu^{k}(\tilde{Q})$ w.r.t. supremum norm, uniformly in $\tilde{Q}\in {\cal F}^{(0)}_M((0,1]^d)$, relying on ${\cal R}(d,J)$ satisfying both (\ref{Ra}) and (\ref{Rb}). 
  

  After the proof of Theorem \ref{mainkthordercdf}, in the last subsection Theorem \ref{mainkthordercdfnew}  generalizes this result to ${\cal R}(d,J)$ only relying on (\ref{Rb}). 
     
 \subsection{Uniform approximation error for linear combination of $k$-th order splines of a $k$-th order primitive function}
  Let's first state our general main theorem.
\begin{theorem}\label{mainkthordercdf}
Let ${\cal R}(d,J)$ be defined as in Corollary \ref{defknots} satisfying both (\ref{Ra}) and (\ref{Rb}).
Let $k\in \{1,2,\ldots\}$ be a given integer.
Let $\mu^k_{J,\beta}(\tilde{Q})\equiv \mu^k(\tilde{Q}_{J,\beta})$ be the approximation of $\mu^k(\tilde{Q})$  implied by plugging in the zero-order spline approximation $\tilde{Q}_{J,\beta}=\sum_{u\in {\cal R}(d,J)}\beta(u)\phi_u^0$ for $\tilde{Q}$ in $\mu^{k}(\tilde{Q})$.

Define the criterion:
\[
R^k_J(\beta)\equiv  J^{-1}\sum_{v\in {\cal R}(d,J)}\left( \mu^{k}_{J,\beta}(\tilde{Q})-\mu^{k}(\tilde{Q})\right)^2(v),\]
Let $\beta_J^{*}=\arg\min_{\beta}R^k_J(\beta)$ be the minimizer, which can also be denoted with $\beta_J^*(\tilde{Q})$ to emphasize that it depends on $\tilde{Q}$.
We have 
\[\sup_{\tilde{Q}\in {\cal F}^{(0)}_M((0,1]^d)} \pl \mu^{k}_{J,\beta_J^{*}(\tilde{Q})}(\tilde{Q})-\mu^{k}(\tilde{Q})\pl_{\infty}=O(C(M) r(d,J)^{k+1}).\]

This proves that
\[
\sup_{Q\in {\cal F}^{(k)}_{M}((0,1]^d)}\inf_{\beta}\pl \sum_{u\in {\cal R}(d,J)}\beta(u)\phi_u^k-Q\pl_{\infty}=O(C(M) r(d,J)^{k+1}).\]

In addition, we also have for $k=2,\ldots$
\[
\sup_{\tilde{Q}\in {\cal F}^{(0)}_M((0,1]^d)} \pl \mu^{k-1}_{J,\beta_J^*(\tilde{Q})}(\tilde{Q})-\mu^{k-1}(\tilde{Q})\pl_{\infty}=O(C(M)r(d,J)^k),\]
while for $k=1$, we have
\[
\sup_{\tilde{Q}\in {\cal F}^{(0)}_M((0,1]^d)} \pl \tilde{Q}_{J,\beta_J^*(\tilde{Q})}-\tilde{Q}\pl_{\mu}=O(C(M) r(d,J)).\]
This shows also convergence in variation norm of $\mu^{k}_{J,\beta_J^{*}}(\tilde{Q})$:
for $k=1,\ldots$
\[ 
\sup_{\tilde{Q}\in {\cal F}^{(0)}_M((0,1]^d)} \pl \mu^{k}_{J,\beta_J^*(\tilde{Q})}(\tilde{Q})-\mu^{k}(\tilde{Q})\pl_{v}=O(C(M) r(d,J)^{k}).\]
\end{theorem}

\subsection{Proof of Theorem \ref{mainkthordercdf} for $k=1$.}
We first state the theorem specifically for $k=1$ and will first prove this special case. However, this proof will essentially just be recursively  applied so that it captures most of the essence. 

\begin{theorem}\label{1storder}
Recall $\tilde{\phi}_u(x)=I(x<u)$, and we use also notation $\tilde{\phi}_x$ for $u\rightarrow\tilde{\phi}_u(x)$.
Define the knot-point set ${\cal R}(d,J)$ as in Corollary \ref{defknots} so that the $J$ knot-points are chosen so  that
\begin{eqnarray*}
\sup_x \inf_{\alpha}\pl \tilde{\phi}_x-\sum_{v\in {\cal R}(d,J)}\alpha(v)\tilde{\phi}_v\pl_{\mu}&=&O(r(d,J))\\
\sup_{\tilde{Q}\in {\cal F}^{(0)}_M((0,1]^d)}\inf_{\alpha}\pl \sum_v\alpha(v)\phi_v-\tilde{Q}\pl_{\mu}&=&O(C(M) r(d,J)).
\end{eqnarray*}
Let $\tilde{Q}_{J,\beta}=\sum_{u\in {\cal R}(d,J)}\beta(u)\phi_u^0$.
Define the criterion:
\[
R^1_J(\beta)\equiv  J^{-1}\sum_{v\in {\cal R}(d,J)}\left( \int_{(0,v]}\left(\sum_{u\in {\cal R}(d,J)}\beta(u)\phi_u^0-\tilde{Q}\right )d\mu \right)^2.\]
Note, 
\[
R^1_J(\beta)=J^{-1}\sum_{v\in {\cal R}(d,J)}\left(\mu^1_{J,\beta}(\tilde{Q})-\mu(\tilde{Q})\right)^2(v),\]
where $\mu^1_{J,\beta}(\tilde{Q})=\mu(\tilde{Q}_{J,\beta})$. 
 
Let $\beta_J^*=\arg\min_{\beta}R^1_J(\beta)$ be the minimizer, which could also be denoted with $\beta_J^*(\tilde{Q})$.
Then,\[
\sup_{\tilde{Q}\in {\cal F}_M^{(0)}((0,1]^d)} \pl \mu^1_{J,\beta_J^*(\tilde{Q})}(\tilde{Q})-\mu(\tilde{Q})\pl_{\infty}=O(C(M) r(d,J)^2).\]
In addition, 
\[
\sup_{\tilde{Q}\in {\cal F}_M^{(0)}((0,1]^d)} \pl \tilde{Q}_{J,\beta_J^*(\tilde{Q})}-\tilde{Q}\pl_{\mu}=O(C(M)r(d,J)).\]
This implies
\[
\sup_{\tilde{Q}\in {\cal F}_M^{(0)}((0,1]^d)} \pl \mu^1_{J,\beta_J^*(\tilde{Q})}(\tilde{Q})-\mu(\tilde{Q})\pl_{v}=O(C(M)r(d,J)).\]
\end{theorem}
{\bf Proof of Theorem \ref{1storder}:}
Let $\beta_J$ be such that $\pl \tilde{Q}_{J,\beta_J}-\tilde{Q}\pl_{\mu}=O(r(d,J))$.
{\bf The steepest descent definition of $\beta_J^*$ and its derivative equations:}
We will apply a steepest descent algorithm, starting at $\beta_J$,  aiming to minimize over $\beta$ the criterion
\[
R(\beta)= 0.5 J^{-1}\sum_{v\in {\cal R}(d,J)}\left( \int_{(0,v]}\left (\sum_{u\in {\cal R}(d,J)}\beta(u)\phi_u^0-\tilde{Q} \right ) d\mu\right)^2,\]
and the resulting update $\beta_J^1$ will be shown to solve the derivative equations \[
H_v(\beta_J^1)\equiv \int_{(0,v]}\{ \tilde{Q}_{J,\beta_J^1(\tilde{Q})}-\tilde{Q}\} d\mu=0 \] for all $v\in {\cal R}(d,J)$.  Note that $R(\beta)=0.5/J \sum_{v\in {\cal R}(d,J)}\{H_v(\beta)\}^2$ so that the equations $H_v(\beta_J^1)=0$ are equations implied by the derivative equations $d/d\beta R(\beta)=0$ at $\beta=\beta_J^1$, as shown in detail below, the same equations solved by the global minimizer  $\beta_J^*$  
 Since $R(\beta)$ is strictly convex it  has only one local minimum, which proves that $\beta_J^1 $ has to be equal to the  global minimizer $\beta_J^*$ of  $R(\beta)$ (by definition of $\beta_J^*$). So we will denote $\beta_J^1$ with $\beta_J^*$.

We need to bound the supremum norm of    $\mu_{J,\beta_J^*}(\tilde{Q})(x)=\int_{(0,x]}\tilde{Q}_{J,\beta_J^*} (y) d\mu(y)- \mu(\tilde{Q})$. 
We will prove that the update $\beta_J^*$ of $\beta_J$ preserves the rate of $\beta_J$ so that (just as known to hold for $\beta_J$) \begin{equation}\label{k1a}
\pl \tilde{Q}_{J,\beta_J^* }-\tilde{Q}\pl_{\mu}=O(r(d,J)).\end{equation} 


Assume for now that $H_v(\beta_J^*)=0$ for all $v\in {\cal R}(d,J)$ and (\ref{k1a}) has been shown to hold. We then proceed as follows.
\newline
{\bf Main proof:}
Note that 
\[
H_v(\beta_J^*)=\int \tilde{\phi}_v(y)(\tilde{Q}_{J,\beta_J^*}-\tilde{Q})  d\mu(y).\]
Since $H_v(\beta_J^*)=0$ for all $v\in {\cal R}(d,J)$, for all vectors $\alpha$
\[
0=\int \sum_{v\in {\cal R}(d,J)}\alpha(v)\tilde{\phi}_v(y)(\tilde{Q}_{J,\beta_J^*}-\tilde{Q}) d\mu(y).\]

 We can then carry out the  following proof. For any vector $\alpha_x$ we have
\begin{eqnarray*}
(\mu_{J,\beta_J^*}(\tilde{Q})-\mu(\tilde{Q}))(x)&=&\int_{(0,x]}(\tilde{Q}_{J,\beta_J^*}-\tilde{Q})) d\mu(y)\\
&=&\int \tilde{\phi}_x(y) (\tilde{Q}_{J,\beta_J^*}-\tilde{Q})(y)d\mu(y)\\
&=&\int (\tilde{\phi}_x(y)-\sum_{v\in {\cal R}(d,J)}\alpha_x(v)\tilde{\phi}_v(y) )    (\tilde{Q}_{J,\beta_J^*}-\tilde{Q})(y)       d\mu(y)\\
&\leq&  \pl \tilde{\phi}_x-\sum_{v\in {\cal R}(d,J)}\alpha_x(v)\tilde{\phi}_v\pl_{\mu} \pl \tilde{Q}_{J,\beta_J^*} -\tilde{Q}\pl_{\mu},
\end{eqnarray*}
where we apply the Cauchy-Schwarz inequality at the last line.
Therefore, we have
\begin{eqnarray*}
(\mu_{J,\beta_J^*}(\tilde{Q})-\mu(\tilde{Q}))(x)&\leq&
\inf_{\alpha} \pl \tilde{\phi}_x-\sum_{v\in {\cal R}(d,J)}\alpha(v)\tilde{\phi}_v\pl_{\mu} \pl \tilde{Q}_{J,\beta_J^*}-\tilde{Q}\pl_{\mu}\\
&=& O(r(d,J)^{2}),
\end{eqnarray*}
by the property (\ref{Ra})  of ${\cal R}(d,J)$ and (\ref{k1a}). 
This bound is uniformly in $x$ and uniformly in all $\tilde{Q}$ with $\pl\tilde{ Q}\pl_{v}<1$. Thus, this proves $\pl \mu_{J,\beta_J^*}(\tilde{Q})-\mu(\tilde{Q})\pl_{\infty}=O(r(d,J)^2)$. To conclude, we have constructed a $\mu_{J,\beta_J^*}(\tilde{Q})$ so that uniformly in $\tilde{Q}$ with $\pl \tilde{Q}\pl_{v}<1$:
1) $\pl \tilde{Q}_{J,\beta_J^*}-\tilde{Q} \pl_{\mu} =O(r(d,J))$; 2) $\pl \mu(\tilde{Q}_{J,\beta_J^*})-\mu(\tilde{Q})\pl_{\infty}=O(r(d,J)^{2})$. This then completes  the proof of the theorem. It remains to prove that 1) the steepest descent algorithm update $\beta_J^*$ of $\beta_J$ preserves the rate of the initial $\beta_J$ so that (\ref{k1a}) holds, and 2), that $\beta_J^*$ solves the equations $H_v(\beta_J^*)=0$ for all $v\in {\cal R}(d,J)$. 



{\bf Why should the steepest descent update $\beta_J^*$ of $\beta_J$ satisfy (\ref{k1a})?}
Before we start the formal proof of these two remaining tasks, let's provide the intuition behind 1).
We already have that $H_v(\beta_J)=\int_{(0,v]}(\tilde{Q}_{J,\beta_J}-\tilde{Q})d\mu=O(r(d,J))$ uniformly in $v\in {\cal R}(d,J)$, so that we should only have to induce a change of $\tilde{Q}_{J,\beta_J}$ in $L^2(\mu)$-norm of size at most  $O(r(d,J))$ to make it solve these equations $H_v(\beta_J^*)=0$ for all $v\in {\cal R}(d,J)$. That would then mean we still have that the modified $\pl \tilde{Q}_{J,\beta_J^*}-\tilde{Q}\pl_{\mu}=O(r(d,J))$.

{\bf Defining the local steepest descent algorithm applied to initial $\beta_J$:}
Note that $R(\beta)=0.5 J^{-1}\sum_{v\in {\cal R}(d,J)}(\mu_{J,\beta}(\tilde{Q}) -\mu(\tilde{Q}))^2$, where $\mu_{J,\beta}(
\tilde{Q})=\mu (\sum_{u\in {\cal R}(d,J)}\beta(u)\phi_u^0)$.
In other words, we are locally optimizing the performance for $\pl \mu_{J,\beta}(\tilde{Q})-\mu(\tilde{Q})\pl_{J}$ w.r.t. sum of squared residuals at the knot-points starting at $\beta_J$, where we define the Euclidean type norm $\pl Q\pl_J\equiv \left(1/J\sum_{v\in {\cal R}(d,J)} Q(v)^2\right)^{1/2}$. We also define the corresponding inner product $\langle h_1,h_2\rangle_J=1/J\sum_{v\in {\cal R}(d,J)}h_1(v)h_2(v)$.

Recall $H_v(\beta)= \int_{(0,v]}(\sum_{u\in {\cal R}(d,J)}\beta(u)\phi_u^0-\mu(Q))  d\mu$. 
Since \[
R(\beta)=\frac{0.5}{J}\sum_{v\in {\cal R}(d,J)}\{H_v(\beta)\}^2,\]
 using that $\phi_u^1(v)=\int_{(0,v]} \phi_u^0(y_1)d\mu(y_1)$, it follows that the pathwise derivative along paths $\beta+\delta h$ of this criterion $R(\beta)$ at $\delta=\delta_0=0$ is given by:
\begin{eqnarray*}
\frac{d}{d\delta_0}R(\beta+\delta_0 h)&=&
J^{-1}\sum_{v\in {\cal R}(d,J)}H_v(\beta) \sum_{u\in {\cal R}(d,J)}h(u)\phi_u^1(v)\\
&=&
J^{-1}\sum_{u\in {\cal R}(d,J)} h(u) \sum_{v\in {\cal R}(d,J)}H_v(\beta)\phi_u^1(v).\end{eqnarray*}
Then, \[
D^*_{\beta}(u)=\sum_{v\in {\cal R}(d,J)}H_v(\beta)\phi_u^1(v)\] is the canonical gradient of this pathwise derivative of $R(\beta)$ at $\beta$ w.r.t. inner product  $\langle h_1,h_2\rangle_J$.  That is, 
\[
\frac{d}{d\delta_0}R(\beta+\delta_0 h)=\langle D^*_{\beta},h\rangle_J .\]
Thus, the local steepest descent path for this criterion $R(\beta)$  through $\beta_J$ is given by \[
\beta_{J,\delta}^{lsd}=\beta_J+\delta D^*_{\beta_J}.\]
Since  $\frac{d}{d\delta_0}R(\beta_J+\delta_0 D^*_{\beta_J})>0$, it follows that minimizing $R(\beta)$ along this path corresponds with moving in direction of $\delta<0$.
The local steepest descent path  defines a universal steepest descent path by recursively tracking the local steepest descent  path (in direction of decreasing $R(\cdot)$) with infinitesimal small steps starting at $\beta_J$, which corresponds with the definition:
\[
\beta_{J,\epsilon}^{usd}=\beta_J-\int_{(0,\epsilon]} D^*_{\beta_{J,x-}^{usd}}  dx.\]
This universal steepest descent path  has the property that at any $\epsilon\geq 0$, we have
\[
\frac{d}{d\epsilon}\beta_{J,\epsilon}^{usd}=-D^*_{\beta_{J,\epsilon}^{usd}}.\]

Along this universal steepest descent path $\{\beta_{J,\epsilon}^{usd}:\epsilon\}$ we have for all $\epsilon\geq 0$ (direction of decreasing $R(\beta)$)
\[
\frac{d}{d\epsilon}R(\beta_{J,\epsilon}^{usd})=-\langle D^*_{\beta_{J,\epsilon}^{usd}},D^*_{\beta_{J,\epsilon}^{usd}}\rangle_J=-\pl D^*_{\beta_{J,\epsilon}^{usd}}\pl_J^2
 .\]
Let $\epsilon_J=\arg\min_{\epsilon}R(\beta_{J,\epsilon}^{usd})$ be the minimizer over $\epsilon$ along this path. This defines our proposed modification $\beta_J^*\equiv \beta_{J,\epsilon_J}^{usd}$ of $\beta_J$ that solves
\[
D^*_{\beta_J^*}(u)=0\mbox{ for all $u\in {\cal R}(d,J)$} .\]
{\bf Showing that the universal steepest descent algorithm  update $\beta_J^*$ of $\beta_J$ solves $H_v(\beta_J^*)=0$ for all $v\in {\cal R}(d,J)$:}
The last equations implies
\[
\sum_{v\in {\cal R}(d,J)}H_v(\beta_J^*)\phi_u^1(v)=0\mbox{ for all $u\in {\cal R}(d,J)$.}\]
The $J$ vectors $(\phi_u^1(v):v\in {\cal R}(d,J))$ indexed by $u\in {\cal R}(d,J))$ are independent in $\openr^J$.  Therefore, the last equation implies that $H_v(\beta_J^*)=0$ for all $v\in {\cal R}(d,J)$ so that the desired constraints are  indeed solved with this update of $\beta_J$.
 
{\bf Showing that the update $\beta_J^*$ preserves rate of $\beta_J$:}
To prove (\ref{k1a}) it  suffices to bound the $L_1$-norm \[
\pl \beta_{J,\epsilon_J}^{usd}-\beta_J\pl_1=\pl \int_{(0,\epsilon_J]} D^*_{\beta_{J,x-}^{usd}}  dx\pl_1 \]
by $O(r(d,J))$.
We note that 
\begin{eqnarray*} \pl \int_{(0,\epsilon_J]} D^*_{\beta_{J,x-}^{usd}}  dx\pl_1 &=&
\sum_{u\in {\cal R}(d,J)}\left | \int_{(0,\epsilon_J]} D^*_{\beta_{J,x-}^{usd}}(u) dx\right |\\
&\leq& \int_{(0,\epsilon_J]}\sum_{u\in {\cal R}(d,J)}\left | D^*_{\beta_{J,x-}^{usd}}(u)\right |  dx\\
&=& \int_{(0,\epsilon_J]}\pl D^*_{\beta_{J,x-}^{usd}}\pl_1 dx.
\end{eqnarray*}
The Euclidean norm $\pl D^*_{\beta_{J,x-}^{usd}}\pl_J$ monotonically decreases as $x$ moves from $0$ to $\epsilon_J$. Moreover, $\max_u \mid D^*_{\beta_{J,x-}^{usd}}\mid(u)$ approaches zero as $x$ goes from zero to $\epsilon_J$ and attains zero at $\epsilon_J$.  Therefore, we can bound the $L_1$-norm $\pl D^*_{\beta_{J,x-}^{usd}}\pl_1$  at $x\in (0,\epsilon_J]$  by  its value at $x=0$ $\pl D^*_{\beta_J}\pl_1$ or a constant times the latter. So
\[ \int_{(0,\epsilon_J]}\pl D^*_{\beta_{J,x-}^{usd}}\pl_1 dx\sim \epsilon_J \pl D^*_{\beta_J}\pl_1 .\]
Therefore, we can focus on bounding
\[
\pl \epsilon_J D^*_{\beta_J}\pl_1.\]


We have
\begin{eqnarray}
R(\beta_J)-R(\beta_{J,\epsilon_J}^{usd})&\approx& \epsilon_J\frac{d}{d\epsilon_0}
R(\beta_{J,\epsilon_0}^{usd})\nonumber \\
&=&\epsilon_J\langle D^*_{\beta_J},D^*_{\beta_J}\rangle_J,\label{preservea}
\end{eqnarray}
where $\epsilon_0=0$.




Based on this equation (\ref{preservea}) we want to understand the $L_1$-norm of $\epsilon_J D^*_{\beta_J}$ in terms of $r(d,J)$.
Firstly, note $R(\beta_J)=O(r(d,J)^{2})$. 
We also  know that $\max_{v\in {\cal R}(d,J)}\mid H_v(\beta_J)\mid =O(r(d,J))$.
The left-hand side of (\ref{preservea}) is bounded by $R(\beta_J)=O(r(d,J)^{2}$.
Thus, we have
\[
\epsilon_J=\frac{O(r(d,J)^{2})}{\pl D^*_{\beta_J}\pl_J^2} .\]

Recall $D^*_{\beta_J}(u)=\sum_{v\in {\cal R}(d,J)}H_v(\beta_J)\phi_u^1(v)$, where $H_v(\beta_J)$  behaves as $\sim r(d,J)$, while $\phi_u^1(v)$ is uniformly bounded by $1$. 
Therefore, under the assumption that $H_v(\beta_J)\sim r(d,J)$, we do not only have that $\max_u\mid D^*_{\beta_J}(u)\mid =O(Jr(d,J))$, but we can also state that $ D^*_{\beta_J} $ behaves as $\sim J r(d,J)$ (i.e., not smaller order). Specifically, we  consider the case that
\[
C(J)\equiv \frac{J r(d,J)}{
\pl D^*_{\beta_J}\pl_J }=O(1).\]
 Then, we have
 \[
 \epsilon_J=\frac{O(r(d,J)^{2})}{\pl D^*_{\beta_J}\pl_J^2} =O\left (C(J) \frac{r(d,J)}{\pl D^*_{\beta_J}\pl_J} J^{-1}\right ) =O(J^{-2}),\]
 where at last equality we use that $r(d,J)/\pl D^*_{\beta_J}\pl_J=O(1/J)$ by $C(J)=O(1)$.
 Combined with $\max_u\mid D^*_{\beta_J}(u)\mid =O(Jr(d,J))$, this gives
 that 
 \[
 \begin{array}{l}
 \pl \epsilon_J D^*_{\beta_J}\pl_1=\epsilon_J\sum_u\mid D^*_{\beta_J}(u)\mid=
 O (J^{-2} J) \max_u\mid D^*_{\beta_J}(u)\mid \\
 =O(J^{-1} J r(d,J))=O(r(d,J)).\end{array}
 \]
 This proves $\pl \beta_J^*-\beta_J\pl_1=O(r(d,J))$, and thus
  \[
 \pl \mu_{J,\beta_J^*}(\tilde{Q})-\mu(\tilde{Q})  \pl_{\mu}=O(r(d,J)).\]
 The above proof relied on $C(J)=O(1)$.
 
 Suppose now that $C(J)$ converges to infinity. Then we can write
 $\epsilon_J D^*_{\beta_J}=\epsilon_{1,J} \tilde{D}^*_{\beta_J}$, where $\epsilon_{1,J}=\epsilon_J C(J)$ and $\tilde{D}^*_{\beta_J}=D^*_{\beta_J}/C(J)$. 
 Now we have
 $(Jr(d,J))/\pl \tilde{D}^*_{\beta_J}\pl_J=O(1)$. We have
  \begin{equation}\label{z1a}
R(\beta_J)-R(\beta_{J,\epsilon_J}^{usd})\approx \epsilon_J \pl D^*_{\beta_J}\pl_J^2=\epsilon_{1,J} C(J)\pl \tilde{D}^*_{\beta_J}\pl_J^2.
\end{equation}
We can now use that $1/\pl \tilde{D}^*_{\beta_J}\pl_J=O(J^{-1}r(d,J)^{-1})$.
Thus.
\begin{eqnarray*}
\epsilon_{1,J}&=&\frac{O(r(d,J)^{2})}{C(J)\pl D^*_{\beta_J}\pl_J^2} \\
&=&\frac{O(r(d,J)^{2}r(d,J)^{-2}J^{-2}}{C(J)}\\
&=&O(C(J)^{-1} J^{-2}).
\end{eqnarray*}
Then, the $L_1$-norm of $\epsilon_{1,J}\tilde{D}^*_{\beta_J}$ is bounded as:
\begin{eqnarray*}
\pl \epsilon_{1,J}\tilde{D}^*_{\beta_J}\pl_1&=&O(C(J)^{-1}J^{-2}) C(J)^{-1}\pl D^*_{\beta_J}\pl_1\\
&=&O(C(J)^{-2}J^{-2})O(J^2 r(d,J))\\
&=&O(C(J)^{-2}r(d,J)).
\end{eqnarray*}
So this shows that our earlier scenario with $C(J)=O(1)$  represents the conservative scenario in the sense  that the rate of $\pl \beta_J^*-\beta_J\pl_1$ is even smaller order than $O(r(d,J))$  if $C(J)\rightarrow\infty$.



We have shown that $\pl \tilde{Q}_{J,\beta_J^*} -\tilde{Q} \pl_{\mu} =O(r(d,J))$ (same rate as $\pl \tilde{Q}_{J,\beta_J}-\tilde{Q} \pl_{\mu}$); $\pl \mu_{J,\beta_J^*}(\tilde{Q})-\mu(\tilde{Q})\pl_{\infty} =O(r(d,J)^2)$.
It is also clear that this proof of the preservation of the rate is uniformly in $\tilde{Q}\in {\cal F}^{(0)}_M((0,1]^d)$.
This completes the proof of the Theorem \ref{1storder}. 
$\Box$


\subsection{Proof of  Theorem \ref{mainkthordercdf} for general $k$.}
We will carry out a proof by induction. 
Above we already proved the Theorem \ref{1storder} for $k=1$. 
Let $k\in \{2,\ldots\}$ be given and we can assume that we already proved the result for $k-1$ so that there exists a $\beta_J$ so that $\pl \mu^{k-1}_{J,\beta_J}(\tilde{Q})-\mu^{k-1}(\tilde{Q})\pl_{\infty}=O(r(d,J)^k)$. 
We now want to prove that $\pl \mu^{k}_{J,\beta_J^*} (\tilde{Q})-\mu^{k}(\tilde{Q})\pl_{\infty}=O(r(d,J)^{k+1})$. The proof essentially repeats the proof of Theorem \ref{1storder} for $k=1$, but is presented here for completeness.


We will update $\beta_J$ into $\beta_J^*$ through a steepest descent algorithm that starts at $\beta_J$ and solves its derivative equations \begin{equation}\label{Hbeta}
H_v(\beta_J^*)\equiv \int_{(0,v]} \left\{ \mu^{k-1}_{J,\beta_J^*}(\tilde{Q})-\mu^{k-1}(\tilde{Q}) \right\} d\mu(y)=0 \end{equation}
 for all $v\in {\cal R}(d,J)$, and we will also show that it preserves the rate of $\mu^{k-1}_{J,\beta_J}$ in the sense that (the result known to hold for $\beta_J$)
\begin{equation}\label{bk1}
\pl \mu_{J,\beta_J^*}^{k-1}(\tilde{Q})-\mu^{k-1}(\tilde{Q}) \pl_{\infty}=O(r(d,J)^k).\end{equation}
The equations $H_v(\beta_J^*)=0$ for all $v\in {\cal R}(d,J)$ follow from fact that $\beta_J^*$ is the minimizer of $R^k_J(\beta)$, as shown in detail below. 
Viewing $\beta_J^*$ as a steepest descent update of $\beta_J$ provides an explicit proof that $\beta_J^*$ preserves the rate of $\beta_J$  so that (\ref{bk1}) holds. 
Our proof will also make clear that the results hold uniformly in $\tilde{Q}\in {\cal F}_M^{(0)}((0,1]^d)$.

{\bf Main proof:}
For now, assume that we have shown (\ref{Hbeta}) and (\ref{bk1}). 
Note that by (\ref{Hbeta})) we have for any vector $\alpha$ that
\[
\int \sum_{v\in {\cal R}(d,J)}\alpha(v)\tilde{\phi}_v(y) \left\{ \mu^{k-1}_{J,\beta_J^*}(\tilde{Q})-\mu^{k-1}(\tilde{Q}) \right\}
d\mu(y)=0.\]
Let $\alpha_x$ be the vector so that 
\[
\sup_x \pl \tilde{\phi}_x-\sum_{v\in {\cal R}(d,J) }\alpha_x(v)\tilde{\phi}_v\pl_{\mu}=O(r(d,J)).\]
 We can now  carry out the same proof as used for $k=1$:
\[
\begin{array}{l}
\mu_{J,\beta_J^*}^{k}(\tilde{Q}) -\mu^{k}(\tilde{Q}))(x)=\int_{(0,x]}\left\{ \mu_{J,\beta_J^*}^{k-1}(\tilde{Q})-\mu^{k-1}(\tilde{Q}) \right\} d\mu(y)\\
=\int \tilde{\phi}_x(y) \left\{ \mu_{J,\beta_J^*}^{k-1}(\tilde{Q})-\mu^{k-1}(\tilde{Q}) \right\} 
d\mu(y)\\
=\int (\tilde{\phi}_x(y)-\sum_{v\in {\cal R}(d,J)}\alpha_x(v)\tilde{\phi}_v(y) ) 
\left\{ \mu_{J,\beta_J^*}^{k-1}(\tilde{Q})-\mu^{k-1}(\tilde{Q}) \right\} 
d\mu(y)\\
\leq\pl \tilde{\phi}_x-\sum_{v\in {\cal R}(d,J)}\alpha_x(v)\tilde{\phi}_v\pl_{\mu}\pl 
 \mu_{J,\beta_J^*}^{k-1}(\tilde{Q})-\mu^{k-1}(\tilde{Q}) 
\pl_{\mu}\\
\leq\sup_x \pl \tilde{\phi}_x-\sum_{v\in {\cal R}(d,J)}\alpha_x(v)\tilde{\phi}_v\pl_{\mu}\pl  \mu_{J,\beta_J^*}^{k-1}(\tilde{Q})-\mu^{k-1}(\tilde{Q})
\pl_{\mu}\\
=O(r(d,J))O(r(d,J)^{k})\\
= O(r(d,J)^{k+1}),
\end{array}
\]
where the  latter bound $O(r(d,J)^{k+1})$ is constant in $x$. 
So this then proves that $\pl \mu_{J,\beta_J^*}^{k}(\tilde{Q}) -\mu^{k}(\tilde{Q}))\pl_{\infty}
=O(r(d,J)^{k+1})$. 


{\bf Why should the rate of $\mu_{J,\beta_J}^{k-1}(\tilde{Q})$ be preserved by the modification $\mu_{J,\beta_J^*}^{k-1}(\tilde{Q})$?}
Due to $\pl \mu_{J,\beta_J}^{k-1}(\tilde{Q})-\mu^{k-1}(\tilde{Q})\pl_{\infty}=O(r(d,J)^k)$, we already have that $\int_{(0,v]}(\mu_{J,\beta_J}^{k-1}(\tilde{Q})-\mu^{k-1}(\tilde{Q}))d\mu=O(r(d,J)^k)$ uniformly in $v\in {\cal R}(d,J)$, so that we should only have to induce a change in $\mu_{J,\beta_J}^{k-1}(\tilde{Q})$ of size  $O(r(d,J)^k)$ w.r.t. $L^1(\mu)$-norm  to make $\beta_J^*$ solve these equations $H_{\beta}(v)=0$ for all $v\in {\cal R}(d,J)$. As a result, we should still have 
$\pl \mu_{J,\beta_J^*}^{k-1}(\tilde{Q})-\mu^{k-1}(\tilde{Q})\pl_{\infty}=O(r(d,J)^k)$. 

{\bf The steepest descent algorithm definition of $\beta_J^*$:}
Let $\beta_J^1$ the the update of $\beta_J$ using a  steepest descent algorithm aiming to minimize over $\beta$
\[
R(\beta)= 0.5 J^{-1}\sum_{v\in {\cal R}(d,J)}\left(\mu_{J,\beta}^{k}(\tilde{Q})-\mu^{k}(\tilde{Q})\right)(v) ^2.\]  This steepest descent algorithm will find a local minimum around $\beta_J$ in which all partial derivatives are zero, which then shows that $H_{\beta_J^1}(v)=0$ for all $v\in {\cal R}(d,J)$.
Since $R(\beta)$ is strictly convex, it follows that  there is only one local minimum so that  $\beta_J^1=\beta_J^*$. Therefore, we will denote $\beta_J^1$ with $\beta_J^*$ although it now is defined as a steepest descent algorithm applied to $\beta_J$ (e.g., like the one-step Newton-Raphson algorithm update).


Note $R(\beta)=0.5 J^{-1}\sum_{v\in {\cal R}(d,J)}(\mu^{k}(\sum_{u\in {\cal R}(d,J)}\beta(u)\phi_u^0) -\mu^k(\tilde{Q}))^2(v)$.
In other words, $\beta_J^*$ minimizes  $\beta\rightarrow \pl \mu^{k}_{J,\beta}(\tilde{Q})-\mu^{k}(\tilde{Q}) \pl_{J}$, where we define the Eucliean type  norm $\pl Q\pl_J=\left(1/J\sum_{u\in {\cal R}(d,J)} Q(u)^2\right)^{1/2}$. We also define the corresponding inner product $\langle h_1,h_2\rangle_J=1/J\sum_{u\in {\cal R}(d,J)}h_1(u)h_2(u)$.

Recall $H_v(\beta)=(\mu^{k}_{J,\beta}(\tilde{Q})-\mu^{k}(\tilde{Q}))(v)$.  Thus \[
R(\beta)=\frac{0.5}{J}\sum_{v\in {\cal R}(d,J)}\{H_v(\beta)\}^2.\]
Recall that $\phi_u^l(v)=\int_{(0,v]} \phi_u^{l-1}(y)d\mu(y)$, $l=1,\ldots,k$.
Note, for $k\in \{2,3,\ldots\}$, we have
\[
H_v(\beta)=\int_{(0,v]}\int_{(0,y_1]}\ldots\int_{(0,y_{k-1}]}\left (\sum_{u\in {\cal R}(d,J)}\beta(u)\phi_u-\tilde{Q}\right )d\mu(y) \prod_{j=1}^{k-1}d\mu(y_l).\]
Consider a path $\beta+\delta h$ and let $\delta_0=0$. Then,
\[
\begin{array}{l}
\frac{d}{d\delta_0}H_v(\beta+\delta_0 h)(v)=\int_{(0,v]}\int_{(0,y_1]}\ldots\int_{(0,y_{k-1}]}\sum_{u\in {\cal R}(d,J)} h(u)\phi_u(y) d\mu(y) \prod_{j=1}^{k-1}d\mu(y_l)\\
=\sum_{u\in {\cal R}(d,J) } h(u) \int_{(0,v]}\int_{(0,y_1]}\ldots\int_{(0,y_{k-1}]}\phi_u(y) d\mu(y) \prod_{j=1}^{k-1}d\mu(y_l)\\
=\sum_{u\in {\cal R}(d,J)}h(u)\phi_u^k(v).
\end{array}
\]
 It follows that the pathwise derivative along paths $\beta+\delta h$ of  criterion $R(\beta)$ at $\delta=\delta_0=0$ is given by:
\begin{eqnarray*}
\frac{d}{d\delta_0}R(\beta+\delta_0 h)&=&
J^{-1}\sum_{v\in {\cal R}(d,J)}H_v(\beta) \sum_{u\in {\cal R}(d,J)}h(u)\phi_u^k(v)\\
&=&
J^{-1}\sum_{u\in {\cal R}(d,J)} h(u) \sum_{v\in {\cal R}(d,J)}H_v(\beta)\phi_u^k(v).\end{eqnarray*}
Therefore it follows that \[
D^*_{\beta}(u)=\sum_{v\in {\cal R}(d,J)}H_v(\beta)\phi_u^k(v)\] is the canonical gradient of this pathwise derivative of $R(\beta)$ at $\beta$ w.r.t. inner product  $\langle h_1,h_2\rangle_J=J^{-1}\sum_{u\in {\cal R}(d,J)} h_1(u)h_2(u)$.  That is, 
\[
\frac{d}{d\delta_0}R(\beta+\delta_0 h)=\langle D^*_{\beta},h\rangle_J .\]
Thus, the local steepest descent path $\{\beta_{J,\delta}^{0,lsd}:\delta\}$ for this criterion $R(\beta)$  through $\beta_J$ is given by \[
\beta_{J,\delta}^{0,lsd}=\beta_J+\delta D^*_{\beta_J},\]
where $\delta<0$ represents direction of descent.
The local steepest descent path  defines a universal steepest descent path
$\{\beta_{J,\epsilon}^{usd}:\epsilon\}$ through $\beta_J$
 recursively defined as: for $\epsilon>0$ 
\[
\beta_{J,\epsilon}^{usd}=\beta_J-\int_{(0,\epsilon]} D^*_{\beta_{J,x-}^{usd}}  dx.\]
Along this universal steepest descent path $\{\beta_{J,\epsilon}^{usd}:\epsilon\}$ we have for all $\epsilon\geq 0$ 
\[
\frac{d}{d\epsilon}R(\beta_{J,\epsilon}^{usd})=-\langle D^*_{\beta_{J,\epsilon}^{usd} },D^*_{\beta_{J,\epsilon}^{usd} } \rangle_J .\]
Let $\epsilon_J=\arg\min_{\epsilon}R(\beta_{J,\epsilon}^{usd})$ be the minimizer along this path. That defines our proposed modification $\beta_J^*\equiv \beta_{J,\epsilon_J}^{usd}$ of $\beta_J$.

{\bf Showing the $\beta_J^*$ solves $H_v(\beta_J^*)=0$ for all $v\in {\cal R}(d,J)$:}
Since $\frac{d}{d\epsilon_J}R(\beta_{J,\epsilon_J}^{usd})=0$, it follows that $\pl D^*_{\beta_J^* } \pl_J^2=0$ and thus
\[
D^*_{\beta_J^*}(u)=0\mbox{ for all $u\in {\cal R}(d,J)$} .\]
In other words
\[
\sum_{v\in {\cal R}(d,J)}H_v(\beta_J^*)\phi_u^k(v)=0\mbox{ for all $u\in {\cal R}(d,J)$.}\]
Thus, $J$-dimensional vector $H_v(\beta_J^*)$ is orthogonal to $J$-dimensional vector $\phi_u^k({\bf v})\equiv (\phi_u^k(v):v\in {\cal R}(J))$ for all $J$ vectors $\phi_u^k({\bf v})$, $u\in {\cal R}(d,J)$.
Due to independence of the vectors $\phi_u^k({\bf v})$ across $u\in {\cal R}(d,J)$, it follows that $H_v(\beta_J^*)=0$ for all $v\in {\cal R}(J)$ so that the desired constraints are  indeed solved by our proposed modification. We already explained that this implies that $\beta_J^*=\arg\min_{\beta} R(\beta)$, due to strict convexity of $R(\beta)$ that $\beta_J^*$ is a local minimum of $R(\beta)$ at which all partial derivatives are equal to zero. 
 
{\bf Showing that $\beta_J^*$ preserves rates of $\beta_J$ (\ref{bk1}):}
For establishing (\ref{bk1}) it suffices to bound the $L_1$-norm \[
\pl \beta_{J,\epsilon_J}^{usd}-\beta_J\pl_1=\pl \int_{(0,\epsilon_J]} D^*_{\beta_{J,x-}^{usd}}  dx\pl_1 \]
by $O(r(d,J)^k)$.
We note that 
\begin{eqnarray*} \pl \int_{(0,\epsilon_J]} D^*_{\beta_{J,x-}^{usd}}  dx\pl_1 &=&
\sum_{u\in {\cal R}(d,J)}\left | \int_{(0,\epsilon_J]} D^*_{\beta_{J,x-}^{usd}}(u) dx\right | \\
&\leq& \int_{(0,\epsilon_J]}\sum_{u\in {\cal R}(d,J)}\left | D^*_{\beta_{J,x-}^{usd}}(u)\right | dx\\
&=& \int_{(0,\epsilon_J]}\pl D^*_{\beta_{J,x-}^{usd}}\pl_1 dx.
\end{eqnarray*}
The Euclidean norm $\pl D^*_{\beta_{J,x-}^{usd}}\pl_J$ monotonically decreases as $x$ moves from $0$ to $\epsilon_J$. Moreover, $\max_u \mid D^*_{\beta_{J,x-}^{usd}}\mid$ approaches zero as $x$ goes from zero to $\epsilon_J$ and attains zero at $\epsilon_J$.  Therefore, we can bound the $L_1$-norm of  $ \pl D^*_{\beta_{J,x-}^{usd}}\pl_1$
by its  $L_1$-norm at $x=0$ given by  $\pl D^*_{\beta_J}\pl_1$, at most up till a constant factor. So
\[ \int_{(0,\epsilon_J]}\pl D^*_{\beta_{J,x-}^{usd}}\pl_1 dx\sim \epsilon_J \pl D^*_{\beta_J}\pl_1 .\]
Therefore, we can focus on bounding
\[
\pl \epsilon_J D^*_{\beta_J}\pl_1.\]

We have
\begin{eqnarray}
R(\beta_J)-R(\beta_{J,\epsilon_J}^{usd})&\approx& \epsilon_J\frac{d}{d\epsilon_0}
R(\beta_{J,\epsilon_0}^{usd})\nonumber \\
&=&\epsilon_J\langle D^*_{\beta_J},D^*_{\beta_J}\rangle_J,\label{preservek}
\end{eqnarray}
where $\epsilon_0=0$.



Firstly, we note that \begin{eqnarray*}
H_v(\beta_J)&=&\mu^{k}_{J,\beta_J}(\tilde{Q})(v) -\mu^{k}(\tilde{Q})(v)\\
&=&\int_{(0,v]}(\mu^{k-1}_{J,\beta_J}(\tilde{Q})-\mu^{k-1}(\tilde{Q}))(y)d\mu(y)\\
&\leq& \pl \mu^{k-1}_{J,\beta_J}(\tilde{Q})-\mu^{k-1}(\tilde{Q})\pl_{\infty}\mu(v)\\
&=&O(r(d,J)^{k}),\end{eqnarray*}
uniformly in $v$ (i.e. $\mu(v)\leq 1$).
This shows that $\max_{v\in {\cal R}(d,J)}\mid H_v(\beta_J)\mid =O(r(d,J)^{k})$ and $R(\beta_J)=O(r(d,J)^{2k})$.
Thus, the left-hand side of (\ref{preservek})  behaves as $R(\beta_J)=O(r(d,J)^{2k})$.
Thus, we have
\[
\epsilon_J=\frac{O(r(d,J)^{2k})}{\pl D^*_{\beta_J}\pl_J^2} .\]

Recall $D^*_{\beta_J}(u)=\sum_{v\in {\cal R}(d,J)}H_v(\beta_J)\phi_u^k(v)$, where $H_v(\beta_J)$ behaves as $\sim r(d,J)^k$, while $\phi_u^k(v)$ is uniformly bounded by $1$.  
Therefore, under the assumption that $H_v(\beta_J)\sim r(d,J)^k$, we do not only have that $\max_u\mid D^*_{\beta_J}(u)\mid =O(Jr(d,J)^k)$, but we can also state that $ D^*_{\beta_J} $ behaves as $\sim J r(d,J)^k$ (i.e., not smaller order). Specifically, we  first consider the case that
\[
C(J)\equiv \frac{J r(d,J)^k}{
\pl D^*_{\beta_J}\pl_J }=O(1).\]
 Then, we have
 \[
 \epsilon_J=\frac{O(r(d,J)^{2k})}{\pl D^*_{\beta_J}\pl_J^2} =O(J^{-2}).\]
 Combined with $\max_u\mid D^*_{\beta_J}(u)\mid =O(Jr(d,J)^k)$, this gives
 that 
 \[
 \begin{array}{l}
 \pl \epsilon_J D^*_{\beta_J}\pl_1=\epsilon_J\sum_u\mid D^*_{\beta_J}(u)\mid=
 O(J^{-2} J) \max_u\mid D^*_{\beta_J}(u)\mid \\
 =O(J^{-1} J r(d,J)^k) =O(r(d,J)^k).\end{array}
 \]
 This proves \[
 \pl \mu^{k-1}_{J,\beta_J^* }(\tilde{Q})  -\mu^{k-1}(\tilde{Q})\pl_{\infty}=O(r(d,J)^k).\]
 
 Suppose now that $C(J)$ converges to infinity. Then we can write
 $\epsilon_J D^*_{\beta_J}=\epsilon_{1,J} \tilde{D}^*_{\beta_J}$, where $\epsilon_{1,J}=\epsilon_J C(J)$ and $\tilde{D}^*_{\beta_J}=D^*_{\beta_J}/C(J)$. 
 Now we have
 $(Jr(d,J)^k)/\pl \tilde{D}^*_{\beta_J}\pl_J=O(1)$. Then,
  \begin{equation}\label{z1ak}
R(\beta_J)-R(\beta_J^*) \approx \epsilon_J \pl D^*_{\beta_J}\pl_J^2=\epsilon_{1,J} C(J)\pl \tilde{D}^*_{\beta_J}\pl_J^2 .
\end{equation}
We can now use that $1/\pl \tilde{D}^*_{\beta_J}\pl_J=O(J^{-1}r(d,J)^{-k})$.
Thus.
\begin{eqnarray*}
\epsilon_{1,J}&=&\frac{O(r(d,J)^{2k})}{C(J)\pl D^*_{\beta_J}\pl_J^2} \\
&=&\frac{O(r(d,J)^{2k}r(d,J)^{-2k}J^{-2}}{C(J)}\\
&=&O(C(J)^{-1} J^{-2}).
\end{eqnarray*}
Then, the $L_1$-norm of $\epsilon_{1,J}\tilde{D}^*_{\beta_J}$ is bounded as:
\begin{eqnarray*}
\pl \epsilon_{1,J}\tilde{D}^*_{\beta_J}\pl_1&=&O(C(J)^{-1}J^{-2}) C(J)^{-1}\pl D^*_{\beta_J}\pl_1\\
&=&O(C(J)^{-2}J^{-2})O(J^2 r(d,J)^k)\\
&=&O(C(J)^{-2}r(d,J)^k).
\end{eqnarray*}
So this shows that our earlier scenario with $C(J)=O(1)$  represents the conservative scenario in the sense  that the rate of $\pl \beta_J^*-\beta_J\pl_1$ is even smaller order than $O(r(d,J)^k)$  if $C(J)\rightarrow\infty$.


This completes the proof of  Theorem \ref{mainkthordercdf}. $\Box$
\subsection{Theorem only requiring (\ref{Rb}) for the knot-point set ${\cal R}(d,J)$.}
Here we will generalize the approximation error result of Theorem \ref{mainkthordercdf}, relying on ${\cal R}(d,J)$ satisfying both (\ref{Ra}) and (\ref{Rb}), to the same result but now only requiring ${\cal R}(d,J)$ to satisfy (\ref{Rb}). The proof is till a large degree a copy of the proof of Theorem \ref{mainkthordercdf} so that we only demonstrate the new part of the proof and borrow the other parts. Before we present the theorem and its proof, let's explain the main idea. 
The proof of Theorem \ref{mainkthordercdf} was relying on approximating the indicator $\tilde{\phi}_x$ with linear combinations of $\tilde{\phi}_v$, $v\in {\cal R}(d,J)$, thereby relying on condition (\ref{Ra}). However,  this was due to working with $\mu( G_J-G)$, where $G_J=\mu^{k-1}_{J,\beta}(\tilde{Q})$; $G=\mu^{k-1}(\tilde{Q})$ and $\mu(Q)=\int_{(0,x]} Qd\mu=\int \tilde{\phi}_x(y)Q(y)d\mu(y)$. However, if we work with $\bar{\mu}(G_J-G)$ with $\bar{\mu}(Q)(x)=\int_{(x,1]} fd\mu=\int \phi_x(y)Q(y)d\mu(y)$, then it becomes all about approximating $\phi_x$ with linear combinations of $\phi_v$, $v\in {\cal R}(d,M)$, which is then covered by condition (\ref{Rb}). This suggest then defining $R(\beta)$ as sum of squared residuals  of $\bar{\mu}\mu^{k-1}_{J,\beta}(\tilde{Q})-\bar{\mu}\mu^{k-1}(\tilde{Q})$ instead and defining $\beta_J^*$ as the minimizer of this new criterion. This then allows us to copy the proof by  induction of Theorem \ref{mainkthordercdf} to establish uniform convergence of $\bar{\mu}\mu^{k-1}_{J,\beta_J^*}(\tilde{Q})$ to $\bar{\mu}\mu^{k-1}_{J,\beta_J^*}(\tilde{Q})$ at rate $O(r(d,J)^{k+1})$, assuming we already established the rate $O(r(d,J)^k)$ for $\mu^{k-1}_{J,\beta_J}(\tilde{Q})-\mu^{k-1}(\tilde{Q})$ for a $\beta_J$. Finally,  by using that $\mu(G_J-G)(x)$ can be expressed in terms of $\bar{\mu}(G_J-G)$ at corners of $(0,x]$, we obtain the desired result for $\mu^k(G_J-G)$ as well, so that the proof by induction follows.

\begin{theorem}\label{mainkthordercdfnew}
Let ${\cal R}(d,J)$ be defined as in Corollary \ref{defknots}  but only satisfying (\ref{Rb}):
\[
\sup_{Q\in {\cal F}^0_M((0,1]^d)}\inf_{\alpha}\pl \sum_{v\in {\cal R}(d,J)}\alpha(v)\phi_v^0-Q\pl_{\mu}=O(C(M)r(d,J)). \]
Let $k\in \{1,2,\ldots\}$ be a given integer.
Let $\bar{\mu}(\tilde{Q})(x)=\int_{(x,1]}\tilde{Q}(y)d\mu(u)$. 
Recall definition of $\mu^k_{J,\beta}(\tilde{Q})=\mu^k(\tilde{Q}_{J,\beta})$ implied by plugging in the zero-order spline approximation $\tilde{Q}_{J,\beta}=\sum_{u\in {\cal R}(d,J)}\beta(u)\phi_u^0$ for $\tilde{Q}$ in $\mu^{k}(\tilde{Q})$.

Define the criterion:
\[
R^k_J(\beta)\equiv  J^{-1}\sum_{v\in {\cal R}(d,J)}\left( \bar{\mu}\mu^{k-1}_{J,\beta}(\tilde{Q})-\bar{\mu}\mu^{k-1}(\tilde{Q})\right)^2(v).\]
Let $\beta_J^{*}=\arg\min_{\beta}R^k_J(\beta)$ be the minimizer.
We have 
\[\sup_{\tilde{Q}\in {\cal F}^{(0)}_M((0,1]^d)} \pl \mu^{k}_{J,\beta_J^{*}}(\tilde{Q})-\mu^{k}(\tilde{Q})\pl_{\infty}=O(C(M) r(d,J)^{k+1}),\]

This proves that
\[
\sup_{Q\in {\cal F}^{(k)}_{M}((0,1]^d)}\inf_{\beta}\pl \sum_{u\in {\cal R}(d,J)}\beta(u)\phi_u^k-Q\pl_{\infty}=O(C(M) r(d,J)^{k+1}).\]

In addition, we also have for $k=2,\ldots$
\[
\sup_{\tilde{Q}\in {\cal F}^{(0)}_M((0,1]^d)} \pl \mu^{k-1}_{J,\beta_J^*}(\tilde{Q})-\mu^{k-1}(\tilde{Q})\pl_{\infty}=O(C(M) r(d,J)^k),\]
while for $k=1$, we have
\[\sup_{\tilde{Q}\in {\cal F}^{(0)}_M((0,1]^d)}\pl \tilde{Q}_{J,\beta_J^*(\tilde{Q})}-\tilde{Q}\pl_{\mu}=O(C(M)r(d,J)).\]
The latter implies  convergence in variation norm of $\mu^{k}_{J,\beta_J^{*}}(\tilde{Q})$: for $k=1,2,\ldots$
\[ 
\sup_{\tilde{Q}\in {\cal F}^{(0)}_M((0,1]^d)} \pl \mu^{k}_{J,\beta_J^*}(\tilde{Q})-\mu^{k}(\tilde{Q})\pl_{v}=O(C(M) r(d,J)^{k}).\]
\end{theorem}


\subsection{Proof of  Theorem \ref{mainkthordercdfnew} for general $k$.}
As in Theorem \ref{mainkthordercdf} we carry out the a proof by induction. 
Suppose that we have proven the theorem for $k=1$, which is shown in completely analogue manner as in the general proof for $k$ below, and is therefore omitted here. Let $k\in \{2,\ldots\}$ be given and we can assume that we already proved the result for $k-1$ so that there exists a $\beta_J$ so that $\pl \mu^{k-1}_{J,\beta_J}(\tilde{Q})-\mu^{k-1}(\tilde{Q})\pl_{\infty}=O(r(d,J)^k)$.
We now want to prove that $\pl \bar{\mu}\mu^{k-1}_{J,\beta_J^*} (\tilde{Q})-\bar{\mu}\mu^{k-1}(\tilde{Q})\pl_{\infty}=O(r(d,J)^{k+1})$. 


We will update $\beta_J$ into $\beta_J^*$ so that \begin{equation}\label{Hbetanew}
H_v(\beta_J^*)\equiv \int_{(v,1]} \left\{ \mu^{k-1}_{J,\beta_J^*}(\tilde{Q})-\mu^{k-1}(\tilde{Q}) \right\} d\mu(y)=0 \end{equation}
 for all $v\in {\cal R}(d,J)$, and in such a way that we still have
\begin{equation}\label{helpz1}
\pl \mu_{J,\beta_J^*}^{k-1}(\tilde{Q})-\mu^{k-1}(\tilde{Q}) \pl_{\infty}=O(r(d,J)^k).\end{equation}
Note that $\beta_J^*$ is minimizing $\sum_v \{H_v(\beta)\}^2$, and $\beta$ has same dimension as number of squared residuals $H_v(\beta)^2$ across $v\in {\cal R}(d,J)$. Therefore we should indeed have that the minimizer $\beta_J^*$  will satisfy $H_v(\beta_J^*)=0$. This is  shown in precisely the same way as in proof of Theorem  \ref{mainkthordercdf}. In addition, the fact that the rate of $\beta_J$ is preserved so that (\ref{helpz1}) holds is also proven in completely same manner as in the proof of Theorem \ref{mainkthordercdf}.
Note that
\[
  \int_{(v,1]} \left\{ \mu^{k-1}_{J,\beta_J^*}(\tilde{Q})-\mu^{k-1}(\tilde{Q}) \right\} d\mu(y)=
  \int \phi_v^0(y)\left\{ \mu^{k-1}_{J,\beta_J^*}(\tilde{Q})-\mu^{k-1}(\tilde{Q}) \right\} d\mu(y),\]
  where $\phi_v^0(y)=I(y\geq v)$. 
  Thus (\ref{Hbetanew}) implies that for any vector $\alpha$ 
\[
\int \sum_{v\in {\cal R}(d,J)}\alpha(v){\phi}_v(y) \left\{ \mu^{k-1}_{J,\beta_J^*}(\tilde{Q})-\mu^{k-1}(\tilde{Q}) \right\}
d\mu(y)=0.\]
Let $\alpha_x$ be the vector so that 
\[
\sup_x \pl {\phi}_x-\sum_{v\in {\cal R}(d,J) }\alpha_x(v){\phi}_v\pl_{\mu}=O(r(d,J)).\]
We know this choice $\alpha_x$ exists due to ${\cal R}(d,J)$ satisfying (\ref{Rb}) and noting that $\phi_x\in {\cal F}^{(0)}_M((0,1]^d)$, thereby avoiding the need for condition (\ref{Ra}).
 Then,
\[
\begin{array}{l}
\bar{\mu}\mu_{J,\beta_J^*}^{k-1}(\tilde{Q}) -\bar{\mu}\mu^{k-1}(\tilde{Q}))(x)=\int_{(x,1]}\left\{ \mu_{J,\beta_J^*}^{k-1}(\tilde{Q})-\mu^{k-1}(\tilde{Q}) \right\} d\mu(y)\\
=\int \phi_x(y) \left\{ \mu_{J,\beta_J^*}^{k-1}(\tilde{Q})-\mu^{k-1}(\tilde{Q}) \right\} 
d\mu(y)\\
=\int (\phi_x(y)-\sum_{v\in {\cal R}(d,J)}\alpha_x(v)\phi_v(y) ) 
\left\{ \mu_{J,\beta_J^*}^{k-1}(\tilde{Q})-\mu^{k-1}(\tilde{Q}) \right\} 
d\mu(y)\\
\leq\pl {\phi}_x-\sum_{v\in {\cal R}(d,J)}\alpha_x(v)\phi_v\pl_{\mu}\pl 
 \mu_{J,\beta_J^*}^{k-1}(\tilde{Q})-\mu^{k-1}(\tilde{Q}) 
\pl_{\mu}\\
\leq\sup_x \pl \phi_x-\sum_{v\in {\cal R}(d,J)}\alpha_x(v)\phi_v\pl_{\mu}\pl  \mu_{J,\beta_J^*}^{k-1}(\tilde{Q})-\mu^{k-1}(\tilde{Q})
\pl_{\mu}\\
=O(r(d,J))O(r(d,J)^{k})\\
= O(r(d,J)^{k+1}),
\end{array}
\]
where the  latter bound $O(r(d,J)^{k+1})$ is constant in $x$. 
So this then proves that $\pl \bar{\mu}\mu_{J,\beta_J^*}^{k-1}(\tilde{Q}) -\bar{\mu}\mu^{k-1}(\tilde{Q}))\pl_{\infty}
=O(r(d,J)^{k+1})$. 

We now need to prove that this implies the same rate for $\mu\mu_{J,\beta_J^*}^{k-1}(\tilde{Q})-\mu^k(\tilde{Q})$.
The cube $(0,x]$ has $2^d$-corners given by $\{(x(s),0(-s)):s\subset\{1,\ldots,d\}\}$.
Moreover, $Q(x)=\int_{(0,x]} dQ$ is a generalized difference over these corners of the function $\bar{F}(x)=\int_{(x,1]}dQ$. For example, consider the case that $d=2$. Then we have
\[
Q(x)=\bar{F}(x_1,x_2)-\bar{F}(x_1,0)-\bar{F}(0,x_2)+\bar{F}(0,0).\]
Therefore, we have
\begin{equation}\label{lhelpa}
\mu(Q)(x)=\int_{(0,x]} Q(y)d\mu(y)
=\sum_{s\subset \{1,\ldots,d\}} (-1)^{\mid s\mid} \bar{\mu}(Q)(x(s),0(-s)).\end{equation}
 Let $G_J\equiv \mu_{J,\beta_J^*}^{k-1}(\tilde{Q})$ and $G=\mu^{k-1}(\tilde{Q})$ where we know $\pl G_J-G\pl_{\infty}=O(r(d,J)^k)$. 
We can apply this result (\ref{lhelpa}) to express  $\mu(G_J)-\mu(G)$ in terms of $\bar{\mu}(G_J-G)$:
\begin{eqnarray*}
(\mu(G_J)-\mu(G))(x)&=&\mu(G_J-G)(x)\\
&=&\sum_{s\subset\{1,\ldots,d\}}(-1)^{\mid s\mid} \{\bar{\mu}(G_J)-\bar{\mu}(G)\}(x(s),0(-s))\\
&\leq &2^d\pl \bar{\mu}(G_J)-\bar{\mu}(G)\pl_{\infty}\\
&=&O (r(d,J)^{k+1}).
\end{eqnarray*}
Thus, this proves that also $\pl \mu(G_J)-\mu(G)\pl_{\infty}=O(r(d,J)^{k+1})$, and thus that
$\pl \mu\mu_{J,\beta_J^*}^{k-1}(\tilde{Q})-\mu\mu^{k-1}(\tilde{Q})\pl_{\infty}=O(r(d,J)^{k+1})$. This completes the proof of the theorem. 
$\Box$

\section{Submodels of $k$-th order smoothness class}\label{AppendixE}
\subsection{Particular submodel of interest assuming higher order derivatives are zero on the $0$-edges}
We want to highlight a particular type of submodel of $D^{(k)}([0,1]^d)$ of interest that heavily reduces the dimension in terms of family of basis functions. 
\nl 
{\bf Submodel of first order spline representation:}
For the sake of demonstration, let's first show this for the first order spline representation. 
Suppose that $Q^{(1)}_{s_1,s_2}=0$ for all strict subsets $s_2\subset s_1$, $\mid s_2\mid<\mid s_1\mid$. 
Then the first order spline representation of $Q\in D^{(1)}([0,1]^d)$ reduces to
\begin{eqnarray*}
Q(x)&=&\sum_{s\subset \{1,\ldots,d\}} \mu^1(\tilde{Q}^{(1)}_{s,s})(x(s))\\
&=&\sum_{s\subset\{1,\ldots,d\}}\int \phi^1_{u(s)}(x(s)) d\tilde{Q}^{(1)}_{s,s}(u(s))\\
&=&\sum_{s\subset\{1,\ldots,d\}}\int \phi^1_{u(s)}(x(s)) Q^{(1)}_s(du(s)),
\end{eqnarray*}
where
\[
\tilde{Q}^{(1)}_{s,s}(u)=\int_{(0(s),u(s)]}Q^{(1)}_{s,s}(dy(s))=\int_{(0(s),u(s)]}Q^{(1)}_{s}(dy(s)),\]
and recall $Q^{(1)}_s(y(s))=\frac{d}{dy(s)}Q(y(s),0(-s))$ is the Radon-Nikodym derivative w.r.t. Lebesgue of the section $y(s)\rightarrow Q(y(s),0(-s))$.

{\bf Submodel of $k$-th order spline representation, assuming first and higher order derivatives are zero or constant  on $0$-edges:}
More generally, if for all $j=1,\ldots,k$, we have that $Q^{(j)}_{\bar{s}(j),s_{j+1}}=0$ for strict subsets $s_{j+1}\subset s_j$, $\mid s_{j+1}\mid<\mid s_j\mid$, then the $k$-th order spline representation reduces to:
\begin{eqnarray*}
Q(x)&=&\sum_{s\subset\{1,\ldots,d\}}\mu^k(\tilde{Q}^{(k)}_{s,\ldots,s})(x(s))\\
&=&\sum_{s\subset\{1,\ldots,d\}}\int \phi_{u(s)}^k(x(s)) d\tilde{Q}^{(k)}_{s,\ldots,s}(u(s))\\
&=&\sum_{s\subset\{1,\ldots,d\}}\int \phi_{u(s)}^k(x(s)) Q^{(k)}_s(du(s)),\end{eqnarray*}
where 
$Q^{(k)}_s(y(s))=\frac{d^k}{dy(s)}Q(y(s),0(-s))$ is the $k$-th order Radon-Nikodym derivative of the section $y(s)\rightarrow Q(y(s),0(-s))$. 

This submodel of $D^{(k)}([0,1]^d)$ corresponds with assuming that the first and higher order (up till $k$-th order) derivatives of the $s$-specific section $y(s)\rightarrow Q(y(s),0(-s))$ are zero at $y(s_1)=0$ for strict subsets $s_1\subset s$, across all $s\subset\{1,\ldots,d\}$. We could slightly extend this model by assuming these derivatives are constant at $y(s_1)=0$ instead of zero. We suggest that this  might be an effective submodel of interest and might also do a good job approximating functions not satisfying that these derivatives are constant (or zero) on the zero edges.

{\bf Submodel of $k$-th order spline representation, assuming $m$-th and higher order derivatives are zero or constant on $0$-edges:}
We could weaken this model by only assuming that the $m_1$-th order derivatives are zero on the $0$-edges for all $m_1=m,m+1,\ldots,k$, The above submodel corresponds with $m=1$. We could start with using the full $m-1$-th order spline representation of $Q\in D^{(m)}([0,1]^d)$ given by $Q(x)=\sum_{\bar{s}(m)}\bar{\phi}_{\bar{s}(m)}(x)\int \phi_{u(s_m)}^{m-1}(x(s_m))Q^{(m-1})_{\bar{s}(m)}(du(s_m))$. 
In that expression we replace $Q^{(m-1)}_{\bar{s}(m)}(du(s_m))$ by $Q^{(m)}_{\bar{s}(m)}(u(s_m))d\mu(u(s_m))$. Then we use the $k-m$-th order spline representation of $Q^{(m)}_{\bar{s}(m)}(u(s_m))$ that now only considers $s_{m+1}=\ldots=s_{k+1}$ all equal to $s_m$. So this makes clear that this model starts with an arbitrary flexible model for $D^{(m)}([0,1]^d)$, but then starts making assumptions about the $s_{m+1}$-sections of the $\geq m$-th order derivatives $Q^{(m)}_{\bar{s}(m)}$ for non-empty $s_m$ being zero across all real subsets $s_{m+1}$ of $s_m$, and similarly for the sections of its derivative $Q^{(m+1)}_{\bar{s}(m+1)}$ with $s_{m+1}=s_m$, and  so on. Note that all the terms with $s_m$ the empty set in the $m$-th order representation are not changing. 
 Equivalently, one takes the $k$-th order representation of $D^{(k)}([0,1]^d)$  in terms of $\bar{\phi}_{\bar{s}(k+1)}\int \phi_{u(s_{k+1})}^k(x(s_{k+1}))dQ^{(k)}_{\bar{s}(k+1)}(u(s_{k+1}))$, and constraints $\bar{s}(k+1)$ to $\bar{s}(m)$ with $s_m=s_{m+1}=\ldots=s_{k+1}$ since the other terms are assumed to be zero. 
So then it is  clear that we just need to understand if these terms simplify due to this structure of $\bar{s}(k+1)$.  
We note that if $s_m=s_{m+1}=\ldots=s_{k+1}$, then 
\[
\bar{\phi}_{\bar{s}(k+1)}(x)=\prod_{j=1}^{m-1}\phi^j_0(x(s_j/s_{j+1})),\]
which thus equals $\bar{\phi}_{\bar{s}(m)}(x)$. 
Note that ${Q}^{(k)}_{\bar{s}(m),s(m+1:k+1)=s_m}(u(s_m))$ 
corresponds with a $k-m$-th order derivative  of $Q^{(m)}_{\bar{s}(m)}$, all w.r.t. $\mu_{s_m}$. 
We will denote this $k-m$-th order derivative of $Q^{(m)}_{\bar{s}(m)}$ with $Q^{(k)}_{\bar{s}(m)}$. 

Then the $k$-th order spline representation reduces to:
\begin{eqnarray*}
Q(x)&=&\sum_{\bar{s}(m)\subset{\cal S}^{m}(d)} \bar{\phi}_{\bar{s}(m) }(x(s_1/s_m))\mu^k(\tilde{Q}^{(k)}_{\bar{s}(m)})(x(s_m))\\
&=& \sum_{\bar{s}(m)\subset{\cal S}^{m}(d)}\bar{\phi}_{\bar{s}(m)}(x(s_1/s_m)) \int \phi^k_{u_k(s_m)}(x(s_m)) Q^{(k)}_{\bar{s}(m)}(du_k(s_m) ).
\end{eqnarray*}
We can still separate out the terms with $s_m$ being the empty set given by $\sum_{\bar{s}(m), s_m=\emptyset}\bar{\phi}_{\bar{s}(m)}Q^{(m)}_{\bar{s}(m)}(0(s_m))$, so that we have
\[
\begin{array}{l}
Q(x)=\sum_{\bar{s}(m), s_m=\emptyset}\bar{\phi}_{\bar{s}(m)}Q^{(m)}_{\bar{s}(m)}(0(s_m))\\
+
\sum_{\bar{s}(m)\subset{\cal S}^{m}(d),\mid s_m\mid>0}\bar{\phi}_{\bar{s}(m)}(x(s_1/s_m)) \int \phi^k_{u_k(s_m)}(x(s_m)) Q^{(k)}_{\bar{s}(m)}(du_k(s_m) ).
\end{array}
\]
It will be of interest to evaluate these different submodels for simulated data to understand the impact of misspecifying the $0$-edge behavior of first and higher order derivatives. Our sense is that higher order splines can make  a real difference in finite samples but that correctly modeling the  $0$-edge behavior of first and higher order derivatives will only become important for large sample sizes, suggesting that the above models could be highly practical.  This could result in practical recommendations, but either way, using a discrete super learner will allow us to let the data decide what submodel to choose.

\subsection{Models allowing for varying degrees of smoothness across different regions in domain and different subsets of its coordinates} In Section \ref{section2} we already clarified that our $k$-th order spline representation and corresponding sieve and HAL-MLEs immediately extend to $D^{(k)}(\{0,1\}^m\times [0,1]^d)$.
However, we can also consider many variations of spaces by assuming $Q_0(b,\cdot)\in D^{(k(b)})([0,1]^d)$ for specified $k(b)\in \{0,\ldots\}$  and apply the $k(b)$-th order spline representation to each $Q_0(b,\cdot)$ across all $b\in \{0,1\}^m$.
Of course, the rate of convergence is now dominated by the rate of convergence for the  lowest level of smoothness $k(b)$, but, in practice, this allows one to model different $b$-specific functions with different degree of smoothness and that might result in important finite sample gains, relative to having to use the same amoung of smoothness across all functions. Similarly,  we could also assume different degree of smoothness of $s$-specific sections $Q_{0,s}$ of $Q_0$, thereby allowing that the target function has various degree of smoothness in different subsets of its coordinates. 

 \section{Generalization of uniform  approximation error results to {\em submodels} $D^{(k)}({\cal R}^k(d))$ of $k$-th order smoothness class $D^{(k)}([0,1]^d)$.}\label{AppendixF}

In Section \ref{section3} we described a $k$-th order smoothness class $D^{(k)}([0,1]^d)$ as linear combinations of basis functions indexed by $\bar{s}(k+1)\in {\cal S}^{k+1}(d)$ and knot-points $u\in (0,1]^{\mid s(k+1)\mid }$, with $(\bar{s}(k+1),u)$ varying over a maximal index set ${\cal R}^k(d)$. We noted that any such function $Q\in D^{(k)}((0,1]^d)$ could be represented as
\[
Q=\sum_{\bar{s}(k+1)}\bar{\phi}_{\bar{s}(k+1)}\int_{u\in {\cal R}(\bar{s}(k+1))} \phi_u^k \tilde{Q}^{(k)}_{\bar{s}(k+1)}(du),\]
where ${\cal R}(\bar{s}(k+1))=(0,1]^{\mid s(k+1)\mid}$. 
Here $\int_{u\in {\cal R}(\bar{s}(k+1))} \phi_u^k \tilde{Q}^{(k)}_{\bar{s}(k+1)}(du)=\mu^k(\tilde{Q}^{(k)}_{\bar{s}(k+1)})$ as well. As a result of this representation as a linear span of basis functions identified by ${\cal R}^k(d)$  we also denote this model with $D^{(k)}({\cal R}^k(d))$. 

Suppose now we consider a submodel of this class by restricting $\bar{s}(k+1)$ and $u\in {\cal R}(\bar{s}(k+1)\subset (0,1]^{\mid s(k+1)\mid}$. Let this restricted set be denoted with ${\cal R}^{k,*}(d)$ and the unrestricted set with ${\cal R}^k(d)$.
This then defines a submodel of ${D}^{(k)}([0,1]^d)$ implied by this set ${\cal R}^{k,*}(d)$, which we could denote with ${D}^{(k)}_{{\cal R}^{k,*}(d)}([0,1]^d) $ or $D^{(k)}({\cal R}^{k,*}(d))$.

Our uniform approximation results trivially extend to restrictions on the set of possible  $\bar{s}(k+1)$ since it just means we need to approximate fewer $k$-th order primitives $\mu^k(\tilde{Q}^{(k)}_{\bar{s}(k+1)})$. Recall that our uniform approximation result in Section \ref{section4} was  just about approximating a function in $\tilde{Q}^{(k)}_{\bar{s}(k+1)}\in {\cal F}^{(0)}_M((0(s_{k+1}),1(s_{k+1})]$  across all $\bar{s}(k+1)$ with a linear combination of zero-order splines, and plugging that function in  $\mu^k(\tilde{Q}_{\bar{s}(k+1)}^{(k)})$.
Therefore, it remains to understand if our uniform approximation results for $k$-th order primitives $\mu^k(\tilde{Q})$ for  a function $\tilde{Q}\in {\cal F}^{(0)}_M((0,1]^d)$ generalize to functions  $\tilde{Q}  \in {\cal F}^{(0)}_M({\cal R})$ for a subset ${\cal R}\subset (0,1]^d$. So we need to generalize our main Theorem \ref{mainkthordercdf} in Appendix \ref{AppendixD} to approximations $\mu^k(\tilde{Q}_{J,\beta})$ of $\mu^k(\tilde{Q})$ with $\tilde{Q}\in {\cal F}^{(0)}_M({\cal R})$. 
We could keep it simple by just noting that ${\cal F}^{(0)}_M({\cal R})\subset{\cal F}^{(0)}_M((0,1]^d)$ so that if we use the same $\tilde{Q}_{J,\beta_J^*}$ as in Theorem \ref{mainkthordercdf}, then this will imply the same uniform approximation error rates $r(d,J)^{k+1}$ of $\mu^k(\tilde{Q}_{J,\beta_J^*})$ w.r.t $\mu^k(\tilde{Q})$. However, if ${\cal R}$ is a strict subset of $(0,1]^d$, then this choice $\tilde{Q}_{J,\beta_J^*}$ is using more knots than needed so that this will hurt the finite sample performance of the resulting sieve and HAL-MLEs.
Therefore, we want to establish these results for a sensible selection $\tilde{Q}_{J,\beta_J}$ that respects the known model ${\cal F}^{(0)}_M({\cal R})$.


Firstly, we need to generalize our requirements on ${\cal R}(d,J)$ as defined in Corollary (\ref{defknots}) that was specific to approximating ${\cal F}^{(0)}((0,1]^d)$ to approximation the smaller set ${\cal F}^{(0)}({\cal R})$. 
The following lemma generalizes these requirements and states the existence of such as set. 
\begin{lemma}\label{defknotsgen}
Let ${\cal R}\subset (0,1]^d$ be such that any $x\in {\cal R}$ is an interior point or a point on its boundary. 
We will refer to such a set as an almost everywhere open set. 
Let $\pl Q\pl_{\mu}=\left( \int_{(0,1]^d} Q^2 d\mu\right)^{1/2}$ be the $L^2(\mu)$-norm.
Let ${\cal F}^{(0)}({\cal R})=\{x\rightarrow \int_{u\in {\cal R}}\phi_u(x)dQ(u):\int_{u\in {\cal R}}\mid dQ(u)\mid <\infty\}$ and notice that this is a linear space spanned by $\{\phi_u:u\in {\cal R}\}$. 
There exists a set of $J$ knot-points ${\cal R}(d,J)\subset {\cal R}$ so that both
\begin{equation}\label{Ragenweak}
\sup_{x\in {\cal R}} \inf_{\alpha}\left |\left | \Pi_{{\cal F}^{(0)}({\cal R})} \left(  \tilde{\phi}_x-\sum_{v\in {\cal R}(d,J)}\alpha(v)\tilde{\phi}_v   \right)\right |\right |  _{\mu}=O(r(d,J)),
\end{equation} where $\Pi_{{\cal F}^{(0)}({\cal R}}$ is the projection operator on the linear space ${\cal F}^{(0)}({\cal R})$ within the Hilbert space $L^2(\mu)$,
and 
\begin{eqnarray}
\sup_{Q\in {\cal F}^0_M({\cal R})}\inf_{\alpha}\pl \sum_{v\in {\cal R}(d,J)}\alpha(v)\phi_v^0-Q\pl_{\mu}&=&O(C(M)r(d,J))\label{Rbgen} .
\end{eqnarray}
\end{lemma}

We can now state the generalization of Theorem \ref{mainkthordercdf}.
\begin{theorem}\label{mainkthordercdfgen}
Let ${\cal R}\subset (0,1]^d$ be an almost everywhere open set. Recall the definition 
${\cal F}^{(0)}_M({\cal R})=\left\{x\rightarrow \int_{u\in {\cal R}}\phi_u(x)dQ(u): \pl Q\pl_v<M\right\}$ and ${\cal F}^{(0)}({\cal R})$ as this same definition but with $M=\infty$.
Let $\tilde{Q}_{J,\beta}=\sum_{u\in {\cal R}(d,J)}\beta(u)\phi_u^0$ for  a knot-point set ${\cal R}(d,J)\subset{\cal R}$ satisfying (\ref{Ragenweak})  and (\ref{Rbgen}).

Let $k\in \{1,2,\ldots\}$ be a given integer.
Recall definition of $\mu^k_{J,\beta}(\tilde{Q})=\mu^k(\tilde{Q}_{J,\beta})$ implied by plugging in the zero-order spline approximation $\tilde{Q}_{J,\beta}=\sum_{u\in {\cal R}(d,J)}\beta(u)\phi_u^0$ for $\tilde{Q}$ in $\mu^{k}(\tilde{Q})$.

Define the criterion:
\[
R^k_J(\beta)\equiv  J^{-1}\sum_{v\in {\cal R}(d,J)}\left( \mu^{k}_{J,\beta}(\tilde{Q})-\mu^{k}(\tilde{Q})\right)^2(v).\]
Let $\beta_J^{*}=\arg\min_{\beta}R^k_J(\beta)$ be the minimizer, which can also be denoted with $\beta_J^*(\tilde{Q})$ to emphasize that it depends on $\tilde{Q}$.
We have 
\[\sup_{\tilde{Q}\in {\cal F}^{(0)}_M({\cal R} )} \pl \mu^{k}_{J,\beta_J^{*}(\tilde{Q})}(\tilde{Q})-\mu^{k}(\tilde{Q})\pl_{\infty}=O(r(d,J)^{k+1}).\]

This proves that
\[
\sup_{Q\in {\cal F}^{(k)}_{M}({\cal R})}\inf_{\beta}\pl \sum_{u\in {\cal R}(d,J)}\beta(u)\phi_u^k-Q\pl_{\infty}=O(r(d,J)^{k+1}).\]

In addition, we also have for $k=2,\ldots$
\[
\sup_{\tilde{Q}\in {\cal F}^{(0)}_M({\cal R} )} \pl \mu^{k-1}_{J,\beta_J^*(\tilde{Q})}(\tilde{Q})-\mu^{k-1}(\tilde{Q})\pl_{\infty}=O(r(d,J)^k),\]
while for $k=1$, we have
\[
\sup_{\tilde{Q}\in {\cal F}^{(0)}_M({\cal R})} \pl \tilde{Q}_{J,\beta_J^*(\tilde{Q})}-\tilde{Q}\pl_{\mu}=O(r(d,J)).\]
This shows also convergence in variation norm of $\mu^{k}_{J,\beta_J^{*}}(\tilde{Q})$:
for $k=1,\ldots$
\[ 
\sup_{\tilde{Q}\in {\cal F}^{(0)}_M({\cal R} )} \pl \mu^{k}_{J,\beta_J^*(\tilde{Q})}(\tilde{Q})-\mu^{k}(\tilde{Q})\pl_{v}=O(r(d,J)^{k}).\]
\end{theorem}



The following theorem restates this theorem specifically for $k=1$.

\begin{theorem}\label{1stordergen} Let ${\cal R}\subset (0,1]^d$ be an almost everywhere open set.  Recall the definition 
Let $\tilde{Q}_{J,\beta}=\sum_{u\in {\cal R}(d,J)}\beta(u)\phi_u^0$ for  a knot-point set ${\cal R}(d,J)\subset{\cal R}$ satisfying (\ref{Ragenweak}) and (\ref{Rbgen}) of Lemma \ref{defknotsgen}.

Define the criterion:
\[
R^1_J(\beta)\equiv  J^{-1}\sum_{v\in {\cal R}(d,J)}\left( \int_{(0,v]}\left(\sum_{u\in {\cal R}(d,J)}\beta(u)\phi_u^0-\tilde{Q}\right )d\mu \right)^2.\]
Note, 
\[
R^1_J(\beta)=J^{-1}\sum_{v\in {\cal R}(d,J)}\left(\mu^1_{J,\beta}(\tilde{Q})-\mu(\tilde{Q})\right)^2(v),\]
where $\mu^1_{J,\beta}(\tilde{Q})=\mu(\tilde{Q}_{J,\beta})$. 
 
Let $\beta_J^*=\arg\min_{\beta}R^1_J(\beta)$ be the minimizer, which could also be denoted with $\beta_J^*(\tilde{Q})$.
Then,\[
\sup_{\tilde{Q}\in {\cal F}_M^{(0)}({\cal R})} \pl \mu^1_{J,\beta_J^*(\tilde{Q})}(\tilde{Q})-\mu(\tilde{Q})\pl_{\infty}=O(r(d,J)^2).\]
In addition, 
\[
\sup_{\tilde{Q}\in {\cal F}_M^{(0)}({\cal R})} \pl \tilde{Q}_{J,\beta_J^*(\tilde{Q})}-\tilde{Q}\pl_{\mu}=O(r(d,J)).\]
This implies
\[
\sup_{\tilde{Q}\in {\cal F}_M^{(0)}({\cal R})} \pl \mu^1_{J,\beta_J^*(\tilde{Q})}(\tilde{Q})-\mu(\tilde{Q})\pl_{v}=O(r(d,J)).\]
\end{theorem}
{\bf Proof of Theorem \ref{1stordergen}:}
Let $\tilde{Q}\in {\cal F}^{(0)}_M({\cal R})$ be given.
Let $\beta_J$ be such that $\pl \tilde{Q}_{J,\beta}-\tilde{Q}\pl_{\mu}=O(r(d,J))$, which exists by definition of ${\cal R}(d,J)$.
We want to analyze   $\mu_{J,\beta_J^*}(\tilde{Q})=\int_{(0,x]}\tilde{Q}_{J,\beta_J^*} (y) d\mu(y)$ relative to $\mu(\tilde{Q})$. In this proof we represent $\beta_J^*$ as a steepest descent algorithm update of  $\beta_J$ applied to $R^1_J(\beta)$    satisfying \[
H_v(\beta_J^*)\equiv \int_{(0,v]}(\mu_{J,\beta_J^*}(\tilde{Q})-\mu(\tilde{Q}))d\mu=0\]
 for all $v\in {\cal R}(d,J)$, and  
\begin{equation}\label{k1agen}
\pl \tilde{Q}_{J,\beta_J^* }-\tilde{Q}\pl_{\mu}=O(r(d,J)).\end{equation} 
Note that $R^1_J(\beta)=1/J \sum_{u\in {\cal R}(d,J)}\{H_v(\beta)\}^2$ so that the equations $H_v(\beta_J^*)=0$ are equations implied by the derivative equations $d/d\beta R^1_J(\beta)=0$ at $\beta=\beta_J^*$ solved by a minimum $\beta_J^*$, as shown in detail below.   Viewing $\beta_J^*$ as an update of a steepest descent algorithm to minimize $R^1_J(\beta)$ allows us to explicitly show that $\beta_J^*$ preserves the rate of $\beta_J$ so that (\ref{k1a}) holds.
This is completely analogue to our proof of Theorem  in Appendix \ref{AppendixD}.

Assume that $H_v(\beta_J^*)=0$ for all $v\in {\cal R}(d,J)\subset{\cal R}$ and (\ref{k1a}) has been shown. So then, for all $v\in {\cal R}(d,J)$
\begin{eqnarray*}
0&=& H_v(\beta_J^*)\\
&=&\int \tilde{\phi}_v(y)(\tilde{Q}_{J,\beta_J^*}-\tilde{Q}) ) d\mu(y)\\
&=&\int \Pi\left (\tilde{\phi}_v\mid {\cal F}^{(0)}({\cal R}) \right)(y) (\tilde{Q}_{J,\beta_J^*}-\tilde{Q}) ) d\mu(y).
\end{eqnarray*}
Since $H_v(\beta_J^*)=0$ for all $v\in {\cal R}(d,J)$, for all vectors $\alpha$ we have
\[
0=\int \sum_{v\in {\cal R}(d,J)}\alpha(v)\tilde{\phi}_v(y) (\tilde{Q}_{J,\beta_J^*}-\tilde{Q}))d\mu(y).\]

 We can then carry out the  following proof. Let $x$ be given. For any vector $\alpha_x$ we have
\begin{eqnarray*}
(\mu_{J,\beta_J^*}(\tilde{Q})-\mu(\tilde{Q}))(x)&=&\int_{(0,x]}(\tilde{Q}_{J,\beta_J^*}-\tilde{Q})) d\mu(y)\\
&=&\int \tilde{\phi}_x(y) (\tilde{Q}_{J,\beta_J^*}-\tilde{Q})(y)d\mu(y)\\
&=&\int \Pi_{{\cal F}^{(0)}({\cal R})} \left (\tilde{\phi}_x \right) (y)(\tilde{Q}_{J,\beta_J^*}-\tilde{Q})(y)d\mu(y)\\
&=&\int \left( \tilde{\phi}_x(y)-\sum_{v\in {\cal R}(d,J)}\alpha_x(v)\tilde{\phi}_v(y) \right)    (\tilde{Q}_{J,\beta_J^*}-\tilde{Q})(y)       d\mu(y)\\
&=& \int \Pi_{{\cal F}^{(0)}({\cal R})}  \left( \tilde{\phi}_x-\sum_{v\in {\cal R}(d,J)}\alpha_x(v)\tilde{\phi}_v  \right)(y)
(\tilde{Q}_{J,\beta_J^*}-\tilde{Q})(y)       d\mu(y).
\end{eqnarray*}
The latter projection representation allows us to use (\ref{Ragenweak}). 
 Under assumption (\ref{Ragenweak}) we can bound it with Cauchy-Schwarz inequality by
\[
\leq  \pl \Pi_{{\cal F}^{(0)}({\cal R})} \left (\tilde{\phi}_x  -\sum_{v\in {\cal R}(d,J)}\alpha_x(v)\tilde{\phi}_v \right)  \pl_{\mu} \pl \tilde{Q}_{J,\beta_J^*} -\tilde{Q}\pl_{\mu}. \]

Therefore, under assumption (\ref{Ragenweak}) we have
\begin{eqnarray*}
(\mu_{J,\beta_J^*}(\tilde{Q})-\mu(\tilde{Q}))(x)&\leq&
\inf_{\alpha}  \pl \Pi_{{\cal F}^{(0)}({\cal R})} \left (\tilde{\phi}_x  -\sum_{v\in {\cal R}(d,J)}\alpha_x(v)\tilde{\phi}_v \right)  \pl_{\mu} 
 \pl \tilde{Q}_{J,\beta_J^*}-\tilde{Q}\pl_{\mu}\\
&=& O(r(d,J)^{2}),
\end{eqnarray*}
by the properties (\ref{Ragenweak}) and (\ref{Rbgen}) of ${\cal R}(d,J)$. 

This bound is uniformly in $x$ and uniformly in all $\tilde{Q}$ with $\pl\tilde{ f}\pl_{v}<1$. Thus, this proves $\pl \mu_{J,\beta_J^*}(\tilde{Q})-\mu(\tilde{Q})\pl_{\infty}=O(r(d,J)^2)$. To conclude, we have constructed a $\mu_{J,\beta_J^*}(\tilde{Q})$ so that uniformly in $\tilde{Q}$ with $\pl \tilde{Q}\pl_{v}<1$
1) $\pl \tilde{Q}_{J,\beta_J^*}-\tilde{Q} \pl_{\mu} =O(r(d,J))$; 2) $\pl \mu(\tilde{Q}_{J,\beta_J^*})-\mu(\tilde{Q})\pl_{\infty}=O(r(d,J)^{2})$. This then completes  the proof of the theorem. 

The remaining part of proof is identical to the proof of Theorem \ref{1storder}.

\subsection{Proof of  Theorem \ref{mainkthordercdfgen} for general $k$.}
We will carry out a proof by induction. 
Above we already proved the Theorem \ref{1stordergen} for $k=1$. 
Let $k\in \{2,\ldots\}$ be given and we can assume that we already proved the result for $k-1$ so that there exists a $\beta_J$ so that $\pl \mu^{k-1}_{J,\beta_J}(\tilde{Q})-\mu^{k-1}(\tilde{Q})\pl_{\infty}=O(r(d,J)^k)$. 
We now want to prove that $\pl \mu^{k}_{J,\beta_J^*} (\tilde{Q})-\mu^{k}(\tilde{Q})\pl_{\infty}=O(r(d,J)^{k+1})$. 


We will update $\beta_J$ into $\beta_J^*$ through a steepest descent algorithm that starts at $\beta_J$ so that \begin{equation}\label{Hbetagen}
H_v(\beta_J^*)\equiv \int_{(0,v]} \left\{ \mu^{k-1}_{J,\beta_J^*}(\tilde{Q})-\mu^{k-1}(\tilde{Q}) \right\} d\mu(y)=0 \end{equation}
 for all $v\in {\cal R}(d,J)$, and in such a way that we still have
\begin{equation}\label{bk1gen}
\pl \mu_{J,\beta_J^*}^{k-1}(\tilde{Q})-\mu^{k-1}(\tilde{Q}) \pl_{\infty}=O(r(d,J)^k).\end{equation}
The equations $H_v(\beta_J^*)=0$ for all $v\in {\cal R}(d,J)$ follow from fact that $\beta_J^*$ is the minimizer of $R^k_J(\beta)$, as shown in detail below. 
Viewing $\beta_J^*$ as a steepest descent update of $\beta_J$ provides an explicit proof that $\beta_J^*$ preserves the rate of $\beta_J$  so that (\ref{bk1}) holds. That part of proof is identical to the proof of Theorem \ref{mainkthordercdf}.

Assume that we have shown (\ref{Hbetagen}) and (\ref{bk1gen}). 
Note that by (\ref{Hbetagen})) we have for any vector $\alpha$ that
\[
\int \sum_{v\in {\cal R}(d,J)}\alpha(v)\tilde{\phi}_v(y) \left\{ \mu^{k-1}_{J,\beta_J^*}(\tilde{Q})-\mu^{k-1}(\tilde{Q}) \right\}
d\mu(y)=0.\]
Let $\alpha_x$ be the vector so that 
\[
\sup_x \pl \tilde{\phi}_x-\sum_{v\in {\cal R}(d,J) }\alpha_x(v)\tilde{\phi}_v\pl_{\mu}=O(r(d,J)).\]
We can now  carry out the same proof as used for $k=1$.
We note that $\left\{ \mu_{J,\beta_J^*}^{k-1}(\tilde{Q})-\mu^{k-1}(\tilde{Q}) \right\} $ is an element of ${\cal F}^{(k)}({\cal R})=\{\mu^k(\tilde{Q}): \tilde{Q}\in {\cal F}^{(0)}({\cal R})\}$.
Thus, we can replace  $(\tilde{\phi}_x(y)-\sum_{v\in {\cal R}(d,J)}\alpha_x(v)\tilde{\phi}_v(y) ) $ in the proof below by its projection onto ${\cal F}^{(k)}({\cal R})$. Since ${\cal F}^{(k)}({\cal R})$ is a subspace of ${\cal F}^{(0)}({\cal R})$ this projection is bounded by the projection onto the bigger subspace ${\cal F}^{(0)}({\cal R})$. Thus it then remains to bound the $L^2(\mu)$-norm of 
\[
\Pi_{{\cal F}^{(0)}({\cal R})} \left( \tilde{\phi}_x-\sum_{v\in {\cal R}(d,J)}\alpha_x(v)\tilde{\phi}_v\right) .\]
Specifically, we have
\[\begin{array}{l}\mu_{J,\beta_J^*}^{k}(\tilde{Q}) -\mu^{k}(\tilde{Q}))(x)=\int_{(0,x]}\left\{ \mu_{J,\beta_J^*}^{k-1}(\tilde{Q})-\mu^{k-1}(\tilde{Q}) \right\} d\mu(y)\\
=\int \tilde{\phi}_x(y) \left\{ \mu_{J,\beta_J^*}^{k-1}(\tilde{Q})-\mu^{k-1}(\tilde{Q}) \right\} d\mu(y)\\
=\int (\tilde{\phi}_x(y)-\sum_{v\in {\cal R}(d,J)}\alpha_x(v)\tilde{\phi}_v(y) ) \left\{ \mu_{J,\beta_J^*}^{k-1}(\tilde{Q})-\mu^{k-1}(\tilde{Q}) \right\} d\mu(y)\\
=\int \Pi_{{\cal F}^{(0)}({\cal R})} \left( 
 \tilde{\phi}_x-\sum_{v\in {\cal R}(d,J)}\alpha_x(v)\tilde{\phi}_v \right)(y) \left\{ \mu_{J,\beta_J^*}^{k-1}(  \tilde{Q})-\mu^{k-1}(\tilde{Q}) \right\} d\mu(y)\\
 \leq \left | \left | \Pi_{{\cal F}^{(0)}({\cal R})} \left( 
 \tilde{\phi}_x-\sum_{v\in {\cal R}(d,J)}\alpha_x(v)\tilde{\phi}_v \right) \right | \right |_{\mu}
 \left | \left | \mu_{J,\beta_J^*}^{k-1}(  \tilde{Q})-\mu^{k-1}(\tilde{Q}) \right | \right |_{\mu}.
 \end{array}
 \] 
 Taking the infimum over $x$ on the right-hand side yields the $O(r(d,J)^2)$ by assumptions (\ref{Ragenweak}) and (\ref{Rbgen}).

The remaining part of proof is identical to the proof of Theorem \ref{mainkthordercdf}.

\section{Sup-norm covering number for $k$-th order primitives and smoothness class $D^{(k)}_M([0,1]^d)$}\label{AppendixG}
Our approximation results for approximating a $k$-th order primitive function with a linear combination of $k$-th order splines translates into sup-norm covering number for the class ${\cal F}^{(k)}_M((0,1]^d)$ of $k$-th order primitives. This then also  implies the same covering number (up till constant factor) for the  $k$-th order smoothness class $D^{(k)}_M([0,1]^d)$.
Our main proof relies on the next assumption which we claim to hold, is intuitively expected to hold easily as discussed in some detail below, including sufficient conditions and a sketch of proof. This assumption allows us to bridge our uniform approximation error $O(r(d,J)^{k+1})$ of the (infinite set) working model $\{\mu^k(\tilde{Q}_{J,\beta}): \beta\}$ with $\tilde{Q}_{J,\beta}=\sum_{u\in {\cal R}(d,J)}\beta(u)\phi_u^0$ consisting  of linear combinations of $\{\phi_u^k:u\in {\cal R}(d,J)\}$ w.r.t. ${\cal F}^{(k)}((0,1]^d)$, based on Theorem \ref{mainkthordercdf}, towards a discrete finite set within this working model representing the $\epsilon$-net (with $\epsilon\sim r(d,J)^{k+1}$) for obtaining a covering number. 

\subsection{Assumption that allows bridging the uniform approximation error approximating ${\cal F}^{(k)}((0,1]^d)$ with $J$-dimensional working model $\{\mu^k(\tilde{Q}_{J,\beta}):\beta\}$ to same uniform approximation error for a finite discretization with $J^J$-elements}

\begin{assumption}\label{assumptionepsnet}
Let ${\cal R}(d,J)$ be given and $M<\infty$ and assume it satisfies conditions of Theorem \ref{mainkthordercdf}.
Consider the discrete set ${\cal E}_{J,d}\subset{\cal E}_J\equiv \{(x(u):u\in {\cal R}(d,J)):\max_u\mid x(u)\mid<M\}$ of size $\sim J^J$ consisting of $x$-vectors of dimension $J$ whose values $x(u)$ are multiples of $1/J$ for all $u\in {\cal R}(d,J)$. Define ${\cal B}_{J,d}\equiv \{\beta: (\tilde{Q}_{J,\beta}(v):v\in {\cal R}(d,J),\pl \beta\pl_1<M)\in {\cal E}_{J,d}\}$ as  the set of $\beta$-vectors for which $\tilde{Q}_{J,\beta}(v)$ is a multiple of $1/J$ for all $v\in {\cal R}(d,J)$.
 Note that each $\beta$-value in this set is uniquely determined by solving $\tilde{Q}_{J,\beta}(v)=e(v)$ for $v\in {\cal R}(d,J)$ for a given vector $e\in {\cal E}_{J,d}$, while excluding it if this solution has $\pl\beta\pl_1>M$,  so that also ${\cal B}_{J,d}$ has maximal $O(J^J)$-values (fewer due to removing ones with $L_1$-norm exceeding $M$). 
 Our assumption is that
 \begin{equation}\label{assumptionepsnet1}
\max_{\beta\in {\cal B}_J} \min_{\beta_J\in {\cal B}_{J,d}}\max_{v\in {\cal R}(d,J)}\mid  \mu^k(\tilde{Q}_{J,\beta_J})-\mu^k(\tilde{Q}_{J,\beta})\mid (v)=O(r(d,J)^{k+1}).\end{equation}
\end{assumption}
Let $\beta_J(\beta)\equiv \arg\min_{\beta_J\in {\cal B}_{J,d}}\max_{v\in {\cal R}(d,J)}\mid  \mu^k(\tilde{Q}_{J,\beta_J})-\mu^k(\tilde{Q}_{J,\beta})\mid (v)$ the $\beta\in {\cal B}_{J,d}$ giving this best approximation of $\mu^k(\tilde{Q}_{J,\beta})$.

\begin{lemma}
Consider the setting of Assumption \ref{assumptionepsnet}. 
The assumption (\ref{assumptionepsnet1}) implies
\begin{equation}\label{implicepsnet1}
\max_{\beta\in {\cal B}_J} \min_{\beta_J\in {\cal B}_{J,d}}\pl \mid  \mu^k(\tilde{Q}_{J,\beta_J})-\mu^k(\tilde{Q}_{J,\beta})\pl_{\infty}=O(r(d,J)^{k+1}).
\end{equation}
\end{lemma}
{\bf Proof:}
The proof of this last lemma follows from the analogue proof of Theorem \ref{mainkthordercdf} (see main outline proof), but instead of using that the derivative equations $H(\beta_J(\beta))(v)\equiv (\mu^k(\tilde{Q}_{J,\beta_J(\beta)})-\mu^k(\tilde{Q}_{J,\beta}))(v)=0$, $v\in {\cal R}(d,J)$, are solved exactly, allowing $\tilde{\phi}_x^0$ to be approximated by $\sum_{u\in {\cal R}(d,J)}\alpha_x(u)\tilde{\phi}_u^0$, they are now solved at $O(r(d,J)^{k+1})$ by (\ref{implicepsnet1}). Nonetheless, this still  allows us to  approximate $\tilde{\phi}_x(y)$ by linear combinations of $\{\tilde{\phi}_v:v\in {\cal R}(d,J)\}$, but adding back this non-zero term implied by the linear combination of $\{\tilde{\phi}_v:v\in {\cal R}(d,J)\}$  yields an extra remainder $O(r(d,J)^{k+1})$ (which thus does not affect the claimed sup-norm rate). 
$\Box$

\paragraph{Discussion of Assumption \ref{assumptionepsnet1}.}
Let  \[
\beta_J(\beta)\equiv \arg\min_{\beta_J\in {\cal B}_{J,d}}\sum_{v\in {\cal R}(d,J)}\left(\mu
^k(\tilde{Q}_{J,\beta_J})-\mu^k(\tilde{Q}_{J,\beta})\right)^2(v),\] which is the complete  analogue of  $\beta_J^*(\tilde{Q}_{J,\beta})$ (in short-hand notation, $\beta_J^*(\beta)$) defined in Theorem \ref{mainkthordercdf}, but restricting the minimizer to be in the discrete finite set ${\cal B}_{J,d}$ instead of the full continuous set ${\cal B}_J$. Of course, in setting of Theorem \ref{mainkthordercdf}$\beta_J^*(\tilde{Q}_{J,\beta})=\beta$ since $\beta$ is then already included in the set ${\cal B}_J$, so that the approximation error equals zero. 

The question is if with this discrete approximation $\beta_J(\beta)$ of $\beta_J^*(\beta)=\beta$ can we still  mostly copy the proof of Theorem \ref{mainkthordercdf}, which proved  that the supremum norm of $H(\beta_J^*(\beta))=\mu^k(\tilde{Q}_{J,\beta_J^*(\beta)})-\mu^k(\tilde{Q})$ is $O(r(d,J)^{k+1})$, uniformly in $\tilde{Q}\in {\cal F}^{(k)}_M((0,1]^d)$, except that we only need this result for the much easier case that $\tilde{Q}=\tilde{Q}_{J,\beta}$.  For example, consider the proof of Theorem \ref{1storder} for $k=1$ with $\mu^1(\tilde{Q}_{J,\beta})(x)=\int \tilde{\phi}_x(u)\tilde{Q}_{J,\beta}(y)d\mu(y)$. That proof relied on the values of this function  $H(\beta_J(\beta))(v)=0$ for $v\in {\cal R}(d,J)$. That allowed us to approximate $\tilde{\phi}_x$ in the expression $\mu^1(\tilde{Q}_{J,\beta_J^*})-\mu^1(\tilde{Q}_{J,\beta})$ at $x$ with a linear combination of $\tilde{\phi}_v$, $v\in {\cal R}(d,J)$, and with that our proof of Theorem \ref{1storder} could be completed. For our proof here we do not need to control the supremum norm of $H(\beta_J(\beta))$, but only the max-norm over ${\cal R}(d,J)$. Therefore, we only need to approximate $\tilde{\phi}_v$, $v\in {\cal R}(d,J)$, (instead of $\tilde{\phi}_x$) by linear combinations of  $\tilde{\phi}_{m_d(v)}$ for which $H(\beta_J(\beta))(m_d(v))\approx 0$ (since we do not have $H(\beta_J(\beta))(v)=0$ exactly).
Therefore, we could approximate $\tilde{\phi}_v$, $v\in {\cal R}(d,J)$, by a linear combination of $\tilde{\phi}_{m_d(v)}$, in $L^2(\mu)$-norm, where $m_d(v)\in (0,1]^d$ is an approximation of $v\in {\cal R}(d,J)$ chosen such that $(\mu^k(\tilde{Q}_{J,\beta_{J}(\beta)})-\mu^k(\tilde{Q}_{J,\beta})(m_v(d))=O(r(d,J)^{k+1})$. For example, one could select $m_d(v)$ so that $\tilde{\phi}_v-\tilde{\phi}_{m_d(v)}$ is minimal in $L^2(\mu)$-norm under the constraint that $H(\beta_J(\beta))(m_d(v))=O(r(d,J)^{k+1})$ uniformly in  $v\in {\cal R}(d,J)$.
We then need that for some $M<\infty$, $\max_{v\in {\cal R}(d,J)}\inf_{\alpha, \pl \alpha\pl_1<M}
\pl \tilde{\phi}_v-\sum_{u\in {\cal R}(d,J)}\alpha(u)\tilde{\phi}_{m_d(u)}\pl_{\mu}=O(r(d,J))$.
Assumption \ref{assumptionepsnet1} for $k=1$ then follows by our main part of proof of Theorem \ref{1storder}. This proof is then generalized to general $k$ just as we did that in proof of Theorem \ref{mainkthordercdf}, bounding the max-norm of $H(\beta_J(\beta))$ over ${\cal R}(d,J)$  by $O(r(d,J)^{k+1})$ instead of bounding its supremum norm. 
This proves the following lemma. 

\begin{lemma}\label{lemmasufficient}
Suppose that there exists a $m_d:{\cal R}(d,J)\rightarrow (0,1]^d$ so that 
1) $\max_{v\in {\cal R}(d,J)}\mid H(\beta_{J}(\beta))(m_d(v))\mid =O(r(d,J)^{k+1})$ and for some $M<\infty$ 2)
\[
\max_{v\in {\cal R}(d,J)}\inf_{\pl \alpha\pl_1<M}\pl \tilde{\phi}_v-\sum_{u\in {\cal R}(d,J)}\alpha(u)\tilde{\phi}_{m_d(u)}\pl_{\mu}=O(r(d,J)).\]
Then, Assumption (\ref{assumptionepsnet1}) holds.
\end{lemma}

We conjecture that the sufficient condition of this lemma holds, but, for now, we just state it as assumption. In fact, we conjecture that $m_d(v)=v$ applies, i.e., $\beta_J(\beta)$ already achieves the desired $O(r(d,J)^{k+1})$ approximation.

\subsection{Main lemma providing bound on sup-norm covering number for ${\cal F}^{(k)}_M((0,1]^d)$}
We will now state the lemma relying on (\ref{assumptionepsnet1}) to hold. 
 Subsequently, we will establish the proof of this assumption.
 \begin{lemma}\label{lemmasupnormcoveringnumber}
Let ${\cal F}^d={\cal F}^{(k)}_M((0,1]^d)$ be the set of $k$-th order primitives $\mu^k(\tilde{Q})$  of a function $\tilde{Q}\in {\cal F}^{(0)}_M((0,1]^d)$ with bounded  variation-norm $\pl\tilde{Q}\pl_{v}\leq M$.  Let $N_{\infty}(\epsilon,{\cal F}^d)$ be the number of spheres of size $\epsilon$ w.r.t. supremum norm that are needed to cover ${\cal F}^d$. Assume  Assumption (\ref{assumptionepsnet}) above, or, the sufficient condition of Lemma \ref{lemmasufficient}. Let $d_1=2(d-1)$.
Then,
\[
\log^{1/2}N_{\infty}(\epsilon,{\cal F}^d) \sim \epsilon^{-1/(2(k+1))}(-\log \epsilon)^{(d_1-1)/2}.\]
The corresponding entropy integral is  given by 
\[
J_{\infty}(\delta,{\cal F}^k_M((0,1]^d)=\int_{(0,\delta]} \epsilon^{-1/(2(k+1))}(-\log \epsilon)^{(d_1+1)/2}
d\epsilon.\]
Using integration by parts shows that this entropy integral behaves as (edge term dominates integral term that is of smaller order than original integral)
\[
J_{\infty}(\delta,{\cal F}^{(k)}_M((0,1]^d))  )\sim \delta^{(2k+1)/(2k+2)} (-\log \delta)^{(d_1+1)/2 -1)}.\]
This result implies that $D^{(k)}_M([0,1]^d)$ has the same covering number and entropy integral up till a constant: (noting $(d_1+1)/2-1=d-3/2$)
\[
J_{\infty}(\delta,D^{(k)}_M([0,1]^d))\sim \delta^{(2k+1)/(2k+2)} (-\log \delta)^{d-3/2}.\]
\end{lemma}
{\bf Proof:}
Consider a $k$-th order primitive  $\mu^k(\tilde{Q})$ of a function $\tilde{Q}$ with $\pl \tilde{Q}\pl_{v}<M$.
Theorem \ref{mainkthordercdf} shows that there is a set of  $J$ support points ${\cal R}(d,J)$ so that for any $\tilde{Q}\in {\cal F}^{(0)}_M((0,1]^d)$, there exists a $\beta_J^*(\tilde{Q})$ such that $\mu^{k}(\tilde{Q}_{J,\beta_J^*(\tilde{Q})})$ provides a sup-norm approximation of $\mu^{k}(\tilde{Q})$ of $O(r(d,J)^{k+1})$, where $\tilde{Q}_{J,\beta}=\sum_{u\in {\cal R}(d,J)}\beta(u)\phi_u^0$.  Moreover, this choice $\beta_J^*(\tilde{Q})$ is uniquely determined as the minimizer of $\sum_{v\in {\cal R}(d,J)}\left(\mu^k(\tilde{Q}_{J,\beta})-\mu^k(\tilde{Q})\right)^2(v)$, which corresponds with solving in $\beta$ the linear equations $\mu^k(\tilde{Q}_{J,\beta})(v)=\mu^k(\tilde{Q})(v)$ for $v\in {\cal R}(d,J)$. 

Consider the set ${\cal E}_{J,d}$ of vectors of dimension $J$ whose values are multiples of $1/J$. We note that this set has $O(J^J)$ values..
Consider now the set of values  ${\cal B}_{J,d}\equiv \{\beta: (\tilde{Q}_{J,\beta})(v):v\in {\cal R}(d,J))\in {\cal E}_{J,d},\pl\beta\pl_1<M\}$, i..e this is the set of $\beta$-vectors with $\pl \beta\pl_1<M$ for which $\tilde{Q}_{J,\beta}(v)$ is a multiple of $1/J$ for all $v\in {\cal R}(d,J)$. Let ${\cal B}_J=\{\beta:\pl\beta\pl_1<M\}$.
 Note that each $\beta$-value in this set is uniquely determined by solving $\tilde{Q}_{J,\beta}(v)=e(v)$ for $v\in {\cal R}(d,J)$ for a given vector $e\in {\cal E}_J$ that yields a solution $\beta(e)$ with $\pl \beta(e)\pl_1<M$, so that also ${\cal B}_{J,d}$ has $O(J^J)$-values. 
 Let ${\cal F}_{J}=\{\mu^k(\tilde{Q}_{J,\beta}):\beta\in {\cal B}_J\}\subset {\cal F}_d$, while
 ${\cal F}_{J,d}=\{\mu^k(\tilde{Q}_{J,\beta}:\beta\in {\cal B}_{J,d}\}\subset {\cal F}_{J}$.

 Theorem \ref{mainkthordercdf} proves that  ${\cal F}_J\subset{\cal F}_d$ 
 yields an $O(r(d,J)^{k+1})$ sup-norm approximation of ${\cal F}_d$, while Assumption \ref{assumptionepsnet} shows that ${\cal F}_{J,d}$ yields an $O(r(d,J)^{k+1})$ sup-norm approximation of ${\cal F}_J$. Thus, the set ${\cal F}_{J,d}$ with $O(J^J)$ elements yields an $O(r(d,J)^{k+1})$ sup-norm approximation of
 ${\cal F}_{d}$. 

 Let $\epsilon =r(d,J)^{k+1}$, and solve for $J=J(\epsilon)$. Then, the covering number behaves as $N_{\infty}(\epsilon)=J(\epsilon)^{J(\epsilon)}$ and thus
$\log N_{\infty}(\epsilon)=J(\epsilon)\log J(\epsilon)$.
We have that $r(d,J)=J^{-1}(\log J)^{d_1}$ for  $d_1=2(d-1)$.
So we need to solve $\epsilon=J^{-{k+1}}(\log J)^{(k+1)d_1}$. 
This yields $J(\epsilon)\sim \epsilon^{-1/(k+1)}(-\log \epsilon)^{d_1}$.
Thus, $\log N_{\infty}(\epsilon)=\epsilon^{-1/(k+1)}(-\log \epsilon)^{d_1}\{ (\log -\epsilon)+\log(-\log \epsilon)\}$ which behaves as $\log N_{\infty}(\epsilon)\sim \epsilon^{-1/(k+1)}(-\log \epsilon)^{d_1+1})$.
Thus $\log^{1/2}N_{\infty}(\epsilon) \sim \epsilon^{-1/(2(k+1))}(-\log \epsilon)^{(d_1+1)/2}$.
The entropy integral is thus given by 
\[
J_{\infty}(\delta,{\cal F}^k_M((0,1]^d))=\int_{(0,\delta]} \epsilon^{-1/(2(k+1))}(-\log \epsilon)^{(d_1+1)/2}
d\epsilon.\] Finally, we can use integration by parts to write it as  a sum of an edge term
\[
\left . \epsilon^{1-1/(2k+2)}(-\log \epsilon)^{(d_1+1)/2}\right |_0^{\delta}\] and integral 
$\int_{(0,\delta]} \epsilon^{1-1/(2k+2)}(-\log \epsilon)^{(d_1+1)/2 -1} d\epsilon$. We now note that the integral term is of smaller order than original term showing that the dominating term is the edge term and that is reported in the lemma. 
$\Box$

 \section{Proof of Theorem \ref{theoremoraclemle}}\label{AppendixI}  
For convenience, let's present the Theorem \ref{theoremoraclemle} here as well.
\begin{theorem}
Let $R_0(\beta)\equiv P_0L(Q_{J,\beta})-P_0L(Q_0)$; $\beta_J=\arg\min_{\beta}R_0(\beta)$; and let $D^*_{\beta}=(\frac{d}{d\beta(u)}R_0(\beta):u\in {\cal R}(d,{\bf J}))$ be the gradient of $R_0(\beta)$ at $\beta$. Note  $D^*_{\beta_J}=0$. Let $\pl \beta\pl_J^2=\sum_{u\in {\cal R}(d,{\bf J})}\beta(u)^2$.
Let $Q_{0,J}=Q_{J,\beta_{J}}$ with $\beta_{J}=\arg\min_{\beta}R_0(\beta)$ be the oracle MLE.
Let $\tilde{\beta}_J$ be so that $\tilde{Q}_{0,J}=Q_{J,\tilde{\beta}_J}$ with $\pl \tilde{Q}_{0,J}-Q_0\pl_{\infty}=O(C(M)r(d,J)^{k+1})$.  Note that under a weak regularity condition,  $R_0(\tilde{\beta}_J)=O(C(M)^2 r(d,J)^{2(k+1)})$.

{\bf Assumptions:}
Assume that $D^*_{\beta}=0$ implies $\beta=\beta_J$, and  (e.g., using that $d_0(Q_{0,J},Q_0)=O(C(M)^2r(d,J)^{2(k+1)})$ and Cauchy-Schwarz inequality)
\begin{equation}\label{z22}
\pl D^*_{\tilde{\beta}_J}\pl_J=O(J^{1/2} C(M)r(d,J)^{(k+1)}).
\end{equation}

Then, $\pl Q_{0,J}-Q_0\pl_{\infty}=O(C(M) r(d,J)^{k+1})$. 
More generally, 
\begin{equation}\label{uniformoraclemle}
\sup_{Q\in D_M^{(k)}([0,1]^d)}\inf_{\beta}\pl Q_{J,\beta}-Q\pl_{\infty}=O(C(M)r(d,J)^{k+1}).\end{equation}

\end{theorem} 
{\bf Proof:}
We have that $Q_{0,J}$ solves the vector-valued score equation $P_0 S_{J,Q_{0,J}}=0$ with $S_{J,f}=\left( \frac{d}{dQ}L(Q)(\phi_u):u\in {\cal R}(d,{\bf J}\right)$. This implies that indeed $D^*_{\beta_J}=0$. We want to show that $\pl Q_{0,J}-Q_0\pl_{\infty}=O(r(d,J)^{k+1})$, just as $\pl \tilde{Q}_{0,J}-Q_0\pl_{\infty}=O(r(d,J)^{k+1})$. It is immediate that $d_0(Q_{0,J},Q_0)=P_0L(Q_{0,J})-P_0L(Q_0)=O(r(d,J)^{2(k+1)})$, due to $d_0(Q_{0,J},Q_0)\leq d_0(\tilde{Q}_{0,J},Q_0)$ and that $\pl \tilde{Q}_{0,J}-Q_0\pl_{\infty}=O(r(d,J)^{k+1})$.  In order to establish the sup-norm result for this loss-based projection $Q_{0,J}$,
we will now follow the proof of Theorem \ref{mainkthordercdf} carried out in Appendix \ref{AppendixD}.

Let $\tilde{\beta}_J$ be so that $\tilde{Q}_{0,J}=Q_{J,\tilde{\beta}_J}$, while $\beta_J$ is such that $Q_{0,J}=Q_{J,\beta_J}$. 
Let $R_0(\beta)\equiv P_0L(Q_{J,\beta})-P_0L(Q_0)$.
Let $\frac{d}{d\delta_0}R_0(\beta+\delta_0 h)$, where $\delta_0=0$, be the pathwise derivative, and we can represent it as $\langle D^*_{\beta},h\rangle_J$, where $\langle h_1,h_2\rangle_J=\sum_j h_1(j)h_2(j)$ is the standard inner product. Note that $D^*_{\beta}$ is the standard gradient  whose components are the partial derivatives of $R_0(\beta)$. 
Consider the steepest descent path $\tilde{\beta}_{J,\delta}^{lsd}=\tilde{\beta}_J+\delta D^*_{\tilde{\beta}_J}$ through $\tilde{\beta}_J$, where we move in direction $\delta\leq 0$. 

This implies a universal steepest descent path 
$\tilde{\beta}_{J,\epsilon}^{usd}=\tilde{\beta_J}-\int_{(0,\epsilon]}D^*_{\tilde{\beta}_{J,x-}^{usd}} dx$. 
This universal steepest descent path  has the property that at any $\epsilon\geq 0$, we have
\[
\frac{d}{d\epsilon}\beta_{J,\epsilon}^{usd}=-D^*_{\tilde{\beta}_{J,\epsilon}^{usd}}.\]
Along this universal steepest descent path $\{\tilde{\beta}_{J,\epsilon}^{usd}:\epsilon\}$ we have for all $\epsilon\geq 0$ (direction of decreasing $R(\beta)$)
\[
\frac{d}{d\epsilon}R(\tilde{\beta}_{J,\epsilon}^{usd})=-\langle D^*_{\tilde{\beta}_{J,\epsilon}^{usd}},D^*_{\tilde{\beta}_{J,\epsilon}^{usd}}\rangle_J=-\pl D^*_{\tilde{\beta}_{J,\epsilon}^{usd}}\pl_J^2
 ,\]
 where $\pl \beta\pl_J^2=\sum_j \beta^2(j)$.

Let $\epsilon_J=\arg\min_{\epsilon}R_0(\tilde{\beta}_{J,\epsilon}^{usd})$. Let $\beta_J^*=\tilde{\beta}_{J,\epsilon_J}^{usd}$. We note that $\beta_J^*$ solves $D^*_{\beta_J^*}=0$ so that, by assumption, $\beta_J^*=\beta_J$.
For preservation of the sup-norm rate of convergence, $\pl Q_{J,\beta_J^*}-Q_0\pl_{\infty}=O(r(d,J)^{k+1})$, it suffices to show that $\pl \beta_J^*-\tilde{\beta}_J\pl_1=O(r(d,J)^{k+1})$. 
We note that 
\begin{eqnarray*}
\pl \beta_J^*-\tilde{\beta}_J\pl_1&=& \pl \int_{(0,\epsilon_J]} D^*_{\tilde{\beta}_{J,x-}^{usd}}  dx\pl_1 \\
&=&
\sum_{u\in {\cal R}(d,J)}\mid \int_{(0,\epsilon_J]} D^*_{\tilde{\beta}_{J,x-}^{usd}}(u) dx\mid\\
&\leq& \int_{(0,\epsilon_J]}\sum_{u\in {\cal R}(d,J)}\mid D^*_{\tilde{\beta}_{J,x-}^{usd}}(u)\mid dx\\
&=& \int_{(0,\epsilon_J]}\pl D^*_{\tilde{\beta}_{J,x-}^{usd}}\pl_1 dx.
\end{eqnarray*}
The Euclidean norm $\pl D^*_{\tilde{\beta}_{J,x-}^{usd}}\pl_J$ monotonically decreases as $x$ moves from $0$ to $\epsilon_J$. Moreover, $\max_u \mid D^*_{\tilde{\beta}_{J,x-}^{usd}}\mid$ approaches zero as $x$ goes from zero to $\epsilon_J$ and attains zero at $\epsilon_J$.  Therefore, we can bound the $L_1$-norm $\pl D^*_{\tilde{\beta}_{J,x-}^{usd}}\pl_1$  at $x$  by $\pl D^*_{\tilde{\beta}_J}\pl_1$ or a constant times the latter. So
\[ \int_{(0,\epsilon_J]}\pl D^*_{\tilde{\beta}_{J,x-}^{usd}}\pl_1 dx\sim \epsilon_J \pl D^*_{\tilde{\beta}_J}\pl_1 .\]
Therefore, it suffices to prove that 
\[
\pl \epsilon_J D^*_{\beta_J}\pl_1=O(r(d,J)^{k+1}).\]

Note that this can be bounded as $\epsilon_J J^{1/2} \pl D^*_{\tilde{\beta}_J}\pl_J$, which equals $\epsilon_J\pl D^*_{\tilde{\beta}_J}\pl_J^2 J^{1/2}/\pl D^*_{\tilde{\beta}_J}\pl_J$. 

We have $R_0(\tilde{\beta}_J)-R_0(\beta_J^*)\sim \epsilon_J \pl D^*_{\tilde{\beta}_J}\pl_J^2$. 
The left-hand side is $O(r(d,J)^{2(k+1)}$. 
Thus, this gives the bound:
\begin{equation}\label{z11} \epsilon_J \pl D^*_{\tilde{\beta}_J}\pl_1\sim r(d,J)^{2(k+1)} \frac{J^{1/2}}{\pl D^*_{\tilde{\beta}_J}\pl_J } .\end{equation}

One expects that 
\[\pl D^*_{\tilde{\beta}_J}\pl_J^2=O(J r(d,J)^{2(k+1)}),
\]
 due to rate of convergence of $R_0(\tilde{\beta}_J)=O(r(d,J)^{2(k+1)})$.  Assume this to be the case. In addition, let's assume that  the rate of $\pl D^*_{\tilde{\beta}_J}\pl_J$ is not faster  so that $C(J)\equiv \frac{J^{1/2} r(d,J)^{k+1}}{\pl D^*_{\tilde{\beta}_J}\pl_J}=O(1)$. Then, (\ref{z11}) gives the bound
\[
\epsilon_J \pl D^*_{\tilde{\beta}_J}\pl_1\sim r(d,J)^{k+1}.
\]

The scenario that $C(J)$ converges to infinity is a more conservative scenario in which the initial estimator solves the gradient at a faster level than above, making it more likely that an MLE update preserves the sup-norm rate of convergence.  Suppose now that $C(J)$ converges to infinity. Then we can write
 $\epsilon_J D^*_{\tilde{\beta}_J}=\epsilon_{1,J} \tilde{D}^*_{\tilde{\beta}_J}$, where $\epsilon_{1,J}=\epsilon_J C(J)$ and $\tilde{D}^*_{\tilde{\beta}_J}=D^*_{\tilde{\beta}_J}/C(J)$. 
 Now we have
 $(J^{1/2}r(d,J)^{k+1})/\pl \tilde{D}^*_{\tilde{\beta}_J}\pl_J=O(1)$. 

We have
  \begin{equation}\label{z1a}
R(\tilde{\beta}_J)-R(\beta_J) \approx \epsilon_J \pl D^*_{\tilde{\beta}_J}\pl_J^2=\epsilon_{1,J} C(J)\pl \tilde{D}^*_{\tilde{\beta}_J}\pl_J^2
\end{equation}
We can now use that $1/\pl \tilde{D}^*_{\beta_J}\pl_J=O(J^{-1/2}r(d,J)^{-(k+1)})$.
Thus.
\begin{eqnarray*}
\epsilon_{1,J}&=&\frac{O(r(d,J)^{2(k+1)})}{C(J)\pl \tilde{D}^*_{\tilde{\beta}_J}\pl_J^2} \\
&=&\frac{O(r(d,J)^{2(k+1)}r(d,J)^{-2(k+1)}J^{-1}}{C(J)}\\
&=&O(C(J)^{-1} J^{-1}).
\end{eqnarray*}
Then, the $L_1$-norm of $\epsilon_{1,J}\tilde{D}^*_{\tilde{\beta}_J}$ is bounded as:
\begin{eqnarray*}
\pl \epsilon_{1,J}\tilde{D}^*_{\tilde{\beta}_J}\pl_1&=&O(C(J)^{-1}J^{-1}) C(J)^{-1}\pl D^*_{\tilde{\beta}_J}\pl_1\\
&=&O(C(J)^{-2}J^{-1}) J^{1/2}\pl D^*_{\tilde{\beta}_J}\pl_J \\
&=&O(C(J)^{-2}J^{-1/2}) O(J^{1/2} r(d,J)^{k+1} )\\
&=&O(C(J)^{-2} r(d,J)^{k+1}),
\end{eqnarray*}
where we used that $\pl D^*_{\tilde{\beta}_J}\pl_J=O( J^{1/2} r(d,J)^{k+1})$, by (\ref{z22}). 
So this shows that our earlier scenario with $C(J)=O(1)$  represents the conservative scenario in the sense  that the rate of $\pl \tilde{\beta}_{J,\epsilon_J}^{usd}-\tilde{\beta}_J\pl_1$ is even smaller order than $O(r(d,J))$  if $C(J)\rightarrow\infty$.
This completes the proof of the theorem. 
$\Box$

\section{Score equations for HAL-MLE}
\label{AppendixJ}

\subsection{$k$-th order spline Highly Adaptive Lasso}
Suppose that we initially select rich enough non-informative index  set ${\cal R}(d,{\bf J}_{max})$. For example, one might set  ${\bf J}(\bar{s}(k+1))=J_{max}$ for a single size $J_{max}$ for the knot-point set ${\cal R}(\mid s_{k+1}\mid, J)$ for approximating $\mu^{(k)}(\tilde{Q}^{(k)}_{\bar{s}(k+1)})$ across $\bar{s}(k+1)$ with $s_{k+1}$ non-empty set. 
Specifically, as shown in Appendix \ref{AppendixB}, in the regression case, we can select ${\cal R}(\bar{s}_{k+1},J=n)=\{X_i(s_{k+1}):i=1,\ldots,n\}$ across all $\bar{s}(k+1)$ with $s_{k+1}$ non-empty, and set ${\cal R}(d,{\bf J}_{max})=\{(\bar{s}(k+1),u): \mid s_{k+1}\mid>0,u\in {\cal R}(\bar{s}(k+1),J=n)\}\cup\{\bar{s}(k+1):\mid s_{k+1}\mid =0\}$. This defines now an initial working model $D({\cal R}(d,{\bf J}_{max}))$. We then run the lasso and select the $L_1$-norm of the coefficients with cross-validation. This results now in selection of a random index set ${\cal R}_n\subset {\cal R}(d,{\bf J}_{max})$. We refer to the resulting estimator $Q_n=Q_n^k$ as the $k$-th order spline Highly Adaptive Lasso MLE. For the random set ${\cal R}_n$ we can define the oracle MLE $Q_{0,n}=Q_{0,{\cal R}_n}$. The random set ${\cal R}_n$ will correspond with a $J_n$ so that each $\mu^{k}(\tilde{Q}^{(k+1)}_{\bar{s}(k+1)})$ is fitted $\sim J_n$ knot-points (asymptotically).

\subsection{Undersmoothing Lepski's method  for selecting $L_1$-norm in HAL-MLE.}
 Our asymptotic normality theorems for the HAL-MLE suggest that some undersmoothing relative to the cross-validation selector of the $L_1$-norm might be needed: 1) to control asymptotic bias relative to standard error; 2) to solve the key score equation $\tilde{r}_n(x)$ at level $o_P((n/d_n)^{-1/2})$. Here we suggest a particular undersmoothing method of interest. 
 Firstly, we  compute the HAL-MLE $Q_{n,\lambda}(x)$ for an in increasing sequence of $L_1$-norms $\lambda$ starting at the cross-validation selector, with corresponding $\delta$-method based variance estimators $\sigma^2_{n,\lambda}(x)$ for the resulting working model $D^{(k)}({\cal R}_n(\lambda))$ where ${\cal R}_n(\lambda)$ is the set of non-zero coefficients for the fit based on $L_1$-norm $\lambda$. 
  We can then use (an undersmoothed version  of) Lepski's method to determine the $L_1$-norm where the plateau in $Q_{n,\lambda}(x)$ is reached (w.r.t. change in standard error) at which point the bias is negligible relative to the standard error, which is where the normal approximation should become valid as well. Specifically, one keeps increasing the $L_1$-norm till the change in neighboring values of $Q_{n,\lambda}(x)$ are equal to a constant like $1.96$ times the change in corresponding standard errors. This also corresponds with optimizing the upper or lower bound of the $0.95$-confidence interval $Q_{n,\lambda}(x)\pm 1.96 \sigma_{n,\lambda}(x)$: if $Q_{n,\lambda}(x)$ is increasing in $\lambda$, then we maximize the lower bound, and if it is decreasing, then we maximize the upper bound. 
  We refer to Chapter 25 in \citep{vanderLaan&Rose18} for a more detailed explanation of this method.

\subsection{Interesting undersmoothing method for the relax HAL-MLE}
A selector of the $L_1$-norm of interest for the relax HAL-MLE is the cross-validation selector $C_{n,cv}$ for the HAL-MLE. Clearly, this choice of $C_{n,cv}$ would exceed the cross-validation selector of the relax HAL-MLE itself, and with this choice we guarantee to solve all the score equations for its non-zero coefficients of the cross-validated HAL-MLE, $P_n \frac{d}{dQ_n}L(Q_n)(\phi_j)=0$,  exactly for all $j$. 
We suggest that this might be an effective undersmoothing method for the relax HAL-MLE, while computationally being trivial.

\subsection{Undersmoothed HAL-MLE solves unpenalized score equations for non-zero coefficients}
The following lemma establishes the rate at which the HAL-MLE solves the regular score equations $P_n S_{j_0}(Q_n)=P_n \frac{d}{dQ_n}L(Q_n)(\phi_{j_0})\approx 0$, $j_0\in {\cal R}_n$, which would have been solved exactly if $Q_n$ is replaced by the regular (non-penalized) MLE for the model $D^{(k)}({\cal R}_n)$ (i.e., for the relax HAL-MLE). In particular, it shows under what conditions this rate makes the score equations negligible relative to the rate of convergence of $(Q_n-Q_{0,n})$. Solving these score equations at the right rate uniformly in all $j_0\in {\cal R}_n$ also shows that it solves the linear span of these score equations $P_n\frac{d}{dQ_n}L(Q_n)(\sum_{j\in {\cal R}_n} \alpha(j)\phi_{j})\approx 0$ at the same rate, uniformly in $\pl \alpha\pl_1<M$ for any $M<\infty$. The latter range of score equations will then suffice to control the score equations uniformly in a fixed  ${\cal R}_{0,n}$ for large enough $L_1$-norm (the score equation $\tilde{r}_n(x)\approx 0$ of our asymptotic normality theorems).  We refer to Lemma \ref{lemmatildernx} for an alternative manner of expressing $P_n \frac{d}{dQ_n}L(Q_n)(\phi_{j_0})=(P_n-\tilde{P}_n)(S_j(Q_n)-\tilde{S}_j(Q_n))$, where $\tilde{S}_j(Q_n)$ is a best approximation of $S_j(Q_n)$ among scores exactly solved by HAL-MLE (i.e., $P_n \tilde{Q}_n(Q_n)=0$), which also directly shows that this score equation is easily solved at level $o_P((n/d_n)^{-1/2})$, typically without a need for undersmoothing. 

\begin{lemma}\label{lemmascoreeqnhal}
Let $\beta_n=\arg\min_{\pl \beta\pl_1\leq C_n}P_n L(\sum_{j\in {\cal R}_{n,max}}\beta(j)\phi_j)$ be the lasso estimator, so that the $k$-th order spline HAL-MLE over ${D}^{(k)}_{C_n}({\cal R}_{n,max})$ is given by $Q_n=\sum_{j\in {\cal R}_{n,max}} \beta_n(j)\phi_j$. 
Let ${\cal R}_{n}=\{j\in {\cal R}_{n,max}:\beta_n(j)\not =0\}$ be the set of indices of size $d_n$ for the non-zero coefficients. Let $\tilde{Q}_n=\arg\min_{Q\in {\cal D}^{(k)}({\cal R}_n)}P_n L(Q)$ be the unconstrained MLE over the linear model implied by these non-zero coefficients.
Let ${Q}_{0,n}=\arg\min_{Q\in D^{(k)}({\cal R}_{n})}P_0 L(Q)$.
Let $j^*=\arg\min_{j\in {\cal R}_n}\pl \phi_{j}\pl_{1,P_0}$; 
$\delta_{1,n}(j^*)\equiv \pl \phi_{j^*}\pl_{1,P_0}=P_0\mid \phi_{j^*}(x)\mid$. Assume 
 $P_0\left\{ \frac{d}{dQ_n}L(Q_n)-\frac{d}{d\tilde{Q}_n}L(\tilde{Q}_n)\right\}(\phi_{j^*})=O\left(\pl Q_n-\tilde{Q}_n\pl_{\infty} \delta_{1,n}(j^*)\right)$.

{\bf Score equations exactly solved:}
A set of score equations solved by $Q_n$, corresponding with paths $(1+\delta h(j))\beta_n(j)$,  is given by $0=P_n \frac{d}{dQ_n}L(Q_n)(\sum_{j\in {\cal R}_{n}}h(j)\beta_n(j)\phi_j)$ for all $h$ satisfying $\sum_{j\in {\cal R}_{n}}h(j)\mid \beta_n(j)\mid =0$. 
For $j_0\in {\cal R}_{n}$, a  ''regular'' score equation for $\beta(j_0)$ solved by the unconstrained MLE $\tilde{Q}_n$ over ${\cal D}^{(k)}({\cal R}_{n})$ is given by $P_n \frac{d}{d\tilde{Q}_n}L(\tilde{Q}_n)(\phi_{j_0})=0$, $j_0\in {\cal R}_{n}$. 

{\bf Regular score equations:} We have the following bound for  $P_n \frac{d}{dQ_n}L(Q_n)(\phi_{j_0})$: for all $j_0\in {\cal R}_{n}$
\[
\begin{array}{l}
\pl P_n \frac{d}{dQ_n}L(Q_n)(\phi_{j_0})\pl \\
=
\beta_n(j_0)^{-1}\pl P_n\left\{ \frac{d}{dQ_n}L(Q_n)-\frac{d}{d\tilde{Q}_n}L(\tilde{Q}_n)\right\}(\phi_{j^*})\pl \\
\leq \beta_n(j_0)^{-1}\pl (P_n-P_0)\left\{ \frac{d}{dQ_n}L(Q_n)-\frac{d}{d\tilde{Q}_n}L(\tilde{Q}_n)\right\}(\phi_{j^*})\pl \\
+\beta_n(j_0)^{-1}\pl P_0 \left\{ \frac{d}{dQ_n}L(Q_n)-\frac{d}{d\tilde{Q}_n}L(\tilde{Q}_n)\right\}(\phi_{j^*})\pl \\
\sim \beta_n(j_0)^{-1} d_n n^{-1}\pl \phi_{j^*}\pl_{P_0}+
 \beta_n(j_0)^{-1}d_n^{1/2}n^{-1/2}\log n \pl \phi_{j^*}\pl_{1,P_0}.\end{array}
\]
Assume $d_n^{1/2}n^{-1/2}\log^{-1}n \pl \phi_{j^*}\pl_{P_0}/\pl \phi_{j^*}\pl_{1,P_0}=O(1)$ so that the second term dominates. For example, if $\pl \phi_{j^*}\pl_{P_0}/\pl \phi_{j^*}\pl_{1,P_0}\sim \pl \phi_{j^*}\pl_{1,P_0}^{1/2}$, then, this holds if 
\[
\frac{d_n n^{-1}\log^{-2}n}{
\pl \phi_{j^*}\pl_{1,P_0}}=o(1).\]
Then, a sufficient condition for this bound on the score equation $P_n \frac{d}{dQ_n}L(Q_n)(\phi_{j_0})$ to be of smaller order than $(d_n/n)^{-1/2}\log n$ (i.e., the rate of $(Q_n-Q_{0,n})(x)$) is given by 
\[
\min_{j\in {\cal R}_n}\pl \phi_{j}\pl_{1,P_0}=O(\beta_n(j_0)n^{-\delta})\mbox{ for some $\delta>0$}.\]
\end{lemma}
{\bf Proof:}
One selects a path $(1+\delta h(j))\beta_n(j)$ so that $h(j)=0$ for $j\not\in \{j_0,j^*\}$, $h(j_0)=1$ and solve for $h(j^*)$ so that $\sum_{j\in {\cal R}_n}h(j)\mid \beta_n(j)\mid =0$.
This results in the first equality.
The second inequality is the triangle inequality. The empirical process term can be viewed as an empirical process $n^{1/2}(P_n-P_0)f_n$ indexed by a $d_n$-dimensional class $f_n\in {\cal F}_n$ of functions, while we know that $\pl Q_n-\tilde{Q}_n\pl_{\infty}=O_P((d_n/n)^{1/2}\log n)$, giving a bound on the $L^2(P_0)$-norm of $f_n$. Bounding the empirical process term with the entropy integral of this finite dimensional class then yields the bound $d_n n^{-1}\pl \phi_{j^*}\pl_{P_0}$. The second term is bounded by using the bounding assumption and the rate of $Q_n-\tilde{Q}_n$. The remaining statements are straightforward. $\Box$


\paragraph{Discussion of assumptions:}
For example, if $\beta_n(j_0)^{-1}\sim d_n$, then, the last displayed  condition corresponds with $\min_{j}\pl \phi_{j}\pl_{1,P_0}=o(d_n^{-1})$ by some polynomial power.
The reason that $\phi_{j^*}$ is small is that its support is small, so $\phi_{j^*}$ is analogue to an indicator of a shrinking set $I(X\geq u)$ for a knot-point $u$: in fact, for the zero order splines, it is exactly an indicator. Therefore, one expects that $\int (\phi_{j^*})^2 d\mu \sim \int \mid \phi_{j^*}\mid d\mu$ so that  $\pl \phi_{j^*}\pl_{P_0}/\pl \phi_{j^*}\pl_{1,P_0}\sim \pl \phi_{j^*}\pl_{1,P_0}^{1/2}$, which is the sufficient assumption in the above lemma for making the second term the dominating term.

\section{Proofs for asymptotic normality  of $k$-th order spline sieve and (relax) HAL-MLE of regression function based on unweighted least squares}\label{AppendixK}

In Section \ref{section8} we established asymptotic normality of sieve and HAL-MLEs based on the weighted least squares regression. In this section we study unweighted least squares. Firstly, we show the analogue analysis which just ends up with a non-closed form asymptotic variance due to not being able to rely on the log-likelihood behavior of the loss function (see Section \ref{section10}). Subsequently, we show that we can still obtain the desired simple asymptotic variance expression under a stronger uniform approximation assumption (still implied by the nonparametric uniform approximation assumption on $D^{(k)}({\cal R}_{0,n})$).

\subsection{Asymptotic normality, when parametrizing regression working model in terms of orthonormal basis in $L^2(P_0)$}


Let ${\cal  R}_{0,n}$ be a fixed set approximating  ${\cal R}_n$ in terms of score equations, as defined formally  below, so that $r_n(j)=P_n \frac{d}{dQ_n}L(Q_n)(\phi_j)\approx 0$, $j\in {\cal R}_{0,n}$, is approximately solved. Recall $Q_{0,n}=\arg\min_{Q\in D^{(k)}({\cal R}_{0,n})}P_0 L(Q)$.
 Let $\{\phi_j^*:j\in {\cal R}_{0,n}\}$ be an orthonormal basis for $\{\phi_j:j\in {\cal R}_{0,n}\}$ in $L^2(P_0)$. This does not change the working model $D^{(k)}({\cal R}_{0,n})$.
 
 Solving the score equations $r_n(j)=P_n \frac{d}{dQ_n}L(Q_n)(\phi_j)\approx 0$, $j\in {\cal R}_{0,n}$, implies  solving 
 $r_n^*(j)\equiv P_n \frac{d}{dQ_n}L(Q_n)(\phi_j^*)\approx 0$, $j\in {\cal R}_{0,n}$.
Define the undersmoothing term:
\begin{equation}\label{tildernxls}
\tilde{r}_n(x)\equiv \sum_{j\in {\cal R}_{0,n}}r_n^*(j)\phi_j^*(x).\end{equation}

As starting point for our analysis of $(Q_n-Q_{0,n})(x)$, using that $P_0 \phi_j^* (Y-Q_{0,n})=0$ and $P_n \phi_j^*(Y-Q_n)=r_n^*(j)$,  we have
\[
P_0 \phi_j^*(Y-{Q_n})-P_0\phi_j^*(Y-{Q_{0,n}})=-(P_n-P_0)\phi_j^*(Y-{Q_n})+r_n^*(j),\]
so that
\[
P_0 \phi_j^*({Q_n}-{Q_{0,n}} )=(P_n-P_0)\phi_j^*(Y-{Q_n})-r_n^*(j).\]
We decompose $Q_n=\Pi(Q_n\mid D^{(k)}({\cal R}_{0,n}))+\{Q_n-\Pi(Q_n\mid D^{(k)}({\cal R}_{0,n}))\}$. 
Let $\Pi_{J_{0,n}}(Q_n)$ be short-hand notation for this projection onto $D^{(k)}({\cal R}_{0,n})$ in $L^2(P_0)$. This creates a term with $(Q_n-\Pi_{J_{0,n}}(Q_n))$.
By our uniform approximation result we know that there exists a $\tilde{Q}_n\in D^{(k)}({\cal R}_{0,n})$ so that $\pl Q_n-\tilde{Q}_n\pl_{\infty}=O_P(C(M)r(d,J_{0,n})^{k+1})$.
We have that $\Pi_{J_{0,n}}(Q_n)=\arg\min_{Q\in D^{(k)}({\cal R}_{0,n})}P_0(Q_n-Q)^2$. We can apply Theorem \ref{theoremoraclemle} with $R_0(\beta)=P_0(Q_{J_{0,n},\beta}-Q_n)^2$  and $\tilde{\beta}_n$ is so that $\tilde{Q}_n=Q_{J_{0,n},\tilde{\beta}_n}$.  The conditions are clearly satisfied as demonstrated with the least squares loss function under Theorem \ref{theoremoraclemle}. Assuming that $\pl Q_n\pl_{v,k}^*<M$,  this proves that $\pl \Pi_{J_{0,n}}(Q_n)-Q_n\pl_{\infty}=O_P(C(M)r(d,J_{0,n})^{k+1})$. 

Since for any $Q\in D^{(k)}({\cal R}_{0,n})$, $Q(x)=\sum_{j\in {\cal R}_{0,n}}\{ P_0 Q\phi_j^*\} \phi_j^*(x)$, we have
\begin{eqnarray*}
(\Pi(Q_n)-Q_{0,n})(x)&=& \sum_{j\in {\cal R}_{0,n}}\{P_0\phi_j^*(\Pi(Q_n)-{Q_{0,n}}) \}\phi_j^*(x)\\
&=&\sum_{j\in {\cal R}_{0,n}}(P_n-P_0)\phi_j^*(Y-{Q_n}) \phi_j^*(x)-\sum_{j\in {\cal R}_{0,n}}r_n^*(j) \phi_j^*(x)\\
&&-\sum_{j\in {\cal R}_{0,n}}P_0\phi_j^*\{(Q_n-\Pi_{J_{0,n}}(Q_n)\}\phi_j^*(x)\\
&=&\sum_{j\in {\cal R}_{0,n}}(P_n-P_0)\phi_j^*(Y-{Q_n}) \phi_j^*(x)-
\tilde{r}_n(x),
\end{eqnarray*}
where we use that the last term  before the last equality equals $\Pi_{J_{0,n}}(Q_n-\Pi_{J_{0,n}}(Q_n))=0$. 
Replacing the left-hand side by $(Q_n-Q_{0,n})(x)$, using our sup-norm bound on $Q_n-\Pi_{J_{0,n}}(Q_n)$, yields
\[
(Q_n-Q_{0,n})(x) 
=\sum_{j\in {\cal R}_{0,n}}(P_n-P_0)\phi_j^*(Y-{Q_n}) \phi_j^*(x)-
\tilde{r}_n(x)+O_P(C_n r(d,J_{0,n})^{k+1}),\]
where $C_n=\pl \beta_n\pl_1$ is the $k$-th order sectional variation norm of $Q_n=\sum_{u\in {\cal R}_n}\beta_n(u)\phi_u$.

Replacing $Q_n$ by $Q_{0,n}$ makes the leading term a sum of independent mean zero random variables
\begin{equation}\label{D0nxls}
D_{Q_{0,n},x}(O)\equiv \sum_{j\in {\cal R}_{0,n}}\phi_j^*(Y-{Q_{0,n}}) \phi_j^*(x).
\end{equation}
By CLT this empirical mean converges to a normal limit distribution at rate $\sim (n/J_{0,n})^{-1/2}\tilde{\sigma}_{0,n})$ with $\tilde{\sigma}_{0,n}^2$ is the variance of $D_{Q_{0,n},x}$. We will then have to assume that the undersmoothing $\tilde{r}_n(x)$ and $O(C_n r(d,J_{0,n})^{k+1})$  is of smaller order than this rate. 


\subsection{Asymptotic normality theorem for unweighted least squares}
This proves  the following asymptotic normality theorem.
Let 
\begin{equation}\label{Enxls}
E_n(x)\equiv \sum_{j\in {\cal R}_{0,n}}(P_n-P_0)\phi_j^*(Q_n-Q_{0,n}) \phi_j^*(x).
\end{equation}
We need the same assumptions as in our Theorem \ref{thasnormal} for the weighted least squares regression, with this definition of $E_n(x)$, and we avoided having another remainder $R_{n,1}(x)$ of Assumption A2 (\ref{A2}).


\begin{theorem}\label{thasnormalB}\ \nl
Consider the $k$-th order sieve MLE, HAL-MLE or relax HAL-MLE $Q_n$.  Let $k^*=k+1$.
Let ${\cal R}_n$ be the set of $d_n$ non-zero coefficients in $Q_n=\sum_{j\in {\cal R}_n}\beta_n(j)\phi_j$. Assume $\inf_{x\in [0,1]^d}\sigma^{2}_{0,n}(x)>0$. 
Recall the definitions of $\tilde{r}_n(x)$ (\ref{tildernxls}, $E_n(x)$ (\ref{Enxls}), $D_{Q_{0,n},x}(O)$ (\ref{D0nxls}) above. 
Assume the Assumptions (\ref{A0}), (\ref{undersmoothm}), (\ref{keyconditiona1m}), (\ref{neglbias})  of Theorem \ref{thasnormal}.

Due to Assumption A0, we have the following expansion:
\begin{eqnarray}
(n/d_{0,n})^{1/2}(Q_n-Q_{0,n})(x)&=& n^{1/2}(P_n-P_0) d_{0,n}^{-1/2}\sum_{j\in {\cal R}_{0,n}}\phi_j^*(Y-{Q_{0,n}}) \phi_j^*(x)\nonumber\\
&&-(n/d_{0,n})^{1/2}E_n(x) -
(n/d_{0,n})^{1/2}\tilde{r}_n(x)
\nonumber \\
&&\hspace*{-3cm}+(n/d_{0,n})^{1/2}O_P(C(M_n)r(d,J_{0,n})^{k+1}),\label{keyexpansion}\end{eqnarray}
where $M_n=\pl Q_n\pl_{v,k}^*$, which equals the $L_1$-norm  of the vector of non-zero coefficients in $Q_n$.

Due to Assumptions A0,A1,A3 and A4 we have
\begin{eqnarray}
\frac{(Q_n-Q_{0,n})(x)}{(n/d_{0,n})^{-1/2}}= n^{1/2}(P_n-P_0) d_{0,n}^{-1/2}\sum_{j\in {\cal R}_{0,n}}\phi_j^*(Y-{Q_{0,n}}) \phi_j^*(x)+o_P(1). \label{keyexpansion} 
\end{eqnarray}

For a given $x$, the leading term is now  a sum of independent mean zero random variables $D_{Q_{0,n},x}(O)$ (\ref{D0nxls}) so that the central limit theorem can be applied. 
The variance is given by:
\begin{equation}\label{sigma2nf}
\tilde{\sigma}^2_{0,n}(x)=
\frac{1}{d_{0,n}}\sum_{j_1,j_2\in {\cal R}_{0,n} }\Sigma_{0,n}(j_1,j_2) \phi_{j_1}^*(x)\phi_{j_2}^*(x),\end{equation}
where
\[
\Sigma_{0,n}(j_1,j_2)=P_0 \{ \phi_{j_1}^*\phi_{j_2}^*(Y-{Q_{0,n}})^2 \}.
\]

{\bf Bounded variance condition:}
Assume $\lambda_n$ is $O(1)$ or grows as a power of $\log n$-factor: for some $m<\infty$
\begin{equation}\label{keyconditionbx}\sup_x \tilde{\sigma}^2_{0,n}(x)=O(\log^m(n)).\end{equation}
By Lemma \ref{lemmaboundedvar}, we have \[
\tilde{\sigma}^2_{0,n}(x)=O\left(\lambda_n d_n^{-1}\sum_{j\in {\cal R}_{0,n}}\phi^{*2}_j(x)\right),\]
where $\lambda_n$ is the maximal eigenvalue of $\Sigma_{0,n}$.
 Thus a sufficient condition for (\ref{keyconditionbx}) is that $d_n^{-1}\sum_{j\in {\cal R}_{0,n}}\phi^{*2}_j(x)=O(1)$ and $\lambda_n=O(\log^m(n))$ for some $m<\infty$.

{\bf Conclusion:}
We have
\[
\tilde{\sigma}^{-1}_{0,n} (n/d_{0,n})^{1/2}({Q_n}-{Q_{0,n}})(x)\Rightarrow_d N(0,1),\]
and
\[
\tilde{\sigma}^{-1}_{0,n}(n/d_{0,n})^{1/2}(Q_n-Q_0)(x)\Rightarrow_d N(0,1).\]
Assuming $M_n=O(\log^{m}n)$ for some $m<\infty$, by choosing $J_{0,n}=n^{1/(2k^*+1)}\log^{m_1} n$ for some  $m_1<\infty$,  it follows that for  some  $m<\infty$, \[
\mid Q_n-Q_0\mid(x)=O_P(n^{-k^*/(2k^*+1)}\log^m n),\]
while Assumption A4 holds too. 
\end{theorem}


\paragraph{Discussion of bounded variance condition:}
The rate of convergence of ${Q_n}-{Q_{0,n}}$ is thus $O_P((d_{0,n}/n)^{1/2}\tilde{\sigma}_{0,n})$. Given that in $L^2$-norm the rate of convergence of ${Q_n}-{Q_{0,n}}$ is $O((d/n)^{1/2}\log n)$, we should expect that at most $\tilde{\sigma}_{0,n}(x)=O(\log n)$ so that the rate of convergence of $({Q_n}-{Q_{0,n}})(x)$ is $(d_n/n)^{1/2}$ up till a power of  $\log n$-factor.
The following argument suggests that in fact  $\tilde{\sigma}_{0,n}(x)=O(1)$.
Suppose that $E((Y-{Q_{0,n}})^2\mid X)$ is constant in $X$. Then, it follows that \[
\tilde{\sigma}^2_{0,n}(x)=
d_n^{-1}\sum_{j\in {\cal R}_{0,n} } \{\phi_{j}^*(x)\}^2 ,
\]
which is $O(1)$. It would be strange that the conditional variance of the residuals, as long  as bounded away from infinity and $0$, would change the rate of convergence relative to the rate when the conditional variance is constant. Similarly, it would be strange if the weighted least squares estimator converges at a faster rate than the unweighted least squares estimators. 
Either way, our assumption (\ref{keyconditionbx}) only assumes that the rate of $\tilde{\sigma}_{0,n}(x)$ is at most a power of $\log n$.

\subsection{Lemma for understanding second order remainder condition (\ref{keyconditiona1m}).}
The following lemma proves the sufficient condition $J_n=O(n^{2/5-\delta})$ for some $\delta>0$ for  the empirical process condition (\ref{keyconditiona1m}) $E_n(x)=o_P((n/d_n)^{-1/2})$.
\begin{lemma}\label{lemmasuffkeyconditiona}
Consider the squared error loss and assume $Y$ is uniformly bounded. In addition, assume that there exists a $C<\infty$ so that with probability tending to 1 $Q_n,Q_{0,n}\in D^{(0)}_C([0,1]^d)$.
Then, $\pl Q_n-Q_{0,n}\pl_{P_0}=O_P((J_n/n)^{1/2} \log n)$.

For $(j_1,j_2)\in {\cal R}_{0,n}^2$, let \[
\Sigma_n(j_1,j_2)\equiv
P_0 \left\{ \phi_{j_1}^*\phi_{j_2}^*({Q_n}-{Q_{0,n}} )^2 \right\} .\]
Let $\lambda_n$ be the maximal eigenvalue of matrix $\Sigma_n$. 
(\ref{keyconditiona1m}) holds if  for some $\delta>0$ $d_n \lambda_n=O(n^{-\delta})$, where we are reminded $d_n\sim J_n$. 

\end{lemma}


{\bf Proof of Lemma \ref{lemmasuffkeyconditiona}:} 
The conditions guarantee that  $\sup_{Q\in D^{(0)}_C([0,1]^d) }\pl L(Q)\pl_{\infty}=O(1)$ and $\sup_{Q\in D^{(0)}_C([0,1]^d)}\frac{P_0(L(Q)-L(Q_{0,n}))^2}{d_0(Q,Q_{0,n})}=O(1)$, while $Q_n,Q_{0,n}\in D^{(0)}_C([0,1]^d)$. These conditions were sufficient for establishing the rate for  $d_0(Q_n,Q_{0,n})$, which equals $\pl Q_n-Q_{0,n}\pl^2_{P_0}$.


We now prove the next statement. The function  $d_n^{-1}\sum_{j\in {\cal R}_{0,n}}\phi_j^*({Q_{n,\beta_n}}-{Q_{n,\beta_{0,n}}})\phi_j^*(x)$ is included in the fixed class (for some $M<\infty$)
\[
{\cal F}_{n,x}=\left\{ d_n^{-1}\sum_{j\in {\cal R}_{0,n}}\phi_j^*({Q_{n,\beta}}-{Q_{n,\beta_{0,n}}} )\phi_j^*(x):\beta\right\}.\]
We can further embed ${\cal F}_{n,x}$ in
\[
{\cal F}_n= \left\{ d_n^{-1}\sum_{j\in {\cal R}_{0,n}}\phi_j^*({Q_{n,\beta}}-{Q_{n,\beta_{0,n}}} )\phi_j^*(x):\beta,x\in [0,1]^d\right\}.\]
Let $g_{\beta,x}$ be one element in ${\cal F}_n$. Then
\[
P_0 g_{\beta,x}^2=d_n^{-2}\sum_{j_1,j_2\in {\cal R}_{0,n}}\Sigma_0(\beta,\beta_{0,n}) (j_1,j_2)\phi_{j_1}^*(x)\phi_{j_2}^*(x),\]
where
\[
\Sigma_0(\beta,\beta_{0,n})(j_1,j_2)=
P_0 \left\{ \phi_{j_1}^*\phi_{j_2}^*({Q_{n,\beta}}-{Q_{n,\beta_{0,n} } })^2 \right\} .\]
Let $\lambda(\beta,\beta_{0,n})$ be the  maximal eigenvalue of $\Sigma_0(\beta,\beta_{0,n})$. Then, the proof of Lemma \ref{lemmaboundedvar} above  establishes
\[
P_0 g_{\beta,x}^2=O(d_n^{-1}\lambda(\beta,\beta_{0,n})).\]
Moreover, ${\cal F}_n$ has  universally  bounded sup-norm,
and 
$J(\delta,{\cal F}_n)\sim -d_n^{1/2}\int_0^{\delta} (\log \epsilon)^{1/2} d\epsilon$.  This can be conservatively bounded by \[ J(\delta,{\cal F}_n)\sim -d_n^{1/2}\int_0^{\delta}\log \epsilon d\epsilon = -d_n^{1/2}\{\delta \log \delta +\delta\}\sim -d_n^{1/2}\delta \log \delta.\]

Let $\lambda_n$ be the rate of $\lambda(\beta_n,\beta_{0,n})$. 
Then we can set $\delta_n=d_n^{-1/2}\lambda_n^{1/2}$.
Applying empirical process inequalities for $E \sup_{Q\in {\cal F}_n,\pl Q\pl_{P_0}\leq \delta_n}\mid n^{1/2}(P_n-P_0)Q\mid \sim J(\delta_n,{\cal F}_n,L^2)$ yields now:
\[
\begin{array}{l}
\sup_x\left | \sum_{j\in {\cal R}_{0,n}}(P_n-P_0)\phi_j^*({Q_n}-{Q_{0,n}}) \phi_j^*(x)\right |\\
=d_n\sup_x\left | d_n^{-1}\sum_{j\in {\cal R}_{0,n}}(P_n-P_0)\phi_j^*({Q_n}-{Q_{0,n}}) \phi_j^*(x)\right | \\
\leq
d_n n^{-1/2} \sup_{g\in {\cal F}_n,\pl g\pl_{P_0}\leq \delta_n} \mid n^{1/2}(P_n-P_0)g\mid \\
\sim  (d_n/n)^{1/2} d_n^{1/2} J(\delta_n,{\cal F}_n,L^2)\\
  \sim (d_n/n)^{1/2}  d_n^{1/2} d_n^{1/2} \delta_n \log \delta_n\\
  = (d_n/n)^{1/2} d_n d_n^{-1/2}\lambda_n^{1/2} \log \delta_n\\
  = (d_n/n)^{1/2} d_n^{1/2} \lambda_n^{1/2} \log \delta_n\\
  = (d_n/n)^{1/2}  (d_n\lambda_n)^{1/2} \log \delta_n.
\end{array}
\]
Thus, for (\ref{keyconditiona1m}) it suffices that $  d_n\lambda_n=O(n^{-\delta})$ for some $\delta>0$. $\Box$

\subsection{Asymptotic normality for unweighted least squares regression with explicit formula for asymptotic variance.}\label{AppendixunweightedLS}
Under a stronger uniform approximation assumption on $D^{(k)}({\cal R}_{0,n})$ we can obtain asymptotic normality with a closed form expression for the asymptotic variance instead of reliance on a maximal eigenvalue as in the previous Theorem \ref{thasnormalB}

Let $dP^*_0(x)=\sigma^2_{0,n}(x)dP_{X,0}(x)$, where $\sigma^2_{0,n}(X)=E_0(Y-Q_{0,n})^2\mid X)$, and let $L^2(P_0^*)$ be the corresponding Hilbert space of functions of $X$ with inner product $\langle f,g\rangle_{P_0^*}=P_0^* fg$. 
For simplicity, we assume that $\inf_{x\in [0,1]^d}\sigma^2_{0,n}(x)>\delta>0$ for some $\delta>0$, although for pointwise convergence at a particular $x_0$ one only needs to assume $\lim\sup \sigma^2_{0,n}(x_0)>0$.
 Let $\{\phi_j^*:j\in {\cal R}_{0,n}\}$ be an orthonormal basis for $\{\phi_j:j\in {\cal R}_{0,n}\}$ in $L^2(P_0^*)$. 
 
{\bf Undersmoothing to solve the target parameter specific efficient score equation:}
 Solving the score equations $r_n(j)=P_n \frac{d}{dQ_n}L(Q_n)(\phi_j)\approx 0$, $j\in {\cal R}_{0,n}$, implies  solving 
 $r_n^*(j)\equiv P_n \frac{d}{dQ_n}L(Q_n)(\phi_j^*)\approx 0$, $j\in {\cal R}_{0,n}$.
Define the undersmoothing term:
\begin{equation}\label{tildernxlsb}
\tilde{r}_n(x)\equiv \sum_{j\in {\cal R}_{0,n}}r_n^*(j)\phi_j^*(x).\end{equation}

As starting point for our analysis of $(Q_n-Q_{0,n})(x)$, using that $P_0 \phi_j^* (Y-Q_{0,n})=0$ and $P_n \phi_j^*(Y-Q_n)=r_n^*(j)$,  we have
\[
P_0 \phi_j^*(Y-{Q_n})-P_0\phi_j^*(Y-{Q_{0,n}})=-(P_n-P_0)\phi_j^*(Y-{Q_n})+r_n^*(j),\]
so that
\[
P_0 \phi_j^*({Q_n}-{Q_{0,n}} )=(P_n-P_0)\phi_j^*(Y-{Q_n})-r_n^*(j).\]

We write this as
\[
P_0^* g_{0,n} \phi_j^*=(P_n-P_0)\phi_j^*(Y-Q_n)-r_n^*(j),\]
where 
\begin{equation}
\label{g0n}
g_{0,n}\equiv \sigma^{-2}_{0,n}(Q_n-Q_{0,n}).
\end{equation} 
We decompose $g_{0,n}=\Pi_{J_{0,n}}(g_{0,n})+(g_{0,n}-\Pi_{J_{0,n}}(g_{0,n}))$, where $\Pi_{J_{0,n}}(g_{0,n})=\arg\min_{g\in D^{(k)}({\cal R}_{0,n})}P_0^*(g-g_{0,n})^2$ is the projection of $g_{0,n}$ on the fixed working model $D^{(k)}({\cal R}_{0,n})$ in $L^2(P_0^*)$.

Since $\{\phi_j^*:j\in {\cal R}_{0,n}\}$, is an orthonormal basis of the linear subspace $D^{(k)}({\cal R}_{0,n})$ of $L^2(P_0^*)$, we have
\begin{eqnarray*}
\Pi_{J_{0,n}}(g_{0,n})(x)&=& \sum_{j\in {\cal R}_{0,n}}\{P_0^*g_{0,n} \phi_j^* \}\phi_j^*(x)\\
&=&\sum_{j\in {\cal R}_{0,n}}(P_n-P_0)\phi_j^*(Y-{Q_n}) \phi_j^*(x)-\sum_{j\in {\cal R}_{0,n}}r_n^*(j) \phi_j^*(x)\\
&=&\sum_{j\in {\cal R}_{0,n}}(P_n-P_0)\phi_j^*(Y-{Q_n}) \phi_j^*(x)-
\tilde{r}_n(x).
\end{eqnarray*}
{\bf Proving that the projection $\Pi_{J_{0,n}}(g_{0,n})$  approximates $g_{0,n}$ uniformly at rate $O(r(d,J_{0,n})^{k+1})$:}
By our nonparametric uniform approximation assumptoin  we know that there exists a $\tilde{g}_{0,n}\in D^{(k)}({\cal R}_{0,n})$ so that $\pl g_{0,n}-\tilde{g}_{0,n}\pl_{\infty}=O_P(C(M)r(d,J_{0,n})^{k+1})$, where $M$ is a bound on $\pl g_{0,n}\pl_{v,k}^*$.  Let $Q_{J_{0,n},\beta}=\sum_{j\in {\cal R}_{0,n}}\beta(j)\phi_j^*$ and $\tilde{\beta}_n$ is so that $\tilde{g}_{0,n}=Q_{J_{0,n},\tilde{\beta}_n}$.

We have that $\Pi_{J_{0,n}}(g_{0,n})=\arg\min_{g\in D^{(k)}({\cal R}_{0,n})}P_0^*(g_{0,n}-g)^2$. We can apply Theorem \ref{theoremoraclemle} with $R_0(\beta)=P_0^*(Q_{J_{0,n},\beta}-g_{0,n})^2$; $\beta_n=\arg\min_{\beta}R_0(\beta)$; $\Pi_{J_{0,n}}(g_{0,n})=Q_{J_{0,n},\beta_n}$; and $\tilde{g}_{0,n}=Q_{J_{0,n},\tilde{\beta}_n}$.  The conditions of Theorem \ref{theoremoraclemle} are  satisfied as demonstrated with the least squares loss function under Theorem \ref{theoremoraclemle}. Applying Theorem \ref{theoremoraclemle} proves that the uniform approximation error of $\tilde{g}_{0,n}$ w.r.t. $g_{0,n}$ carries over to the minimizer $Q_{J_{0,n},\beta_n}$ of the squared error risk $R_0(\beta$. Specifically, it proves that if $\pl g_{0,n}\pl_{v,k}^*<M$,  then $\pl \Pi_{J_{0,n}}(g_{0,n})-g_{0,n}\pl_{\infty}=O_P(C(M)r(d,J_{0,n})^{k+1})$. 


So this proves that
\[
 \sigma^{-2}_{0,n}(Q_n-Q_{0,n})(x)
=\sum_{j\in {\cal R}_{0,n}}(P_n-P_0)\phi_j^*(Y-{Q_n}) \phi_j^*(x)-
\tilde{r}_n(x)+O_P(C(M_n) r(d,J_{0,n})^{k+1}),\]
where $M_n=\pl \sigma^{-2}_{0,n}(Q_n-Q_{0,n})\pl_{v,k}^*$ is the $k$-th order sectional variation norm. The latter can be bounded in terms of $\pl \sigma^{-2}_{0,n}\pl_{v,k}^*$ and the $L_1$-norms of the vectors  $\beta_n$ and $\beta_{0,n}$ of non-zero coefficients that define $Q_n\in D^{(k)}({\cal R}_n)$ and $Q_{0,n}\in D^{(k)}({\cal R}_{0,n})$. 

{\bf Expansion for $(Q_n-Q_{0,n})(x)$:}
Thus, 
\begin{eqnarray*}
(Q_n-Q_{0,n})(x)&=&\sigma^2_{0,n}(x)\sum_{j\in {\cal R}_{0,n}}(P_n-P_0)\phi_j^*(Y-{Q_n}) \phi_j^*(x)-
\sigma^2_{0,n}(x)\tilde{r}_n(x)\\
&&+O_P(\sigma^2_{0,n}(x)C(M_n) r(d,J_{0,n})^{k+1}).\end{eqnarray*}

Assuming the second order remainder condition $E_n(x)=o_P((n/d_n)^{-1/2})$ with 
\begin{equation}\label{Enxlsb}
E_n(x)\equiv \sum_{j\in {\cal R}_{0,n}}(P_n-P_0)\phi_j^*(Q_n-Q_{0,n}) \phi_j^*(x),
\end{equation}
 this gives
\begin{eqnarray*}
(n/d_{0,n})^{1/2}(Q_n-Q_{0,n})(x)&=&\sigma^2_{0,n}(x)n^{1/2} (P_n-P_0)d_{0,n}^{-1/2}\sum_{j\in {\cal R}_{0,n}}\phi_j^*(Y-{Q_{0,n}}) \phi_j^*(x)\\
&&\hspace*{-3cm}-
(n/d_{0,n})^{1/2}\sigma^2_{0,n}(x)\tilde{r}_n(x)+O_P((n/d_{0,n})^{1/2}\sigma^2_{0,n}(x)C(M_n) r(d,J_{0,n})^{k+1}).\end{eqnarray*}
The last two terms are remainders controlled by undersmoothing, so that we can focus on the leading term.

The leading term is a sum of independent mean zero random variables
\begin{equation}\label{D0nxlsb}
D_{Q_{0,n},x}(O)\equiv d_{0,n}^{-1/2}\sum_{j\in {\cal R}_{0,n}}\phi_j^*(Y-{Q_{0,n}}) \phi_j^*(x).
\end{equation}
The variance of $D_{Q_{0,n},x}$  is given by:
\begin{equation}\label{sigma2ng}
\tilde{\sigma}^2_{0,n}(x)=\sigma^4_{0,n}(x)
\frac{1}{d_{0,n}}\sum_{j_1,j_2\in {\cal R}_{0,n} }\Sigma_{0,n}(j_1,j_2) \phi_{j_1}^*(x)\phi_{j_2}^*(x),\end{equation}
where
\[
\Sigma_{0,n}(j_1,j_2)=P_0 \{ \phi_{j_1}^*\phi_{j_2}^*\sigma^2_{0,n} \}=P_0^*\phi_{j_1}^*\phi_{j_2}^*.
\]
Thus, by the orthogonality of $\phi_j^*$ in $L^2(P_0^*)$ we have
\[
\tilde{\sigma}^2_{0,n}(x)=\sigma^4_{0,n}(x) \frac{1}{d_{0,n}}\sum_{j\in {\cal R}_{0,n}}\{\phi_j^*(x)\}^2.\]
This simplified expression for $\tilde{\sigma}_{0,n}^2(x)$ motivated our choice for an orthonormal basis in $L^2(P_0^*)$ (e.g., instead of in $L^2(P_0)$).
Lemma \ref{lemmaexplicit}  demonstrates that  $d_{0,n}^{-1}\sum_{j\in {\cal R}_{0,n}}\{\phi_j^*(x)\}^2=O(1)$, and   that  this expression approximates $1/(\sigma^2_{0,n}(x)p_0(x))$ under the nonparametric uniform approximation condition on $D^{(k)}({\cal R}_{0,n})$. The latter then implies that $\tilde{\sigma}^2_{0,n}(x)\rightarrow_p \sigma^2_0(x)/p_0(x)$, where $\sigma^2_0(x)=E_0( (Y-Q_0)^2\mid X=x)$.

\subsection{Asymptotic normality theorem for unweighted least squares with closed form asymptotic variance}
Recall the definitions of $\tilde{r}_n(x)$ (\ref{tildernxlsb}), $E_n(x)$ (\ref{Enxlsb}) and $g_{0,n}$ (\ref{g0n}). 
We make the same assumptions as in previous theorem except that we n eed a slightly stronger uniform approximation condition A0n on $D^{(k)}({\cal R}_{0,n})$ stated below.
\newline 
{\bf Assumption A0n:}
Let $Q_{0,n}=Q_{{\cal R}_{0,n}}=\arg\min_{Q\in D^{(k)}({\cal R}_{0,n})}P_0 L(Q)$ and
$g_{0,n}=\sigma^{-2}_{0,n}(Q_n-Q_{0,n})$.
We assume that  $\inf_{Q\in D^{(k)}({\cal R}_{0,n})}\pl Q-Q_0\pl_{\infty}=O_P(C(M_n)r(d,J_{0,n})^{k+1})$ and $\inf_{Q\in D^{(k)}({\cal R}_{0,n})}\pl g_{0,n}-Q\pl_{\infty}=O_P(C(M_n)r(d,J_{0,n})^{k+1})$. This implies  that $\pl Q_{0,n}-Q_0\pl_{\infty}=O_P(C(M_n)r(d,J_{0,n})^{k+1})$ and
$\pl g_{0,n}-\Pi_{J_{0,n}}(g_{0,n})\pl_{\infty}=O_P(C(M_n)r(d,J_{0,n})^{k+1})$.

\begin{theorem}\label{thasnormal2}\ \nl
Consider the $k$-th order sieve MLE, HAL-MLE or relax HAL-MLE $Q_n$.  Let $k^*=k+1$.
Let ${\cal R}_n$ be the set of $d_n$ non-zero coefficients in $Q_n=\sum_{j\in {\cal R}_n}\beta_n(j)\phi_j$. Assume $\inf_{x\in [0,1]^d}\sigma^{2}_{0,n}(x)>0$. 
Recall the definitions of $\tilde{r}_n(x)$ (\ref{tildernxlsb}, $E_n(x)$ (\ref{Enxlsb}), $D_{Q_{0,n},x}(O)$ (\ref{D0nxlsb}) above. 
Assume the Assumptions A0n above and conditions (\ref{undersmoothm}), (\ref{keyconditiona1m}), (\ref{neglbias})  of Theorem \ref{thasnormal}.

We have the following expansion: 
\begin{eqnarray}
(Q_n-Q_{0,n})(x)&=&\sigma^2_{0,n}(x)\sum_{j\in {\cal R}_{0,n}}(P_n-P_0)\phi_j^*(Y-{Q_n}) \phi_j^*(x)\nonumber\\
&&-
\sigma^2_{0,n}(x)\tilde{r}_n(x)+O_P(\sigma^2_{0,n}(x)C(M_n) r(d,J_{0,n})^{k+1}),\label{keyexpansion}\end{eqnarray}
where $M_n=\pl \sigma^{-2}_{0,n}(Q_n-Q_{0,n})\pl_{v,k}^*$ is the $k$-th order sectional variation norm. 
The latter can be bounded in terms of  $\pl \sigma^{-2}_{0,n}\pl_{v,k}^*$ and the $L_1$-norms $\pl \beta_n\pl_1$ and $\pl \beta_{0,n}\pl_1$ of the vectors of non-zero coefficients that define $Q_n$ and $Q_{0,n}$. 
Due to Assumptions A0n,A1,A3 and A4 we have
\begin{eqnarray}
\frac{(Q_n-Q_{0,n})(x)}{(d_{0,n}/n)^{1/2}}&=&\sigma^2_{0,n}(x)n^{1/2} (P_n-P_0)d_{0,n}^{-1/2}\sum_{j\in {\cal R}_{0,n}}\phi_j^*(Y-{Q_{0,n}}) \phi_j^*(x)+o_P(1)\label{keyexpansion} 
\end{eqnarray}

For a given $x$, the leading term is now  a sum of independent mean zero random variables $D_{Q_{0,n},x}$ so that the central limit theorem can be applied. 
The variance is given by:
\begin{equation}\label{sigma2nh}
\tilde{\sigma}^2_{0,n}(x)=\sigma^4_{0,n}(x) \frac{1}{d_{0,n}}\sum_{j\in {\cal R}_{0,n}}\{\phi_j^*(x)\}^2.\end{equation}
By Lemma \ref{lemmaexplicit} it follows that
$\tilde{\sigma}^2_{0,n}(x)=O(1)$, and, under the nonparametric uniform approximation assumption on $D^{(k)}({\cal R}_{0,n})$ we have
\[
\tilde{\sigma}^2_{0,n}(x)=\frac{\sigma^4_{0,n}(x)}{ p_0^*(x)}+o_P(1)=\frac{\sigma^2_{0,n}(x)}{p_0(x)}+o_P(1).\]

{\bf Conclusion:}
We have
\[
\tilde{\sigma}^{-1}_{0,n} (n/d_{0,n})^{1/2}({Q_n}-{Q_{0,n}})(x)\Rightarrow_d N(0,1),\]
and
\[
\tilde{\sigma}^{-1}_{0,n}(n/d_{0,n})^{1/2}(Q_n-Q_0)(x)\Rightarrow_d N(0,1).\]
Assuming $M_n=O(\log^{m}n)$ for some $m<\infty$, by choosing $J_{0,n}=n^{1/(2k^*+1)}\log^{m_1} n$ for some  $m_1<\infty$,  it follows that for  some  $m<\infty$, \[
\mid Q_n-Q_0\mid(x)=O_P(n^{-k^*/(2k^*+1)}\log^m n),\]
while Assumption A4 holds too.

\end{theorem}

\section{Asymptotic normality of arbitrary (non-pathwise differentiable) functions of sieve and HAL-MLEs, building on Section \ref{section10}}\label{AppendixK1}

Theorem \ref{theoremasnormgen} establishes asymptotic linearity $(Q_n-Q_{0,n})(x) =P_n D_{Q_{0,n},x}+R_n(x)$ for some remainder $R_n(x)=o_P((n/d_{0,n})^{-1/2})+O_P(C(M_n) r(d,J_{0,n})^{k+1})$, and thereby asymptotic normality at rate $(n/d_{0,n})^{-1/2}$ when $d_{0,n}$ is chosen so that this rate dominates the bias $C(M_n)r(d,J_{0,n})^{k+1}$. We assume the latter to be true. 
The influence curve is given by  \[
D_{Q_{0,n},x}\equiv \sum_{j\in {\cal R}_{0,n}}  S_{Q_{0,n}}({\phi}_j^*)\phi_j^*(x).\]
As in our proof of the theorem, the same asymptotic linearity applies to $(\Pi_{J_{0,n}}(Q_n)-Q_{0,n})(x)$, where $\Pi_{J_{0,n}}(Q_n)\in D^{(k)}({\cal R}_{0,n})$ is known to satisfy $\pl Q_n-\Pi_{J_{0,n}}(Q_n)\pl_{\infty}=O_P(C(M_n) r(d,J_{0,n})^{k+1})$, and for notational convenience, let's denote it with $\tilde{Q}_n=\Pi_{J_{0,n}}(Q_n)$. 

 Expression (\ref{tempidb}) also showed that
\begin{equation}\label{tempidbn}
\langle \phi_j^*,\Pi_{J_{0,n}}(Q_n)-Q_{0,n}\rangle_{J_{0,n}}=(P_n-P_0) S_{Q_{0,n}}({\phi}_j^*)-r_n({\phi}_j^*)+E_n({\phi}_j^*)+\sum_{m=1}^3R_{m,n}({\phi}_j^*),
\end{equation}
where $E_n({\phi}_j^*)\equiv (P_n-P_0)(S_{Q_n}-S_{Q_{0,n}})({\phi}_j^*)$.
Let's denote the total remainder in last formula with $R_n({\phi}_j^*)\equiv r_n({\phi}_j^*)+\sum_{m=1}^3 R_{m,n}({\phi}_j^*) +E_n({\phi}_j^*)$. We can represent $\tilde{Q}_n=\sum_j\beta_n^*(j)\phi_j^*$ and $Q_{0,n}=\sum_j \beta_{0,n}(j)\phi_j^*$, so that (\ref{tempidbn}) states 
\begin{equation}\label{diffcoeff}
(\beta_n^*-\beta_{0,n})(j)=P_n S_{Q_{0,n}}({\phi}_j^*)+R_n({\phi}_j^*).\end{equation}


Let's now consider a general functional $\Phi$ that maps $Q\in D^{(k)}([0,1]^d)$ into a function $\Phi(Q)$. We want to establish asymptotic normality of $(n/d_{0,n})^{1/2}(\Phi(Q_n)(x)-\Phi(Q_{0,n})(x))$.  For this to hold one excludes the case that $\Phi(Q)$ is a pathwise differentiable functional of $P_0)$ since then the natural rate of convergence would be $n^{1/2}$. The latter is addressed in the next section.


Let $R_{\Phi,n,x}$ be defined as the exact remainder for the first order Tailor expansion of $\Phi({f}_n)$ at $Q_{0,n}$.
\[
\Phi({f}_n)(x)-\Phi(Q_{0,n})(x)=d\Phi(Q_{0,n})({f}_n-Q_{0,n})(x)+R_{\Phi,n,x},\]
where $d\Phi(Q_{0,n})(h)=d/d \delta_0\Phi(Q_{0,n}+\delta_0 h)$ at $\delta_0=0$ is the directional derivative of $\Phi$ at $Q_{0,n}$ in direction $h$.  
 Plugging in the asymptotic linearity result above for $({f}_n-Q_{0,n})(x)$ yields
\[
\Phi({f}_n)(x)-\Phi(Q_{0,n})(x)= P_n d\Phi(Q_{0,n})(D_{Q_{0,n},\cdot} )(x)+
d\Phi(Q_{0,n})(R_n)(x)+R_{\Phi,n,x}.\]
Given we assumed $R_n(x)=o_P((n/d_{0,n})^{1/2})$ for all $x$ (or uniformly in $x$) and we already have a rate of convergence for $(Q_n-Q_{0,n})(x)=O_P((n/d_{0,n})^{1/2})$ for all $x$ (or uniformly),  it is reasonable to assume that $R_{\Phi,n,x}=o_P((n/d_{0,n})^{-1/2})$ and $d\Phi(Q_{0,n})(R_n)(x)=o_P((n/d_{0,n})^{-1/2})$.  We then obtained the desired asymptotic linearity
\[
\Phi({f}_n)(x)-\Phi(Q_{0,n})(x)= P_n D_{\Phi,Q_{0,n},x}+o_P((n/d_{0,n})^{-1/2}),\]
where
\[
 D_{\Phi,Q_{0,n},x}=d\Phi(Q_{0,n})(D_{Q_{0,n},\cdot} )(x).\]
For the sake of studying the expression for the asymptotic variance, it  is helpful to note that we can repeat the same analysis for $\Phi(\tilde{Q}_n)-\Phi(Q_{0,n})$, so that for the sake of asymptotic variance we can focus on the latter, having the advantage that $\tilde{Q}_n,Q_{0,n}\in D^{(k)}({\cal R}_{0,n})$.
To be concrete, 
let $R_{\Phi,n,x,1}$ be defined as the exact remainder for the first order Tailor expansion of $\Phi(\tilde{Q}_n)$ at $Q_{0,n}$.
\[
\Phi(\tilde{Q}_n)(x)-\Phi(Q_{0,n})(x)=d\Phi(Q_{0,n})(\tilde{Q}_n-Q_{0,n})(x)+R_{\Phi,n,x,1}.\]
 Plugging in the same asymptotic linearity result above for $(\tilde{Q}_n-Q_{0,n})(x)$, assuming $R_{\Phi,n,x,1}=o_P((n/d_{0,n})^{-1/2})$, yields now the same asymptotic linearity: 
\[
\Phi(\tilde{Q}_n)(x)-\Phi(Q_{0,n})(x)= P_n D_{\Phi,Q_{0,n},x}+o_P((n/d_{0,n})^{-1/2}).\]

{\bf Expression for asymptotic variance of $\Phi(\tilde{Q}_n)-\Phi(Q_{0,n})$:}
We can represent $\tilde{Q}_n=\sum_j\beta_n^*(j)\phi_j^*$ and $Q_{0,n}=\sum_j \beta_{0,n}(j)\phi_j^*$. 
 Let's dive a little deeper into the form of $D_{\Phi,Q_{0,n},x}$,  Firstly we notice that \[
 \Phi(\tilde{Q}_n)-\Phi(Q_{0,n})=\Phi\left(\sum_j\beta_n^*(j)\phi_j^*\right)-\Phi\left(\sum_j \beta_{0,n}(j)\phi_j^*\right)\equiv \Phi_1(\beta_n^*)-\Phi_1(\beta_{0,n}),\]
  so that we can represent it as $\Phi_1(\beta_n^*)-\Phi_1(\beta_{0,n})$. 
 Assuming again differentiability  of $\Phi_1$  yields
 \[
 \Phi_1(\beta_n^*)-\Phi_1(\beta_{0,n})=\sum_{j\in {\cal R}_{0,n}}\frac{d}{d\beta_{0,n}(j)}\Phi_1(\beta_{0,n}) (\beta_n^*-\beta_{0,n})(j)+R_{\Phi,n}.\]
 Let $c_{0,n}(j)\equiv \frac{d}{d\beta_{0,n}(j)}\Phi_1(\beta_{0,n})$ so that
 \[
 \Phi_1(\beta_n^*)-\Phi_1(\beta_{0,n})=\sum_{j\in {\cal R}_{0,n}}(\beta_n^*-\beta_{0,n})(j)c_{0,n}(j)+R_{\Phi,n}.\]
 We can plug-in  (\ref{diffcoeff}) $(\beta_n^*-\beta_{0,n})(j)=P_n S_{Q_{0,n}}({\phi}_j^*)+R_n({\phi}_j^*)$ so that we obtain
 \[
  \Phi_1(\beta_n^*)-\Phi_1(\beta_{0,n})=\sum_{j\in {\cal R}_{0,n}}S_{Q_{0,n}}({\phi}_j^*)c_{0,n}(j)+\sum_{j\in {\cal R}_{0,n}}R_n(\phi_j^*) c_{0,n}(j). \]
Our assumptions on second order remainders corresponds with assuming 
$ \sum_{j\in {\cal R}_{0,n}}R_n(\phi_j^*) c_{0,n}(j)=o_P((n/d_{0,n})^{-1/2})$. 
Thus the influence curve is given by
\[
D_{\Phi,Q_{0,n},x}=\sum_{j\in {\cal R}_{0,n}}S_{Q_{0,n}}({\phi}_j^*)c_{0,n}(j,x).\]

Suppose now that the loss function behaves as  log-likelihood so that $\{S_{Q_{0,n}}({\phi}_j^*):j\in {\cal R}_{0,n}\}$ is  an orthonormal basis in $L^2(P_0)$. 
The asymptotic variance of $(n/d_{0,n})^{1/2}(\Phi_1(\beta_n^*)-\Phi_1(\beta_{0,n}))$ is given by 
\[
\tilde{\sigma}_{0,n}^2(x)=\frac{1}{d_{0,n}}\sum_{j\in {\cal R}_{0,n}}\{c_{0,n}(j,x)\}^2 ,\]
due to the scores $S_{Q_{0,n}}({\phi}_j^*)$ being orthonormal in $L^2(P_0)$.
 We proved the following theorem.

\begin{theorem}\label{theoremasnormgenphi}
Consider the setting of Theorem \ref{theoremasnormgen} so that
$(Q_n-Q_{0,n})(x)= P_n D_{Q_{0,n},x}+R_n(x)+O_P(C(M_n) r(d,J_{0,n})^{k+1})$ for a specified remainder $R_n(x)=o_P((n/d_{0,n})^{-1/2})$.
Let $R_n({\phi}_j^*)\equiv r_n({\phi}_j^*)+\sum_{m=1}^3 R_{m,n}({\phi}_j^*)+E_n({\phi}_j^*)$ so that $R_n(x)=\sum_{j\in {\cal R}_{0,n}}R_n({\phi}_j^*)\phi_j^*(x)$.
Let $\Phi$ be a functional and define $R_{\Phi,n,x}$ as follows:
\[
\Phi({f}_n)(x)-\Phi(Q_{0,n})(x)=d\Phi(Q_{0,n})(\tilde{Q}_n-Q_{0,n})(x)+R_{\Phi,n,x},\] 
and similarly define $R_{\Phi,n,x,1}$ with $Q_n$ replaced by $\tilde{Q}_n$.

{\bf Assumptions:}
$C(M_n) r(d,J_{0,n})^{k+1}=o((n/d_{0,n})^{-1/2})$;  $R_{\Phi,n,x}=o_P((n/d_{0,n})^{-1/2})$; and $d\Phi(Q_{0,n})(R_n)(x)=o_P((n/d_{0,n})^{-1/2})$.

Then,  \[
\Phi(Q_n)(x)-\Phi(Q_{0,n})(x)= P_n D_{\Phi,Q_{0,n},x}+o_P((n/d_{0,n})^{-1/2}),\]
where
\[
 D_{\Phi,Q_{0,n},x}=d\Phi(Q_{0,n})(D_{Q_{0,n},\cdot} )(x).\]
 Let $\tilde{\sigma}^2_{0,n}(x)=P_0 D_{\Phi,Q_{0,n},x}^2$.
 Then, $(n/d_{0,n})^{1/2}(\Phi(Q_n)-\Phi(Q_{0,n}))(x)/\tilde{\sigma}_{0,n}(x)\Rightarrow_d N(0,1)$.
 
 {\bf Expression for influence curve:}
Assume also  $R_{\Phi,n,x,1}=o_P((n/d_{0,n})^{-1/2})$; 
 Let $\Phi_1(\beta)=\Phi(\sum_{j\in {\cal R}_{0,n}}\beta(j)\phi_j^*)$.
 Let $c_{0,n}(j)\equiv \frac{d}{d\beta_{0,n}(j)}\Phi_1(\beta_{0,n})$ 
 Then, we can represent the influence curve as
 \[
D_{\Phi,Q_{0,n},x}=\sum_{j\in {\cal R}_{0,n}}S_{Q_{0,n}}({\phi}_j^*)c_{0,n}(j,x).\]
Moreover, if the scores $S_{Q_{0,n}}({\phi}_j^*)$, $j\in {\cal R}_{0,n}$, are orthonormal in $L^2(P_0)$, as we showed to be the case for a log-likelihood loss function, then \[
\tilde{\sigma}_{0,n}^2(x)=\frac{1}{d_{0,n}}\sum_{j\in {\cal R}_{0,n}}\{c_{0,n}(j,x)\}^2 .\]
 \end{theorem}

\begin{corollary}
Assume the conditions of the above theorem so that $\Phi(Q_n)-\Phi(Q_{0,n})=P_n D_{\Phi,Q_{0,n},x}+o_P((n/d_{0,n})^{-1/2})$. Assume also that $\tilde{\sigma}_{0,n}(x)=O(1)$. 
Assume that a special analysis establishes that $\pl \Phi(Q_{0,n})-\Phi(Q_0)\pl_{\infty}=O_P(r(J_{0,n},k,\Phi))$ for some rate $r(J_{0,n},k,\Phi)$, where one can use that $\pl Q_{0,n}-Q_0\pl_{\infty}=O(C(M_n) r(d,J_{0,n})^{k+1})$. 
Suppose that $J_{0,n}$ is chosen so that $r(J_{0,n},k,\Phi)=o((J_{0,n}/n)^{1/2})$. Then, 
$(n/d_{0,n})^{1/2}(\Phi(Q_n)-\Phi(Q_0))(x)/\tilde{\sigma}_{0,n}(x)\Rightarrow_d N(0,1)$.
By choosing $J_{0,n}$ so that $(\log n) r(J_{0,n},k,\Phi)\sim (J_{0,n}/n)^{1/2}$ one optimizes the rate of convergence, while keeping the bias of smaller order (by $\log n$-factor) than the standard error. 
\end{corollary}


\subsection{Uniform rate of convergence of $\Phi(Q_n)-\Phi(Q_{0,n})$ and $\Phi(Q_n)-\Phi(Q_0)$ for general functionals $\Phi$.}

We can generalize the pointwise results in the previous subsection to uniform results. 

\begin{theorem}\label{theoremasnormgenphiunif}
Assume the setting of previous theorem and assume  $(\log n) C(M_n) r(d,J_{0,n})^{k+1}=o_P((n/d_{0,n})^{-1/2})$;  $\sup_x\mid R_{\Phi,n,x}\mid =o_P((n/d_{0,n})^{-1/2})$; $\sup_x\mid R_{\Phi,n,x,1}\mid =o_P((n/d_{0,n})^{-1/2})$, and $\sup_x\mid d\Phi(Q_{0,n})(R_n)(x)\mid=o_P((n/d_{0,n})^{-1/2})$.

Then, 
\[ \Phi({f}_n)(x)-\Phi(Q_{0,n})(x)= P_n D_{\Phi,Q_{0,n},x}+R_n(x),\]
where $\pl R_n\pl_{\infty}=o_P((n/d_{0,n})^{-1/2})$, and
\[ D_{\Phi,Q_{0,n},x}=\sum_{j\in {\cal R}_{0,n}}S_{Q_{0,n}}(\tilde{\phi}_j)c_{0,n}(j,x).\]

Let $\tilde{\Sigma}_{0,n}$ be covariance matrix  of $(S_{Q_{0,n}}(\tilde{\phi}_j): j\in {\cal R}_{0,n})$,
and let $\lambda_{0,n}$ be the  maximal eigenvalue of $\tilde{\Sigma}_{0,n}$. We assume that $\lambda_{0,n}=O_P((\log n)^m) $ for some $m<\infty$. In the special case of the log-likelihood loss we showed that this covariance matrix is the identity matrix so that $\lambda_{0,n}=1$.


We have 
\[
\sup_{x\in [0,1]^d}\mid \Phi(Q_n)(x)-\Phi(Q_{0,n})(x)\mid = o_P((d_{0,n}/n)^{1/2} \lambda_{0,n}\log n).
\]
Assume that a special analysis establishes that $\pl \Phi(Q_{0,n})-\Phi(Q_0)\pl_{\infty}=O_P(r(J_{0,n},k,\Phi))$ for some rate $r(J_{0,n},k,\Phi)$, where one can use that $\pl Q_{0,n}-Q_0\pl_{\infty}=O(C(M_n) r(d,J_{0,n})^{k+1})$.

{\bf Uniform convergence of $k$-th order sieve or (relax) HAL-MLE to target function $Q_0$:}
Suppose that $J_{0,n}$ is chosen so that $(\log n) C(M_n)r(J_{0,n},k,\Phi)=o_P((J_{0,n}/n)^{1/2}\lambda_{0,n})$. Then, 
\[
\sup_{x\in [0,1]^d}\mid \Phi(Q_n)-\Phi(Q_0)\mid(x)=o_P((d_{0,n}/n)^{1/2}\lambda_{0,n}\log n).\]
By choosing $J_{0,n}$ so that $r(J_{0,n},k,\Phi)\sim (J_{0,n}/n)^{1/2}\lambda_{0,n}$ one optimizes the rate of convergence while keeping bias $\Phi(Q_{0,n})-\Phi(Q_0)$ asymptotically negligible.
\end{theorem}

\section{Rate of convergence of the first and higher order spline sieve and HAL-MLE to target function w.r.t.  the sectional variation norm}\label{AppendixL}
In this section we study the rate of convergence of the sectional variation norm of the $k$-th order spline sieve MLE approximation error $Q^k_{{\bf J},\beta_n}-Q_0$ w.r.t $Q_0\in D^{(k)}([0,1]^d)$ with $k=1,2,\ldots$. We then discuss the implication of these results for the $k$-th order spline HAL-MLE.

\subsection{Expressing sectional variation norm of $Q\in D^{(k)}_M([0,1]^d)$.}
For notational convenience, we use $\beta^*$ for $\beta_{0,n}$ and $Q$ for $Q_0$.
We recall the $k$-th order spline representation of $Q\in D^{(k)}([0,1]^d)$ given by 
\[
\begin{array}{l}
Q(x)=\sum_{\bar{s}(k+1)}\bar{\phi}_{\bar{s}(k+1)}(x)\mu^k(\tilde{Q}^{(k)}_{\bar{s}(k+1)})\\
=\sum_{\bar{s}(k+1),\mid s_{k+1}\mid =0}\beta^*(\bar{s}(k+1))\bar{\phi}_{\bar{s}(k+1)}+
\sum_{\bar{s}(k+1),\mid s_{k+1}\mid>0}\bar{\phi}_{\bar{s}(k+1)}\mu^k(\tilde{Q}^{(k)}_{\bar{s}(k+1)}).
\end{array}
\]
Here for $\mid s_{k+1}\mid>0$,  $\bar{\phi}_{\bar{s}(k+1)}=\prod_{j=1}^k \phi_0^j(x(s_j/s_{j+1}))$ is a function of $x(s_1/s_{k+1})$ while $\mu^k(\tilde{Q}^{(k)}_{\bar{s}(k+1)})$ is a function of $x(s_{k+1})$. For $\mid s_{k+1}\mid =0$, $\bar{\phi}_{\bar{s}(k+1)}=\bar{\phi}_{\bar{s}(m+1)}=\prod_{j=1}^m\phi_0^j(x(s_j/s_{j+1}))$ with $m=m(\bar{s}(k+1))$ so that $s_{m+1}$ is first empty set in $\bar{s}(k+1)$.

By the zero-order spline representation for a function $Q\in D^{(0)}([0,1]^d)$, we  have that $Q=\sum_{s}\tilde{Q}_{s}$, where $\tilde{Q}_{s}(x(s))=\int_{(0(s),x(s)]}Q(du(s),0(-s))$, and $\pl Q\pl_v^*=\sum_{s}\pl \tilde{Q}_{s}\pl_v$.  We want to evaluate this $\pl Q\pl_v^*$ in terms of the $k$-th order spline representation of $Q\in D^{(k)}([0,1]^d)$. Firstly, it follows that
\[
\tilde{Q}_{s}(x(s))=\sum_{\bar{s}(k+1),s_1=s}\bar{\phi}_{\bar{s}(k+1)}(x(s/s_{{k+1}}))\mu^k(\tilde{Q}^{(k)}_{\bar{s}(k+1)} ) (x(s_{k+1})).\]
We now want to compute $\tilde{Q}_s(du(s))=Q(du(s),0(-s))$. We can do this by noting that this equals
$\tilde{Q}_s((du(j):j\in s/s_{k+1}),(du(j): j\in s_{k+1}))$. 
Represent $\tilde{Q}_s=\sum_{\bar{s}(k+1),s_1=s}\tilde{Q}_{s,\bar{s}(k+1)}$. Since for $\mid s_{k+1}\mid >0$,  $\tilde{Q}_{s,\bar{s}(k+1)}$ is a product of a function that only depends on $x(s/s_{k+1})$ times a function that only depends on $x_{s_{k+1}}$, this implies that for $\mid s_{k+1}\mid >0$
\[
\tilde{Q}_{s,\bar{s}(k+1)}(du(s))=\bar{\phi}_{\bar{s}(k+1)}(du(j):j\in s/s_{k+1})\mu^k(\tilde{Q}^{(k)}_{\bar{s}(k+1)})(du(j):j\in s_{k+1}),\]
while for $\mid s_{k+1}\mid =0$, (noting that $\bar{\phi}_{\bar{s}(k+1)}(x)$ is a function of all of $x(s)$ when $s_1=s$) we have
\[
\tilde{Q}_{s,\bar{s}(k+1)}(du(s))=\beta(\bar{s}(k+1))\bar{\phi}_{\bar{s}(m+1)}(du(j): j\in s)=\beta^*(\bar{s}(k+1)\prod_{j=1}^m\phi_0^j(du(s_j/s_{j+1})) .\]
Let's first focus on $\mid s_{k+1}\mid>0$.
We have
\[ \mu^k(\tilde{Q}^{(k)}_{\bar{s}(k+1)})(du(j):j\in s_{k+1})=\mu^{k-1}(\tilde{Q}^{(k)}_{\bar{s}(k+1)})(u(j):j\in s_{k+1})d\mu(u(s_{k+1})),\]
where $d\mu(u(s_{k+1}))= \prod_{j, j\in s_{k+1}}du(j)$. 
We note that
$\bar{\phi}_{\bar{s}(k+1)}=\prod_{j=1}^k \phi_0^j(x(s_j/s_{j+1}))$  is a product of $k$ functions where each function $\phi_0^j(x(s_j/s_{j+1}))$ depends on disjoint set of coordinates $(x(j):j\in s_j/s_{j+1})$. Moreover, each $\phi_0^j(x(s_j/s_{j+1}))=\prod_{l\in s_j/s_{j+1}}\phi_0^j(x_l)$ is a product univariate spline basis functions of order $j$.  
Therefore, we have that
\[
\bar{\phi}_{\bar{s}(k+1)}(du(j):j\in s_1/s_{k+1})=\prod_{j=1}^k \prod_{l\in s_j/s_{j+1}}\phi_0^j(du(l)).\]
However, by Lemma \ref{lemmamuphij} we have that $\phi_0^j(x)=\mu (\phi_0^{j-1})(x)=\int_{(0,x]}\phi_0^{j-1}(u) du$ so that it follows that $\phi_0^j(du(l))=\phi_0^{j-1}(u(l))du(l)$.
So we have
\begin{eqnarray*}
\bar{\phi}_{\bar{s}(k+1)}(du(j):j\in s_1/s_{k+1})&=&\prod_{j=1}^k \prod_{l\in s_j/s_{j+1}}\phi_0^{j-1}(u(l)) du(l)\\
&=&\prod_{j=1}^k \phi_0^{j-1}(u(s_j/s_{j+1}))d\mu(u(s_1/s_{k+1}).\end{eqnarray*}
Similarly, for $\mid s_{k+1}\mid =0$ (and $s_1=s$), we have
\[
\bar{\phi}_{\bar{s}(k+1)}(du(j):j\in s)=\prod_{j=1}^m\phi_0^{j-1}(u(s_j/s_{j+1}))d\mu(u(s)).\]
Let's define
\[
\bar{\phi}_{\bar{s}(k+1)}^{k-1}(x(s_1/s_{k+1})\equiv \prod_{j=1}^k \phi_0^{j-1}(x(s_j/s_{j+1})),\]
where, if $\mid s_{k+1}\mid =0$, this reduces to $\prod_{j=1}^m\phi_0^{j-1}(x(s_j/s_{j+1}))$.

Combining $d\mu(u(s_1/s_{k+1})d\mu(u(s_{k+1}))=d\mu(u(s_1))$, gives us, for $\mid s_{k+1}\mid>0$, 
\[
\tilde{Q}_{s,\bar{s}(k+1)}(du(s))=\bar{\phi}_{\bar{s}(k+1)}^{k-1}(u(s_1/s_{k+1})) \mu^{k-1}(\tilde{Q}^{(k)}_{\bar{s}(k+1)})(u(s_{k+1}))  d\mu(u(s)),\]
while for $\mid s_{k+1}\mid =0$, 
\[
\tilde{Q}_{s,\bar{s}(k+1)}(du(s))=\beta^*(\bar{s}(k+1))\bar{\phi}_{\bar{s}(k+1)}^{k-1}(u(s)) d\mu(u(s))=
\prod_{j=1}^{m(\bar{s}(k+1))}\phi_0^{j-1}(u(s_j/s_{j+1}))d\mu(u(s)).\]
Thus,
\[
\begin{array}{l}
\tilde{Q}_s(du(s))=\left\{ \sum_{\bar{s}(k+1),s_1=s} \bar{\phi}_{\bar{s}(k+1)}^{k-1}(u(s_1/s_{k+1})) \mu^{k-1}(\tilde{Q}^{(k)}_{\bar{s}(k+1)})(u(s_{k+1})) \right\}  d\mu(u(s))\\
=\left\{ 
\sum_{\bar{s}(k+1),\mid s_{k+1}\mid =0,s_1=s}
\beta^*(\bar{s}(k+1))\bar{\phi}_{\bar{s}(k+1)}^{k-1}(u(s))+\right .  \\
\left . \sum_{\bar{s}(k+1),\mid s_{k+1}\mid>0,s_1=s}\bar{\phi}_{\bar{s}(k+1)}^{k-1}(u(s_1/s_{k+1})) \mu^{k-1}(\tilde{Q}^{(k)}_{\bar{s}(k+1)} (u(s_{k+1})) \right\}  d\mu(u(s))\\
\equiv
\tilde{Q}_s^{(1)}(u(s))d\mu(u(s)).
\end{array}
\]
This proves that 
\begin{eqnarray*}
\pl Q\pl_v^*&=&\sum_s \int | \tilde{Q}_s^{(1)}(u(s))| d\mu(u(s))\\
&=& \sum_s \pl \tilde{Q}_s^{(1)}\pl_{1,\mu}.
\end{eqnarray*}
This proves the following lemma. 
\begin{lemma}
Let $Q\in D^{(k)}([0,1]^d)$.
Let's define
\[
\bar{\phi}_{\bar{s}(k+1)}^{k-1}(x(s_1/s_{k+1})\equiv \prod_{j=1}^k \phi_0^{j-1}(x(s_j/s_{j+1})),\]
where, if $\mid s_{k+1}\mid =0$, this reduces to $\prod_{j=1}^{m(\bar{s}(k+1)}\phi_0^{j-1}(x(s_j/s_{j+1}))$.
Let
\[
\tilde{Q}_s^{(1)}\equiv \sum_{\bar{s}(k+1),s_1=s}\bar{\phi}^{k-1}_{\bar{s}(k+1)}\mu^{k-1}(\tilde{Q}^{(k)}_{\bar{s}(k+1)}),\]
where for $\mid s_{k+1}\mid =0$ the terms within the sum reduces to $\beta(Q)(\bar{s}(k+1))\bar{\phi}^{k-1}_{\bar{s}(k+1)}$.
We have that the (zero-order) sectional variation norm $\pl Q\pl_v^*=\sum_s\pl Q_s\pl_v$ with $Q_s(x(s))=Q(x(s),0(-s))$ being the $s$-specific section of $Q$ can be expressed as follows:
\[
\pl Q\pl_v^*=\sum_s \pl \tilde{Q}_s^{(1)}\pl_{1,\mu} .\]
\end{lemma}
\subsection{Sectional variation norm of $k$-th order spline  finite dimensional approximation error w.r.t. $k$-th order smooth target function.}
We now  want to obtain an analogue expression  of the sectional variation norm $\pl Q_{{\bf J},\beta}^k-Q\pl_v^*$.
For notational convenience, let $j$ represent  $\bar{s}(k+1)$ and define $j(1)=s_1$ and $j(k+1)=s_{k+1}$ for a given $j=\bar{s}(k+1)$.
 Then, using our conventions for the case that $\mid j(k+1)\mid =0$, we can write 
\[
Q(x)=\sum_{j}\bar{\phi}_j(x)\mu^k(\tilde{Q}^{(k)}_j),\]
and 
\[
\tilde{Q}_s^{(1)}=\sum_{j, j(1)=s}\bar{\phi}^{k-1}_{j}\mu^{k-1}(\tilde{Q}^{(k)}_j).\]
The $k$-th order spline approximation implied by ${\bf J}$ and its knot-point set ${\cal R}(d,{\bf J})=\cup_{j,\mid j(k+1)\mid >0} {\cal R}(j,{\bf J}(j))\cup\{j:\mid j(k+1)\mid =0\}$ can then be denoted as
\[
Q_{{\bf J},\beta}^k(x)=\sum_{j}\bar{\phi}_j(x)\mu^k_{{\bf J},\beta_j}(\tilde{Q}^{(k)}_j),\]
where $\mu^k_{{\bf J},\beta_j}(\tilde{Q}
^{(k)}_j)=\mu^{k}(\tilde{Q}_{{\bf J},\beta_j})$ and $\tilde{Q}_{{\bf J},\beta_j}=\sum_{u\in {\cal R}(j,{\bf J}(j))}\beta(j,u)\phi_u^0$. 
Analogue to above  for general $Q\in D^{(k)}([0,1]^d)$ we have $Q^k_{{\bf J},\beta}=\sum_s Q^k_{{\bf J},\beta,s}$ with 
\[
{Q}^k_{{\bf J},\beta,s}(x(s))=\sum_{j,j(1)=s}\bar{\phi}_j\mu^k_{{\bf J},\beta_j} (\tilde{Q}^{(k)}_j), \]
and
\[
\pl {Q}^k_{{\bf J},\beta,s}\pl_v=\pl \sum_{j,j(1)=s}\bar{\phi}_j^{k-1}\mu^{k-1}_{{\bf J},\beta_j}(\tilde{Q}^{(k)}_j) \pl_{1,\mu} .\]
Thus, we also have the following lemma for representing $\pl Q-Q^k_{{\bf J},\beta}\pl_v^*$ for $Q\in D^{(k)}([0,1]^d)$.
\begin{lemma}\label{lemmasectvarexpression}
We have
\begin{eqnarray*}
\pl Q-Q^k_{{\bf J},\beta}\pl_v^*&=&\sum_s \pl \sum_{j,j(1)=s}\bar{\phi}_j^{k-1}\left\{\mu^{k-1}_{{\bf J},\beta_j}(\tilde{Q}^{(k)}_j)-\mu^{k-1}(\tilde{Q}^{(k)}_j)\right\}\pl_{1,\mu},
\end{eqnarray*}
where the inner sum within the $L_1$-norm can be further written as:
\[
\begin{array}{l}
\sum_{j,j(1)=s,\mid j(k+1)\mid =0}(\beta_j-\beta_j(Q))\bar{\phi}_j^{k-1} +\sum_{j,j(1)=s,\mid j(k+1)\mid >0}\bar{\phi}_j^{k-1}\left\{\mu^{k-1}_{{\bf J},\beta_j}(\tilde{Q}^{(k)}_j)-\mu^{k-1}(\tilde{Q}^{(k)}_j)\right\}.
\end{array}
\]

\end{lemma}
This lemma shows that we can bound $\pl Q-Q^k_{{\bf J},\beta}\pl_v^*$ in terms of the $L_1(\mu)$-norm of $\mu^{k-1}(\tilde{Q}_{J,\beta})-\mu^{k-1}(\tilde{Q})$ with $\tilde{Q}_{J,\beta}=\sum_{u\in {\cal R}(J)}\beta(u)\phi_u^0$ being the zero-order spline approximation, plus
$\sum_s\sum_{j,j(1)=s,\mid j(k+1)\mid =0}\mid \beta_j-\beta_j(Q)\mid$. 
In the our main approximation result Theorem \ref{mainkthordercdf} for approximating $\mu^k(\tilde{Q})$ with $\mu^k(\tilde{Q}_{J,\beta})$ we showed that 
there is a choice $\beta $ (which, in particular, selects $\beta_j=\beta_j(Q)$ for $j$ with $\mid j(k+1)\mid =0$, so that the sup-norm and variation norm of $\mu^k(\tilde{Q}_{J,\beta})-\mu^k(\tilde{Q})$ converge at rate $r(d,J)^{k+1}$ and $r(d,J)^k$, respectively. 
Applying this result to each $\mu^k(\tilde{Q}^{(k)}_j)$ in the above expression for $\pl Q-Q^k_{{\bf J},\beta}\pl_v^*$ proves then that there is an overall stacked vector $\beta$ so that both $\pl Q-Q^k_{{\bf J},\beta}\pl_{\infty}=O(r(d,J)^{k+1})$ and $\pl Q-Q^k_{{\bf J},\beta}\pl_v^*=O(r(d,J)^k)$, $k=1,2,\ldots$. 

\begin{theorem}\label{theoremconvsectvarnorma}
Let $Q\in D^{(k)}([0,1]^d)$ for given $k=1,2,\ldots$. 
We have 
\[
\sup_{Q\in D^{(k)}_{M}([0,1]^d)}
\inf_{\beta}\pl Q-Q^k_{{\bf J},\beta}\pl_v^*= O(r(d,J)^k).\]
Moreover, there is a $\beta=\beta(Q)$ so that both $\pl Q^k_{{\bf J},\beta}-Q\pl_{\infty}=O(r(d,J)^{k+1})$ and $\pl Q^k_{{\bf J},\beta}-Q\pl_v^*=O(r(d,J)^k)$, and, again, this holds uniformly in $Q\in D^{(k)}_M([0,1]^d)$.

Let $\beta_{0,n}=\arg\min_{\beta}P_0 L(Q^k_{{\bf J},\beta})$ and $Q_0=\arg\min_{Q}P_0L(Q)$ and assume $Q_0\in D^{(k)}_M([0,1]^d)$. Then
we also have
\[
\pl Q^k_{{\bf J},\beta_{0,n}}-Q_0\pl_v^*=O(r(d,J)^k).\]
\end{theorem}
{\bf Proof:}
Above we obtained expression for $\pl Q_{{\bf J},\beta}-Q\pl_v^*$. From Theorems \ref{mainkthordercdf} and \ref{mainkthordercdfnew} (for weaker condition on knot-point set) we know that there exists a $\beta^*$ so that 
 for all $j$
\[
\pl \mu^k_{{\bf J},\beta^*_j}(\tilde{Q}^{(k)}_j)-\mu^k(\tilde{Q}^{(k)}_j)\pl_{\infty}=O(r(d,J)^{k+1)}\]
and 
\begin{equation}\label{preserve}
\pl \mu^{k-1}_{{\bf J},\beta^*_j}(\tilde{Q}^{(k)}_j)-\mu^{k-1}(\tilde{Q}^{(k)}_j)\pl_{\infty}=O(r(d,J)^k).\end{equation}
That is the approximation for the $k$-th order primitive function was such that it preserved the sup-norm approximation for the $k-1$-th order primitive function, and as a consequence, we also obtained
\[\pl \mu^{k}_{{\bf J},\beta^*_j}(\tilde{Q}^{(k)}_j)-\mu^{k}(\tilde{Q}^{(k)}_j)\pl_{v}=O(r(d,J)^k).\]
This proves that for this $\beta^*$ for which $\pl Q-Q^k_{{\bf J},\beta^*}\pl_{\infty}=O(r(d,J)^{k+1})$ we also have
\[
\pl Q-Q^k_{{\bf J},\beta^*}\pl_v^*= O(r(d,J)^k).\]
In addition, in Section \ref{section7} we also showed that these same approximation results  apply to the loss-based projection  $Q^k_{{\bf J},\beta_{0,n}}$ of $Q_0$. This completes the proof. $\Box$

\subsection{Sectional variation norm of approximation error of  $k$-th order spline working-model specific MLE w.r.t. oracle MLE.}
Having obtained a bound for $\pl Q^k_{{\bf J},\beta_{0,n}}-Q_0\pl_v^*$ in previous subsection, where $Q^k_{{\bf J},\beta_{0,n}}=\arg\min_{Q\in D^{(k)}({\cal R}(d,{\bf J}))}P_0L(Q)($ is the oracle MLE of $Q_0$ for the given working model $D^{(k)}({\cal R}(d,{\bf J}))$, we now proceed with the MLE analysis $\pl Q^k_{{\bf J},\beta_n}-Q^k_{{\bf J},\beta_{0,n}}\pl_v^*$. 
Let $\beta_n=\arg\min_{\beta}P_n L(Q^k_{{\bf J},\beta})$ be the MLE.
By Lemma \ref{lemmasectvarexpression}  we have
\begin{eqnarray*}
\pl Q_{{\bf J},\beta_n}^k-Q_{{\bf J},\beta_{0,n}}^k\pl_v^*&=&
\sum_s \pl \sum_{j,j(1)=s}\bar{\phi}_j^{k-1}\left\{\mu^{k-1}_{{\bf J},\beta_{n,j}-\beta_{0,n,j}} (\tilde{Q}^{(k)}_j)\right\} \pl_{1,\mu}\\
&=&\sum_s \pl \sum_{j,j(1)=s}\bar{\phi}_j^{k-1}\left\{\mu^{k-1}\left(\tilde{Q}_{{\bf J},\beta_{n,j}-\beta_{0,n,j} } \right )  \right\} \pl_{1,\mu},
\end{eqnarray*}
where $\tilde{Q}_{{\bf J},\beta_{n,j}}=\sum_{u\in {\cal R}(j,{\bf J})}\beta_{n,j}(u)\phi_u^0$ and $\tilde{Q}_{{\bf J},\beta_{0,n,j}}=\sum_{u\in {\cal R}(j,{\bf J})}\beta_{0,n,j}(u)\phi_u^0$ represents the MLE and oracle MLE zero-order spline approximation of $Q^{(k)}_{j}$. 

\begin{lemma}\label{lemma3a}
For each $s\subset\{1,\ldots,d\}$, let \[
\Phi_s(\beta)\equiv \sum_{j,j(1)=s}\bar{\phi}_j^{k-1}\mu^{k-1}_{{\bf J},\beta_{n,j}-\beta_{0,n,j}}(\tilde{Q}^{(k)}_j).\]
For $j$ with $\mid j(k+1)\mid =0$, the $j$-specific terms reduce to $\bar{\phi}_j^{k-1}(\beta_{n,j}-\beta_{0,n,j})$.
Assume the conditions of 
Theorem \ref{theoremasnormgenphiunif} for pointwise asymptotic normality and uniform rate convergence of $\Phi_s(\beta_n)-\Phi_s(\beta_{0,n})$ so that  for some $m<\infty$
\[
\pl \Phi_s(\beta_n)-\Phi_s(\beta_{0,n})\pl_{\infty}=O_P((d_n/n)^{1/2}\lambda_{0,n}\log n).\]
Then, we have
\[
\pl Q_{{\bf J},\beta_n}^k-Q_{{\bf J},\beta_{0,n}}^k\pl_v^*=O_P((d_n/n)^{1/2}\lambda_{0,n}\log n).\]
\end{lemma}
As stated in this theorem, for log-likelihood behaving loss functions, $\lambda_{0,n}=1$ and generally one expects $\lambda_{0,n}=O(1)$ and at minimal $O(\log^m n)$ for some $m<\infty$.
\newline
{\bf Proof:}
This is just an application of the general asymptotic normality and uniform convergence theorems for analyzing $\Phi(\beta_n)-\Phi(\beta_{0,n})$ for these particular choices of $\Phi$. We only need $L_1$-norm convergence instead of sup-norm convergence, so the conditions stated that suffice for uniform convergence could be weakened. $\Box$

So we have proven the desired theorem.

\begin{theorem}\label{theoremconvsectvarnorm}
Let $Q_0\in D^{(k)}_M([0,1]^d)$ and  let $k\in \{1,\ldots\}$ be given.
Let $D^{(k)}({\cal R}(d,{\bf J}))=\{Q^k_{{\bf J},\beta}:\beta\}$ be given working model; $\beta_{0,n}=\arg\min_{\beta}P_0L(Q^k_{{\bf J},\beta})$ and $\beta_n=\arg\min_{\beta}P_n L(Q^k_{{\bf J},\beta})$ be the oracle MLE and MLE, respectively. 
We have $Q_{{\bf J},\beta_n}^k-Q_0=Q_{{\bf J},\beta_n-\beta_{0,n}}^k+Q_{{\bf J},\beta_{0,n}}^k-Q_0$. We have
\[
\pl Q_{{\bf J},\beta_{0,n}}^k-Q_0\pl_v^*=O(r(d,J)^k).\]
In addition, under conditions of Lemma \ref{lemma3a} (i.e., conditions in Theorem \ref{theoremasnormgenphiunif} needed for uniform convergence of $\Phi_s(\beta_n)-\Phi_s(\beta_{0,n})$), we have
\[
\pl Q_{{\bf J},\beta_n-\beta_{0,n}}^k\pl_v^*=O_P((J/n)^{1/2}\lambda_{0,n}\log n).\]
Then,
\[
\pl Q_{{\bf J},\beta_n}^k-Q_0\pl_v^*=O_P( r(d,J)^k+(J/n)^{1/2}\lambda_{0,n}\log n).\]
\end{theorem}
We can now apply this theorem to the $k$-th order spline sieve MLE that uses a data adaptive selector ${\bf J}_n$, for example, restricted to values ${\bf J}$ that already guarantee the desired trade-off between bias and variance. The above theorem provides the desired trade-off and optimal rate for $J_n$ for optimizing the rate of the sectional variation norm of   $\pl Q_{{\bf J}_n,\beta_n}-Q_0\pl_v^*$, while for supremum norm convergence the rate $J_n$ would be different balancing 
$r(d,J)^{k+1}$ and $(n/J)^{-1/2}\log n$ as stated in our theorems. 
For example, suppose $k=1$ and we choose a $J$ that optimizes $O( r^2(d,J) +(J/n)^{1/2}\log n)$ to make $Q^k_{{\bf J},\beta_n}$ achieve the best pointwise and sup-norm rate for $Q$ itself (up till possible $\log n$-factor).
Let's ignore $\log n$-factors. This means that $J$ is chosen so that $1/J^2\sim (J/n)^{1/2}$ and thus $J\sim n^{1/5}$ providing a sup-norm rate of $n^{-2/5}$ for $Q_{{\bf J},\beta_n}^1-Q_0$.
Then, the resulting rate for the sectional variation norm of $Q_{{\bf J},\beta_n}^1-Q_0$ is given by $r(d,J)^1+(J/n)^{1/2}\sim J^{-1}+(J/n)^{1/2}\sim n^{-1/5}+n^{-2/5}\sim n^{-1/5}$. However, by selecting $J\sim n^{1/3}$ we would have
\[
\pl Q_{{\bf J},\beta_n}^1-Q\pl_v^*\sim n^{-1/3}\]
up till power of $\log n$-factor. So by undersmoothing we obtain the desired optimal rate in sectional variation norm. The rate for $\pl Q_n-Q\pl_{\infty}$ then becomes
$J^{-2}+(J/n)^{1/2}\sim n^{-1/3}$ as well, thereby sacrificing the sup-norm rate of $n^{-2/5}$ by having it lowered to $n^{-1/3}$.
As another example, consider the case that $k$ is large such as $k=4$. 
For optimizing the pointwise or uniform rate we balance $1/J^5$ with $(J/n)^{1/2}$, giving $J\sim n^{1/11}$ and rate of convergence $n^{-5/11}$, an almost parametric rate of convergence.
For that rate $J\sim n^{1/11}$, the sectional variation norm behaves as $1/J^4+(J/n)^{1/2}\sim 
n^{-4/11}+n^{-5/11}\sim n^{-4/11}$. If the goal is to optimize the sectional variation norm one balances $1/J^4$ with $(J/n)^{1/2}$ giving $J\sim n^{1/9}$ and the sectional variation norm then converges at rate $n^{-4/9}$.  As one observes, as $k$ gets large, even when optimizing convergence  in sup-norm as the cross-validation selector, one ends up with very good rates of convergence of $Q_{{\bf J}_n,\beta_n}-Q_0$ for both the sup-norm as well as the sectional variation norm. Moreover, we observe that for optimizing the convergence rate in sectional variation norm one wants to undersmooth relative to the cross-validation selector.

\subsection{Implications for sectional variation norm of  $k$-th order spline HAL-MLE minus target function} In the definition of $k$-th order spline HAL and relax HAL the working model $D^{(k)}({\cal R}(d,{\bf J}_{0,n})$  (with ${\bf J}_{0,n}$ the sieve-MLE selector ${\bf J}_n$ but applied to independent sample $P_n^{\#}$) for the sieve MLE is replaced by a working model $D^{(k)}({\cal R}_{0,n})$, where ${\cal R}_{0,n}={\cal R}(P_n^{\#})$ is the independent version of the  data adaptively determined set ${\cal R}_n={\cal R}(P_n)$ of basis functions within a bigger working model $D^{(k)}({\cal R}(d,{\bf J}_{max,n})$ obtained by  using $L_1$-norm regularization of the MLE for the initial working model. We know that this set ${\cal R}_{0,n}$ 	is of the form ${\cal R}_{0,n}={\cal R}(d,{\bf J}_{0,n,hal})$ with the only twist that the knot-point sets used for modeling $\mu^k(Q^{(k)}_{\bar{s}(k+1)})$  are data adaptively determined (but based on independent data $P_n^{\#}$). However,  and as we argued in Appendix \ref{AppendixB}, HAL-regularization is able to select knot-point sets ${\cal R}(d,{\bf J}_{0,n}(\bar{s}(k+1)))$ for modeling each $\mu^k(Q^{(k)}_{\bar{s}(k+1)})$ that still satisfy the $O(r(d,J))$-$L^2$-approximation error property of ${\cal R}(d,{\bf J}(\bar{s}(k+1)))$ (since it achieves optimal rate of convergence as proven in Section \ref{section7}). So one can interpret and analyze HAL as a sieve MLE with a particular selector ${\bf J}_{0,n,hal}$ relative to the (e.g.) cross-validation selector ${\cal J}_{0,n}$ used by the sieve-MLE. In that manner, our theorems above generalize to HAL and relax HAL, as long as one makes the assumption that the HAL-selector ${\cal R}_{0,n}$ does yield knot-point sets (for zero order splines) ${\cal R}(d,{\bf J}_{0,n}(\bar{s}(k+1))$ with the $O(r(d,J))$-$L^2$-approximation error (which potentially might require a little undersmoothing).

In addition, generally speaking, the HAL and relax HAL should be an improvement on using the initial (unregularized HAL) sieve MLE implied by working model $D^{(k)}({\cal R}(d,{\bf J}_{max,n}))$. If one  uses the cross-validation selector for the $L_1$-norm it will be asymptotically better or equivalent with the choice of no-regularization, so that the HAL-MLE will achieve the optimal rate of convergence (pointwise and supremum norm) of the sieve-MLE implied by  initial working model. If one uses undersmoothing when selecting the $L_1$-norm in HAL-MLE and the initial ${\bf J}_{max,n}$, then that provides a way to push towards optimization of the rate of convergence w.r.t. sectional variation norm.

\section{Behavior of sectional variation norm of  zero-order spline sieve and relax HAL-MLE}\label{AppendixM}
In the previous section we studied the rate of convergence of the sectional variation norm of the approximation error of the first and higher order spline approximations  w.r.t. a first and higher order smooth target function. In this section we investigate the behavior of the sectional variation norm of the zero-order spline approximation error of a cadlag function with finite sectional variation norm.  In other words, here we consider the case $k=0$, while the previous section studied the case $k=1,2,\ldots$.
 
Let $Q_n^0=\sum_{u\in {\cal R}_{n}}\beta_n(u)\phi_u^0$ and $Q_{0,n}^0=\sum_{u\in {\cal R}_{n}}\beta_{0,n}(u)\phi_u^0$ be the MLE and oracle MLE for working model $D^{(0)}({\cal R}_{n})$, where ${\cal R}_{n}={\cal R}(P_n)$ has size $d_{n}=O(J_{n})$.
We like to understand the rates of $\pl Q_n^0\pl_v^*=\sum_{u\in {\cal R}_n}\mid \beta_n(u)\mid$ and $\pl Q_n^0-Q_{0,n}^0\pl_v^*=\sum_{u\in {\cal R}_n}\mid \beta_n-\beta_{0,n}\mid(u)$ as function of $J_n$ and sample size $n$. Of course , if we use the HAL-MLE with $L_1$-norm $C_n$, then we would have $\pl Q_n^0\pl_v^*=C_n$, so we are focussing on the case that we do not use $L_1$-norm regularization as in the sieve or relax HAL-MLE. For example, Theorem  \ref{theoremasnormgenphiunif} (or our theorems for $(Q_n^0-Q_{0,n}^0)$ itself) that if we replace ${\cal R}_n$ by its independent version ${\cal R}_{0,n}={\cal R}(P_n^{\#})$, 
 $\pl Q_n^0-Q_{0,n}^0\pl_{\infty}=O(\log n \lambda_{0,n} n(J_n/n)^{1/2})$, where $\lambda_{0,n}$ is the maximal eigenvalue of the covariance matrix of the score vector for an orthonormalized basis for $D^{(0)}({\cal R}_{0,n})$, known to be $1$ for log-likelihood behaving loss functions. We will assume the bound $\pl Q_n^0-Q_{0,n}^0\pl_{\infty}= O(\log n (J_n/n)^{1/2})$ as starting point of our analysis of the sectional variation norm of $Q_n^0$.

The following lemma bounds the sectional variation norm of $Q_n^0-Q_{0,n}^0$ in terms of its supremum norm.
\begin{lemma}
Recall for $Q_{\beta}^0\in D^{(0)}({\cal R}_n)$ we have $\pl Q_{\beta}^0\pl_v^*=\pl {\beta}\pl_1=\sum_{u\in {\cal R}_n}\mid \beta(u)\mid$. 

Suppose that $\pl Q_n^0-Q_{0,n}^0\pl_{\infty}=O_P(\log n(J_n/n)^{1/2})$. Then we also have that 
$\max_{u\in {\cal R}_n}\mid \beta_n(u)-\beta_{0,n}(u)\mid =O_P(\log n(J_n/n)^{1/2})$, and thus
\[
\pl Q_n^0-Q_{0,n}^0\pl_v^*=\sum_{u\in {\cal R}_n}\mid \beta_n-\beta_{0,n}\mid (u)=O_P( \log n J_n^{3/2}n^{-1/2}).\]
Thus if $J_n/( n^{1/3}(\log n)^{2/3})=O(1)$, then $\pl Q_n^0-Q_{0,n}^0\pl_v^*=O_P(1)$. 
In addition, $\pl Q_n^0\pl_v^*-\pl Q_{0,n}^0\pl_v^*=O_P(J_n^{3/2}n^{-1/2})$.

\end{lemma}
Note that $J_n\sim n^{1/3}$ is the rate optimizing the rate of convergence of $Q_n^0-Q_0$ providing the rate of convergence $d_0(Q_n,Q_0)=O_P(n^{-2/3})$ up till $\log n$-factors. 
Thus, if we select this choice, the zero order spline sieve MLE or relax HAL-MLE will have bounded sectional  variation norm. In addition, if one undersmooths up till $\log n$-factors, then this sectional variation norm will not grow faster to infinity than a $\log n$-factor (e.g., slow enough for establishing asymptotic equicontinuity conditions in proofs of asymptotic efficiency for smooth features).
\newline
{\bf Proof:}
Let the supremum norm of $Q_n^0-Q_{0,n}^0$ be $O_P(r_1(n))$.
For any knot-point $u\in {\cal R}_n$, for $u_1$ close enough to $u$ we have
 $\mid ( Q_n^0-Q_{0,n}^0)((u_1,u])\mid =\mid (\beta_n-\beta_{0,n})(u)\mid $. 
 From this identity it follows that  $\max_{u\in {\cal R}_n}\mid (\beta_n-\beta_{0,n})(u)\mid =O_P(r_1(n))$ and thereby $\pl Q_n^0-Q_{0,n}^0\pl_v^*=O_P(r_1(n)J_n)$.

We also know that for each $u$
 \[
 \mid \mid \beta_n\mid -\mid \beta_{0,n}\mid \mid  (u)\leq \mid \beta_n-\beta_{0,n}\mid (u).\]
 Thus, also
 \[
 \mid \pl \beta_n\pl_1-\pl \beta_{0,n}\pl_1\mid \leq \pl \beta_n-\beta_{0,n}\pl_1 =O_P(r_1(n)J_n).\Box\]

\section{Special results for zero-order spline linear and logistic regression when $X$ discrete on growing set of knot-points}
\label{AppendixN}
We consider the regression case where $Q_0(X)=E_0(Y\mid X)$.
Let ${\cal R}_n$ be a random set of knot-points with corresponding working model $D^{(0)}({\cal R}_n)=\{Q_{n,\beta}=\sum_{u\in {\cal R}_n}\beta(u)\phi_u:\beta\}$.  If $Y$ is continuous, then  $L(Q)=(Y-Q(X))^2$ and, if $Y\in \{0,1\}$, then $L(Q)=-Y\log m_Q(X)-(1-Y)\log(1-m_Q(X))$ with $m_Q(x)=1/(1+\exp(-Q(x)))$. Let $m_Q(X)=X$ when $Y$ is continuous. 
Let $\beta_n$ be the MLE and $\beta_{0,n}$ the oracle MLE. 
Let $Q_n=Q_{n,\beta_n}$ and $Q_{0,n}=Q_{n,\beta_{0,n}}$.
The score of $\beta(u)$  is given by $S_{\beta,u}=\phi_u (Y-m_{Q_{n,\beta_n}}(X))$, so that we have $P_n S_{\beta_n}=0=P_0 S_{\beta_{0,n}}$. In certain designs, one might be able to control $X$ as part of the experiment. 
In that case, we could, for example, select a (possible random and informative) set ${\cal R}_n$ and arrange that  $P(X=u)=1/d_n$ for each $u\in {\cal R}_n$,, or use a discrete non-uniform distribution still satisfying that 
\begin{equation}\label{keyconditiondiscretea}
\min_{u\in {\cal R}_n}P(X=u)>\delta /d_n\mbox{ for some $\delta>0$.}\end{equation} 
The discreteness of $X$ allows us to immediately obtain an expression \[
(m_{Q_n}-m_{Q_{0,n}})(x)=(P_n-P_0)I_x/p_0(x) (Y-m_{Q_n}),\]
 where $I_x(X)=I(X=x)$ as shown by following lemma.
\begin{lemma}
Consider the setting described above. 
We note that for all $x\in {\cal R}_n$, $m_{Q_{0,n}}(x)=E(Y|X=x)$ due to the discreteness of $X$. 
 For any orthonormal basis $\{\phi_u^*:u\in {\cal R}_n\}$ of $\{\phi_u:u\in {\cal R}_n\}$ in  Hilbert space $L^2(\mu_0)$ for functions of $X$, such as $L^2(P_0)$, with inner product $\langle h_1,h_2\rangle_{\mu_0}=\int h_1h_2d\mu_0$, we have for each $u\in {\cal R}_n$:
\[
P_0(m_{Q_n}-m_{Q_{0,n}})\phi_u^*=(P_n-P_0)\phi_u^*(Y-m_{Q_n}).\]
Due to $X$ being discrete on ${\cal R}_n$, we can select as orthonormal basis $\phi_u^*(X)=I_u(X)/p_0^{1/2}(u)$, $u\in {\cal R}_n$, where $p_0(u)=P_0(X=u)$ and $I_u(X)=I(X=u)$.
Then, this implies: for each $x\in {\cal R}_n$, we have
\[
(m_{Q_n}-m_{Q_{0,n}})(x)p_0^{1/2}(x)=(P_n-P_0)\phi_x^*(Y-m_{Q_n}),\]
and thus
\[
(m_{Q_n}-m_{Q_{0,n}})(x)= (P_n-P_0)\frac{I_x}{p_0(x)} (Y-m_{Q_n}).\]
\end{lemma}
The next lemma succeed in replacing  $Q_n$ by $Q_{0,n}$ to obtain the desired asymptotic linearity with a remainder that is $o_P(n/J_n)^{-1/2})$.
  
\begin{lemma}
Assume (\ref{keyconditiondiscretea}).
We have for $x\in {\cal R}_n$
\[
(m_{Q_n}-m_{Q_{0,n}})(x)=(P_n-P_0)\frac{I_x}{p_0(x)}(Y-m_{Q_{0,n}})+O_P(J_n/n).\]
So if $J_n/n\rightarrow 0$, then for $x\in {\cal R}_n$
\[
(m_{Q_n}-m_{Q_{0,n}})(x)=(P_n-P_0)\frac{I_x}{p_0(x)}(Y-m_{Q_{0,n}})+o_P((n/J_n)^{-1/2}).\]
Using that $Q_{0,n}(x)=Q_0(x)$ for $x\in {\cal R}_n$, it follows that for $x\in {\cal R}_n$:
\[
(m_{Q_n}-m_{Q_0})(x)=(P_n-P_0)\frac{I_x}{p_0(x)}(Y-m_{Q_0}(x))+o_P((n/J_n)^{-1/2}).\]
\end{lemma}
{\bf Proof:}
Using standard empirical process theory on a $O(J_n)$ dimensional class of functions, gives
 $\max_u \mid (P_n-P_0)(X=u)\mid =o_P(\log n (n/J_n)^{-1/2})$.
 For a given $x$, we have $\mid (P_n-P_0)(X=x)\mid =O_P((n/J_n)^{-1/2})$.
Therefore,
\begin{eqnarray*}
 (P_n-P_0)I_x/p_0(x) (m_{Q_n}-m_{Q_{0,n}})& =&(m_{Q_n}-m_{Q_{0,n}})(x) p_0(x)^{-1}(P_n-P_0)(X=x)\\
 &=&(m_{Q_n}-m_{Q_{0,n}})(x)O(J_n)o_P(n^{-1/2}J_n^{-1/2})\\
 &=&(m_{Q_n}-m_{Q_{0,n}})(x)O_P( n^{-1/2}J_n^{1/2})  .
 \end{eqnarray*}
 Thus, \[
 (m_{Q_n}-m_{Q_{0,n}})(x)=(P_n-P_0)I_x/p_0(x)(Y-m_{Q_{0,n}})+(m_{Q_n}-m_{Q_{0,n}})(x)O_P(n^{-1/2}J_n^{1/2}),\]
which implies 
\[(m_{Q_n}-m_{Q_{0,n}}))(x)=O_P(J_n^{1/2}n^{-1/2}))+(m_{Q_n}-m_{Q_{0,n}})(x)O_P(n^{-1/2}J_n^{1/2}).\]
If $n/J_n\rightarrow \infty$, then it follows that $(m_{Q_n}-m_{Q_{0,n}})(x)=O_P((n/J_n)^{-1/2})$.
So then, the remainder in last displayed expression for $(m_{Q_n}-m_{Q_{0,n}})(x)$ is $O_P(n^{-1}J_n)$. This proves that for each $x\in {\cal R}_n$
 \[ (m_{Q_n}-m_{Q_{0,n}})(x)=(P_n-P_0)I_x/p_0(x)(Y-m_{Q_{0,n}})+O_P((n/J_n)^{-1}) .\Box\]

{\bf Asymptotic normality:}
So for $x\in {\cal R}_n$ we have
\[ 
(n/d_n)^{1/2}(m_{Q_n}-m_{Q_0})(x)= (n/d_n)^{1/2}(P_n-P_0)\frac{I_x}{p_0(x)}(Y-m_{Q_0}(x))
+o_P(1).\]
The leading term on right-hand side is now a sum of mean zero independent random variables. The variance of 
$(n/d_n)^{1/2}(P_n-P_0)I_x/p_0(x) (Y-m_{Q_0}(x)) $ is given by:
\begin{eqnarray*}
\tilde{\sigma}_{0,n}^2(x)&\equiv &d_n^{-1}p_0(x)^{-2} P_0 \sigma^2_{0} I_x\\
&=&
d_n^{-1}p_0(x)^{-1}\sigma^2_{0}(x),
\end{eqnarray*}
where $\sigma^2_{0}(x)=E((Y-Q_0(x))^2\mid X=x)$  if $Y$ is continuous and  $\sigma^2_{0}(x)=Q_0(x)(1-Q_0(x))$ if $Y$ is binary. We know that $p_0(x)$ depends on $n$, which is suppressed in our notation. 
By assumption (\ref{keyconditiondiscretea}) we have $p_0(x)\geq \delta 1/d_n$ so that $d_n^{-1}/p_0(x)=O(1)$ or might converge to a fixed limit $\mu_0(x)$ of $d_np_{0,n}(x)$
Thus, we have that $\tilde{\sigma}_{0,n}^2(x)=O(1)$ or it converges to $\mu_0(x)\sigma^2_0(x)$.
The final issue is that the above relied on $x\in {\cal R}_n$ and we allowed ${\cal R}_n$ to be informative and random. Therefore we can only apply it to $x$ for which we know that $P(x\in {\cal R}_n)=1$ for $n>N(x)$ for some $N(x)<\infty$.  
\begin{theorem}\label{theorempointwisediscrete}
Consider the setting described above and 
assume  (\ref{keyconditiondiscretea}) and that $n/d_n\rightarrow\infty$.

For each $x$ for which $P(x\in {\cal R}_n)=1$ for $n>N(x)$ for some $N(x)<\infty$, 
 we have
 \[
\tilde{\sigma}_{0,n}^{-1}(x)(n/d_n)^{1/2}(m_{Q_n}-m_{Q_0})(x)\Rightarrow_d N(0,1).\]
For any $(x_1,\ldots,x_m)\in {\cal R}_n^m$ with $x_{j_1}\not =x_{j_2}$ for $j_1\not =j_2$, we have
that $(\tilde{\sigma}_{0,n}^{-1}(x_j)(n/d_n)^{1/2}(m_{Q_n}-m_{Q_0})(x_j):j=1,\ldots,m)$
converges in distribution to a multivariate normal with mean zero and covariance matrix being the (diagonal) identity matrix. 
 \end{theorem}
We note that by selecting $d_n=n/\log n$, we still obtain the desired asymptotic normality but the rate of convergence is as slow as $1/\log n$. By selecting $d_n=O(1)$, one obtains the parametric rate $n^{-1/2}$. So there is a trade-off between the size of the set ${\cal R}_n$ at which $x$ we want to estimate $Q_0(x)$ and the  rate of convergence  at which $Q_n(x)-Q_0(x)$ converges to zero. 

\subsection{Rate of convergence w.r..t supremum norm}
We have
\[(m_{Q_n}-m_{Q_{0,n}})(x)= (P_n-P_0)\frac{I_x}{p_0(x)} (Y-m_{Q_n}),\]
and
\[
\sup_{x\in [0,1]^d} \left | (P_n-P_0)\frac{I_x}{p_0(x)}(m_{Q_n}-m_{Q_{0,n}})\right | =o_P(\pl Q_n-Q_{0,n}\pl_{\infty}\log n n^{-1/2}J_n^{1/2}),
\]
using that $\sup_x\mid (P_n-P_0)(X=x)\mid =o_P(\log n (n/J_n)^{-1/2})$.
Therefore, if $\log n n^{-1/2}J_n^{1/2}=O(1)$, then
\[
\sup_x\left |  m_{Q_n}-m_{Q_{0,n}}\right |  (x)=\sup_x \left | (P_n-P_0)I_x/p_0(x)(Y-m_{Q_{0,n}})\right | +o_P(\pl m_{Q_n}-m_{Q_{0,n}}\pl_{\infty}).\]
We define the following class of functions of $O$
\[
{\cal F}_n=\left\{I_x/p_0^{1/2}(x)(Y-m_{Q_{0,n}}): x\in [0,1]^d\right\},\]
which are standardized to have uniformly bounded variance, by our assumption on $\sup_x \tilde{\sigma}_{0,n}^2(x)=O_P(1)$.

We note  that $J_n^{-1/2}{\cal F}_n=\{I_x(Y-m_{Q_{0,n}}):x\in [0,1]^d\}$  has an envelop $Q_n$ with $\pl Q_n\pl_{P_0}<M$ for some $M<\infty$.
The entropy integral $J(\delta,J_n^{-1/2}{\cal F}_n)\equiv \int_0^{\delta}\sup_Q\sqrt{\log N(\epsilon,J_n^{-1/2}{\cal F}_n,L^2(Q))} d\epsilon$ of $J_n^{-1/2} {\cal F}_n $ is bounded by $-\int_0^{\delta} (\log \epsilon)^{1/2} d\epsilon$.  This can be conservatively bounded by $-\int_0^{\delta}\log \epsilon d\epsilon = -\delta \log \delta +\delta\sim -\delta \log \delta$.
Therefore
\begin{eqnarray*}
(n/d_n)^{1/2}\pl m_{Q_n}-m_{Q_{0,n}}\pl_{\infty}&=&
\sup_{Q\in {\cal F}_n}\mid n^{1/2}(P_n-P_0)Q\mid  \\
&=&
d_n^{1/2}\sup_{Q\in d_n^{-1/2}{\cal F}_n,\pl Q\pl_{P_0}\leq d_n^{-1/2}}\mid n^{1/2}(P_n-P_0)Q\mid \\
&\sim&d_n^{1/2} J(d_n^{-1/2}{\cal F}_n,\delta= d_n^{-1/2})\\
&= &o(d_n^{1/2} d_n^{-1/2} \log d_n)\\
&=&o(\log d_n).
\end{eqnarray*}
This proves the following theorem.
\begin{theorem}
Consider the setting of Theorem \ref{theorempointwisediscrete} and assume that $P_{X,0}$ is discrete with support ${\cal R}_n$  satisfying (\ref{keyconditiondiscretea}). 
In addition, assume $(\log n) n^{-1/2}d_n^{1/2}=O(1)$.
Then,
\[
\sup_{x\in [0,1]^d}\mid Q_n-Q_{0,n}\mid (x)= o_P(n^{-1/2}d_n^{1/2}\log d_n).
\]
\end{theorem}

\end{document}